\documentclass[12pt]{amsart}
\usepackage{latexsym,amsxtra}

\calclayout

\renewcommand{\:}{\colon}
\newcommand{\+}{\nobreakdash-} 
\renewcommand{\.}{\text{$\mskip .5\thinmuskip$}}
\renewcommand{\;}{,\medspace}

\global\let\le\undefined
\global\let\ge\undefined
\DeclareMathSymbol{\le}{\mathrel}{AMSa}{"36}      
\DeclareMathSymbol{\ge}{\mathrel}{AMSa}{"3E}      
\DeclareMathSymbol{\birarrow}{\mathrel}{AMSa}{"13}  

\newcommand{\bu}{{\text{\smaller\smaller$\scriptstyle\bullet$}}}

\newcommand{\oc}{\mathbin{\text{\smaller$\square$}}}

\newcommand{\ocn}{\odot}
\newcommand{\up}{\uparrow}
\newcommand{\down}{\downarrow}

\newcommand{\lrarrow}{\.\relbar\joinrel\relbar\joinrel\rightarrow\.}

\newcommand{\rarrow}{\longrightarrow}
\newcommand{\larrow}{\longleftarrow}
\newcommand{\mpsto}{\longmapsto}
\newcommand{\ot}{\otimes}
\newcommand{\dsb}{\dotsb}
\renewcommand{\d}{\partial}
\newcommand{\eps}{\varepsilon}

\newcommand{\g}{{\mathfrak g}}
\newcommand{\oA}{{\overline{A}}}
\newcommand{\oC}{{\overline{C}}}

\newcommand{\sD}{{\mathsf D}}
\newcommand{\sH}{{\mathsf H}}
\newcommand{\sA}{{\mathsf A}}
\newcommand{\sB}{{\mathsf B}}
\newcommand{\sC}{{\mathsf C}}
\newcommand{\sF}{{\mathsf F}}
\newcommand{\sE}{{\mathsf E}}
\newcommand{\sM}{{\mathsf M}}

\newcommand{\sDG}{\mathsf{DG}}
\newcommand{\Hot}{\mathsf{Hot}}
\newcommand{\Acycl}{\mathsf{Acycl}}
\newcommand{\Quis}{\mathsf{Quis}}
\newcommand{\FQuis}{\mathsf{FQuis}}

\newcommand{\inj}{{\mathsf{inj}}}
\newcommand{\proj}{{\mathsf{proj}}}
\newcommand{\fl}{{\mathsf{fl}}}
\newcommand{\co}{{\mathsf{co}}}
\newcommand{\ctr}{{\mathsf{ctr}}}
\renewcommand{\b}{{\mathsf{b}}}
\newcommand{\abs}{{\mathsf{abs}}}
\newcommand{\sgr}{{\mathsf{gr}}}
\newcommand{\gr}{{\mathrm{gr}}}
\newcommand{\op}{{\mathsf{op}}}
\newcommand{\rop}{{\mathrm{op}}}
\renewcommand{\ss}{{\mathrm{ss}}}
\newcommand{\dg}{{\mathsf{dg}}}
\newcommand{\cdg}{{\mathsf{cdg}}}
\newcommand{\conilp}{{\mathsf{conilp}}}
\newcommand{\aug}{{\mathsf{aug}}}
\newcommand{\coaug}{{\mathsf{coaug}}}
\newcommand{\fin}{{\mathsf{fin}}}
\newcommand{\fg}{{\mathsf{fg}}}
\newcommand{\fgp}{{\mathsf{fgp}}}
\newcommand{\fd}{{\mathsf{fd}}}
\newcommand{\fpid}{{\mathsf{fpid}}}
\newcommand{\fil}{{\mathsf{fil}}}

\newcommand{\vect}{{\operatorname{\mathsf{--vect}}}}
\newcommand{\modl}{{\operatorname{\mathsf{--mod}}}}
\newcommand{\modr}{{\operatorname{\mathsf{mod--}}}}
\newcommand{\modrfl}{{\operatorname{\mathsf{mod_{fl}--}}}}
\newcommand{\modrcoh}{{\operatorname{\mathsf{mod^{coh}--}}}}
\newcommand{\modrfc}{{\operatorname{\mathsf{mod_{fil}^{coh}--}}}}
\newcommand{\modrgc}{{\operatorname{\mathsf{mod_{gr}^{coh}--}}}}
\newcommand{\comodl}{{\operatorname{\mathsf{--comod}}}}
\newcommand{\comodr}{{\operatorname{\mathsf{comod--}}}}
\newcommand{\comodrinj}{{\operatorname{\mathsf{comod_{inj}--}}}}
\newcommand{\contra}{{\operatorname{\mathsf{--contra}}}}
\newcommand{\alg}{{\operatorname{\mathsf{--alg}}}}
\newcommand{\coalg}{{\operatorname{\mathsf{--coalg}}}}
\newcommand{\modrdown}{{\operatorname{\mathsf{mod}}}^\down
        {\operatorname{\mathsf{--}}}}
\newcommand{\modrup}{{\operatorname{\mathsf{mod}}}^\up
        {\operatorname{\mathsf{--}}}}
\newcommand{\modrfu}{{\operatorname{\mathsf{mod}}}
        _{\operatorname{\mathsf{fil}}}^\up
        {\operatorname{\mathsf{--}}}}

\newcommand{\boZ}{{\mathbb Z}}
\newcommand{\boQ}{{\mathbb Q}}
\newcommand{\boL}{{\mathbb L}}
\newcommand{\boR}{{\mathbb R}}
\newcommand{\boI}{{\mathbb I}}

\DeclareMathOperator{\Hom}{Hom}
\DeclareMathOperator{\End}{End}
\DeclareMathOperator{\Tor}{Tor}
\DeclareMathOperator{\Ext}{Ext}
\DeclareMathOperator{\Cotor}{Cotor}
\DeclareMathOperator{\Cohom}{Cohom}
\DeclareMathOperator{\Coext}{Coext}
\DeclareMathOperator{\Ctrtor}{Ctrtor}
\DeclareMathOperator{\Tot}{Tot}
\DeclareMathOperator{\Br}{Bar}
\DeclareMathOperator{\Cb}{Cob}

\DeclareMathOperator{\id}{id}
\DeclareMathOperator{\im}{im}
\DeclareMathOperator{\cone}{cone}
\DeclareMathOperator{\coker}{coker}

\renewcommand{\O}{{\mathcal O}}
\newcommand{\B}{{\mathcal B}}
\newcommand{\A}{{\mathcal A}}
\newcommand{\K}{{\mathcal K}}
\renewcommand{\L}{{\mathcal L}}
\newcommand{\M}{{\mathcal M}}
\newcommand{\N}{{\mathcal N}}
\renewcommand{\P}{{\mathcal P}}
\newcommand{\Q}{{\mathcal Q}}
\newcommand{\D}{{\mathcal D}}
\newcommand{\E}{{\mathcal E}}
\newcommand{\F}{{\mathcal F}}
\newcommand{\T}{{\mathcal T}}
\newcommand{\END}{{\mathcal E}nd}
\newcommand{\HOM}{{\mathcal H}om}

\theoremstyle{plain}
\newtheorem*{thm}{Theorem}
\newtheorem*{thm1}{Theorem 1}
\newtheorem*{thm2}{Theorem 2}
\newtheorem*{lem}{Lemma}
\newtheorem*{lem1}{Lemma 1}
\newtheorem*{lem2}{Lemma 2}
\newtheorem*{lem3}{Lemma 3}
\newtheorem*{prop}{Proposition}
\newtheorem*{cor}{Corollary}
\newtheorem*{cor1}{Corollary 1}
\newtheorem*{cor2}{Corollary 2}
\theoremstyle{definition}
\newtheorem*{rem}{Remark}
\newtheorem*{rem1}{Remark 1}
\newtheorem*{rem2}{Remark 2}
\newtheorem*{rem3}{Remark 3}
\newtheorem*{ex}{Example}
\newtheorem*{exs}{Examples}
\newtheorem*{qst}{Question}

\newcommand{\Ainfty}{\mathrm{A}_\infty}

\newcommand{\Section}[1]{\bigskip\section{#1}\medskip}
\begin{document}

\title{Two kinds of derived categories, Koszul duality, \\
and comodule-contramodule correspondence}
\author{Leonid Positselski}

\address{Sector of Algebra and Number Theory, Institute for
Information Transmission Problems, Bolshoy Karetny per.~19 str.~1,
Moscow 127994, Russia}
\email{posic@mccme.ru}


\maketitle

\tableofcontents

\section*{Introduction}
\medskip

\subsection{{}}  \label{introduction-two-completions}
 A common wisdom says that difficulties arise in Koszul duality
because important spectral sequences diverge.
 What really happens here is that one considers the spectral
sequence of a complex endowed with, typically, a decreasing
filtration which is not complete.
 Indeed, the spectral sequence of a complete and cocomplete
filtered complex always converges in the relevant sense~\cite{EM1}.
 The solution to the problem, therefore, is to either replace
the complex with its completion, or choose a different filtration.
 In this paper, we mostly follow the second path.
 This involves elaboration of the distinction between two kinds
of derived categories, as we will see below.

 The first conclusion is that one has to pay attention to completions
if one wants one's spectral sequences to converge.
 What this means in the case of the spectral sequence related to
a bicomplex is that the familiar picture of two spectral sequences
converging to the same limit splits in two halves when the bicomplex
becomes infinite enough.
 The two spectral sequences essentially converge to the cohomology
of two different total complexes.
 To obtains those, one takes infinite products in the ``positive''
direction along the diagonals and infinite direct sums in
the ``negative'' direction (like in Laurent series).
 The two  possible choices of the ``positive'' and ``negative''
directions give rise to the two completions.
 The word ``essentially'' here is to be understood as ``ignoring
the delicate, but often manageable issues related to nonexactness
of the inverse limit''. 

\subsection{{}}  \label{introduction-classical-two-kinds}
 This alternative between taking infinite direct sums and infinite
products when constructing the total complex leads to the classical
distinction between differential derived functors of the first and
the second kind~\cite[section~I.4]{HMS}.
 Roughly speaking, one can consider a DG-module either as
a deformation of its cohomology or as a deformation of itself
considered with zero differential; the spectral sequences related
to the former and the latter kind of deformations essentially
converge to the cohomology of the differential derived functors
of the first and the second kind, respectively.

 \emph{Derived categories of the first} and \emph{the second kind}
are intended to serve as the domains of the differential derived
functors of the first and the second kind.
 This does not always work as smoothly as one wishes; one
discovers that, for technical reasons, it is better to consider
derived categories of the first kind for \emph{algebras} and
derived categories of the second kind for \emph{coalgebras}.
 The distinction between the derived functors/categories of the first
and the second kind is only relevant when certain finiteness
conditions no longer hold; this happens when one considers either
unbounded complexes, or differential graded modules.

 Let us discuss the story of two derived categories in more detail.
 When the finiteness conditions do hold, the derived category can
be represented in two simple ways.
 It is both the quotient category of the homotopy category by
the thick subcategory of complexes with zero cohomology and
the triangulated subcategory of the homotopy category formed by
the complexes of projective or injective objects.
 In the general case, this simple picture splits in two halves.
 The derived category of the first kind is still defined as
the quotient category of the homotopy category by the thick
subcategory of complexes (DG\+modules,~\dots)\ with zero cohomology.
 It can be also obtained as a full subcategory of the homotopy
category, but the description of this subcategory is more
complicated~\cite{Spal,Kel,BL}.
 On the other hand, the derived category of the second kind is defined
as the quotient category of the homotopy category by a thick
subcategory with a rather complicated description.
 At the same time, it is equivalent to the full subcategory of
the homotopy category formed by complexes (DG\+comodules,
DG\+contramodules,~\dots)\ which become injective or projective when
considered without the differential.

\subsection{{}}
 The time has come to mention that there exist two kinds of module
categories for a coalgebra: besides the familiar \emph{comodules},
there are also \emph{contramodules}~\cite{EM2}.
 Comodules can be thought of as discrete modules which are unions
of their finite-dimensional subcomodules, while contramodules are
modules where certain infinite summation operations are defined.
 For example, the space of linear maps from a comodule to any vector
space has a natural contramodule structure.

 The derived category of the first kind is what is known as just
\emph{the derived category}: the unbounded derived category,
the derived category of DG\+modules, etc.
 The derived category of the second kind comes in two dual versions:
the \emph{coderived} and the \emph{contraderived category}.
 The coderived category works well for comodules, while
the contraderived category is useful for contramodules.
 The classical notion of a DG\+(co)algebra itself can be generalized
in two ways; the derived category of the first kind is well-defined
for an $\Ainfty$\+algebra, while the derived category of the second
kind makes perfect sense for a \emph{CDG\+coalgebra}.

 Other situations exist when derived categories of the second kind
are well-behaved.
 One of them is that of a CDG\+ring whose underlying graded ring has
a finite homological dimension.
 In this case, the coderived and contraderived categories coincide.
 In particular, this includes the case of a CDG\+algebra whose
underlying graded algebra is free.
 Such CDG\+algebras can be thought of as strictly counital
\emph{curved $\Ainfty$\+coalgebras}; CDG\+modules over the former
with free and cofree underlying graded modules correspond to strictly
counital curved $\Ainfty$\+comodules and $\Ainfty$\+contramodules
over the latter.
 For a \emph{cofibrant} associative DG\+algebra, the derived,
coderived, and contraderived categories of DG\+modules coincide.
 Since any DG\+algebra is quasi-isomorphic to a cofibrant one, it
follows that the derived category of DG\+modules over any DG\+algebra
can be also presented as a coderived and contraderived category.

 The functors of forgetting the differentials, assigning graded
(co/contra)modules to CDG\+(co/contra)modules, play a crucial role
in the whole theory of derived categories of the second kind.
 So it is helpful to have versions of these functors defined for
arbitrary DG\+categories.
 An attempt to obtain such forgetful functors leads to a nice
construction of an \emph{almost involution} on the category of
DG\+categories.
 The related constructions for CDG\+rings and CDG\+coalgebras
are important for the nonhomogeneous quadratic duality theory,
particularly in the relative case~\cite{P}.

\subsection{{}}
 Now let us turn to (derived) Koszul duality.
 This subject originates from the classical Bernstein--Gelfand--Gelfand
duality (equivalence) between the bounded derived categories of 
finitely generated graded modules over the symmetric and exterior
algebras with dual vector spaces of generators~\cite{BGG}.
 Attempting to generalize this straightforwardly to arbitrary algebras,
one discovers that many restricting conditions have to be imposed:
it is important here that one works with algebras over a field, that
the algebras and modules are graded, that the algebras are Koszul,
that one of them is finite-dimensional, while the other is Noetherian
(or at least coherent) and has a finite homological dimension.

 The standard contemporary source is~\cite{BGSoe}, where many of these
restrictions are eliminated, but it is still assumed that
everything happens over a semisimple base ring, that the algebras
and modules are graded, and that the algebras are Koszul.
 In~\cite{BGS}, Koszulity is not assumed, but positive grading and
semisimplicity of the base ring still are.
 The main goal of this paper is to work out the Koszul duality for
ungraded algebras and coalgebras over a field, and more generally,
differential graded algebras and coalgebras.
 In this setting, the Koszulity condition is less important, although
it allows to obtain certain generalizations of the duality results.
 As to the duality over a base more general than a field, in this paper
we only consider the special case of $\D$--$\Omega$ duality, i.~e.,
the duality between complexes of modules over the ring of differential
operators and (C)DG\+modules over the de Rham (C)DG\+algebra of
differential forms (see~\ref{introd-d-omega}).
 The ring of functions (or sections of the bundle of endomorphisms of
a vector bundle) is the base ring in this case.
 For a more general treatment of the relative situation, we
refer the reader to~\cite[Chapter~11]{P}, where a version of Koszul
duality is obtained for a base coring over a base ring.

 The thematic example of nonhomogeneous Koszul duality over a field is
the relation between complexes of modules over a Lie algebra~$\g$ and
DG\+comodules over its standard homological complex.
 Here one discovers that, when $\g$ is reductive, the standard
homological complex with coefficients in a nontrivial irreducible
$\g$\+module has zero cohomology---even though it is not contractible,
and becomes an injective graded comodule when one forgets
the differential.
 So one has to consider a version of derived category of DG\+comodules
where certain acyclic DG\+comodules survive if one wishes this 
category to be equivalent to the derived category of $\g$\+modules.
 That is how derived categories of the second kind appear in Koszul
duality~\cite{Fl,Lef,Kel2}.

\subsection{{}}
 Let us say a few words about the homogeneous case.
 In the generality of DG-(co)algebras, the homogeneous situation is
distinguished by the presence of an additional positive grading
preserved by the differentials.
 Such a grading is well-known to force convergence of the spectral
sequences, so there is no difference between the derived categories
of the first and the second kind in the homogeneous case.
 It is very essential here that the grading be indeed positive
(or negative) on the DG\+(co/contra)modules as well as
the DG\+(co)algebras, as one can see already in the example of
the duality between the symmetric and the exterior algebras in
one variable, $S=k[x]$ and $\Lambda=k[\eps]/\eps^2$.
 The graded $S$\+module $M=k[x,x^{-1}]$ corresponds to the acyclic
complex of $\Lambda$\+modules $K=(\dsb\rarrow\Lambda\rarrow\Lambda
\rarrow\dsb)$ whose every term is $\Lambda$ and every differential
is the multiplication with~$\eps$.

 The acyclic, but not contractible complex $K$ of projective and
injective $\Lambda$\+modules provides the simplest way to distinguish
between the derived categories of the first and the second kind.
 In derived categories of the first kind, it represents the zero
object and is not adjusted to various derived functors, while
in derived categories of the second kind, it is adjusted to
derived functors and represents a nonzero object.
 So the $S$\+module $M$ has to be excluded from the category of
modules under consideration for a duality between conventional
derived categories of $S$\+modules and $\Lambda$\+modules to hold.
 The positivity condition on the internal grading accomplishes that
much in the homogeneous case.
 All the stronger conditions on the gradings considered
in~\cite{BGSoe,BGS} are unnecessary for the purposes of establishing
derived Koszul duality.

\subsection{{}} \label{introd-d-omega}
 Several attempts have been made in the literature~\cite{Kap,BD2} to
obtain an equivalence between the derived category of modules over
the ring/sheaf of differential operators and an appropriately defined
version of derived category of DG\+modules over the de Rham complex
of differential forms.
 More generally, let $X$ be a smooth algebraic variety and $\E$ be
a vector bundle on $X$ with a global connection~$\nabla$.
 Let $\Omega(X,\END(\E))$ be the sheaf of graded algebras of
differential forms with coefficients in the vector bundle $\END(\E)$ of
endomorphisms of~$\E$, \ $d_\nabla$~be the de Rham differential in
$\Omega(X,\END(\E))$ depending on the connection~$\nabla$, and
$h_\nabla\in\Omega^2(X,\END(\E))$ be the curvature of~$\nabla$.
 Then the triple consisting of the sheaf of graded rings
$\Omega(X,\END(\E))$, its derivation~$d_\nabla$, and
the section~$h_\nabla$ is a sheaf of CDG\+rings over~$X$.
 The derived category of modules over the sheaf of rings $\D_{X,\E}$ of
differential operators acting in the sections of the vector bundle $\E$
on~$X$ turns out to be equivalent to the coderived category (and
the contraderived category, when $X$ is affine) of quasi-coherent
CDG\+modules over the sheaf of CDG\+rings $\Omega(X,\END(\E))$.
 The assumption about the existence of a global connection in~$\E$
can be dropped (see subsection~\ref{d-omega-duality-functors}
for details).

\subsection{{}}
 Yet another very good reason for considering derived categories of 
the second kind is that in their terms a certain relation between
comodules and contramodules can be established.
 Namely, the coderived category of CDG\+comodules and the contraderived
category of CDG\+contramodules over a given CDG\+coalgebra are
naturally equivalent.
 We call this phenomenon the \emph{comodule-contramodule
correspondence}; it appears to be almost as important as
the Koszul duality.

 One can generalize the comodule-contramodule correspondence to
the case of strictly counital curved $\Ainfty$\+comodules and
$\Ainfty$\+contramodules over a curved $\Ainfty$\+coal\-gebra by
considering the derived category of the second kind for CDG\+modules
over a CDG\+algebra whose underlying graded algebra is a free
associative algebra.

\subsection{{}}
 This paper can be thought of as an extended introduction to
the monograph~\cite{P}, as indeed, its key ideas precede those
of~\cite{P} both historically and logically.
 It would had been all but impossible to invent the use of exotic
derived categories for the purposes of~\cite{P} if these were not
previously discovered in the work presented below.
 Nevertheless, most results of this paper are not covered by~\cite{P}, 
since it is written in the generality of DG- and CDG\+modules,
comodules, and contramodules, while~\cite{P} deals with nondifferential
semi(contra)modules most of the time.
 This paper is focused on Koszul duality, while the goal of~\cite{P}
is the semi-infinite cohomology.

 The fact that exotic derived categories arise in Koszul duality was
essentially discovered by Hinich~\cite{Hin2}, whose ideas were
developed by Lef\`evre-Hasegawa~\cite[Chapitres~1 and~2]{Lef};
see also  Fl\o ystad~\cite{Fl}, Huebschmann~\cite{Hue}, and
Nicol\'as~\cite{Nic}.
 The terminology of ``coderived categories'' was introduced in Keller's
exposition~\cite{Kel2}.
 However, the definition of coderived categories in~\cite{Lef,Kel2} was
not entirely satisfactory, in our view, in that the right hand side of
the purported duality is to a certain extent defined in terms of
the left hand side (the approach to $\D$\+$\Omega$ duality developed
in~\cite[section~7.2]{BD2} had the same problem).
 This defect is corrected in the present paper.

 The analogous problem is present in the definitions of model category
structures on the categories of DG\+coalgebras by Hinich~\cite{Hin2}
and Lef\`evre-Hasegawa~\cite{Lef}.
 This proved harder to do away with: we obtain various explicit
descriptions of the distinguished classes of morphisms of
CDG\+coalgebras independent of the Koszul duality, but the duality
functors are still used in the proofs.

 In addition, we emphasize contramodules and CDG\+coalgebras, whose
role in the derived categories of the second kind and derived Koszul
duality business does not seem to have been appreciated enough.

\subsection{{}}
 Now let us describe the content of this paper in more detail.
 In Section~1 we obtain two semiorthogonal decompositions of
the homotopy category of DG\+modules over a DG\+ring, providing
injective and projective resolutions for the derived category of
DG\+modules.
 We also consider flat resolutions and use them to define
the derived functor $\Tor$ for a DG\+ring.
 Besides, we construct a t-structure on the derived category
of DG\+modules over an arbitrary DG\+ring.
 In anticipation of the forthcoming paper~\cite{P2}, we also discuss
silly filtrations on the same derived category.
 This section contains no new results; it is included for the sake
of completeness of the exposition.

 The derived categories of DG\+comodules and DG\+contramodules and
the differential derived functors $\Cotor^{C,I}$ and $\Coext_C^I$
of the first kind for a DG\+coalgebra $C$ are briefly discussed
in Section~2, the proofs of the main results of this section being
postponed to Sections~5 and~7.
 Partial results about injective and projective resolutions for
the coderived and contraderived categories of a CDG\+ring are
obtained in Section~3.
 The finite homological dimension, Noetherian, coherent, and Gorenstein
cases are considered; in the former situation, a natural definition of
the differential derived functor $\Tor^{B,I\!I}$ of the second kind
for a CDG\+ring~$B$ is given.
 In addition, we construct an ``almost involution'' on the category
of DG\+categories.

 In Section~4 we construct semiorthogonal decompositions of
the homotopy categories of CDG\+comodules and CDG\+contramodules
over a CDG\+coalgebra, providing injective and projective
resolutions for the coderived category of CDG\+comodules and
the contraderived category of CDG\+contramodules.
 We also define the differential derived functors $\Cotor$, $\Coext$,
and $\Ctrtor$ for a CDG\+coalgebra, and give a sufficient condition
for a morphism of CDG\+coalgebras to induce equivalences of
the coderived and contraderived categories.
 The comodule-contramodule correspondence for a CDG\+coalgebra is
obtained in Section~5.

 Koszul duality (or ``triality'', as there are actually two module
categories on the coalgebra side) is studied in Section~6.
 Two versions of the duality theorem for (C)DG\+modules, CDG\+comodules,
and CDG\+contramodules are obtained, one valid for conilpotent
CDG\+coalgebras only and one applicable in the general case.
 We also construct an equivalence between natural localizations of
the categories of DG\+algebras (with nonzero units) and conilpotent
CDG\+coalgebras.

 We discuss the derived categories of $\Ainfty$\+modules and
the co/contraderived categories of curved $\Ainfty$\+co/contramodules
in Section~7.
 We explain the relation between strictly unital $\Ainfty$\+algebras
and coaugmented CDG\+coalgebra structures on graded tensor coalgebras,
and use it to prove the standard results about strictly unital
$\Ainfty$\+modules. 
 The similar approach to strictly counital curved $\Ainfty$\+coalgebras
yields the comodule-contramodule correspondence in the $\Ainfty$ case.

 Model category structures (of the first kind) for DG\+modules over
a DG\+ring and model category structures (of the second kind) for
CDG\+comodules and CDG\+contramodules over a CDG\+coalgebra are
constructed in Section~8.
 We also obtain model category structures of the first kind for
DG\+comodules and DG\+contra\-modules over a DG\+coalgebra, and model
category structures of the second kind for CDG\+modules over
a CDG\+ring in the finite homological dimension, Noetherian,
coherent, and Gorenstein cases.
 Quillen equivalences related to the comodule-contramodule
correspondence and Koszul duality are discussed.

 We consider the model categories of DG\+algebras and conilpotent
CDG\+coalgebras in Section~9.
 More precisely, it turns out that the latter category has to be
``finalized'' in order to make it a model category.
 We also discuss DG\+modules over cofibrant DG\+algebras.
 Conilpotent curved $\Ainfty$\+coalgebras and co/contranilpotent
curved $\Ainfty$\+co/contramodules over them are introduced.

 Homogeneous Koszul duality is worked out in Appendix~A\hbox{}.
 The (more general) covariant and the (more symmetric) contravariant
versions of the duality are considered separately.
 The equivalence between the derived category of modules over
the ring/sheaf of differential operations acting in a vector bundle
and the coderived category of CDG\+modules over the corresponding
de Rham CDG\+algebra is constructed in Appendix~B\hbox{}.
 A desription of the bounded derived category of coherent
$\D$\+modules in terms of coherent CDG\+modules is also obtained. 

\subsection*{Acknowledgement}
 The author is grateful to Michael Finkelberg for posing the problem
of constructing derived nonhomogeneous Koszul duality.
 I want to express my gratitude to Vladimir Voevodsky for very
stimulating discussions and encouragement, without which this work
would probably never have been done.
 I also benefited from discussions with Joseph Bernstein, Victor
Ginzburg, Amnon Neeman, Alexander Beilinson, Henning Krause,
Maxim Kontsevich, Tony Pantev, Alexander Polishchuk,
Alexander Kuznetsov, Lars W.~Christensen, Jan \v S\v tov\'\i\v cek,
Pedro Nicol\'as, and Alexander Efimov.
 I am grateful to Ivan Mirkovic, who always urged me to write down
the material presented below.
 I want to thank Dmitry Arinkin, who communicated the proof of
Theorem~\ref{fin-gen-cdg-mod}.2 to me and gave me the permission to
include it in this paper.
 Most of the content of this paper was worked out when I was a Member
of the Institute for Advanced Study, which I wish to thank for its
hospitality.
 I am also indebted to the participants of an informal seminar at IAS,
where I first presented these results in the Spring of 1999.
 The author was partially supported by P.~Deligne 2004 Balzan prize,
an INTAS grant, and an RFBR grant while writing the paper up.
 Parts of this paper have been written when I was visiting
the Institut des Hautes \'Etudes Scientifiques, which I wish to
thank for the excellent working conditions.

\Section{Derived Category of DG-Modules}

\subsection{DG-rings and DG-modules}  \label{dg-rings-modules}
 A \emph{DG\+ring} $A=(A,d)$ is a pair consisting of an associative
graded ring $A=\bigoplus_{i\in\boZ} A^i$ and an odd derivation
$d\:A\rarrow A$ of degree~$1$ such that $d^2=0$.
 In other words, it is supposed that $d(A^i)\subset A^{i+1}$ and
$d(ab)=d(a)b+(-1)^{|a|}ad(b)$ for $a$, $b\in A$, where $|a|$
denotes the degree of a homogeneous element, i.~e., $a\in A^{|a|}$.

 A \emph{left DG\+module} $(M,d_M)$ over a DG\+ring $A$ is a graded
left $A$\+module $M=\bigoplus_{i\in\boZ}M^i$ endowed with
a differential $d_M\:M\rarrow M$ of degree~$1$ compatible with
the derivation of~$A$ and such that $d_M^2=0$.
 The compatibility means that the equation $d_M(ax) = d(a)x +
(-1)^{|a|}ad_M(x)$ holds for all $a\in A$ and $x\in M$.

 A \emph{right DG\+module} $(N,d_N)$ over $A$ is a graded right
$A$\+module $N$ endowed with a differential $d_N$ of degree~$1$
satisfying the equations $d_N(xa)=d_N(x)a+(-1)^{|x|}xd(a)$ and
$d_N^2=0$, where $x\in N^{|x|}$.

 Let $L$ and $M$ be left DG\+modules over~$A$.
 The \emph{complex of homomorphisms} $\Hom_A(L,M)$ from $L$ to $M$
over~$A$ is constructed as follows.
 The component $\Hom_A^i(L,M)$ consists of all homogeneous maps
$f\:L\rarrow M$ of degree~$i$ such that $f(ax)=(-)^{i|a|}af(x)$
for all $a\in A$ and $x\in L$.
 The differential in the complex $\Hom_A(L,M)$ is given by
the formula $d(f)(x)=d_M(f(x))-(-1)^{|f|}f(d_L(x))$.
 Clearly, one has $d^2(f)=0$; for any composable morphisms of
left DG\+modules $f$ and $g$ one has $d(fg)=d(f)g+(-1)^{|f|}fd(g)$.

 For any two right DG\+modules $R$ and $N$ over $A$, the complex
of homomorphisms $\Hom_A(R,N)$ is defined by the same formulas as
above and satisfies the same properties, with the only difference
that a homogeneous map $f\:R\rarrow N$ belonging to $\Hom_A(R,N)$
must satisfy the equation $f(xa)=f(x)a$ for $a\in A$ and $x\in R$.

 Let $N$ be a right DG\+module and $M$ be a left DG\+module over~$A$.
 The \emph{tensor product complex} $N\ot_A M$ is defined as
the graded quotient group of the graded abelian group $N\ot_\boZ M$
by the relations $xa\ot y =x\ot ay$ for $x\in N$, \ $a\in A$, \ 
$y\in M$, endowed with the differential given by the formula
$d(x\ot y)= d(x)\ot y + (-1)^{|x|} x\ot d(y)$.
 For any two right DG\+modules $R$ and $N$ and any two left
DG\+modules $L$ and $M$ the natural map of complexes $\Hom_A(R,N)
\ot_\boZ\Hom_A(L,M)\rarrow\Hom_\boZ(R\ot_AL\;N\ot_AM)$ is defined
by the formula $(f\ot g)(x\ot y)=(-1)^{|g||x|}f(x)\ot g(y)$.
 Here $\boZ$ is considered as a DG\+ring concentrated in degree~$0$.
 
 For any DG\+ring $A$, its cohomology $H(A)=H_d(A)$, defined as
the quotient of the kernel of~$d$ by its image, has a natural
structure of graded ring.
 For a left DG\+module $M$ over~$A$, its cohomology $H(M)$ is
a graded module over $H(A)$; for a right DG\+module $N$, its
cohomology $H(N)$ is a right graded module over $H(A)$.

 A \emph{DG\+algebra} $A$ over a commutative ring~$k$ is
a DG\+ring endowed with DG\+ring homomorphism $k\rarrow A^0$
whose image is contained in the center of the algebra $A$,
where $k$ is considered as a DG\+ring concentrated in degree~$0$;
equivalently, a DG\+algebra is a complex of $k$\+modules with
a $k$\+linear DG\+ring structure.

\begin{rem}
 One can consider DG\+algebras and DG\+modules graded by an abelian
group~$\Gamma$ different from $\boZ$, provided that $\Gamma$ is
endowed with a parity homomorphism $\Gamma\rarrow\boZ/2$ and
an odd element $\boldsymbol{1}\in\Gamma$, so that the differentials
would have degree~$\boldsymbol{1}$.
 In particular, one can take $\Gamma=\boZ/2$, that is have gradings
reduced to parities, or consider fractional gradings by using
some subgroup of $\boQ$ consisting of rationals with odd
denominators in the role of~$\Gamma$.
 Even more generally, one can replace the parity function with 
a symmetric bilinear form $\sigma\:\Gamma\times\Gamma\rarrow\boZ/2$,
to be used in the super sign rule in place of the product of parities;
one just has to assume that $\sigma(\boldsymbol{1},\boldsymbol{1})
= 1\bmod 2$.
 All the most important results of this paper remain valid in
such settings.
 The only exceptions are the results of
subsections \ref{dg-ring-bounded-cases}
and~\ref{dg-coalgebra-bounded-cases}, where we consider
bounded grading.
\end{rem}

\subsection{DG\+categories}  \label{dg-categories}
 A \emph{DG\+category} is a category whose sets of morphisms are
complexes and compositions are biadditive maps compatible with
the gradings and the differentials.
 In other words, a DG\+category $\sDG$ consists of a class of
objects, complexes of abelian groups $\Hom_{\sDG}(X,Y)$, called
the complexes of morphisms from $X$ to~$Y$, defined for any two
objects $X$ and $Y$, and morphisms of complexes $\Hom_{\sDG}(Y,Z)
\ot_{\boZ}\Hom_{\sDG}(X,Y)\rarrow\Hom_{\sDG}(X,Z)$, called
the composition maps, defined for any three objects $X$, $Y$,
and~$Z$.
 The compositions must be associative and unit elements $\id_X\in
\Hom_{\sDG}(X,X)$ must exist; the equations $d(\id_X)=0$ then hold
automatically.

 For example, left DG\+modules over a DG\+ring $A$ form
a DG\+category, which we will denote by $\sDG(A\modl)$.
 The DG\+category of right DG\+modules over $A$ will be
denoted by $\sDG(\modr A)$.

 A \emph{covariant DG\+functor} $\sDG'\rarrow\sDG''$ consists of
a map between the classes of objects and (closed) morphisms between
the complexes of morphisms compatible with the compositions.
 A \emph{contravariant DG\+functor} is defined in the same way,
except that one has to take into account the natural isomorphism
of complexes $V\ot W\simeq W\ot V$ for complexes of abelian
groups $V$ and $W$ that is given by the formula $v\ot w\mpsto
(-1)^{|v||w|}w\ot v$.
 (Covariant or contravariant) DG\+functors between $\sDG'$ and
$\sDG''$ form a DG\+category themselves.
 The complex of morphisms between DG\+functors $F$ and $G$ is
a subcomplex of the product of the complexes of morphisms from
$F(X)$ to $G(X)$ in $\sDG''$ taken over all objects $X\in\sDG'$;
the desired subcomplex is formed by all the systems of morphisms
compatible with all morphisms $X\rarrow Y$ in $\sDG'$.

 For example, a DG\+ring $A$ can be considered as a DG\+category
with a single object; covariant DG\+functors from this
DG\+category to the DG\+category of complexes of abelian groups
are left DG\+modules over~$A$, while contravariant DG\+functors
between the same DG\+categories can be identified with right
DG\+modules over~$A$.

 A \emph{closed morphism} $f\:X\rarrow Y$ in a DG\+category
$\sDG$ is an element of $\Hom^0_{\sDG}(X,Y)$ such that $d(f)=0$.
 The category whose objects are the objects of $\sDG$ and whose
morphisms are closed morphisms in $\sDG$ is denoted by $Z^0(\sDG)$.

 An object $Y$ is called the \emph{product} of a family of objects
$X_\alpha$ (notation: $Y=\prod_\alpha X_\alpha$) if a closed
isomorphism of contravariant DG\+functors $\Hom_{\sDG}({-},Y)\simeq
\prod_\alpha\Hom_{\sDG}({-},X_\alpha)$ is fixed.
 An object $Y$ is called the \emph{direct sum} of a family of
objects $X_\alpha$ (notation: $Y=\bigoplus_\alpha X_\alpha$) if
a closed isomorphism of covariant DG\+functors $\Hom_{\sDG}(Y,{-})
\simeq\prod_\alpha\Hom_{\sDG}(X_\alpha,{-})$ is fixed.

 An object $Y$ is called the \emph{shift} of an object $X$
by an integer~$i$ (notation: $Y=X[i]$) if a closed isomorphism of
contravariant DG\+functors $\Hom_{\sDG}({-},Y)\simeq
\Hom_{\sDG}({-},X)[i]$ is fixed, or equivalently, a closed
isomorphism of covariant DG\+functors $\Hom_{\sDG}\allowbreak(Y,{-})
\simeq\Hom_{\sDG}(X,{-})[-i]$ is fixed.

 An object $Z$ is called the \emph{cone} of a closed morphism
$f\:X\rarrow Y$ (notation: $Z=\cone(f)$) if a closed isomorphism of
contravariant DG\+functors $\Hom_{\sDG}({-},Z)\simeq\cone(f_*)$,
where $f_*\:\Hom_{\sDG}({-},X)\rarrow\Hom_{\sDG}({-},Y)$,
is fixed, or equivalently, a closed isomorphism of covariant
DG\+functors $\Hom_{\sDG}(Z,{-})\simeq\cone(f^*)[-1]$, where
$f^*\:\Hom_{\sDG}(Y,{-})\rarrow\Hom_{\sDG}(X,{-})$, is fixed.

 Let $V$ be a complex of abelian groups and $p\:V\rarrow V$ be
an endomorphism of degree~$1$ satisfying the Maurer--Cartan equation
$d(p)+p^2=0$.
 Then one can define a new differential on $V$ by setting $d'=d+p$;
let us denote the complex so obtained by $V(p)$.
 Let $q\in\Hom_{\sDG}^1(X,X)$ be an endomorphism of degree~$1$
satisfying the equation $d(q)+q^2=0$.
 An object $Y$ is called the \emph{twist} of the object $X$ with
respect to~$q$ if a closed isomorphism of contravariant DG\+functors
$\Hom_{\sDG}({-},X)\simeq\Hom_{\sDG}({-},Y)(q_*)$ is fixed,
where $q_*(g)=q\circ g$ for any morphism~$g$ whose target is~$X$,
or equivalently, a closed isomorphism of covariant DG\+functors
$\Hom_{\sDG}(Y,{-})\simeq \Hom_{\sDG}(X,{-})(-q^*)$ is fixed,
where $q^*(g)=(-1)^{|g|}g\circ q$ for any morphism~$g$ whose source
is~$Y$.

 As any representing objects of DG\+functors, all direct sums,
products, shifts, cones, and twists are defined uniquely up to
a unique closed isomorphism.
 The direct sum of a finite set of objects is naturally also
their product, and vice versa.
 Finite direct sums, products, shifts, cones, and twists are
preserved by any DG\+functors.
 One can express the cone of a closed morphism $f\:X\rarrow Y$ as
the twist of the direct sum $Y\oplus X[1]$ with respect to
the endomorphism~$q$ induced by~$f$.

 Here is another way to think about cones of closed morphisms
in DG\+categories.
 Let $\sDG^\#$ denote the category whose objects are the objects
of $\sDG$ and morphisms are the (not necessarily closed)
morphisms in $\sDG$ of degree~$0$.
 Let $X'\rarrow X\rarrow X''$ be a triple of objects in $\sDG$
with closed morphisms between them that is split exact in
$\sDG^\#$.
 Then $X$ is the cone of a closed morphism $X''[-1]\rarrow X'$.
 Conversely, for any closed morphism $X\rarrow Y$ in $\sDG$
with the cone $Z$ there is a natural triple of objects and closed
morphisms $Y\rarrow Z\rarrow X[1]$, which is split exact in $\sDG^\#$.

 Let $\sDG$ be a DG\+category with shifts, twists, and infinite
direct sums.
 Let $\dsb\rarrow X_n\rarrow X_{n-1}\rarrow \dsb$ be a complex of
objects of $\sDG$ with closed differentials~$\d_n$.
 Then the differentials~$\d_n$ induce an endomorphism~$q$ of
degree~$1$ on the direct sum $\bigoplus_n X_n[n]$ satisfying
the equations $d(q)=0=q^2$.
 The twist of this direct sum with respect to this endomorphism
is called the \emph{total object} of the complex $X_\bu$
\emph{formed by taking infinite direct sums} and denoted by
$\Tot^{\oplus}(X_\bu)$.
 For a DG\+category $\sDG$ with shifts, twists, and infinite
products, one can consider the analogous construction with
the infinite direct sum replaced by the infinite product
$\prod_n X_n[n]$.
 Thus one obtains the definition of the \emph{total object formed
by taking infinite products} $\Tot^{\sqcap}(X_\bu)$.

 For a finite complex $X^\bu$, the two total objects coincide
and are denoted simply by $\Tot(X_\bu)$; this total object only
requires existence of finite direct sums/products for its
construction.
 Alternatively, the total objects $\Tot$, $\Tot^{\oplus}$, and
$\Tot^{\sqcap}$ can be defined as certain representing objects of
DG\+functors.
 The finite total object $\Tot$ can be also expressed in terms of
iterated cones, so it is well-defined whenever cones exist in
a DG\+category $\sDG$, and it is preserved by any DG\+functors.
 
 A DG\+functor $\sDG'\rarrow\sDG''$ is said to be \emph{fully
faithful} if it induces isomorphisms of the complexes of morphisms.
 A DG\+functor is said to be an \emph{equivalence of DG\+categories}
if it is fully faithful and every object of $\sDG''$ admits a closed
isomorphism with an object coming from $\sDG'$.
 This is equivalent to existence of a DG\+functor in the opposite
direction for which both the compositions admit closed isomorphisms
to the identity DG\+functors.
 DG\+functors $F\:\sDG'\rarrow\sDG''$ and $G\:\sDG''\rarrow\sDG'$
are said to be \emph{adjoint} if for every objects $X\in\sDG'$ and
$Y\in\sDG''$ a closed isomorphism of complexes $\Hom_{\sDG''}(F(X),Y)
\simeq\Hom_{\sDG'}(X,G(Y))$ is given such that these isomorphisms
commute with the (not necessarily closed) morphisms induced by
morphisms in $\sDG'$ and $\sDG''$.

 Let $\sDG$ be a DG\+category where (a zero object and) all
shifts and cones exist.
 Then the \emph{homotopy category} $H^0(\sDG)$ is the additive
category with the same class of objects as $\sDG$ and groups of
morphisms given by $\Hom_{H^0(\sDG)}(X,Y)=H^0(\Hom_\sDG(X,Y))$.
 The homotopy category is a triangulated category~\cite{BK}.
 Shifts of objects and cones of closed morphisms in $\sDG$ become
shifts of objects and cones of morphisms in the triangulated
category $H^0(\sDG)$.
 Any direct sums and products of objects of a DG\+category are also
their directs sums and products in the homotopy category.
 Adjoint functors between DG\+categories induce adjoint functors
between the corresponding categories of closed morphisms and
homotopy categories.

 Two closed morphisms $f$, $g\:X\rarrow Y$ in a DG\+category $\sDG$
are called \emph{homotopic} if their images coincide in $H^0(\sDG)$.
 A closed morphism in $\sDG$ is called a \emph{homotopy equivalence}
if it becomes an isomorphism in $H^0(\sDG)$.
 An object of $\sDG$ is called \emph{contractible} if it vanishes
in $H^0(\sDG)$.

 All shifts, twists, infinite direct sums, and infinite direct
products exist in the DG\+categories of DG\+modules.
 The homotopy category of (the DG\+category of) left DG\+modules
over a DG\+ring $A$ is denoted by $\Hot(A\modl)=H^0\sDG(A\modl)$;
the homotopy category of right DG\+modules over~$A$ is denoted by
$\Hot(\modr A)=H^0\sDG(\modr A)$. 

\subsection{Semiorthogonal decompositions}  \label{semiorthogonal}
 Let $\sH$ be a triangulated category and $\sA\subset\sH$ be
a full triangulated subcategory.
 Then the quotient category $\sH/\sA$ is defined as
the localization of $\sH$ with respect to the multiplicative
system of morphisms whose cones belong to~$\sA$.
 The subcategory $\sA$ is called \emph{thick} if it coincides
with the full subcategory formed by all the objects of $\sH$ whose
images in $\sH/\sA$ vanish.
 A triangulated subcategory $\sA\subset\sH$ is thick if and
only if it is closed under direct summands in~$\sH$~\cite{Ver2,Neem}.
 The following Lemma is essentially due to Verdier~\cite{Ver1};
see also~\cite{BBD,Bon}.

\begin{lem}
 Let\/ $\sH$ be a triangulated category and\/ $\sB$, $\sC\subset\sH$
be its full triangulated subcategories such that\/
$\Hom_{\sH}(B,C)=0$ for all $B\in\sB$ and $C\in\sC$.
 Then the natural maps\/ $\Hom_{\sH}(B,X)\rarrow\Hom_{\sH/\sC}(B,X)$
and\/ $\Hom_{\sH}(X,C)\rarrow\Hom_{\sH/\sB}(X,C)$ are isomorphisms
for any objects $B\in\sB$, \ $C\in\sC$, and $X\in\sH$.
 In particular, the functors\/ $\sB\rarrow\sH/\sC$ and\/ $\sC\rarrow
\sH/\sB$ are fully faithful.
 Furthermore, the following conditions are equivalent:
\begin{itemize}
 \item[(a)] $\sB$ is a thick subcategory in\/ $\sH$ and
the functor\/ $\sC\rarrow\sH/\sB$ is an equivalence of triangulated
categories;
 \item[(b)] $\sC$ is a thick subcategory in\/ $\sH$ and
the functor\/ $\sB\rarrow\sH/\sC$ is an equivalence of triangulated
categories;
 \item[(c)] $\sB$ and\/ $\sC$ generate\/ $\sH$ as a triangulated
category, i.~e., any object of\/ $\sH$ can be obtained from objects
of\/ $\sB$ and\/ $\sC$ by iterating the operations of shift and cone;
 \item[(d)] for any object $X\in\sH$ there exists a distinguished
triangle $B\rarrow X\rarrow C\rarrow B[1]$ with $B\in\sB$ and
$C\in\sC$ (and in this case for any morphism $X'\rarrow X''$ in\/ $\sH$
there exists a unique morphism between any distinguished triangles
of the above form for $X'$ and $X''$, so this triangle is unique up
to a unique isomorphism and depends functorially on~$X$);
 \item[(e)] $\sC$ is the full subcategory of\/~$\sH$ formed by
all the objects $C\in\sH$ such that\/ $\Hom_{\sH}(B,C)=0$ for all
$B\in\sB$, and the embedding functor\/ $\sB\rarrow\sH$ has
a right adjoint functor (which can be then identified with
the localization functor $\sH\rarrow\sH/\sC\simeq\sB$);
 \item[(f)] $\sC$ is the full subcategory of\/~$\sH$ formed by
all the objects $C\in\sH$ such that\/ $\Hom_{\sH}(B,C)=0$ for all
$B\in\sB$, \ $\sB$ is a thick subcategory in\/~$\sH$, and
the localization functor\/ $\sH\rarrow\sH/\sB$ has a right adjoint
functor;
 \item[(g)] $\sB$ is the full subcategory of\/~$\sH$ formed by
all the objects $B\in\sH$ such that\/ $\Hom_{\sH}(B,C)=0$ for all
$C\in\sC$, and the embedding functor\/ $\sC\rarrow\sH$ has
a left adjoint functor (which can be then identified with
the localization functor $\sH\rarrow\sH/\sB\simeq\sC$);
 \item[(h)] $\sB$ is the full subcategory of\/~$\sH$ formed by
all the objects $B\in\sH$ such that\/ $\Hom_{\sH}(B,C)=0$ for all
$C\in\sC$, \ $\sC$ is a thick subcategory in\/~$\sH$, and
the localization functor\/ $\sH\rarrow\sH/\sC$ has a left adjoint
functor.
\qed
\end{itemize}
\end{lem}

\subsection{Projective resolutions}  \label{projective-dg-modules}
 A DG\+module $M$ is said to be acyclic if it is acyclic as
a complex of abelian groups, i.~e., $H(M)=0$.
 The thick subcategory of the homotopy category
$\Hot(A\modl)$ formed by the acyclic DG\+modules is denoted
by $\Acycl(A\modl)$.
 The \emph{derived category} of left DG\+modules over~$A$ is
defined as the quotient category $\sD(A\modl)=\Hot(A\modl)/
\Acycl(A\modl)$.

 A left DG\+module $L$ over a DG\+ring $A$ is called
\emph{projective} if for any acyclic left DG\+module $M$
over~$A$ the complex $\Hom_A(L,M)$ is acyclic.
 The full triangulated subcategory of $\Hot(A\modl)$ formed by
the projective DG\+modules is denoted by $\Hot(A\modl)_\proj$.
 The following Theorem says, in particular, that the homotopy category
$\sH=\Hot(A\modl)$ and its subcategories $\sB=\Hot(A\modl)_\proj$
and $\sC=\Acycl(A\modl)$ satisfy the equivalent conditions of
Lemma~\ref{semiorthogonal}, and so describes the derived category
$\sD(A\modl)$.

\begin{thm}
\textup{(a)} The category\/ $\Hot(A\modl)_\proj$ is the minimal
triangulated subcategory of\/ $\Hot(A\modl)$ containing the
DG\+module $A$ and closed under infinite direct sums. \par
\textup{(b)} The composition of functors\/
$\Hot(A\modl)_\proj\rarrow\Hot(A\modl)\rarrow\sD(A\modl)$
is an equivalence of triangulated categories.
\end{thm}

\begin{proof}
 First notice that the category $\Hot(A\modl)_\proj$ is closed
under infinite direct sums.
 It contains the DG\+module $A$, since for any DG\+module $M$
over~$A$ there is a natural isomorphism of complexes of abelian
groups $\Hom_A(A,M)\simeq M$.
 According to Lemma~\ref{semiorthogonal}, it remains to construct
for any DG\+module $M$ a morphism $f\:F\rarrow M$ in the homotopy
category of DG\+modules over~$A$ such that the DG\+module $F$
belongs to the minimal triangulated subcategory containing
the DG\+module $A$ and closed under infinite direct sums, while
the cone of the morphism~$f$ is an acyclic DG\+module.
 When $A$ is a DG\+algebra over a field~$k$, it suffices to consider
the bar-resolution of a DG\+module $M$.
 It is a complex of DG\+modules over~$A$, and its total DG\+module
formed by taking infinite direct sums provides the desired
DG\+module~$F$.

 Let us give a detailed construction in the general case.
 Let $M$ be a DG\+module over~$A$.
 Choose a complex of free abelian groups $M'$ together with
a surjective morphism of complexes $M'\rarrow M$ such that
the cohomology $H(M')$ is also a free graded abelian group
and the induced morphism of cohomology $H(M')\rarrow H(M)$
is also surjective.
 For example, one can take $M'$ to be the graded abelian group
with the components freely generated by nonzero elements of
the components of $M$, endowed with the induced differential.
 Set $F_0=A\ot_{\boZ}M'$; then there is a natural closed
surjective morphism $F_0\rarrow M$ of DG\+modules over~$A$ and
the induced morphism of cohomology $H(F_0)\rarrow H(M)$ is
also surjective.
 Let $K$ be the kernel of the morphism $F_0\rarrow M$ (taken in
the abelian category $Z^0\sDG(A\modl)$ of DG\+modules and closed
morphisms between them).
 Applying the same construction to the DG\+module $K$ in
place of $M$, we obtain the DG\+module $F_1$, etc.
 Let $F$ be the total DG\+module of the complex of DG\+modules
$\dsb \rarrow F_1\rarrow F_0$ formed by taking infinite
direct sums.
 One can easily check that the cone of the morphism $F\rarrow M$
is acyclic, since the complex $\dsb\rarrow H(F_1)\rarrow H(F_0)
\rarrow H(M)\rarrow0$ is acyclic (it suffices to apply
the result of~\cite{EM1} to the increasing filtration of the total
complex of $\dsb\rarrow F_1\rarrow F_0\rarrow M$ coming from
the silly filtration of this complex of complexes).

 It remains to show that the DG\+module $F$ as an object of
the homotopy category can be obtained from the DG\+module $A$
by iterating the operations of shift, cone, and infinite
direct sum.
 Every DG\+module $F_n$ is a direct sum of shifts of
the DG\+module $A$ and shifts of the cone of the identity
endomorphism of the DG\+module~$A$.
 Denote by $X_n$ the total DG\+module of the finite complex
of DG\+modules $F_n\rarrow\dsb\rarrow F_0$.
 Then we have $F=\varinjlim X_n$ in the abelian category
$Z^0\sDG(A\modl)$.
 So there is an exact triple of DG\+modules and closed morphisms
$0\rarrow\bigoplus X_n\rarrow \bigoplus X_n\rarrow F\rarrow0$.
 Since the embeddings $X_n\rarrow X_{n+1}$ split in $\sDG(A\modl)^\#$,
the above exact triple also splits in this additive category.
 Thus $F$ is a cone of the morphism $\bigoplus X_n\rarrow
\bigoplus X_n$ in the triangulated category $\Hot(A\modl)$.
\end{proof}

\subsection{Injective resolutions}  \label{injective-dg-modules}
 A left DG\+module $M$ over a DG\+ring $A$ is said to be
\emph{injective} if for any acyclic DG\+module $L$ over~$A$
the complex $\Hom_A(L,M)$ is acyclic.
 The full triangulated subcategory of $\Hot(A\modl)$ formed by
the injective DG\+modules is denoted by $\Hot(A\modl)_\inj$.

 For any right DG\+module $N$ over $A$ and any complex of abelian
groups $V$ the complex $\Hom_\boZ(N,V)$ has a natural structure
of left DG\+module over $A$ with the graded $A$\+module structure
given by the formula $(af)(n)=(-1)^{|a|(|f|+|n|)}f(na)$.

 The following Theorem provides another semiorthogonal decomposition
of the homotopy category $\Hot(A\modl)$ and another description of
the derived category $\sD(A\modl)$.

\begin{thm}
\textup{(a)} The category\/ $\Hot(A\modl)_\inj$ is the minimal
triangulated subcategory of\/ $\Hot(A\modl)$ containing the
DG\+module\/ $\Hom_\boZ(A,\boQ/\boZ)$ and closed under infinite
products. \par
\textup{(b)} The composition of functors\/
$\Hot(A\modl)_\inj\rarrow\Hot(A\modl)\rarrow\sD(A\modl)$
is an equivalence of triangulated categories.
\end{thm}

\begin{proof}
 The proof is analogous to that of Theorem~\ref{projective-dg-modules}.
 Clearly, the category $\Hot(A\modl)_\inj$ is closed under infinite
products.
 It contains the DG\+module $\Hom_\boZ(A,\boQ/\boZ)$, since
the complex $\Hom_A(L,\Hom_\boZ(A,\boQ/\boZ))\simeq
\Hom_\boZ(L,\boQ/\boZ)$ is acyclic whenever the DG\+module $L$~is.
 To construct an injective resolution of a DG\+module $M$, one
can embed in into a complex of injective abelian groups $M'$
so that the cohomology $H(M')$ is also injective and $H(M)$ also
embeds into $H(M')$.
 For example, one can take the components of $M'$ to be
the products of $\boQ/\boZ$ over all nonzero homomorphisms of
abelian groups from the components of $M$ to $\boQ/\boZ$.
 Take $J_0=\Hom_\boZ(A,M')$ and consider the induced injective
morphism of DG\+modules $M\rarrow J_0$.
 Set $K=J_0/M$, \ $J_{-1}=\Hom_\boZ(A,K')$, etc., and
$J=\Tot^{\sqcap}(J_\bu)$.
 Then the morphism of DG\+modules $M\rarrow J$ has an acyclic cone
and the DG\+module $J$ is isomorphic in $\Hot(A\modl)$ to a DG\+module
obtained from $\Hom_\boZ(A,\boQ/\boZ)$ by iterating the operations
of shift, cone, and infinite product.
\end{proof}

\subsection{Flat resolutions}  \label{flat-dg-resolutions}
 A right DG\+module $N$ over a DG\+ring $A$ is said to be \emph{flat}
if for any acyclic left DG\+module $M$ over~$A$ the complex
$N\ot_AM$ is acyclic.
 Flat left DG\+modules over~$A$ are defined in the analogous way.
 The full triangulated subcategory of $\Hot(A\modl)$ formed by
flat DG\+modules is denoted by $\Hot(A\modl)_\fl$.

 We denote the thick subcategory of acyclic right $A$\+modules by
$\Acycl(\modr A)\subset\Hot(\modr A)$.
 The quotient category $\Hot(\modr A)/\Acycl(\modr A)$ is called
the derived category of right DG\+modules over~$A$ and denoted by
$\sD(\modr A)$.
 The full triangulated subcategory of flat right DG\+modules is
denoted by $\Hot(\modr A)_\fl\subset\Hot(\modr A)$.

 It follows from Theorems~\ref{projective-dg-modules}--%
\ref{injective-dg-modules} and Lemma~\ref{semiorthogonal} that
one can compute the right derived functor $\Ext_A(L,M)=
\Hom_{\sD(A\modl)}(L,M)$ for left DG\+modules $L$ and $M$ over
a DG\+ring $A$ in terms of projective or injective resolutions.
 Namely, one has $\Ext_A(L,M)\simeq H(\Hom_A(L,M))$ whenever $L$ is
a projective DG\+module or $M$ is an injective DG\+module over~$A$.
 The following Theorem allows to define a left derived functor
$\Tor^A(N,M)$ for a right DG\+module $N$ and a left DG\+module $M$
over~$A$ so that it could be computed in terms of flat resolutions.

\begin{thm}
\textup{(a)} The functor\/ $\Hot(A\modl)_\fl/(\Acycl(A\modl)\cap
\Hot(A\modl)_\fl)\rarrow\sD(A\modl)$ induced by the embedding\/
$\Hot(A\modl)_\fl\rarrow\Hot(A\modl)$ is an equivalence of
triangulated categories. \par
\textup{(b)} The functor\/ $\Hot(\modr A)_\fl/(\Acycl(\modr A)\cap
\Hot(\modr A)_\fl)\rarrow\sD(\modr A)$ induced by the embedding\/
$\Hot(\modr A)_\fl\rarrow\Hot(\modr A)$ is an equivalence of
triangulated categories.
\end{thm}

 The proof of Theorem is based on the following Lemma.

\begin{lem}
 Let\/ $\sH$ be a triangulated category and\/ $\sA$, $\sF\subset\sH$
be full triangulated subcategories.
 Then the natural functor\/ $\sF/\sA\cap\sF\rarrow\sH/\sA$ is
an equivalence of triangulated categories whenever one of
the following two conditions holds:
\begin{itemize}
\item[(a)] for any object $X\in\sH$ there exists an object
$F\in\sF$ together with a morphism $F\rarrow X$ in\/~$\sH$ such
that a cone of that morphism belongs to\/~$\sA$, or
\item[(b)] for any object $Y\in\sH$ there exists an object
$F\in\sF$ together with a morphism $Y\rarrow F$ in\/~$\sH$ such
that a cone of that morphism belongs to\/~$\sA$.
\end{itemize}
\end{lem}

\begin{proof}[Proof of Lemma]
 It is clear that the functor $\sF/\sA\cap\sF\rarrow\sH/\sA$
is surjective on the isomorphism classes of objects under either
of the assumptions (a) or~(b).
 To prove that it is bijective on morphisms, represent morphisms
in both quotient categories by fractions of the form
$X\larrow X'\rarrow Y$ in the case~(a) and by fractions of the form
$X\rarrow Y'\larrow Y$ in the case~(b).
\end{proof}

\begin{proof}[Proof of Theorem]
 Part~(a): first notice that any projective left DG\+module $M$ over
a DG\+ring~$A$ is flat.
 Indeed, one has $\Hom_\boZ(N\ot_A M\;\boQ/\boZ)\simeq
\Hom_A(M,\Hom_\boZ(N,\boQ/\boZ))$ for any right DG\+module $N$
over~$A$, so whenever $N$ is acyclic, and consequently
$\Hom_\boZ(N,\boQ/\boZ)$ is acyclic, the left hand side of this
isomorphism is acyclic, too, and therefore $N\ot_A M$ is acyclic.
 So it remains to use Theorem~\ref{projective-dg-modules} together
with Lemma~\ref{semiorthogonal} and the above Lemma. 
 To prove part~(b), switch the left and right sides by passing to
the DG\+ring $A^\rop$ defined as follows.
 As a complex, $A^\rop$ is identified with $A$, while
the multiplication in $A^\rop$ is given by the formula
$a^\rop b^\rop = (-1)^{|a||b|}(ba)^\rop$.
 Then right DG\+modules over $A$ are left DG\+modules over $A^\rop$
and vice versa.
\end{proof}

 Now let us define the derived functor
$$
 \Tor^A\:\sD(\modr A)\times\sD(A\modl)\lrarrow k\modl^\sgr
$$
for a DG\+algebra $A$ over a commutative ring~$k$, where
$k\modl^\sgr$ denotes the category of graded $k$\+modules.
 For this purpose, restrict the functor of tensor product
$\ot_A\:\Hot(\modr A)\times\Hot(A\modl)\rarrow\Hot(k\modl)$ to either
of the full subcategories $\Hot(\modr A)_\fl\times \Hot(A\modl)$ or
$\Hot(\modr A)\times\Hot(A\modl)_\fl$ and compose it with
the cohomology functor $H\:\Hot(k\modl)\rarrow k\modl^\sgr$.
 The functors so obtained factorize through the localizations
$\sD(\modr A)\times\sD(A\modl)$ and the two induced derived functors
$\sD(\modr A)\times\sD(A\modl)\rarrow k\modl^\sgr$ are naturally
isomorphic to each other.

 Indeed, the tensor product $N\ot_A M$ by the definition is acyclic
whenever one of the DG\+modules $N$ and $M$ is acyclic, while
the other one is flat.
 Let us check that the complex $N\ot_A M$ is acyclic whenever
either of the DG\+modules $N$ and $M$ is simultaneously acyclic
and flat.
 Assume that $N$ is acyclic and flat; choose a flat left DG\+module
$F$ over~$A$ together with a morphism of DG\+modules $F\rarrow A$
with an acyclic cone.
 Then the complex $N\ot_A F$ is acyclic, since $N$ is acyclic;
while the morphism $N\ot_A F\rarrow N\ot_A M$ is a quasi-isomorphism,
since $N$ is flat.

 To construct an isomorphism of the two induced derived functors,
it suffices to notice that both of them are isomorpic to the derived
functor obtained by restricting the functor $\ot_A$ to the full
subcategory $\Hot(\modr A)_\fl\times\Hot(A\modl)_\fl$.
 In other words, suppose that $G\rarrow N$ and $F\rarrow M$ are
morphisms of DG\+modules with acyclic cones, where the right
DG\+module $G$ and the left DG\+module $F$ are flat.
 Then there are natural quasi-isomorphisms $G\ot_A M\larrow
G\ot_A F\rarrow N\ot_A F$.
 
\subsection{Restriction and extention of scalars}
\label{dg-mod-scalars}
 Let $f\:A\rarrow B$ be a morphism of DG\+algebras, i.~e.,
a closed morphism of complexes preserving the multiplication.
 Then any DG\+module over $B$ can be also considered as
a DG\+module over $A$, which defines the restriction-of-scalars
functor $R_f\:\Hot(B\modl)\rarrow\Hot(A\modl)$.
 This functor has a left adjoint functor $E_f$ given by the formula
$E_f(M)=B\ot_AM$ and a right adjoint functor $E^f$ given by
the formula $E^f(M)=\Hom_A(B,M)$ (where the DG\+module structure on
$\Hom_A(B,M)$ is defined so that $\Hom_A(B,M)\rarrow\Hom_\boZ(B,M)$
is a closed injective morphism of DG\+modules).

 The functor $R_f$ obviously maps acyclic DG\+modules to acyclic
DG\+modules, and so induces a functor $\sD(B\modl)\rarrow
\sD(A\modl)$, which we will denote by~$\boI R_f$.
 The functor $E_f$ has a left derived functor $\boL E_f$ obtained
by restricting~$E_f$ to either of the full subcategories
$\Hot(A\modl)_\proj$ or $\Hot(A\modl)_\fl\subset\Hot(A\modl)$
and composing it with the localization functor $\Hot(B\modl)\rarrow
\sD(B\modl)$.
 The functor $E^f$ has a right derived functor $\boR E^f$ obtained
by restricting $E^f$ to the full subcategory $\Hot(A\modl)_\inj
\subset\Hot(A\modl)$ and composing it with the localization functor
$\Hot(B\modl)\rarrow\sD(B\modl)$.
 The functor $\boL E_f$ is left adjoint to the functor~$\boI R_f$
and the functor $\boR E^f$ is right adjoint to the functor~$\boI R_f$.

\begin{thm}
 The functors $\boI R_f$, \ $\boL E_f$, \ $\boR E^f$ are equivalences
of triangulated categories if and only if the morphism~$f$ induces
an isomorphism $H(A)\simeq H(B)$.
\end{thm}

\begin{proof}
 Morphisms in $\sD(A\modl)$ between shifts of the DG\+module $A$
recover the cohomology $H(A)$ and analogously for the DG\+algebra~$B$,
so the ``only if'' assertion follows from
the isomorphism $\boL E_f(A)\simeq B$.
 To prove the ``if'' part, we will show that the adjunction
morphisms $\boL E_f\allowbreak(\boI R_f(N))\rarrow N$ and $M\rarrow
\boI R_f(\boL E_f(M))$ are isomorphisms for any left DG\+modules $M$
over~$A$ and $N$ over~$B$.
 The former morphism is represented by the composition $B\ot_A G
\rarrow B\ot_A N\rarrow N$ for any flat DG\+module $G$ over~$A$
endowed with a quasi-isomorphism $G\rarrow N$ of DG\+modules over~$A$.
 This composition is a quasi-isomorphism, since the morphisms
$B\ot_A G\larrow A\ot_A G\rarrow A\ot_A N\simeq N$ are
quasi-isomorphisms.
 The latter morphism is represented by the fraction $M\larrow F
\rarrow B\ot_AF$ for any flat DG\+module $F$ over~$A$ endowed with
a quasi-isomorphism $F\rarrow M$ of DG\+modules over~$A$.
 The morphism $F\simeq A\ot_A F\rarrow B\ot_A F$ is
a quasi-isomorphism.
\end{proof}

\subsection{DG\+module t-structure}  \label{dg-module-t-structure}
 An object $Y$ of a triangulated category $\sD$ is called
an \emph{extension} of objects $Z$ and $X$ if there is
a distinguished triangle $X\rarrow Y\rarrow Z\rarrow X[1]$.
 Let $\sD=\sD(A\modl)$ denote the derived category of left
DG\+modules over a DG\+ring~$A$.
 Let $\sD^{\ge0}\subset\sD$ denote the full subcategory formed
by all DG\+modules $M$ over~$A$ such that $H^i(M)=0$ for $i<0$
and $\sD^{\le0}\subset\sD$ denote the minimal full subcategory
of $\sD(A\modl)$ containing the DG\+modules $A[i]$ for $i\ge0$
and closed under extensions and infinite direct sums.

\begin{thm}
\textup{(a)} The pair of subcategories\/ $(\sD^{\le0},\sD^{\ge0})$
defines a t-structure~\cite{BBD} on the derived category\/
$\sD(A\modl)$. \par
\textup{(b)} The subcategory $\sD^{\le0}\subset\sD$ coincides with
the full subcategory formed by all DG\+modules $M$ over~$A$ such
such that $H^i(M)=0$ for $i>0$ if and only if
$H^i(A)=0$ for all $i>0$.
\end{thm}

\begin{proof}
 Part~(a): clearly, one has $\sD^{\le0}[1]\subset\sD^{\le0}$, \
$\sD^{\ge0}[-1]\subset\sD^{\ge0}$, and $\Hom_{\sD}(\sD^{\le0},
\allowbreak\sD^{\ge0}[-1])=0$.
 It remains to construct for any DG\+module $M$ over~$A$ a closed
morphism of DG\+modules $F\rarrow M$ inducing a monomorphism on
$H^1$ and an isomorphism on $H^i$ for all $i\le0$ such that $F$
can be obtained from the DG\+modules $A[i]$ with $i\ge0$ by
iterated extensions and infinite direct sums in the homotopy
category of DG\+modules.
 This construction is similar to that of the proof of
Theorem~\ref{projective-dg-modules}, with the following changes.
 One chooses a surjective morphism $M'\rarrow M$ onto $M$ from
a complex of free abelian groups $M'$ with free abelian groups
of cohomology so that $H^i(M')=0$ for $i>0$ and the maps
$H^i(M')\rarrow H^i(M)$ are surjective for all $i\le 0$.
 Then for $F_0=A\ot_\boZ M'$ and $K=\ker(F_0\to M)$ one chooses
a surjective morphism $K'\rarrow K$ onto $K$ from a complex of free
abelian groups $K'$ with free abelian groups of cohomology so that
$H^i(K')=0$ for $i>1$ and the maps $H^i(K')\rarrow H^i(K)$ are
surjective for all $i\le 1$ in order to put $F_1=A\ot_\boZ K'$, etc.
 The DG\+module $F$ is constructed as the total DG\+module of
the complex $\dsb \rarrow F_1\rarrow F_0$ formed by taking infinite
direct sums. 
 The ``only if'' assertion in part~(b) is clear.
 To prove ``if'', replace $A$ with its quasi-isomorphic DG\+subring
$\tau_{\le0}A$ with the components $(\tau_{\le0}A)^i=A^i$ for $i<0$, \
$(\tau_{\le0}A)^0=\ker(A^0\to A^1)$, and $(\tau_{\le0}A)^i=0$
for $i>0$; then notice that the canonical filtrations on DG\+modules
over $\tau_{\le0}A$ considered as complexes of abelian groups are
compatible with the action of the ring $\tau_{\le0}A$.
\end{proof}

\begin{rem1}
 The t-structure described in part~(a) of Theorem can well be
degenerate, though it is clearly nondegenerate under the assumptions
of part~(b).
 Namely, one can have $\bigcap_i\sD^{\le0}[i]\ne0$.
 For example, take $A=k[x]$ to be the graded algebra of polynomials
with one generator $x$ of degree~$1$ over a field~$k$ and endow it
with the zero differential.
 Then the graded A\+module $k[x,x^{-1}]$ considered as a DG\+module
with zero differential belongs to the above intersection, since it
can be presented as the inductive limit of the DG\+modules
$x^{-j}k[x]$.
 Moreover, take $A=k[x,x^{-1}]$, where $\deg x=1$ and $d(x)=0$;
then $\sD^{\ge0}=0$ and $\sD^{\le0}=\sD$.
\end{rem1}

\begin{rem2}
 One might wish to define a dual version of the above t-structure
on $\sD(A\modl)$ where $\sD^{\ge0}$ would be the minimal full
subcategory of $\sD$ containing the DG\+modules
$\Hom_\boZ(A,\boQ/\boZ)[i]$ for $i\le0$ and closed under extensions
and infinite products, while $\sD^{\le0}$ would consist of all
DG\+modules $M$ with $H^i(M)=0$ for $i>0$.
 The dual version of the above proof does not seem to work in this
case, however, because of a problem related to nonexactness of
the countable inverse limit.
\end{rem2}

\begin{rem3}
 The above construction of the DG\+module t-structure can be
generalized in the following way (cf.~\cite{Sou}).
 Let $\sD$ be a triangulated category with infinite direct sums.
 An object $C\in\sD$ is said to be compact if the functor
$\Hom_{\sD}(C,{-})$ preserves infinite direct sums.
 Let $\sC\subset\sD$ be a subset of objects of $\sD$ consisting
of compact objects and such that $\sC[1]\subset\sC$.
 Let $\sD^{\ge 0}$ be the full subcategory of $\sD$ formed by
all objects $X$ such that $\Hom_{\sD}(C,X[-1])=0$ for all
$C\in\sC$, and let $\sD^{\le0}$ be the minimal full subcategory
of $\sD$ containing $\sC$ and closed under extensions and
infinite direct sums.
 Then $(\sD^{\le0},\sD^{\ge0})$ is a t-structure on~$\sD$.
 Indeed, let $X$ be an object of~$\sD$.
 Consider the natural map into $X$ from the direct sum of objects
from $\sC$ indexed by morphisms from objects of $\sC$ to $X$;
let $X_1$ be the cone of this map.
 Applying the same construction to the object $X_1$ in place of $X$,
we obtain the object $X_2$, etc.
 Let $Y$ be the homotopy inductive limit of $X_n$, i.~e., the cone
of the natural map $\bigoplus_nX_n\rarrow\bigoplus_nX_n$.
 Then $Y\in\sD^{\ge0}[-1]$ and $\cone(X\to Y)[-1]\in\sD^{\le0}$.
\end{rem3}

\subsection{Silly filtrations}
 Let $A$ be a DG\+ring and $\sD=\sD(A\modl)$ denote the derived
category of left DG\+modules over it.
 Denote by $\sD^{\le 0}\subset\sD$ the full subcategory formed
by all the DG\+modules $M$ such that $H^i(M)=0$ for $i>0$ and
by $\sD^{\ge 1}\subset\sD$ the full subcategory of all 
the DG\+modules $M$ such that $H^i(M)=0$ for $i\le 0$.

 We refer the reader to~\cite[Introduction and Appendices B--C]{P2}
for the general discussion of silly filtrations.
 For the purposes of the present exposition it suffices to say that
the directions of arrows in the distinguished triangle in the next
Theorem are opposite to the ones in distinguished triangles related
to t\+structures.

\begin{thm}
 One has $H^i(A)=0$ for all $i<0$ if and only if any object
$X\in\sD$ can be included into a distinguished triangle
$Y\rarrow X\rarrow Z\rarrow Y[1]$ in $\sD$ with $Y\in\sD^{\ge 1}$
and $Z\in\sD^{\le 0}$.
\end{thm}

\begin{proof}
 ``Only if'': suppose that the left DG\+module $X=A$ over $A$ can
be included into a distinguished triangle $Y\rarrow X\rarrow Z
\rarrow Y[1]$ in $\sD$ with $Y\in\sD^{\ge 0}$ and $Z\in\sD^{\le -1}$
(in the obvious notation).
 Then $\Hom_\sD(X,Z) \simeq H^0(Z)=0$, hence $X$ is a direct summand
of $Y$ in $\sD$, so $H^i(X)=0$ for $i<0$.

 ``If'': let $M$ be a DG\+module representing the object~$X$.
 Let $P$ be the direct sum of left DG\+modules $A[-i]$
over $A$ taken over all the elements $x\in H^i(M)$ for all $i\ge1$.
 There is a natural closed morphism of DG\+modules $P\rarrow M$;
set $M_1$ to be its cone.
 One has $H^i(P)=0$ for all $i\le 0$ and the maps $H^i(P)
\rarrow H^i(M)$ are surjective for all $i\ge 1$.
 Hence the maps $H^i(M)\rarrow H^i(M_1)$ are isomorphisms for
all $i<0$, a monomorphism for $i=0$, and zero maps for $i>0$.
 Repeating the same procedure for the DG\+module $M_1$, etc., we
obtain a sequence of DG\+module morphisms $M\rarrow M_1\rarrow M_2
\rarrow \dsb$, each morphism having the properties described above.
 Set $N=\varinjlim M_i$; then there is a closed morphism of
DG\+modules $M\rarrow N$ over $A$ such that the maps $H^i(M)\rarrow
H^i(N)$ are isomorphisms for all $i<0$ and a monomorphism for $i=0$,
and one has $H^i(N)=0$ for all $i>0$.
 Set $Q=\cone(M\to N)[-1]$; then one has $H^i(Q)=0$ for all $i\le 0$.
 This provides the desired distinguished triangle.
\end{proof}

\begin{rem1}
 There is a much simpler proof of the above Theorem applicable in
the case when $A^i=0$ for all $i<0$.
 However, unlike in the situation of (the proof of)
Theorem~\ref{dg-module-t-structure}(b), it is \emph{not} true that
any DG\+ring with zero cohomology in the negative degrees can be
connected by a chain of quasi-isomorphisms with a DG\+ring with zero
components in the negative degrees.
 For a counterexample, consider the free associative algebra $A$
over a field~$k$ generated by the elements $x$, $y$, and~$\eta$ in
degree~$0$, \ $z$ in degree~$1$, and $\xi$ in degree~$-1$, with
the differential given by $d(x)=d(y)=d(z)=0$, \ $d(\xi)=xy$, and
$d(\eta)=yz$.
 One can check that $H^i(A)=0$ for $i<0$.
 The nontrivial Massey product $\xi z-x\eta$ of the cohomology
classes $x$, $y$, and~$z$ provides the obstruction.
\end{rem1}

\textit{Erratum added five years later:}
 The counterexample in Remark~1 does not work.
 The Massey product of $x$, $y$, and~$z$ is indeed defined and does
not contain zero in $H^*(A)$, but it does not obstruct existence of
a DG\+algebra $B$ quasi-isomorphic to $A$ with $B^i=0$ for $i<0$.
 In fact, the quotient algebra of $A$ by the two-sided ideal
generated by $\xi$ and $d(\xi)$ can be used in the role of~$B$.
 To the best of the author's knowledge, the question whether any
DG\+ring with zero cohomology in the negative degrees can be
connected by a chain of quasi-isomorphisms with a DG\+ring with
zero components in the negative degrees remains open.

\begin{rem2}
 The construction of the above Theorem can be extended to arbitrary
DG\+rings in the following way.
 Given a DG\+ring $A$, denote by $\sD^{\le0}\subset\sD=\sD(A\modl)$
the full subcategory formed by all the DG\+modules $M$ over $A$
such that $H^i(M)=0$ for $i<0$ and by $\sD^{\ge1}\subset\sD$
the minimal full subcategory of $\sD$ containing the DG\+modules
$A[-i]$ with $i\ge1$ and closed under infinite direct sums and
the following operation of countably iterated extension.
 Given a sequence of objects $X_i\in\sD^{\ge 1}$ and a sequence of
distinguished triangles $Y_i\rarrow Y_{i+1}\rarrow X_{i+1}\rarrow
Y_i[1]$ in $\sD$ with $Y_0=X_0$, the homotopy colimit
$\cone(\bigoplus_i Y_i\to \bigoplus_i Y_i)$ of $Y_i$ should
also belong to $\sD^{\ge 1}$.
 Then the same construction as in the above proof provides for
any $X\in\sD$ a distinguished triangle $Y\rarrow X\rarrow Z\rarrow
Y[1]$ with $Y\in\sD^{\ge 1}$ and $Z\in\sD^{\le 0}$.
 This can be generalized even further in the spirit of
Remark~\ref{dg-module-t-structure}.3, by considering an arbitrary
triangulated category $\sD$ admitting infinite direct sums, and
a set of compact objects $\sC\subset\sD$ such that $\sC[-1]\subset
\sC$ in the role of the DG\+modules $A[-i]$ with $i\ge 1$.
 The subcategory $\sD^{\ge 1}$ is then generated by $\sC$ using
the operations of infinite direct sum and countably iterated
extension, and the subcategory $\sD^{\le 0}$ consists of all
object $X\in\sD$ such that $\Hom_{\sD}(C,X)=0$ for all $C\in\sC$.
 If one does not insist on the subcategories $\sD^{\le0}$ and
$\sD^{\ge 1}$ being closed under shifts in the respective
directions, the condition that $\sC[-1]\subset\sC$ can be dropped.
 Moreover, the subcategories $\sD^{\le0}$ and $\sD^{\ge1}$
described in this remark have a semiorthogonality property,
$\Hom_{\sD}(Y,Z)=0$ for all $Y\in\sD^{\ge1}$ and $Z\in\sD^{\le0}$
(cf.~\cite{Pauk}).
 Indeed, given an object $Z$ in a triangulated category $\sD$,
the class of all objects $Y\in\sD$ such that $\Hom_\sD(Y,Z)=0$
is closed under infinite direct sums and countably iterated
extensions, as one can see using the fact that the first derived
functor of projective limit of a sequence of surjective maps of
abelian groups vanishes.
\end{rem2}

\Section{Derived Categories of DG-Comodules and DG-Contramodules}

\subsection{Graded comodules}
 Let $k$ be a fixed ground field.
 A \emph{graded coalgebra} $C$ over~$k$ is a graded $k$\+vector space
$C=\bigoplus_{i\in\boZ}C^i$ endowed with a comultiplication map
$C\rarrow C\ot_kC$ and a counit map $C\rarrow k$, which must be
homogeneous linear maps of degree~$0$ satisfying the coassociativity
and counity equations.
 Namely, the comultiplication map must have equal compositions
with the two maps $C\ot_kC\birarrow C\ot_kC\ot_kC$ induced by
the comultiplication map, while the compositions of
the comultiplication map with the two maps $C\ot_kC\rarrow C$
induced by the counit map must coincide with the identity
endomorphism of~$C$.

 A \emph{graded left comodule} $M$ over~$C$ is a graded $k$\+vector
space $M=\bigoplus_{i\in\boZ}M^i$ endowed with a left coaction map
$M\rarrow C\ot_kM$, which must be a homogeneous linear map of
degree~$0$ satisfying the coassociativity and counity equations.
 Namely, the coaction map must have equal compositions with
the two maps $C\ot_kM\birarrow C\ot_kC\ot_kM$ induced by
the comultiplication map and the coaction map, while the composition
of the coaction map with the map $C\ot_kM\rarrow M$ induced by
the counit map must coincide with the identity endomorphism of~$M$.
 A \emph{graded right comodule} $N$ over~$C$ is a graded vector
space endowed with a right coaction map $N\rarrow N\ot_kC$
satisfying the analogous linearity, homegeneity, coassociativity,
and counity equations.

 The \emph{cotensor product} of a graded right $C$\+comodule $N$ and
a graded left $C$\+comodule $M$ is the graded vector space $N\oc_CM$
defined as the kernel of the pair of linear maps
$N\ot_kM\birarrow N\ot_kC\ot_kM$, one of which is induced by
the right coaction map and the other by the left coaction map.
 There are natural isomorphisms $C\oc_CM\simeq M$ and $N\oc_CC
\simeq N$ for any graded left $C$\+comodule $M$ and graded right
$C$\+comodule $N$.

 Graded left $C$\+comodules of the form $C\ot_k V$, where $V$ is
a graded vector space, are called \emph{cofree} graded left
$C$\+comodules; analogously for graded right $C$\+comodules.
 The category of graded $C$\+comodules is an abelian category
with enough injectives; injective graded $C$\+comodules are
exactly the direct summands of cofree graded $C$\+comodules.

 For any graded left $C$\+comodules $L$ and $M$, the graded vector
space $\Hom_C(L,M)$ consists of homogeneous linear maps
$f\:L\rarrow M$ satisfying the condition that the coaction maps of
$L$ and $M$ form a commutative diagram together with the map~$f$
and the map $f_*\:C\ot_kL\rarrow C\ot_kM$ given by the formula
$f_*(c\ot x)=(-1)^{|f||c|}c\ot f(x)$.
 For any graded right $C$\+comodules $R$ and $N$, the graded vector
space $\Hom_C(R,N)$ consists of homogeneous linear maps $f\:R\rarrow N$
such that the coaction maps of $R$ and $N$ form a commutative diagram
together with the map~$f$ and the map $f_*\:R\ot_kC\rarrow N\ot_kC$
given by the formula $f_*(x\ot c)=f_*(x)\ot c$.
 For any left $C$\+comodule $L$ and any graded vector space $V$
there is a natural isomorphism $\Hom_C(L\;C\ot_kV)\simeq
\Hom_k(L,V)$; analogously in the right comodule case.

\subsection{Graded contramodules}  \label{graded-contramodules}
 A \emph{graded left contramodule} $P$ over a graded coalgebra~$C$
is a graded $k$\+vector space $P=\bigoplus_{i\in\boZ}P^i$
endowed with the following structure.
 Let $\Hom_k(C,P)$ be the graded vector space of homogeneous linear
maps $C\rarrow P$; then a homogeneous linear map
$\Hom_k(C,P)\rarrow P$ of degree~$0$, called the left contraaction map,
must be given and the following contraassociativity and counity
equations must be satisfied.
 For any graded vector spaces $V$, \ $W$, and $P$, define the natural
isomorphism $\Hom_k(V\ot_kW\;P)\simeq\Hom_k(W,\Hom_k(V,P))$ by
the formula $f(w)(v)=(-1)^{|w||v|}f(v\ot w)$.
 The comultiplication and the contraaction maps induce a pair of
maps $\Hom_k(C\ot_kC\;P)\simeq\Hom_k(C,\Hom_k(C,P))
\birarrow\Hom_k(C,P)$.
 These maps must have equal compositions with the contraaction map;
besides, the composition of the map $P\rarrow\Hom_k(C,P)$ induced by
the counit map with the contraaction map must coincide with
the identity endomorphism of~$P$.

 The graded vector space of \emph{cohomomorphisms} $\Cohom_C(M,P)$
from a graded left $C$\+comodule $M$ to a graded left $C$\+contramodule
$P$ is defined as the cokernel of the pair of linear maps
$\Hom_k(C\ot_kM\;P)\rarrow\Hom_k(M,P)$, one of which is induced by
the left coaction map and the other by the left contraaction map.
 For any graded $C$\+contramodule $P$ there is a natural isomorphism
$\Cohom_C(C,P)\simeq P$.
 
 For any graded right $C$\+comodule $N$ and any graded vector space $V$
there is a natural graded left $C$\+contramodule structure on
the graded vector space of homogeneous linear maps $\Hom_k(N,V)$ 
given by the left contraaction map $\Hom_k(C,\Hom_k(N,V))\simeq
\Hom_k(N\ot_kC\;V)\rarrow\Hom_k(N,V)$ induced by the right coaction map.
 For any graded left $C$\+comodule $M$, graded right $C$\+comodule $N$,
and graded vector space $V$, there is a natural isomorphism
$\Hom_k(N\oc_CM\;V)\simeq\Cohom_C(M,\Hom_k(N,V))$.

 Graded left $C$\+contramodules of the form $\Hom_k(C,V)$ are called
\emph{free} graded left $C$\+contramodules.
 The category of graded left $C$\+contramodules is an abelian category
with enough projectives; projective graded left $C$\+contramodules
are exactly the direct summands of free graded left $C$\+contramodules.

 The \emph{contratensor product} of a graded right $C$\+comodule $N$
and a graded left $C$\+contramodule $P$ is the graded vector space
$N\ocn_CP$ defined as the cokernel of the pair of linear maps
$N\ot_k\Hom_k(C,P)\rarrow N\ot_kP$, one of which is induced by
the left contraaction map, while the other one is obtained as
the composition of the map induced by the right coaction map and
the map induced by the evaluation map $C\ot_k\Hom_k(C,P)\rarrow P$
given by the formula $c\ot f\mpsto (-1)^{|c||f|}f(c)$.
 For any graded right $C$\+comodule $N$ and any graded vector space
$V$ there is a natural isomorphism $N\ocn_C\Hom_k(C,V)\simeq N\ot_kV$.

 For any graded left $C$\+contramodules $P$ and $Q$, the graded vector
space $\Hom^C(P,Q)$ consists of all homogeneous linear maps
$f\:P\rarrow Q$ satisfying the condition that the contraaction maps
of $P$ and $Q$ form a commutative diagram together with the map~$f$
and the map $f_*\:\Hom_k(C,P)\rarrow\Hom_k(C,Q)$ given by the formula
$f_*(g)=f\circ g$.
 For any graded left $C$\+contramodule $Q$ and any graded vector
space $V$ there is a natural isomorphism $\Hom^C(\Hom_k(C,V),Q)\simeq
\Hom_k(V,Q)$.
 For any right $C$\+comodule $N$, any graded left $C$\+contramodule $P$,
and any graded vector space $V$, there is a natural isomorphism
$\Hom_k(N\ocn_CP\;V)\simeq\Hom^C(P,\Hom_k(N,V))$.

 The proofs of the results of this subsection are not difficult;
some details can be found in~\cite{P}.
 The assertions stated in the last paragraph can be used to deduce
the assertions of the preceding two paragraphs.

\begin{rem}
 Ungraded contramodules over ungraded coalgebras can be simply
defined as graded contramodules concentrated in degree~$0$ over
graded coalgebras concentrated in degree~$0$.
 One might wish to have a forgetful functor assigning ungraded
contramodules over ungraded coalgebras to graded contramodules
over graded coalgebras.
 The construction of such a functor is delicate in two ways.
 Firstly, to assign an ungraded contramodule to a graded
contramodule~$P$, one has to take the direct product of its
grading components $\prod_{i\in\boZ}P^i$ rather than the direct
sum, while the ungraded coalgebra corresponding to a graded
coalgebra~$C$ is still constructed as the direct sum
$\bigoplus_{i\in\boZ}C^i$.
 Analogously, to assign an ungraded comodule to a graded
comodule~$M$ one takes the direct sum $\bigoplus_{i\in\boZ}M^i$,
to assign an ungraded ring to a graded ring $A$ one takes
the direct sum $\bigoplus_{i\in\boZ}A^i$, while to assign
an ungraded module to a graded module $M$ one can take
\emph{either} the direct sum $\bigoplus_{i\in\boZ}M^i$,
\emph{or} the direct product $\prod_{i\in\boZ}M^i$.
 Secondly, there is a problem of signs in the contraasociativity
equation, which is unique to graded contramodules (no signs
are present in the definitions of graded algebras, modules,
coalgebras, or comodules); it is resolved as follows.
 A morphism between graded vector spaces $f\:V\rarrow W$, when
it is not necessarily even, can be thought of either as a
\emph{left} or as a \emph{right} morphism.
 The left and the right morphisms correspond to each other
according to the sign rule $f(x)=(-1)^{|f||x|}(x)f$, where $f(x)$
is the notation for the left morphisms and $(x)f$ for the right
morphisms.
 The above definition of graded left contramodules~$P$ is given
in terms of \emph{left} morphisms $C\rarrow P$; to define
the functor of forgetting the grading, one has to reinterpret
it in terms of \emph{right} morphisms.
 The exposition in~\cite[Chapter~11]{P} presumes \emph{right}
morphisms in this definition (even if the notation $f(x)$ is
being used from time to time).
\end{rem}

\subsection{DG-comodules and contramodules}
 A \emph{DG\+coalgebra} $C$ over a field~$k$ is a graded coalgebra
endowed with a differential $d\:C\rarrow C$ of degree~$1$ with
$d^2=0$ such that the comultiplication map $C\rarrow C\ot_kC$ and
the counit map $C\rarrow k$ are morphisms of complexes.
 Here the differential on $C\ot_kC$ is defined as on the tensor
product of two copies of the complex~$C$, while the differential
on~$k$ is trivial.

 A \emph{left DG\+comodule} $M$ over a DG\+coalgebra $C$ is a graded
comodule $M$ over the graded coalgebra $C$ together with a differential
$d\:M\rarrow M$ of degree~$1$ with $d^2=0$ such that the left coaction
map $M\rarrow C\ot_kM$ is a morphism of complexes.
 Here $C\ot_kM$ is considered as the tensor product of the complexes
$C$ and $M$ over~$k$.
 \emph{Right DG\+comodules} are defined in the analogous way.
 A \emph{left DG\+contramodule} $P$ over $C$ is a graded contramodule
$P$ over $C$ endowed with a differential $d\:P\rarrow P$ of degree~$1$
with $d^2=0$ such that the left contraaction map $\Hom_k(C,P)\rarrow P$
is a morphism of complexes.
 Here $\Hom_k(C,P)$ is endowed with the differential of the complex
of homomorphisms from the complex $C$ to the complex $P$ over~$k$.

 Whenever $N$ is a right DG\+comodule and $M$ is a left DG\+comodule
over a DG\+coalgebra $C$, the cotensor product $N\oc_CM$ of
the graded comodules $N$ and $M$ over the graded coalgebra $C$ is
endowed with the differential of the subcomplex of the tensor product
complex $N\ot_kM$.
 Whenever $M$ is a left DG\+comodule and $P$ is a left DG\+contramodule
over a DG\+coalgebra $C$, the graded vector space of cohomomorphisms
$\Cohom_C(M,P)$ is endowed with the differential of the quotient
complex of the complex of homomorphisms $\Hom_k(M,P)$.

 Whenever $N$ is a right DG\+comodule and $P$ is a left
DG\+contramodule over a DG\+coalgebra $C$, the contratensor product
$N\ocn_CP$ of the graded comodule $N$ and the graded contramodule $P$
over the graded coalgebra $C$ is endowed with the differential of
the quotient complex of the tensor product complex $N\ot_kP$.

 For any left DG\+comodules $L$ and $M$ over a DG\+coalgebra $C$,
the graded vector space of homomorphisms $\Hom_C(L,M)$ between
the graded comodules $L$ and $M$ over the graded coalgebra $C$ is
endowed with the differential of the subcomplex of the complex
of homomorphisms $\Hom_k(L,M)$.
 Differentials on the graded vector spaces of homomorphisms
$\Hom_C(R,N)$ and $\Hom^C(P,Q)$ for right DG\+comodules $R$, $N$
and left DG\+contramodules $P$, $Q$ over a DG\+coalgebra $C$
are constructed in the completely analogous way.
 These constructions define the DG\+categories $\sDG(C\comodl)$, \
$\sDG(\comodr C)$, and $\sDG(C\contra)$ of left DG\+comodules,
right DG\+comodules, and left DG\+contramodules over~$C$, respectively.

 All shifts, twists, infinite direct sums, and infinite direct products
exist in the DG\+categories of DG\+comodules and DG\+contramodules.
 The homotopy category of (the DG\+category of) left DG\+comodules
over $C$ is denoted by $\Hot(C\comodl)$, the homotopy category of
right DG\+comodules over $C$ is denoted by $\Hot(\comodr C)$, and
the homotopy category of left DG\+contramodules over $C$ is
denoted by $\Hot(C\contra)$.

\subsection{Injective and projective resolutions}
\label{dg-coalgebra-resolutions}
 A DG\+comodule $M$ or a DG\+contra\-module $P$ is said to be acyclic
if it is acyclic as a complex of vector spaces, i.~e., $H(M)=0$ or
$H(P)=0$, respectively.
 The classes of acyclic DG\+comodules over
a DG\+coalgebra $C$  are closed under shifts, cones, and infinite
direct sums, while the class of acyclic DG\+contramodules is closed
under shifts, cones, and infinite products.
 The thick subcategories of the homotopy categories $\Hot(C\comodl)$,
\ $\Hot(\comodr C)$, and $\Hot(C\contra)$ formed by the acyclic
DG\+comodules and DG\+contramodules over~$C$ are denoted by
$\Acycl(C\comodl)$, \ $\Acycl(\comodr C)$, and $\Acycl(C\contra)$,
respectively.
 The \emph{derived categories} of left DG\+comodules, right
DG\+comodules, and left DG\+contramodules over $C$ are defined as
the quotient categories $\sD(C\comodl)=\Hot(C\comodl)/\Acycl(C\comodl)$,
\ $\sD(\comodr C)=\Hot(\comodr C)/\Acycl(\comodr C)$, and
$\sD(C\contra)=\Hot(C\contra)/\Acycl(C\contra)$.

 A left DG\+comodule $M$ over a DG\+coalgebra $C$ is called
\emph{injective} if for any acyclic left DG\+comodule $L$ over $C$
the complex $\Hom_C(L,M)$ is acyclic.
 The full triangulated subcategory of $\Hot(C\comodl)$ formed by
the injective DG\+comodules is denoted by $\Hot(C\comodl)_\inj$.
 A left DG\+contramodule $P$ over a DG\+coalgebra $C$ is called
\emph{projective} if for any acyclic left DG\+contramodule $Q$
over $C$ the complex $\Hom^C(P,Q)$ is acyclic.
 The full triangulated subcategory of $\Hot(C\contra)$ formed by
the projective DG\+contramodules is denoted by $\Hot(C\contra)_\proj$.

\begin{thm}
\textup{(a)} The composition of functors\/ $\Hot(C\comodl)_\inj\rarrow
\Hot(C\comodl)\rarrow\sD(C\comodl)$ is an equivalence of triangulated
categories. \par
\textup{(b)} The composition of functors\/ $\Hot(C\contra)_\proj\rarrow
\Hot(C\contra)\rarrow\sD(C\contra)$ is an equivalence of triangulated
categories.
\end{thm}

Proof will be given in subsection~\ref{dg-coalgebra-resolutions-proof}.

\begin{rem}
 The analogue of Theorem~\ref{dg-mod-scalars} does not hold for
DG-coalgebras.
 More precisely, given a morphism of DG\+coalgebras $f\:C\rarrow D$,
any DG\+comodule over $C$ can be considered as a DG\+comodule over $D$
and any DG\+contramodule over $C$ can be considered as
a DG\+contramodule over $D$, so there are restriction-of-scalars
functors $\boI R_f\:\sD(C\comodl)\rarrow\sD(D\comodl)$ and
$\boI R^f\:\sD(C\contra)\rarrow\sD(D\contra)$.
 Using the above Theorem, one can construct the functor $\boR E_f$
right adjoint to $\boI R_f$ and the functor $\boL E^f$ left adjoint
to $\boI R^f$ (cf.~\ref{cdg-coalgebra-scalars}).
 Just as in the proof of Theorem~\ref{dg-mod-scalars}, one can show
that the morphism $f$ induces an isomorphism of cohomology whenever
either of the functors $\boI R_f$ or $\boI R^f$ is an equivalence of
categories.
 However, there exist quasi-isomorphisms of DG\+coalgebras
$f\:C\rarrow D$ for which the functors $\boI R_f$ and $\boI R^f$
are not equivalences.
 The following counterexample is due to D.~Kaledin~\cite{Kal}.
 Let $E$ and $F$ be DG\+coalgebras and $K$ be a DG\+bicomodule over
$E$ and $F$, i.~e., a complex of $k$\+vector spaces with commuting
structures of left DG\+comodule over $E$ and right DG\+comodule
over~$F$.
 Then there are natural DG\+coalgebra structures on the direct sums
$C=E\oplus K\oplus F$ and $D=E\oplus F$; there are also morphisms of
DG\+coalgebras $D\rarrow C\rarrow D$.
 When the complex $K$ is acyclic, these morphisms are
quasi-isomorphisms.
 A left DG\+comodule over $C$ is the same that a pair $(M,N)$ of left
DG\+comodules $M$ over $E$ and $N$ over $F$ together with a closed
morphism $M\rarrow K\oc_F N$ of left DG\+comodules over~$E$.
 For a DG\+comodule $M$ over $E$ and a DG\+comodule $N$ over $F$,
morphisms $(M,0)\rarrow (0,N[1])$ in the derived category
$\sD(C\comodl)$ are represented by morphisms $M'\rarrow K\oc_FN'$ in
the homotopy category $\Hot(E\comodl)$ together with quasi-isomorphisms
of DG\+comodules $M'\rarrow M$ over $E$ and $N\rarrow N'$ over~$F$.
 In particular, when $M$ is a projective DG\+comodule (i.~e., there are
no nonzero morphisms in $\Hot(E\comodl)$ from $M$ into acyclic
DG\+comodules) and $N$ is an injective DG\+comodule, the vector space
of morphisms $(M,0)\rarrow (0,N[1])$ in $\sD(C\comodl)$ is isomorphic
to the vector space of morphisms $M\rarrow K\oc_FN$ in $\Hot(E\comodl)$.
 This vector space can be nonzero even when $K$ is acyclic, because
the cotensor product of an acyclic DG\+comodule and an injective
DG\+comodule over $F$ can have nonzero cohomology.
 It suffices to take $E$ to be the field~$k$ and $F$ to be the coalgebra
dual to the exterior algebra in one variable $k[\eps]/\eps^2$,
concentrated in degree~$0$ and endowed with zero differential.
 Analogously, a left DG\+contramodule over $C$ is the same that a pair
$(P,Q)$ of left DG\+contramodules $P$ over $E$ and $Q$ over $F$ together
with a closed morphism $\Cohom_E(K,P)\rarrow Q$ of left
DG\+contramodules over~$F$.
 The vector space of morphisms $(P,0)\rarrow (0,Q[1])$ in
$\sD(C\contra)$ can be nonzero even when $K$ is acyclic, because
the complex $\Cohom$ from an acyclic DG\+comodule to a projective
DG\+contramodule over $E$ can have nonzero cohomology.
\end{rem}

\subsection{Cotor and Coext of the first kind}
\label{cotor-coext-first-kind}
 Let $N$ be a right DG\+comodule and $M$ be a left DG\+comodule over
a DG\+coalgebra~$C$.
 Consider the cobar bicomplex $N\ot_k M\rarrow N\ot_kC\ot_kM\rarrow
N\ot_kC\ot_kC\ot_kM\rarrow\dsb$ and construct its total complex by
taking infinite products.
 Let $\Cotor^{C,I}(N,M)$ denote the cohomology of this total complex.

 Let $M$ be a left DG\+comodule and $P$ be a left DG\+contramodule
over a DG\+coalgebra~$C$.
 Consider the bar bicomplex $\dsb\rarrow\Hom_k(M\ot_kC\ot_kC\;P)
\rarrow\Hom_k(M\ot_kC\;P)\rarrow\Hom_k(C,P)$ and construct its total
complex by taking infinite direct sums.
 Let $\Coext_C^I(M,P)$ denote the cohomology of this total complex.

 Let $k\vect$ denote the category of vector spaces over~$k$ and
$k\vect^\sgr$ denote the category of graded vector spaces over~$k$.

\begin{thm}
\textup{(a)}
 The functor $\Cotor^{C,I}$ factorizes through the Cartesian product
of the derived categories of right and left DG\+comodules over $C$, so
there is a well-defined functor $\Cotor^{C,I}\:\sD(\comodr C)\times
\sD(C\comodl)\rarrow k\vect^\sgr$. \par
\textup{(b)}
 The functor $\Coext_C^I$ factorizes through the Cartesian product
of the derived categories of left DG\+comodules and left
DG\+contramodules over $C$, so there is a well-defined functor
$\Coext_C^I\:\sD(C\comodl)\times\sD(C\contra)\rarrow k\vect^\sgr$.
\end{thm}

\begin{proof}
 This follows~\cite[I.3--4]{HMS} from the fact that a complete and
cocomplete filtered complex is acyclic whenever the associated graded
complex is acyclic~\cite{EM1}.
\end{proof}

 Let $N$ be a right DG\+comodule and $P$ be a left DG\+contramodule
over a DG\+coalgebra~$C$.
 Consider the bar bicomplex $\dsb\rarrow N\ot_k\Hom_k(C\ot_kC\;P)
\rarrow N\ot_k\Hom_k(C,P)\rarrow N\ot_k P$ and construct its total
complex by taking infinite direct sums.
 Let $\Ctrtor^{C,I}(N,P)$ denote the cohomology of this total complex.

 Let $L$ and $M$ be left DG\+comodules over a DG\+coalgebra~$C$.
 Consider the cobar bicomplex $\Hom_k(L,M)\rarrow\Hom_k(L\;C\ot_kM)
\rarrow\Hom_k(L\;C\ot_kC\ot_kM)\rarrow\dsb$ and construct its total
complex by taking infinite products.
 Let $\Ext_C^I(L,M)$ denote the cohomology of this total complex.

 Let $P$ and $Q$ be left DG\+contramodules over~$C$.
 Consider the cobar bicomplex $\Hom_k(P,Q)\rarrow\Hom_k(\Hom_k(C,P),Q)
\rarrow\Hom_k(\Hom_k(C\ot_kC\;P)\;Q)\rarrow\dsb$ and construct its
total complex by taking infinite products.
 Let $\Ext^{C,I}(P,Q)$ denote the cohomology of this total complex.

 Just as in the above Theorem, the functors $\Ctrtor^{C,I}$, \
$\Ext_C^I$, and $\Ext^{C,I}$ factorize through the derived categories
of DG\+comodules and DG\+contramodules.

 All of these constructions can be extended to $\Ainfty$\+comodules
and $\Ainfty$\+contra\-modules over $\Ainfty$\+coalgebras; see
Remark~\ref{curved-ainfty-co-contra-derived}.

\begin{rem}
 Another approach to defining derived functors of cotensor product,
cohomomorphisms, contratensor product, etc., whose domains would be
Cartesian products of the derived categories of DG\+comodules and
DG\+contramodules consists in restricting these functors to
the full subcategories of injective DG\+comodules and projective
DG\+contramodules in the homotopy categories.
 To obtain versions of derived functors $\Cotor^C$ and $\Coext_C$
in this way one would have to restrict the functors of cotensor
product and cohomomorphisms to the homotopy categories of injective
DG\+comodules and projective DG\+contramodules in \emph{both}
arguments; resolving only one of the arguments does not provide
a functor factorizable through the derived category in the other
argument~\cite[section~0.2.3]{P}.
 To construct version of derived functors $\Ctrtor^C$, \ $\Ext_C$,
and $\Ext^C$, on the other hand, it suffices to resolve just
\emph{one} of the arguments (the second one, the second one, and
the first one, respectively).
 The versions of $\Ext_C$ and $\Ext^C$ so obtained coincide with
the functors $\Hom$ in the derived categories.
 It looks unlikely that the derived functors defined in the way of
this Remark should agree with the derived functors defined above in
this subsection.
\end{rem}

\Section{Coderived and Contraderived Categories of CDG-Modules}
\label{cdg-module-section}

\subsection{CDG-rings and CDG-modules}
 A \emph{CDG\+ring} (curved differential graded ring) $B=(B,d,h)$ is
a triple consisting of an associative graded ring
$B=\bigoplus_{i\in\boZ}B^i$, an odd derivation $d\:B\rarrow B$
of degree~$1$, and an element $h\in B^2$ satisfying the equations
$d^2(x)=[h,x]$ for all $x\in B$ and $d(h)=0$.
 A \emph{morphism} of CDG\+rings $f\:B\rarrow A$ is a pair $f=(f,a)$
consisting of a morphism of graded rings $f\:B\rarrow A$ and
an element $a\in A^1$ satisfying the equations $f(d_B(x))=
d_A(f(x))+[a,x]$ and $f(h_B)=h_A+d_A(a)+a^2$ for all $x\in B$,
where $B=(B,d_B,h_B)$ and $A=(A,d_A,h_A)$, while the bracket
$[y,z]$ denotes the supercommutator of $y$ and~$z$.
 The composition of morphisms is defined by the rule $(f,a)\circ
(g,b)=(f\circ g\;a+f(b))$.
 Identity morphisms are the morphisms $(\id,0)$.

 The element $h\in B^2$ is called the \emph{curvature element} of
a CDG\+ring $B$.
 The element $a\in A^1$ is called
the \emph{change-of-connection element} of a CDG\+ring morphism~$f$.

 To any DG\+ring structure on a graded ring $A$ one can assign
a CDG\+ring structure on the same graded ring by setting $h=0$.
 This defines a functor from the category of DG\+rings to
the category of CDG\+rings.
 This functor is faithful, but not fully faithful, as non-isomorphic
DG\+rings may become isomorphic as CDG\+rings.

 A \emph{left CDG\+module} $(M,d_M)$ over a CDG\+ring $B$ is
a graded left $B$\+module $M=\bigoplus_{i\in\boZ}M^i$ endowed
with a derivation $d_M\:M\rarrow M$ compatible with the derivation~$d_B$
of~$B$ and such that $d_M^2(x)=hx$ for any $x\in M$.
 A \emph{right CDG\+module} $(N,d_N)$ over a CDG\+ring $B$ is
a graded right $B$\+module $N=\bigoplus_{i\in\boZ}N^i$ endowed
with a derivation $d_N\:N\rarrow N$ compatible with~$d_B$ and
such that $d_N^2(x)=-xh$ for any $x\in N$.

 Notice that there is no natural way to define a left or right
CDG\+module structure on the free left or right graded module $B$
over a CDG\+ring~$B$.
 At the same time, any CDG\+ring $B$ is naturally
a \emph{CDG\+bimodule} over itself, in the obvious sense
(see~\ref{finite-over-gorenstein} for the detailed definition).

 Let $f=(f,a)\:B\rarrow A$ be a morphism of CDG\+rings and $(M,d_M)$
be a left CDG\+module over~$A$.
 Then the left CDG\+module $R_fM$ over~$B$ is defined as the graded
abelian group $M$ with the graded $B$\+module structure obtained
by the restriction of scalars via~$f$, endowed with the differential
$d'_M(x)=d_M(x)+ax$ for $x\in M$.
 Analogously, let $(N,d_N)$ be a right CDG\+module over~$A$.
 Then the right CDG\+module $R_fN$ over~$B$ is defined as the graded
module $N$ over~$B$ with the graded module structure induced by~$f$
endowed with the differental $d'_N(x)=d_N(x)-(-1)^{|x|}xa$.

 For any left CDG\+modules $L$ and $M$ over~$B$, the \emph{complex
of homomorphisms} $\Hom_B(L,M)$ from $L$ to $M$ over~$B$ is constructed
using exactly the same formulas as in~\ref{dg-rings-modules}.
 It turns out that these formulas still define a complex in
the CDG\+module case, as the $h$-related terms cancel each other.
 The same applies to the definitions in the subsequent two paragraphs
in~\ref{dg-rings-modules}, which all remain applicable in
the CDG\+module case, including the definitions of the complex of
homomorphisms between right CDG\+modules, the \emph{tensor product
complex} of a left and a right CDG\+module, etc.

 The cohomology of CDG\+rings and CDG\+modules is not defined, though,
as their differentials may have nonzero squares.
 A \emph{CDG\+algebra} over a commutative ring~$k$ is a graded
$k$\+module with a $k$\+linear CDG\+ring structure. 

 So left and right CDG\+modules over a given CDG\+ring $B$ form
DG\+categories, which we denote, just as the DG\+categories of
DG\+modules, by $\sDG(B\modl)$ and $\sDG(\modr B)$, respectively.
 All shifts, twists, infinite direct sums, and infinite direct products
exist in the DG\+categories of CDG\+modules.
 The corresponding homotopy categories are denoted by $\Hot(B\modl)$
and $\Hot(\modr B)$.

 Notice that there is no obvious way to define derived categories
of CDG\+modules, since it is not clear what should be meant by
an acyclic CDG\+module.
 Moreover, the functors of restriction of scalars $R_f$ related to
CDG\+isomorphisms $f$ between DG\+rings may well transform acyclic
DG\+modules to non-acyclic ones.

 For a CDG\+ring $B$, we will sometimes denote by $B^\#$ the graded
ring $B$ considered without its differential and curvature element
(or with the zero differential and curvature element).
 For a left CDG\+module $M$ and a right CDG\+module $N$ over~$B$,
we denote by $M^\#$ and $N^\#$ the corresponding graded modules
(or CDG\+modules with zero differential) over~$B^\#$.

\subsection{Some constructions for DG\+categories}
\label{dg-categories-constr}
 The reader will easily recover the details of the constructions
sketched below.

 Let $\sDG$ be a DG\+category.
 Define the DG\+category $\sDG^\natural$ by the following
construction.
 An object of $\sDG^\natural$ is a pair $(Z,t)$, where $Z$
in an object of $\sDG$ and $t\in\Hom_{\sDG}^{-1}(Z,Z)$ is
a contracting homotopy with zero square, i.~e., $d(t)=\id_X$
and $t^2=0$.
 Morphisms $(Z',t')\rarrow(Z'',t'')$ of degree~$n$ in
$\sDG^\natural$ are morphisms $f\:Z'\rarrow Z''$ of degree~$-n$
in $\sDG$ such that $d(f)=0$ in $\sDG$.
 The differential on the complex of morphisms in $\sDG^\natural$
is given by the supercommutator with~$t$, i.~e., $d^\natural(f) =
t''f - (-1)^{|f|}ft'$.

 Obviously, all twists of objects by their Maurer--Cartan
endomorphisms (see~\ref{dg-categories}) exist in the DG\+category
$\sDG^\natural$.
 Shifts, (finite or infinite) direct sums, or direct products
exist in $\sDG^\natural$ whenever they exist in $\sDG$.

 Let $B = (B,d,h)$ be a CDG\+ring.
 Construct the DG\+ring $B\sptilde=(B\sptilde,\d)$ as follows.
 The graded ring $B\sptilde$ is obtained by changing the sign of
the grading in the ring $B[\delta]$, which is in turn constructed
by adjoining to $B$ an element $\delta$ of degree~$1$ with
the relations $[\delta,x]=d(x)$ for $x\in B$ and $\delta^2=h$.
 The differential $\d=\d/\d\delta$ is defined by the rules
$\d(\delta)=1$ and $\d(x)=0$ for all $x\in B$.
 This construction can be extended to an equivalence between
the categories of CDG\+rings and acyclic
DG\+rings~\cite[section~0.4.4]{P} (here a DG\+ring is called
acyclic if its cohomology is the zero ring).
 There is a natural isomorphism of DG\+categories
$\sDG(B\modl)^\natural\simeq\sDG(B\sptilde\modl)$.

 Let $\sDG$ be a DG\+category with shifts and cones.
 Denote by $\natural\:X\mpsto X^\natural$ the functor $Z^0(\sDG)
\rarrow Z^0(\sDG^\natural)$ assigning to an object $X$ the object
$\cone(\id_X)[-1]$ with its standard contracting homotopy~$t$.
 This functor can be extended in a natural way to a fully
faithful functor $\sDG^\#\rarrow Z^0(\sDG^\natural)$, since not
necessarily closed morphisms of degree~$0$ also induce closed
morphisms of the cones of identity endomorphisms commuting with
the standard contracting homotopies.

 The functor $\natural$ has left and right adjoint functors $G^+$,
$G^-\:Z^0(\sDG^\natural)\rarrow Z^0(\sDG)$, which are given by
the rules $G^+(Z,t)=Z$ and $G^-(Z,t)=Z[1]$, so $G^+$ and $G^-$
only differ by a shift.
 Whenever all infinite direct sums (products) exist in
the DG\+category $\sDG$, the functors $G^+$, $G^-$, and $\natural$
preserve them.

 In particular, when $\sDG=\sDG(B\modl)$, the category
$Z^0(\sDG^\natural)$ can be identified with the category of graded
left $B^\#$\+modules in such a way that the functor~$\natural$
becomes the functor $M\mpsto M^\#$ of forgetting the differential.
 The category $\sDG(B\modl)^\#$ is then identified with the full
subcategory consisting of all graded $B^\#$\+modules that admit
a structure of CDG\+module over~$B$.

 For any DG\+category $\sDG$, objects of the DG\+category
$\sDG^{\natural\natural}$ are triples $(W,t,s)$, where $W$ is
an object of $\sDG$ and $t$, $s\:W\rarrow W$ are endomorphisms
of degree $-1$ and~$1$, respectively, satisfying the equations
$t^2=0=s^2$, \ $ts+st=\id_W$, \ $d(t)=\id_W$, and $d(s)=0$.
 Assuming that shifts and cones exist in~$\sDG$, there is a natural
fully faithful functor $\sDG\rarrow\sDG^{\natural\natural}$ given
by the formula $W = \cone(\id_X)[-1]$.

 This functor is an equivalence of DG\+categories whenever
all twists of objects exist in $\sDG$ and all images of
idempotent endomorphisms exist in $Z^0(\sDG)$.
 Indeed, to recover the object $X$ from the object $W$, it
suffices to take the image of the closed idempotent
endomorphism $ts$ of the twisted object $W(-s)$.

 In particular, there is a natural equivalence of DG\+categories
$\sDG(B\modl)\simeq\sDG(B\sptilde{}\sptilde\modl)$.
 So the DG\+category of CDG\+modules over an arbitrary CDG\+ring
is equivalent to the DG\+category of DG\+modules over a certain
acyclic DG\+ring.

 The ``almost involution'' $\sDG\mpsto\sDG^\natural$ is not defined
on the level of homotopy categories.
 Indeed, if $\sDG$ is the DG\+category of complexes over an additive
category~$\sA$ containing images of its idempotent endomorphisms,
then all objects of the DG\+category $\sDG^\natural$ are contractible,
while the DG\+category $\sDG^{\natural\natural}$ is again equivalent
to the DG\+category of complexes over~$\sA$.

\subsection{Coderived and contraderived categories}
\label{cdg-modules-co-contra-derived}
 Let $B$ be a CDG\+ring.
 Then the category $Z^0\sDG(B\modl)$ of left CDG\+modules and
closed morphisms between them is an abelian category, so one can
speak about exact triples of CDG\+modules.
 We presume that morphisms constituting an exact triple are closed.
 An exact triple of CDG\+modules can be also viewed as a finite
complex of CDG\+modules, so its total CDG\+module can be
assigned to it.

 A left CDG\+module $L$ over~$B$ is called \emph{absolutely acyclic}
if it belongs to the minimal thick subcategory of the homotopy
category $\Hot(B\modl)$ containing the total CDG\+modules of exact
triples of left CDG\+modules over~$B$.
 The thick subcategory of absolutely acyclic CDG\+modules is denoted
by $\Acycl^\abs(B\modl)\subset\Hot(B\modl)$.
 The quotient category $\sD^\abs(B\modl)=\Hot(B\modl)/
\Acycl^\abs(B\modl)$ is called the \emph{absolute derived category}
of left CDG\+modules over~$B$.

 The thick subcategory $\Acycl^\abs(B\modl)$ is often too small and
some ways of enlarging it turn out to be useful.
 A left CDG\+module over $B$ is called \emph{coacyclic} if it
belongs to the minimal triangulated subcategory of the homotopy
category $\Hot(B\modl)$ containing the total CDG\+modules of exact
triples of left CDG\+modules over~$B$ and closed under infinite
direct sums.
 The coacyclic CDG\+modules form a thick subcategory of the homotopy
category, since a triangulated category with infinite direct sums
contains images of its idempotent endomorphisms~\cite{BN}.
 This thick subcategory is denoted by $\Acycl^\co(B\modl)\subset
\Hot(B\modl)$.
 It is the minimal thick subcategory of $\Hot(B\modl)$ containing
$\Acycl^\abs(B\modl)$ and closed under infinite direct sums.
 The \emph{coderived category} of left CDG\+modules over $B$ is
defined as the quotient category $\sD^\co(B\modl)=\Hot(B\modl)/
\Acycl^\co(B\modl)$.

 Analogously, a left CDG\+module over $B$ is called
\emph{contraacyclic} if it belongs to the minimal triangulated
subcategory of $\Hot(B\modl)$ containing the total CDG\+modules
of exact triples of left CDG\+modules over~$B$ and closed under
infinite products.
 The thick subcategory formed by all contraacyclic CDG\+modules
is denoted by $\Acycl^\ctr(B\modl)\subset\Hot(B\modl)$.
 It is the minimal thick subcategory of $\Hot(B\modl)$ containing
$\Acycl^\abs(B\modl)$ and closed under infinite products.
 The \emph{contraderived category} of left CDG\+modules over $B$ is
defined as the quotient category $\sD^\ctr(B\modl)=\Hot(B\modl)/
\Acycl^\ctr(B\modl)$.

 All the above definitions can be repeated verbatim for right
CDG\+modules, so there are thick subcategories $\Acycl^\co(\modr B)$,
\ $\Acycl^\ctr(\modr B)$, and $\Acycl^\abs(\modr B)$ in $\Hot(\modr B)$
with the corresponding quotient categories $\sD^\co(\modr B)$,
$\sD^\ctr(\modr B)$, and $\sD^\abs(\modr B)$.

\begin{exs}
 When $B$ is a DG\+ring, the coderived and contraderived categories
of DG\+modules over $B$ still differ from the derived category of
DG\+modules and between each other, in general.
 Indeed, they can even all differ when $B$ is simply a ring considered
as a DG\+ring concentrated in degree~$0$.
 For example, let $\Lambda=k[\eps]/\eps^2$ be the exterior algebra in
one variable over a field~$k$.
 Then there is an infinite in both directions, acyclic, noncontractible
complex of free and cofree $\Lambda$\+modules $\dsb\rarrow \Lambda
\rarrow\Lambda\rarrow\dsb$, where the differentials are given by
the action of~$\eps$.
 This complex of $\Lambda$\+modules is neither coacyclic, nor
contraacyclic.
 Furthermore, let $\dsb\rarrow\Lambda\rarrow\Lambda\rarrow k\rarrow0$
and $0\rarrow k\rarrow \Lambda\rarrow \Lambda\rarrow\dsb$ be
the complexes of canonical truncation of the above doubly infinite
complex.
 Then the former of these two complexes is contraacyclic and the latter
is coacyclic, but not the other way.
 There also exist finite-dimensional DG\+modules (over
finite-dimensional DG\+algebras over fields) that are acyclic, but
neither coacyclic, nor contraacyclic.
 The simplest example of this kind is that of the DG\+algebra with
zero differential $B=k[\eps]/\eps^2$, where $\deg\eps=1$, and
the DG\+module $M$ over $B$ constructed as the free graded $B$\+module
with one homogeneous generator $m$ for which $d(m)=\eps m$.
 All these assertions follow from the results
of~\ref{dg-ring-bounded-cases}--\ref{cdg-mod-orthogonality}.
\end{exs}

\begin{rem}
 Much more generally, to define the coderived (contraderived) category,
it suffices to have a DG\+category $\sDG$ with shifts, cones, and
arbitrary infinite direct sums (products), for which the additive
category $Z^0(\sDG)$ is endowed with an exact category structure.
 Examples of such a situation include not only the categories of
CDG\+modules, but also, e.~g., the category of complexes over
an exact category~\cite{P}.
 Then one considers the total objects of exact triples in $Z^0(\sDG)$
as objects of the homotopy category $H^0(\sDG)$ and takes the quotient
category of $H^0(\sDG)$ by the minimal triangulated subcategory
containing all such objects and closed under infinite direct sums
(products).
 It may be advisable to require the class of exact triples in
$Z^0(\sDG)$ to be closed with respect to infinite direct sums
(products) when working with this construction.
 A deeper notion of an \emph{exact DG\+category} is discussed in
Remark~\ref{cdg-mod-orthogonality} below, where some results
provable in this setting are formulated.
\end{rem}

\subsection{Bounded cases}  \label{dg-ring-bounded-cases}
 Let $A$ a DG\+ring.
 Denote by $\Hot^+(A\modl)$ and $\Hot^-(A\modl)$ the homotopy
categories of DG\+modules over~$A$ bounded from below and
from above, respectively.
 That is, $M\in\Hot^+(A\modl)$ iff $M^i=0$ for all $i\ll0$
and $M\in\Hot^-(A\modl)$ iff $M^i=0$ for all $i\gg0$.
 Set $\Acycl^\pm(A\modl) = \Acycl(A\modl)\cap\Hot^\pm(A\modl)$,
and analogously for $\Acycl^{\co,\pm}(A\modl)$ and
$\Acycl^{\ctr,\pm}(A\modl)$.

 Clearly, the thick subcategories $\Acycl^\co(A\modl)$ and
$\Acycl^\ctr(A\modl)$ are contained in the thick subcategory
$\Acycl(A\modl)$ for any DG\+ring~$A$.

\begin{thm1}
 Assume that $A^i=0$ for all $i>0$. Then \par
\textup{(a)} $\Acycl^{\co,+}(A\modl)=\Acycl^+(A\modl)$ and\/
$\Acycl^{\ctr,-}(A\modl)=\Acycl^-(A\modl)$;
{\emergencystretch=0em\hfuzz=1.2pt\par}
\textup{(b)} the natural functors\/ $\Hot^\pm(A\modl)/
\Acycl^\pm(A\modl)\rarrow\sD(A\modl)$ are fully faithful; \par
\textup{(c)} the natural functors\/ $\Hot^\pm(A\modl)/
\Acycl^{\co,\pm}(A\modl)\rarrow\sD^\co(A\modl)$ and\/
$\Hot^\pm(A\modl)/\Acycl^{\ctr,\pm}(A\modl)\rarrow\sD^\ctr(A\modl)$
are fully faithful; \par
\textup{(d)} the triangulated subcategories $\Acycl^\co(A\modl)$ and
$\Acycl^\ctr(A\modl)$ generate the triangulated subcategory
$\Acycl(A\modl)$.
\end{thm1}

\begin{proof}
 For a DG\+module $M$ over $A$, denote by $\tau_{\le n} M$ 
the subcomplexes of canonical filtration of $M$ considered as
a complex of abelian groups.
 Due to the condition on~$A$, these are DG\+submodules.
 Notice that the quotient DG\+modules $\tau_{\le n+1}M/\tau_{\le n} M$
are contractible for any acyclic DG\+module $M$.
 Let $M\in\Acycl^{\co,+}(A\modl)$.
 Then one has $\tau_{\le n}M=0$ for $n$ small enough, hence
$\tau_{\le n}M$ is coacyclic for all~$n$.
 It remains to use the exact triple $\bigoplus_n\tau_{\le n}M
\rarrow \bigoplus_n\tau_{\le n} M\rarrow M$ in order to show
that $M$ is coacyclic.
 The proof of the second assertion of~(a) is analogous.
 To check (b) and~(c), it suffices to notice that whenever a DG\+module
$M$ is acyclic, coacyclic, or contraacyclic, the DG\+modules
$\tau_{\le n}M$ and $M/\tau_{\le n}M$ belong to the same class.
 The same observation allows to deduce~(d) from~(a).
\end{proof}

\begin{thm2}
 Assume that $A^i=0$ for all $i<0$, the ring $A^0$ is semisimple,
and $A^1=0$. Then \par
\textup{(a)} $\Acycl^{\co,-}(A\modl)=\Acycl^-(A\modl)$ and\/
$\Acycl^{\ctr,+}(A\modl)=\Acycl^+(A\modl)$;
{\emergencystretch=0em\hfuzz=1.2pt\par}
\textup{(b)} the natural functors\/ $\Hot^\pm(A\modl)/
\Acycl^\pm(A\modl)\rarrow\sD(A\modl)$ are fully faithful; \par
\textup{(c)} the natural functors\/ $\Hot^\pm(A\modl)/
\Acycl^{\co,\pm}(A\modl)\rarrow\sD^\co(A\modl)$ and\/
$\Hot^\pm(A\modl)/\Acycl^{\ctr,\pm}(A\modl)\rarrow\sD^\ctr(A\modl)$
are fully faithful; \par
\textup{(d)} the triangulated subcategories $\Acycl^\co(A\modl)$ and
$\Acycl^\ctr(A\modl)$ generate the triangulated subcategory
$\Acycl(A\modl)$.
\end{thm2}

\begin{proof}
 Analogous to the proof of Theorem~1, with the only change that
instead of the DG\+submodules $\tau_{\le n}M$ one uses
the (nonfunctorial) DG\+submodules $\sigma_{\ge n}M\subset M$, which
are constructed as follows.
 For any DG\+module $M$ over $A$ and an integer~$n$, choose
a complementary $A^0$\+submodule $K\subset M^n$ to the submodule
$\im(d^{n-1}\:M^{n-1}\to M^n)\subset M^n$.
 Set $(\sigma_{\ge n}M)^i=0$ for $i<n$, \ $(\sigma_{\ge n} M)^n = K$,
and $(\sigma_{\ge n}M)^i = M^i$ for $i>n$.
 Then $\sigma_{\ge n} M$ is a DG\+submodule of $M$ and the quotient
DG\+modules $\sigma_{\ge n-1}M/\sigma_{\ge n}M$ are contractible for
any acyclic DG\+module $M$ over~$A$.
\end{proof}

\begin{rem}
 The assertion~(b) of Theorem~2 holds under somewhat weaker
assumptions: the condition that $A^1=0$ can be replaced with
the condition that $d(A^0)=0$.
 The proof remains the same, except that the quotient DG\+modules
$\sigma_{\ge n-1}M/\sigma_{\ge n}M$ no longer have to be
contractible when $M$ is acyclic.
\end{rem}

\subsection{Semiorthogonality}  \label{cdg-mod-orthogonality}
 Let $B$ be a CDG\+ring.
 Denote by $\Hot(B\modl_\inj)$ the full subcategory of the homotopy
category of left CDG\+modules over $B$ formed by all the CDG\+modules
$M$ for which the graded module $M^\#$ over the graded ring $B^\#$
is injective.
 Analogously, denote by $\Hot(B\modl_\proj)$ the full subcategory
of the homotopy category $\Hot(B\modl)$ formed by all the CDG\+modules
$L$ for which the graded $B^\#$\+module $L^\#$ is projective. 

\begin{thm}
\textup{(a)} For any CDG\+modules $L\in\Acycl^\co(B\modl)$ and
$M\in\Hot(B\modl_\inj)$, the complex\/ $\Hom_B(L,M)$ is acyclic. \par
\textup{(b)} For any CDG\+modules $L\in\Hot(B\modl_\proj)$ and
$M\in\Acycl^\ctr(B\modl)$, the complex\/ $\Hom_B(L,M)$ is acyclic. 
\end{thm}

\begin{proof}
 Part~(a): since the functor $\Hom_B({-},M)$ transforms shifts,
cones, and infinite direct sums into shifts, shifted cones, and
infinite products, it suffices to consider the case when $L$ is
the total CDG\+module of an exact triple of CDG\+modules
$0\rarrow{}'\!K\rarrow K\rarrow{}''\!K\rarrow0$.
 Then $\Hom_B(L,M)$ is the total complex of the bicomplex with
three rows $0\rarrow\Hom_B({}''\!K,M)\rarrow\Hom_B(K,M)\rarrow
\Hom_B({}'K,M)\rarrow0$.
 Since $M^\#$ is an injective graded $B^\#$\+module, this
bicomplex is actually a short exact sequence of complexes, hence
its total complex is acyclic.
 The proof of part~(b) is completely analogous.
\end{proof}

\begin{rem}
 The assertions of Theorem can be extended to the much more general
setting of \emph{exact DG\+categories}.
 Namely, let $\sDG$ be a DG\+category with shifts and cones.
 The category $\sDG$ is said to be an exact DG\+category if an exact
category structure is defined on $Z^0(\sDG^\natural)$ such that
the following conditions formulated in terms of the functor $\natural$
are satisfied.
 Firstly, morphisms in $Z^0(\sDG)$ whose images are admissible
monomorphisms in $Z^0(\sDG^\natural)$ must admit cokernel morphisms
in $Z^0(\sDG^\natural)$ coming from morphisms in $Z^0(\sDG)$.
 Analogously, morphisms in $Z^0(\sDG)$ whose images are
admissible epimorphisms in $Z^0(\sDG^\natural)$ must admit kernel
morphisms in $Z^0(\sDG^\natural)$ coming from morphisms in $Z^0(\sDG)$.
 Secondly, the natural triples $Z\rarrow G^+(Z)^\natural\rarrow Z[-1]$
in $Z^0(\sDG^\natural)$ must be exact for all objects
$Z\in\sDG^\natural$.
 In this case the category $Z^0(\sDG)$ itself acquires an exact
category structure in which a triple is exact if and only if its
image is exact in $Z^0(\sDG^\natural)$; moreover, the functors $G^+$
and $G^-$ preserve and reflect exactness of triples.
 Assume further that all infinite direct sums exist in $\sDG$ and
the class of exact triples in $Z^0(\sDG^\natural)$ is closed under
infinite direct sums.
 Then the coderived category of $\sDG$ is defined as the quotient
category of the homotopy category $H^0(\sDG)$ by the minimal
triangulated subcategory containing the total objects of exact
triples in $Z^0(\sDG)$ and closed under infinite direct sums;
the objects belonging to the latter subcategory are called coacyclic.
 The complex of homomorphisms from any coacyclic object to any
object of $\sDG$ whose image is injective with respect to the exact
category $Z^0(\sDG^\natural)$ is acyclic.
 Analogously, assume that all infinite products exist in $\sDG$ and
the class of exact triples in $Z^0(\sDG^\natural)$ is closed under
infinite products.
 Then the contraderived category of $\sDG$ is defined as the quotient
category of the homotopy category $H^0(\sDG)$ by the minimal
triangulated subcategory containing the total objects of exact triples
in $Z^0(\sDG)$ and closed under infinite products; the objects
belonging to the latter subcategory are called contraacyclic.
 The complex of homomorphisms from any object of $\sDG$ whose image
is projective with respect to the exact category $Z^0(\sDG^\natural)$
into any contraacyclic object is acyclic.
\end{rem}

\subsection{Finite homological dimension case}
\label{finite-homol-dim-cdg-ring}
 Let $B$ be a CDG\+ring.
 Assume that the graded ring $B^\#$ has a finite left homological
dimension (i.~e., the homological dimension of the category of
graded left $B^\#$\+modules is finite).
 The next Theorem identifies the coderived, contraderived, and
absolute derived categories of left CDG\+modules over $B$ and
describes them in terms of projective and injective resolutions.

{\emergencystretch=0em \hfuzz=6pt
\begin{thm}
\textup{(a)} The three thick subcategories\/ $\Acycl^\co(B\modl)$,
\ $\Acycl^\ctr(B\modl)$, and\/ $\Acycl^\abs(B\modl)$ in the homotopy
category\/ $\Hot(B\modl)$ coincide. \par
\textup{(b)} The compositions of functors\/ 
$\Hot(B\modl_\inj)\rarrow\Hot(B\modl)\rarrow\sD^\abs(B\modl)$
and\/ $\Hot(B\modl_\proj)\rarrow \Hot(B\modl)\rarrow\sD^\abs(B\modl)$
are both equivalences of triangulated categories.
\end{thm}}

\begin{proof}
 We will show that the minimal triangulated subcategory of
$\Hot(B\modl)$ containing the total CDG\+modules of exact triples
of CDG\+modules and the triangulated subcategory $\Hot(B\modl)_\inj$
form a semiorthogonal decomposition of $\Hot(B\modl)$, as do
the triangulated subcategory $\Hot(B\modl)_\proj$ and the same
minimal triangulated subcategory.
 This implies both~(b) and the assertion that this minimal
triangulated subcategory is closed under infinite direct sums
and products, which is even stronger than~(a).
 By Lemma~\ref{semiorthogonal}, it suffices to construct for any
CDG\+module $M$ over~$B$ closed CDG\+module morphisms $F\rarrow M
\rarrow J$ whose cones belong to the mentioned minimal triangulated
subcategory, while the graded $B^\#$\+modules $F^\#$ and $J^\#$
are projective and injective, respectively.

 Choose a surjection $P_0\rarrow M^\#$ onto the graded $B^\#$\+module
$M^\#$ from a projective graded $B^\#$\+module $P_0$.
 For any graded left $B^\#$\+module $L$, denote by $G^+(L)$
the CDG\+module over $B$ freely generated by~$L$.
 The elements of $G^+(L)$ are formal expressions of the form $x+dy$
with $x$, $y\in L$.
 The action of $B$ and the differential~$d$ on $G^+(L)$ are given by
the obvious rules $b(x+dy) = bx - (-1)^{|b|}d(b)y + (-1)^{|b|}d(by)$
for $b\in B$ and $d(x+dy)=hy+dx$ implied by the equations in
the definition of a CDG\+module.
 There is a bijective correspondence between morphisms of graded
$B^\#$\+modules $L\rarrow M$ and closed morphisms of CDG\+modules
$G^+(L)\rarrow M$.
 There is also an exact triple of graded $B^\#$\+modules $L\rarrow
G^+(L)^\#\rarrow L[-1]$.
 So, in particular, we have a closed surjective morphism of
CDG\+modules $G^+(P_0)\rarrow M$, where the graded $B^\#$\+module
$G^+(P_0)^\#$ is projective.

 Let $K$ be the kernel of the surjective morphism $G^+(P_0)\rarrow M$
(taken in the abelian category $Z^0\sDG(B\modl)$ of CDG\+modules and
closed morphisms between them).
 Applying the same procedure to the CDG\+module $K$ in place of $M$,
we obtain the CDG\+module $G^+(P_1)$, etc.
 Since the graded left homological dimension of $B^\#$ is finite,
there exists a nonnegative integer~$d$ such that the image $Z$ of
the morphism $G^+(P_d)\rarrow G^+(P_{d-1})$ taken in the abelian
category $Z^0\sDG(B\modl)$ is projective as a graded $B$\+module.
 Set $F$ to be the total CDG\+module of the finite complex of
CDG\+modules $Z\rarrow G^+(P_{d-1})\rarrow\dsb\rarrow G^+(P_0)$.
 Then the graded $B^\#$\+module $F^\#$ is projective and the cone of
the closed morphism $F\rarrow M$, being the total CDG\+module of
a finite exact complex of CDG\+modules, is homotopy equivalent to
a CDG\+module obtained from the total CDG\+modules of exact triples
of CDG\+modules using the operation of cone repeatedly.
 
 Analogously, to obtain the desired morphism $M\rarrow J$ we start
with the construction of the CDG\+module $G^-(L)$ cofreely cogenerated
by a graded $B^\#$\+module $L$.
 Explicitly, $G^-(L)$ as a graded abelian group consists of all
formal expressions of the form $d^{-1}x+py$, where $x$, $y\in L$
and $\deg d^{-1}x=\deg x-1$, \ $\deg py=\deg y$.
 The differential on $G^-(L)$ is given by the formulas $d(d^{-1}x)
=px$ and $d(py)=d^{-1}(hy)$.
 The action of $B$ is given by the formulas $b(d^{-1}x) = (-1)^{|b|}
d^{-1}(bx)$ and $b(px) = p(bx) + d^{-1}(d(b)x)$.
 There is a bijective correspondence between morphisms of graded
$B^\#$\+modules $f\:M^\#\rarrow L$ and closed morphisms of CDG\+modules
$g\:M\rarrow G^-(L)$ which is described by the formula $g(z)=
d^{-1}(f(dz))+pf(z)$.
 There is also a natural closed isomorphism of CDG\+modules
$G^-(L)\simeq G^+(L)[1]$.
 
 Arguing as above, we construct an exact complex of CDG\+modules and
closed morphisms between them $0\rarrow M\rarrow G^-(I_0)\rarrow
G^-(I_{-1}) \rarrow\dsb$, where the graded $B^\#$\+modules
$G^-(I_n)^\#$ are injective.
 Since the graded left homological dimension of $B^\#$ is finite,
there exists a nonnegative integer~$d$ such that the image $Z$ of
the morphism $G^-(I_{-d+1})\rarrow G^-(I_{-d})$ taken in the category
$Z^0\sDG(B\modl)$ is injective as a graded $B$\+module.
 It remains to take $J$ to be the total CDG\+module of the finite
complex of CDG\+modules $G^-(I_0)\rarrow\dsb\rarrow G^-(I_{-d-1})
\rarrow Z$.

 When $B$ is a CDG\+algebra over a field or a DG\+ring, there is
an alternative proof analogous to the proof of
Theorem~\ref{cdg-coalgebra-inj-proj-resolutions} below.
\end{proof}

\begin{rem}
 The results of Theorem can be generalized to exact DG\+categories
in the following way.
 Let $\sDG$ be an exact DG\+category (as defined in
Remark~\ref{cdg-mod-orthogonality}); assume that the exact category
$Z^0(\sDG^\natural)$ has a finite homological dimension and
enough injectives.
 Then the minimal triangulated subcategory of $H^0(\sDG)$
containing the total objects of exact triples in $Z^0(\sDG)$ and
the full triangulated subcategory of $H^0(\sDG)$ consisting of all
the objects whose images are injective in $Z^0(\sDG^\natural)$
form a semiorthogonal decomposition of $H^0(\sDG)$.
 In particular, this minimal triangulated subcategory is closed
under infinite direct sums.
 It would be interesting to deduce the latter conclusion without
assuming existence of injectives, but only finite homological
dimension of the exact category $Z^0(\sDG^\natural)$ together with
existence of infinite direct sums in $\sDG$ and their exactness in
$Z^0(\sDG^\natural)$.
\end{rem}

 So we have $\sD^\co(B\modl)=\sD^\ctr(B\modl)$ for any CDG\+ring $B$
such that the graded ring $B^\#$ has a finite left homological
dimension.
 There are also other situations when the coderived and
contraderived categories of CDG\+modules over a given CDG\+ring, or
over two different CDG\+rings, are naturally equivalent.
 Several results in this direction can be found
in~\ref{gorenstein-cdg-ring-case}--\ref{finite-over-gorenstein}
and~\cite{KI}.

 For a discussion of the assumptions under which the derived and
absolute derived categories $\sD(A\modl)$ and $\sD^\abs(A\modl)$
of DG\+modules over a DG\+ring $A$ coincide,
see Corollaries~\ref{nonconilpotent-duality}
and~\ref{koszul-cogenerators}.2,
Remark~\ref{co-algebra-bar-duality}, and
subsection~\ref{cofibrant-dg-alg}.

\subsection{Noetherian case}  \label{noetherian-cdg-ring-case}
 Let $B$ be a CDG\+ring.
 We will need to consider the following condition on the graded
ring~$B^\#$:
\begin{itemize}
\item[($*$)] Any countable direct sum of injective graded left
$B^\#$\+modules has a finite injective dimension as a graded
$B^\#$\+module.
\end{itemize}
 The condition~($*$) holds, in particular, when the ring $B^\#$ is
graded left Noetherian, i.~e., satisfies the ascending chain condition
for graded left ideals.
 Indeed, a direct sum of injective graded left $B^\#$\+modules
is an injective graded $B^\#$\+module in this case.
 Of course, the condition~($*$) also holds when $B^\#$ has a finite
graded left homological dimension.
 The next Theorem provides a semiorthogonal decomposition of
the homotopy category $\Hot(B\modl)$ and describes the coderived
category $\sD^\co(B\modl)$ in terms of injective resolutions.

\begin{thm}
 Whenever the condition~\textup{($*$)} is satisfied, the composition of
functors\/ $\Hot(B\modl_\inj)\rarrow\Hot(B\modl)\rarrow\sD^\co(B\modl)$
is an equivalence of triangulated categories.
\end{thm}

\begin{proof}
 It suffices to construct for any CDG\+module $M$ a morphism
$M\rarrow J$ in the homotopy category of CDG\+modules over $B$ such
that the graded $B^\#$\+module $J^\#$ is injective and the cone of
the morphism $M\rarrow J$ is coacyclic.
 To do so, we start as in the proof of
Theorem~\ref{finite-homol-dim-cdg-ring}, constructing an exact
complex of CDG\+modules $0\rarrow M\rarrow G^-(I_0)\rarrow G^-(I_{-1})
\rarrow\dsb$ with injective graded $B^\#$\+modules $G^-(I_{-n})^\#$.
 Let $I$ be the total CDG\+module of the complex of CDG\+modules
$G^-(I_0)\rarrow G^-(I_{-1})\rarrow\dsb$ formed by taking infinite
direct sums.
 Let us show that the cone of the closed morphism $M\rarrow I$
is coacyclic.

 Indeed, there is a general fact that the total CDG\+module of an exact
complex of CDG\+modules $0\rarrow E_0\rarrow E_{-1}\rarrow\dsb$
bounded from below formed by taking infinite direct sums is coacyclic.
 To prove this, notice that our total CDG\+module $E$ is
the inductive limit of the total CDG\+modules $X_n$ of the finite
exact complexes of canonical truncation $0\rarrow E_0\rarrow\dsb
\rarrow E_{-n}\rarrow K_{-n}\rarrow 0$.
 So there is an exact triple of CDG\+modules and closed morphisms
$0\rarrow\bigoplus_n X_n\rarrow\bigoplus_n X_n\rarrow E\rarrow0$.
 Clearly, the CDG\+modules $X_n$ are coacyclic.
 Now it remains to notice either that the total CDG\+module of this
exact triple is absolutely acyclic by definition, or that this exact
sequence splits in $\sDG(B\modl)^\#$, and consequently this total
CDG\+module is even contractible.

 When the graded ring $B^\#$ is left Noetherian, the class of injective
graded $B^\#$\+modules is closed with respect to infinite direct sums,
so the graded $B^\#$\+module $I^\#$ is injective and we are done.
 In general, we need to repeat the same procedure for the second time,
using the assumption that the graded $B^\#$\+module $I^\#$ has a finite
injective dimension and arguing as in the proof of
Theorem~\ref{finite-homol-dim-cdg-ring}.
 Construct an exact complex of CDG\+modules $0\rarrow I\rarrow
G^-(J_0)\rarrow G^-(J_{-1})\rarrow\dsb$ with injective graded
$B^\#$\+modules $G^-(J_{-n})^\#$.
 There exists a nonnegative integer~$d$ such that the image $K$ of
the morphism $G^-(J_{-d+1})\rarrow G^-(J_{-d})$ taken in the abelian
category $Z^0\sDG(B\modl)$ is injective as a graded $B$\+module.
 Set $J$ to be the total CDG\+module of the finite complex of
CDG\+modules $G^-(J_0)\rarrow\dsb\rarrow G^-(J_{-d+1})\rarrow K$;
then the cone of the morphism $I\rarrow J$ is absolutely acyclic
and the graded $B^\#$\+module $J^\#$ is injective.
 The composition of closed morphisms $M\rarrow I\rarrow J$ now provides
the desired morphism of CDG\+modules.

 When $B$ is a CDG\+algebra over a field or a DG\+ring, one can
also use a version of the construction from the proof of
Theorem~\ref{cdg-coalgebra-inj-proj-resolutions} below in place of
the above construction with the functor~$G^-$.
\end{proof}

 Notice the difference between the proofs of the above Theorem and
Theorem~\ref{injective-dg-modules}.
 While injective resolutions for the derived category of DG\+modules
are constructed by totalizing bicomplexes by taking infinite products
along the diagonals, to construct injective resolutions for
the coderived category of CDG\+modules one has to totalize by taking
infinite direct sums along the diagonals.
 There is a similar difference for projective resolutions: to
construct them for the derived category, one has to totalize
bicomplexes by taking infinite direct sums (see
Theorem~\ref{projective-dg-modules}), while for the contraderived
category one has to use infinite products (see
Theorem~\ref{coherent-cdg-ring-case} below).
 Cf.\ the proof of Theorem~\ref{cdg-coalgebra-inj-proj-resolutions}.
 For a discussion of this phenomenon,
see~\ref{introduction-two-completions}--%
\ref{introduction-classical-two-kinds}.

\begin{rem}
 Let $\sDG$ be an exact DG\+category in the sense of
Remark~\ref{cdg-mod-orthogonality}.
 Assume that all infinite direct sums exist in $\sDG$ and the class of
exact triples in $Z^0(\sDG^\natural)$ is closed under infinite
direct sums.
 Whenever there are enough injectives in the exact category
$Z^0(\sDG^\natural)$ and countable direct sums of injectives have
finite injective dimensions in $Z^0(\sDG^\natural)$, the full
triangulated subcategory of the homotopy category $H^0(\sDG)$
consisting of all the objects whose images are injective in
$Z^0(\sDG^\natural)$ is equivalent to the coderived
category of~$\sDG$.
\end{rem}

\subsection{Coherent case}  \label{coherent-cdg-ring-case}
 Let $B$ be a CDG\+ring.
 Consider the following condition on the graded ring~$B^\#$:
\begin{itemize}
\item[($**$)] Any countable product of projective graded left
$B^\#$\+modules has a finite projective dimension as a graded
$B^\#$\+module.
\end{itemize}
 In particular, the condition ($**$) holds whenever the graded ring
$B^\#$ is graded right coherent (i.~e., all its finitely generated
graded right ideals are finitely presented) and all flat graded left
$B^\#$\+modules have finite projective dimensions over~$B^\#$.
 Indeed, a product of flat graded left modules over a graded right
coherent ring is flat~\cite{Ch}.

 The condition that all flat graded left modules have finite
projective dimensions is satisfied in many cases.
 This includes graded left perfect rings, for which all flat graded
left modules are projective~\cite{Bas} (so in particular
the condition~($**$) holds for graded right Artinian rings $B^\#$,
i.~e., whenever $B^\#$ satisfies the descending chain condition for
graded right ideals). 
 This also includes~\cite[Th\'eor\`eme~7.10]{GJ} all graded rings
$B^\#$ of cardinality not exceeding $\aleph_n$ for a finite
integer~$n$, and some other important cases (see~\cite[section~3]{KI}
and references therein).

 The next Theorem provides another semiorthogonal decomposition of
the homotopy category $\Hot(B\modl)$ and describes the contraderived
category $\sD^\ctr(B\modl)$ in terms of projective resolutions.

\begin{thm}
 Whenever the condition~\textup{($**$)} is satisfied, the composition
of functors\/ $\Hot(B\modl_\proj)\rarrow\Hot(B\modl)\rarrow
\sD^\ctr(B\modl)$ is an equivalence of triangulated categories.
\end{thm}

\begin{proof}
 Completely analogous to the proof of
Theorem~\ref{noetherian-cdg-ring-case}.
 Let us only point out that one can use, e.~g., the Mittag-Leffler
condition for the vanishing of the derived functor of projective
limit (of abelian groups) in order to show that the total CDG\+module
of an exact complex of CDG\+modules $\dsb\rarrow E_1\rarrow E_0\rarrow0$
bounded from above formed by taking infinite products is contraacyclic.
\end{proof}

 More generally, denote by $\Hot(B\modl_\fl)$ the full subcategory of
the homotopy category of left CDG\+modules over $B$ formed by all
the CDG\+modules $L$ for which the graded $B^\#$\+module $L^\#$ is flat.
 Assume only that the graded ring $B^\#$ is graded right coherent, so
that the class of flat graded left $B^\#$\+modules is closed under
infinite products.
 Then by the same argument and Lemma~\ref{flat-dg-resolutions} we can
conclude that the functor $\Hot(B\modl_\fl)/\Acycl^\ctr(B\modl)\cap
\Hot(B\modl_\fl)\rarrow\sD^\ctr(B\modl)$ is an equivalence of
triangulated categories.
 On the other hand, assume only that all flat graded left
$B^\#$\+modules have finite projective dimensions.
 Then an argument similar to the proof of
Theorem~\ref{finite-homol-dim-cdg-ring} shows that any contraacyclic
CDG\+module over $B$ that is flat as a graded $B$\+module is homotopy
equivalent to a CDG\+module obtained from the total CDG\+modules of
exact triples of CDG\+modules that are flat as graded $B$\+modules
using the operation of cone repeatedly.

\begin{qst}
 Let us call a left CDG\+module over a CDG\+ring $B$ completely
acyclic if it belongs to the minimal triangulated subcategory
of $\Hot(B\modl)$ containing the total CDG\+modules of exact
triples of CDG\+modules and closed under both infinite direct
sums and infinite products.
 The quotient category of $\Hot(B\modl)$ by the thick subcategory
of completely acyclic CDG\+modules can be called the complete
derived category of left CDG\+modules over~$B$.
 It was noticed in~\cite{KLN} that the complex $\Hom_B(L,M)$ is
acyclic for any completely acyclic left CDG\+module $M$ over $B$
and any left CDG\+module $L$ over $B$ for which the graded
$B^\#$\+module $L^\#$ is projective and finitely generated.
 Therefore, the minimal triangulated subcategory of $\Hot(B\modl)$
containing all the left CDG\+modules over $B$ that are projective and
finitely generated as graded $B$\+modules and closed under infinite
direct sums is equivalent to a full triangulated subcategory of
the complete derived category of left CDG\+modules.
 Under what assumptions on a CDG\+ring $B$ does this full triangulated
subcategory coincide with the whole complete derived category?
 Since CDG\+modules $L$ over $B$ for which the graded $B^\#$\+module
$L^\#$ is (projective and) finitely generated are compact objects
of $\Hot(B\modl)$, this question is equivalent to the following one.
 Under what assumptions on a CDG\+ring $B$ is every left CDG\+module
$M$ over $B$ such that the complex $\Hom_B(L,M)$ is acyclic for every
$L$ as above completely acyclic?
 CDG\+modules $L$ are compact objects of the complete derived category;
the question is whether they are its compact generators.
 (Cf.~\cite{Blk}.)
 The answer to this question is certainly positive when $B$ is
a DG\+ring and the complete derived category of left DG\+modules over
$B$ coincides with their derived category.
 This is so for all DG\+rings $B$ satisfying the conditions
of either Theorem~\ref{dg-ring-bounded-cases}.1 or
Theorem~\ref{dg-ring-bounded-cases}.2, and also in the cases listed
in~\ref{cofibrant-dg-alg}, when an even stronger assertion holds.
 Another situation when the answer is positive is that of
Corollary~\ref{koszul-cogenerators}.1.
 Finally, when the graded ring $B^\#$ is left Noetherian and has
a finite left homological dimension, CDG\+modules $L$ that are
finitely generated and projective as graded $B$\+modules are
compact generators of $\sD^\abs(B\modl)$ by the results of
subsection~\ref{fin-gen-cdg-mod}.
\end{qst}

\subsection{Gorenstein case}  \label{gorenstein-cdg-ring-case}
 Let $B$ be a CDG\+ring.
 Assume that the graded ring $B^\#$ is left Gorenstein, in the sense
that the classes of graded left $B^\#$\+modules of finite injective
dimension and finite projective dimension coincide.
 In particular, if $B^\#$ is a graded quasi-Frobenius ring, i.~e.,
the classes of projective graded left $B^\#$\+modules and injective
graded left $B^\#$\+modules coincide, then $B^\#$ is left Gorenstein.
 Of course, any graded ring of finite graded left homological
dimension is left Gorenstein.

\begin{thm}
\textup{(a)} The compositions of functors\/ 
$\Hot(B\modl_\inj)\rarrow\Hot(B\modl)\rarrow\sD^\co(B\modl)$
and\/ $\Hot(B\modl_\proj)\rarrow \Hot(B\modl)\rarrow\sD^\ctr(B\modl)$
are both equivalences of triangulated categories. \par
\textup{(b)} There is a natural equivalence of triangulated categories\/
$\sD^\co(B\modl)\simeq\sD^\ctr(B\modl)$.
\end{thm}

\begin{proof}
 Part~(a) follows from Theorems~\ref{noetherian-cdg-ring-case}--%
\ref{coherent-cdg-ring-case}, since any left Gorenstein ring satisfies
the conditions ($*$) and~($**$).
 Indeed, the projective (injective) dimensions of modules of finite
projective (injective) dimensions over a Gorenstein ring are bounded
by a constant, so any infinite direct sum of injective modules has
a finite projective dimension and any any infinite product of
projective modules has a finite injective dimension.
 To prove part~(b), denote by $\Hot(B\modl_\fpid)$ the homotopy
category of left CDG\+modules over $B$ whose underlying graded modules
have a finite projective (injective) dimension over $B^\#$, and
consider the quotient category $\sD^\abs(B\modl_\fpid)$ of
$\Hot(B\modl_\fpid)$ by the minimal thick subcategory containing
the total CDG\+modules of exact triples of CDG\+modules over $B$
whose underlying graded modules have a finite projective (injective)
dimension over~$B^\#$.
 Then the argument from the proof of
Theorem~\ref{finite-homol-dim-cdg-ring} shows that both compositions
of functors $\Hot(B\modl_\inj)\rarrow\Hot(B\modl_\fpid)\rarrow
\sD^\abs(B\modl_\fpid)$ and $\Hot(B\modl_\proj)\rarrow
\Hot(B\modl_\fpid)\rarrow\sD^\abs(B\modl_\fpid)$ are equivalences
of triangulated categories.
 Alternatively, both functors $\sD^\abs(B\modl_\fpid)\rarrow
\sD^\co(B\modl)$ and $\sD^\abs(B\modl_\fpid)\rarrow\sD^\ctr(B\modl)$
are equivalences of triangulated categories.
\end{proof}

\begin{rem}
 Let $\sDG$ be an exact DG\+category such that all infinite direct
sums and products exist in $\sDG$ and the class of exact triples
in $Z^0(\sDG^\natural)$ is closed under infinite direct sums and
products.
 Assume further that there are enough injectives and projectives 
in $Z^0(\sDG^\natural)$ and the classes of objects of finite
injective and projective dimensions coincide in this exact category.
 Then full triangulated subcategory of $H^0(\sDG)$ consisting of all
objects whose images are injective in $Z^0(\sDG^\natural)$ is equivalent
to the coderived category of $\sDG$, the full triangulated subcategory
of $H^0(\sDG)$ consisting of all objects whose images are projective
in $Z^0(\sDG^\natural)$ is equivalent to the contraderived category of
$\sDG$, and the coderived and contraderived categories of $\sDG$ are
naturally equivalent to each other.
\end{rem}

 Let $A$ be a DG\+ring for which the graded ring $A^\#$ is left
Gorenstein.
 Then the equivalence of categories $\sD^\co(A\modl)\simeq
\sD^\ctr(A\modl)$ makes a commutative diagram with the localization
functors $\sD^\co(A\modl)\rarrow\sD(A\modl)$ and $\sD^\ctr(A\modl)
\rarrow\sD(A\modl)$.
 Consequently, the localization functor $\sD^\co(A\modl)\simeq
\sD^\ctr(A\modl)\rarrow\sD(A\modl)$ has both a left and a right adjoint
functors, which are given by the projective and injective (resolutions
of) DG\+modules in the sense of~\ref{projective-dg-modules}--%
\ref{injective-dg-modules} (cf.~\cite{Kra1}).
 The thick subcategory of $\sD^\co(A\modl)\simeq\sD^\ctr(A\modl)$
annihilated by this localization functor can be thought of as
the domain of the Tate cohomology functor associated with~$A$.
 More precisely, there is the exact functor from $\sD(A\modl)$
to the subcategory annihilated by the localization functor that
assigns to an object of $\sD(A\modl)$ the cone of the natural
morphism between its images under the functors left and right
adjoint to the localization functor.
 The Tate cohomology as a functor of two arguments in $\sD(A\modl)$
is the composition of this functor with the functor $\Hom$ in
the triangulated subcategory annihilated by the localization functor.

\subsection{Finite-over-Gorenstein case} \label{finite-over-gorenstein}
 Let $D$ and $E$ be CDG\+rings.
 A \emph{CDG\+bimodule} $K$ over $D$ and $E$ is a graded
$D$\+$E$\+bimodule endowed with a differential~$d$ of degree~$1$
that is compatible with both the differentials of $D$ and $E$ and
satisfies the equation $d^2(x)=h_Dx-xh_E$ for all $x\in K$.
 Given a CDG\+bimodule $K$ over $D$ and $E$ and a left CDG\+module
$L$ over $E$, there is a natural structure of left CDG\+module over $D$
on the tensor product of graded modules $K\ot_EL$.
 The differential on $K\ot_EL$ is defined by the usual formula.
 Analogously, given a CDG\+bimodule $K$ over $D$ and $E$ and
a left CDG\+module $M$ over $D$, there is a natural structure of
left CDG\+module over $E$ on the graded module $\Hom_D(K,M)$.
 Both the graded left $E$\+module $\Hom_D(K,M)$ and the differential
on it are defined by the usual formulas.
 So a CDG\+bimodule $K$ defines a pair of adjoint functors
$K\ot_E{-}\:\Hot(E\modl)\rarrow\Hot(D\modl)$ and
$\Hom_D(K,{-})\:\Hot(D\modl)\rarrow\Hot(E\modl)$.

 Let $B\rarrow A'$ and $B\rarrow A''$ be two morphisms of
CDG\+rings and $C$ be a CDG\+bimodule over $A'$ and $A''$ with
the following properties.
 The graded ring $A'$ is a finitely generated projective graded
left module over the graded ring~$B$.
 The graded ring $A''$ is a finitely generated projective graded
right module over~$B$.
 The graded $A'$\+$A''$\+bimodule $C$ is a finitely generated
projective graded left and right $B$\+module.
 Finally, the adjoint functors $C\ot_{A''}{-}$ and
$\Hom_{A'}(C,{-})$ induce an equivalence between the categories
of graded left $A''$\+modules induced from $B$ and graded left
$A'$\+modules coinduced from $B$, transforming the graded left
$A''$\+module $A''\ot_BL$ into the graded left $A'$\+module
$\Hom_B(A',L)$ and vice versa, for any graded left $B$\+module~$L$.

 Given a morphism of CDG\+modules $B\rarrow A'$, one can recover
the CDG\+bimodule $C$ as the module of left $B$\+module
homomorphisms $A'\rarrow B$ with its natural CDG\+bimodule structure
and the CDG\+ring $A''$ as the ring opposite to the ring of left
$A'$\+module endomorphisms of $C$ with its natural CDG\+ring
structure, assuming that $C=\Hom_B(A',B)$ is a finitely generated
projective graded left $B$\+module.
 Analogously, given a morphism of CDG\+modules $B\rarrow A''$, one
can recover the CDG\+bimodule $C$ as the module of right $B$\+module
homomorphisms $A''\rarrow B$ and the CDG\+ring $A'$ as the ring
right $A''$\+module endomorphisms of $C$, assuming that
$C=\Hom_{B^\rop}(A'',B)$ is a finitely generated projective
graded right $B$\+module.

 Assume that the graded ring $B^\#$ is left Gorenstein.

\begin{thm}
\textup{(a)} The compositions of functors\/ 
$\Hot(A'\modl_\inj)\rarrow\Hot(A'\modl)\rarrow\sD^\co(A'\modl)$
and\/ $\Hot(A''\modl_\proj)\rarrow \Hot(A''\modl)\rarrow\sD^\ctr
(A''\modl)$ are both equivalences of triangulated categories. \par
\textup{(b)} There is a natural equivalence of triangulated categories\/
$\sD^\co(A'\modl)\simeq\sD^\ctr(A''\modl)$.
\end{thm}

\begin{proof}
 Given a morphism of graded rings $B^\#\rarrow A''{}^\#$ such that
$A''{}^\#$ is a finitely generated projective graded right
$B^\#$\+module, the graded ring $A''{}^\#$ satisfies
the condition~($**$) provided that the graded ring $B^\#$
satisfies~($**$).
 Analogously, given a morphism of graded rings $B^\#\rarrow A'{}^\#$
such that $A'{}^\#$ is a finitely generated projective graded left
$B^\#$\+module, the graded ring $A'{}^\#$ satisfies
the condition~($*$) provided that the graded ring $B^\#$
satisfies~($*$).
 This proves part~(a).
 To prove~(b), notice that the contraderived category
$\sD^\ctr(A''\modl)$ is equivalent to the quotient category of
the homotopy category of CDG\+modules over $A''$ which as graded
$A''{}^\#$\+modules are induced from graded $B^\#$\+modules of
finite projective (injective) dimension by its minimal
triangulated subcategory containing the total CDG\+modules of
exact triples of CDG\+modules over $A''$ which as exact triples
of graded $A''{}^\#$\+modules are induced from exact triples of
graded $B^\#$\+modules of finite projective (injective) dimension.
 The proof of this assertion follows the ideas of the proof of
Theorems~\ref{coherent-cdg-ring-case}
and~\ref{finite-homol-dim-cdg-ring}.
 Analogously, the coderived category $\sD^\co(A'\modl)$ is
equivalent to the quotient category of the homotopy category of
CDG\+modules over $A'$ which as graded $A'{}^\#$\+modules are
induced from graded $B^\#$\+modules of finite projective
(injective) dimension by its minimal triangulated subcategory
containing the total CDG\+modules of exact triples of
CDG\+modules over $A'$ which as exact triples of graded
$A'{}^\#$\+modules are induced from exact triples of graded
$B^\#$\+modules of finite projective (injective) dimension.
 Now the pair of adjoint functors related to the CDG\+bimodule
$C$ induces an equivalence between these triangulated categories.
 (Cf.~\cite[sections~5.4--5.5]{P}.)
\end{proof}

\subsection{Finitely generated CDG\+modules}  \label{fin-gen-cdg-mod}
 Let $B$ be a CDG\+ring.
 Assume that the graded ring $B^\#$ is left Noetherian.
 Denote by $\Hot(B\modl_\fg)$ the homotopy category of
left CDG\+modules over $B$ that are finitely generated as graded
$B^\#$\+modules, or equivalently, finitely generated as
CDG\+modules over~$B$.
 Let $\Acycl^\abs(B\modl_\fg)$ denote the minimal thick subcategory
of $\Hot(B\modl_\fg)$ containing the total CDG\+modules of exact
triples of finitely generated left CDG\+modules over~$B$.
 The quotient category $\sD^\abs(B\modl_\fg) =
\Hot(B\modl_\fg)/\Acycl^\abs(B\modl_\fg)$ is called the \emph{absolute
derived category of finitely generated left CDG\+modules} over~$B$.

\begin{thm1}
 The natural functor\/ $\sD^\abs(B\modl_\fg)\rarrow\sD^\co(B\modl)$
is fully faithful.
 In particular, any object of\/ $\Acycl^\co(B\modl)$ homotopy
equivalent to an object of\/ $\Hot(B\modl_\fg)$ is homotopy equivalent
to an object of\/ $\Acycl^\abs(B\modl_\fg)$.
\end{thm1}

\begin{proof}
 It suffices to show that any morphism $L\rarrow M$ in $\Hot(B\modl)$
between objects $L\in\Hot(B\modl_\fg)$ and $M\in\Acycl^\co(B\modl_\fg)$
factorizes through an object of $\Acycl^\abs(B\modl_\fg)$.
 Indeed, any left CDG\+module over $B$ that can be obtained from
total CDG\+modules of exact triples of CDG\+modules over $B$ using
the operations of cone and infinite direct sum is a filtered inductive
limit of its CDG\+submodules that can be obtained from the total
CDG\+modules of exact triples of finitely generated CDG\+modules
using the operation of cone.
\end{proof}

 It follows from Theorem~\ref{noetherian-cdg-ring-case} that
the objects of $\sD^\abs(B\modl_\fg)$ are compact in $\sD^\co(B\modl)$.
 The following result (cf.~\ref{dg-coalgebra-resolutions-proof} and
Question~\ref{coherent-cdg-ring-case}) is due to Dmitry Arinkin and
is reproduced here with his kind permission.

\begin{thm2}
 The objects of\/ $\sD^\abs(B\modl_\fg)$ form a set of compact
generators of the triangulated category\/ $\sD^\co(B\modl)$.
\end{thm2}

\begin{proof}
 Let $J$ be a left CDG\+module over $B$ such that the graded
$B^\#$\+module $J^\#$ is injective.
 Suppose that the complex $\Hom_B(L,J)$ is acyclic for any
finitely generated left CDG\+module $L$ over~$B$.
 We have to prove that $J$ is contractible.
 Apply Zorn's Lemma to the ordered set of all pairs $(M,h)$, where
$M$ is a CDG\+submodule in $J$ and $h\:M\rarrow J$ is a contracting
homolopy for the identity embedding $M\rarrow J$.
 It suffices to check that whenever $M\ne J$ there exists
a CDG\+submodule $M\subset M'\subset J$, \ $M\ne M'$ and a contracting
homotopy $h'\:M'\rarrow J$ for the identity embedding $M'\rarrow J$
that agrees with~$h$ on~$M$.
 Let $M'$ be any CDG\+submodule in $J$ properly containing $M$ such
that the quotient CDG\+module $M'/M$ is finitely generated.
 Since $J^\#$ is injective, the graded $B^\#$\+module morphism
$h\:M\rarrow J$ of degree~$-1$ can be extended to a graded
$B^\#$\+module morphism $h''\:M'\rarrow J$ of the same degree.
 Let $\iota\:M\rarrow J$ and $\iota'\:M'\rarrow J$ denote
the identity embeddings.
 The map $\iota'-d(h'')$ is a closed morphism of CDG\+modules
$M'\rarrow J$ vanishing in the restriction to~$M$, so it induces
a closed morphism of CDG\+modules $f\:M'/M\rarrow J$.
 By our assumption, there exists a contracting homotopy
$c\:M'/M\rarrow J$ for~$f$.
 Denote the related map of graded $B^\#$\+modules $M'\rarrow J$
of degree~$-1$ also by~$c$.
 Then $h'=h''+c\:M'\rarrow J$ is a contracting homotopy for~$\iota'$
extending~$h$.
\end{proof}

 When the graded ring $B^\#$ is left Noetherian and has a finite
graded left homological dimension, the homotopy category 
$\Hot(B\modl_\fgp)$ of left CDG\+modules over $B$ that are projective
and finitely generated as graded $B$\+modules is equivalent to
the absolute derived category of finitely generated left
CDG\+modules $\sD^\abs(B\modl_\fg)$ over~$B$.

\begin{ex}
 The following simple example of a $\boZ/2$\+graded CDG\+ring received
some attention in the recent literature.
 Let $A$ be a commutative ring and $w\in A$ be an element (often
required to be a nonzero-divisor).
 Set $B^0=A$, \ $B^1=0$, \ $d_B=0$, and $h_B=w$.
 Then CDG\+modules over $(B,d_B,h_B)$ that are free and finitely
generated as $\boZ/2$\+graded $A$\+modules are known classically as
``matrix factorizations''~\cite{Eis} and CDG\+modules over $B$ that
are projective and finitely generated as $\boZ/2$\+graded $A$\+modules
are known as ``B\+branes in the Landau--Ginzburg model''~\cite{Orl}.
 According to the above, the homotopy category of
CDG\+modules over $B$ that are projective and finitely generated
as $\boZ/2$\+graded $A$\+modules is equivalent to the absolute
derived category $\sD^\abs(B\modl_\fg)$ and compactly generates
the coderived category $\sD^\co(B\modl)$ whenever $A$ is a regular
Noetherian ring of finite Krull dimension. 
\end{ex}

\subsection{Tor and Ext of the second kind}  \label{second-tor-ext}
 Let $B$ be a CDG\+algebra over a commutative ring~$k$.
 Our goal is to define differential derived functors of
the second kind
\begin{align*}
 \Tor^{B,I\!I}\:\sD^\abs(\modr B)\times\sD^\abs(B\modl)&\lrarrow
k\modl^\sgr \\
 \Ext_B^{I\!I}\:\sD^\abs(B\modl)^\op\times\sD^\abs(B\modl)&\lrarrow
k\modl^\sgr
\end{align*}
and give a simple categorical interpretation of these definitions in
the finite homological dimension case.
 First we notice that there are enough projective and injective
objects in the category $Z^0\sDG(B\modl)$, and these objects remain
projective and injective in the category of graded $B^\#$\+modules.
 To construct these projectives and injectives, it suffices to
apply the functors $G^+$ and $G^-$ to projective and injective
graded $B^\#$\+modules (see the proof of
Theorem~\ref{finite-homol-dim-cdg-ring}).
 
 Let $N$ and $M$ be a right and a left CDG\+module over~$B$.
 We will consider CDG\+module resolutions $\dsb\rarrow Q_1
\rarrow Q_0\rarrow N\rarrow 0$ and $\dsb\rarrow P_1\rarrow P_0
\rarrow M\rarrow0$, i.~e., exact sequences of this form in
the abelian categories $Z^0\sDG(\modr B)$ and $Z^0\sDG(B\modl)$.
 To any such pair of resolutions we assign the total complex
$T=\Tot^\sqcap(Q_\bu\ot_B P_\bu)$ of the tricomplex $Q_n\ot_B P_m$
formed by taking infinite products along the diagonal planes.
 Whenever all the graded right $B^\#$\+modules $Q_n^\#$ are flat,
the cohomology of the total complex $T$ does not depend of
the choice of the resolution $P_\bu$ and vice versa.
 Indeed, whenever all graded modules $Q_n^\#$ are flat, the natural
map $\Tot^\sqcap(Q_\bu\ot_B P_\bu)\rarrow\Tor^\sqcap (Q_\bu\ot_BM)$ is
a quasi-isomorphism, since it is a quasi-isomorphism on the quotient
complexes by the components of the complete decreasing filtration
induced by the canonical filtration of the complex~$P_\bu$.

 Using existence of projective resolutions in the categories
$Z^0\sDG(\modr B)$ and $Z^0\sDG(B\modl)$, one can see that
the assignment according to which the cohomology $H(T)$ of the total
complex $T$ corresponds to a pair $(N,M)$ whenever either all
the graded right $B^\#$\+modules $Q_n^\#$, or all the graded left
$B^\#$\+modules $P_n^\#$ are flat defines a functor on the category
$Z^0\sDG(\modr B)\times Z^0\sDG(B\modl)$.
 It factorizes through the Cartesian product of the homotopy
categories, defining a triangulated functor of two variables
$H^0\sDG(\modr B)\times H^0\sDG(B\modl)\rarrow k\modl^\sgr$.
 The latter functor factorizes through the Cartesian product
of the absolute derived categories, hence the functor which
we denote by $\Tor^{B,I\!I}$.

 Analogously, let $L$ and $M$ be left CDG\+modules over~$B$.
 Consider CDG\+module resolutions $\dsb\rarrow P_1\rarrow P_0\rarrow L
\rarrow 0$ and $0\rarrow M\rarrow R_0\rarrow R_{-1}\rarrow\dsb$,
i.~e., exact sequences of this form in the category $Z^0\sDG(B\modl)$.
 To any such pair of resolutions we assign the total complex $T =
\Tot^\oplus(\Hom_B(P_\bu,R_\bu))$ of the tricomplex $\Hom_B(P_n,R_m)$
formed by taking infinite direct sums along the diagonal planes.
 The rule according to which the cohomology $H(T)$ corresponds to
a pair $(L,M)$ whenever either all the graded left $B^\#$\+modules
$P_n^\#$ are projective, or all the graded left $B^\#$\+modules
$R_m^\#$ are injective defines a functor on the category
$Z^0\sDG(B\modl)^\op\times Z^0\sDG(B\modl)$.
 This functor factorizes through the Cartesian product of the homotopy
categories, defining a triangulated functor of two variables, which
in turn factorizes through the Cartesian product of the absolute
derived categories.
 Hence the functor which we denote by $\Ext_B^{I\!I}$.
\medskip

 Alternatively, assume that $B$ is a flat graded $k$\+module, and
so is one of the graded modules $N$ and $M$.
 Then one can compute $\Tor^{B,I\!I}(N,M)$ in terms of a CDG version
of the bar-complex of $B$ with coefficients in $N$ and $M$.
 Namely, the graded complex $\dsb\rarrow N\ot_k B\ot_k B\ot_k M\rarrow
N\ot_k B\ot_k M\rarrow N\ot_kM$ is endowed with two additional
differentials, one coming from the differentials on $B$, \ $N$,
and $M$, and the other one from the curvature element of $B$
(see the proof of Theorem~\ref{cdg-coalgebra-inj-proj-resolutions},
subsections~\ref{bar-cobar-constr}--\ref{twisting-cochains-subsect},
or~\cite[Section~7 of Chapter~5]{PP}).
 Then the total complex of this ``bicomplex with three differentials'',
constructed by taking infinite products along the diagonals, computes
the desired modules $\Tor$ of the second kind.

 To check this, consider the bar-complex of $B$ with coefficients in
the complexes of CDG\+modules $Q_\bu$ and $P_\bu$, construct its total
complex by taking infinite products along the diagonal hyperplanes,
and check that the resulting complex is quasi-isomorphic to both
the total complex of $Q_\bu\ot_B P_\bu$ and the total complex of
the bar-complex of $B$ with coefficients in $N$ and~$M$.
 One proves this using the canonical filtrations of the complexes
$Q_\bu$ and $P_\bu$ and the bar-complex.

 Analogously, if $B$ is a projective graded $k$\+module, and either
$L$ is a projecitve graded $k$\+module, or $M$ is an injective graded
$k$\+module, then one can compute $\Ext_B^{I\!I}(L,M)$ in terms of
the total complex of the CDG version of the cobar-complex of $B$ with
coefficients in $L$ and $M$, constructed by taking infinite direct
sums along the diagonals.

\medskip
 Now assume that the graded ring $B^\#$ has a finite weak homological
dimension, i.~e., the homological dimension of the functor $\Tor$
between graded right and left $B^\#$\+modules is finite.
 Using Lemma~\ref{flat-dg-resolutions} and the construction from
the proof of Theorem~\ref{finite-homol-dim-cdg-ring}, one can show
that the natural functor $\Hot(B\modl_\fl)/\Acycl^\abs(B\modl)\cap
\Hot(B\modl_\fl)\rarrow\sD^\abs(B\modl)$ is an equivalence of
triangulated categories.
 To prove the analogous assertion for the homotopy category
$\Hot(\modrfl B)$ of right CDG\+modules over~$B$ that are flat as
graded $B^\#$\+modules, it suffices to pass to the opposite CDG\+ring 
$B^\rop=(B^\rop,d_{B^\rop},h_{B^\rop})$.
 The latter coincides with $B$ as a graded abelian group and has
the multiplication, differential, and curvature element defined by
the formulas $a^\rop b^\rop = (-1)^{|a||b|}(ba)^\rop$, \
$d_{B^\rop}(b^\rop) = d_B(b)^\rop$, and $h_{B^\rop} = -h_B^\rop$,
where $b^\rop$ denotes the element of $B^\rop$ corresponding to
an element $b\in B$.

 The tensor product of CDG\+modules $N\ot_B M$ over~$B$ is acyclic
whenever one of the CDG\+modules $N$ and $M$ is coacyclic and
another one is flat as a graded $B^\#$\+module.
 As in~\ref{flat-dg-resolutions}, it follows that $N\ot_B M$ is
also acyclic whenever one of the CDG\+modules $N$ and $M$ is
simultaneously coacyclic and flat as a graded $B^\#$\+module.
 So restricting the functor of tensor product over~$B$ to
either of the Cartesian products $\Hot(\modrfl B)\times\Hot(B\modl)$
or $\Hot(\modr B)\times\Hot(B\modl_\fl)$, one can construct
the derived functor of tensor product of CDG\+modules, which is defined
on the Cartesian product of absolute derived categories and factorizes
through the Cartesian product of coderived categories.
 Thus we get a functor
 $$
  \sD^\co(\modr B)\times\sD^\co(B\modl)\lrarrow k\modl^\sgr.
 $$
 This derived functor coincides with the above-defined functor
$\Tor^{B,I\!I}$, since one can use finite resolutions $P_\bu$ and
$Q_\bu$ in the construction of the latter functor in the finite
weak homological dimension case.

 Analogously, whenever the graded ring $B^\#$ has a finite left
homological dimension, the functor $\Hom_{\sD^\abs(B\modl)}(L,M)$
of homomorphisms in the absolute derived category coincides with
the above-defined functor $\Ext_B^{I\!I}(L,M)$.

\Section{Coderived Category of CDG-Comodules \\*
and Contraderived Category of CDG-Contramodules}

\subsection{CDG\+comodules and CDG\+contramodules}
\label{cdg-co-contra-modules}
 Let $k$ be a fixed ground field.
 We will consider graded vector spaces $V$ over $k$ endowed with
homogeneous endomorphisms~$d$ of degree~$1$ with \emph{not}
necessarily zero squares.
 The endomorphism~$d$ will be called ``the differential''.
 Given two graded vector spaces $V$ and $W$ with the differentials~$d$,
the differential on the graded tensor product $V\ot_kW$ is defined by
the usual formula $d(v\ot w) = d(v)\ot w + (-1)^{|v|}v\ot d(w)$ and
the differential on the graded vector space of homogeneous
homomorphisms $\Hom_k(V,W)$ is defined by the usual formula $(df)(v)=
d(f(v)) - (-1)^{|f|}f(dv)$.
 The graded vector space~$k$ is endowed with the zero differential.

 Using a version of Sweedler's notation, we will denote symbolically
the comultiplication in a graded coalgebra $C$ by $c\mpsto c_{(1)}
\ot c_{(2)}$.
 The coaction in a graded left comodule $M$ over $C$ will be denoted
by $x\mpsto x_{(-1)}\ot x_{(0)}$, while the coaction in a graded
right comodule $N$ over $C$ will be denoted by $y\mpsto y_{(0)}\ot
y_{(1)}$.
 Here $x_{(0)}\in M$, \ $y_{(0)}\in N$, and $x_{(-1)}$, $y_{(1)}\in C$.
 The contraaction map $\Hom_k(C,P)\rarrow P$ of a graded left
contramodule $P$ over~$C$ will be denoted by~$\pi_P$.

 The graded dual vector space $C^*=\Hom_k(C,k)$ to a graded
coalgebra~$C$ is a graded algebra with the multiplication given
by the formula $(\phi*\psi)(c)=\phi(c_{(2)})\psi(c_{(1)})$.
 Any graded left comodule $M$ over $C$ has a natural structure of
a graded left $C^*$\+module given by the rule $\phi*x=\phi(x_{(-1)})
x_{(0)}$, while any graded right comodule $N$ over $C$ has a natural
structure of a graded right $C^*$\+module given by
$y*\phi=(-1)^{|\phi|}\phi(y_{(1)})y_{(0)}$.
 In particular, the graded coalgebra $C$ itself is a graded
$C^*$\+bimodule.
 Any graded left contramodule $P$ over $C$ has a natural structure
of a graded left $C^*$\+module given by $\phi*p=\pi_P(c\mapsto
(-1)^{|\phi||p|}\phi(c)p)$.

 A \emph{CDG\+coalgebra} over~$k$ is a graded coalgebra~$C$ endowed
with a homogeneous endomorphism $d$ of degree~$1$ (with a not
necessarily zero square) and a homogeneous linear function
$h\:C\rarrow k$ of degree~$2$ (that is $h$ vanishes on all
the components of $C$ except perhaps $C^{-2}$) satisfying
the following equations.
 Firstly, the comultiplication map $C\rarrow C\ot_k C$ and the counit
map $C\rarrow k$ must commute with the differentials on $C$, \ $k$, and
$C\ot_kC$, where the latter two differentials are given by the above
rules.
 Secondly, one must have $d^2(c)=h*c-c*h$ and $h(d(c))=0$ for 
all $c\in C$.
 A homogeneous endomorphism $d$ of degree~$1$ acting on a graded
coalgebra~$C$ is called an \emph{odd coderivation} of degree~$1$ if
it satisfies the first of these two conditions.

 A morphism of CDG\+coalgebras $C\rarrow D$ is a pair $(f,a)$, where
$f\:C\rarrow D$ is a morphism of graded coalgebras and $a\:C\rarrow k$
is a homogeneous linear function of degree~$1$ such that the equations
$d_D(f(c)) = f(d_C(c)) + f(a*c) - (-1)^{|c|}f(c*a)$ and $h_D(f(c)) =
h_C(c) + a(d_C(c)) + a^2(c)$ hold for all $c\in C$.
 The composition of morphisms is defined by the rule $(g,b)\circ(f,a)
= (g\circ f\;b\circ f+a)$.
 Identity morphisms are the morphisms $(\id,0)$.
 So the \emph{category of CDG\+coalgebras} is defined.

 A \emph{left CDG\+comodule} over $C$ is a graded left $C$\+comodule $M$
endowed with a homogeneous linear endomorphism~$d$ of degree~$1$
(with a not necessarily zero square) satisfying the following equations.
 Firstly, the coaction map $M\rarrow C\ot_kM$ must commute with
the differentials on $M$ and $C\ot_kM$.
 Secondly, one must have $d^2(x)=h*x$ for all $x\in M$.
 A \emph{right CDG\+comodule} over $C$ is a graded right $C$\+comodule
$N$ endowed with a homogeneous linear endomorphism~$d$ of degree~$1$
such that the coaction map $N\rarrow N\ot_kC$ commutes with
the differentials and the equation $d^2(y)=-y*h$ holds.
 A \emph{left CDG\+contramodule} over $C$ is a graded left
$C$\+contramodule $P$ endowed with a homogeneous linear endomorphism~$d$
of degree~$1$ such that the contraaction map $\Hom_k(C,P)\rarrow P$
commutes with the differentials on $\Hom_k(C,P)$ and $P$,
and the equation $d^2(p)=h*p$ holds.
 In each of the above three situations, a homogeneous $k$\+linear
endomorphism $d\:M\rarrow M$ or $d\:N\rarrow N$ of degree~$1$ is
called an \emph{odd coderivation} of degree~$1$ \emph{compatible with}
an odd coderivation $d\:C\rarrow C$ of degree~$1$ or a homogeneous
$k$\+linear endomorphism $d\:P\rarrow P$ of degree~$1$ is called
an \emph{odd contraderivation} of degree~$1$ \emph{compatible
with} an odd coderivation $d\:C\rarrow C$ of degree~$1$ if the first
of the two conditions is satisfied.

 For any morphism of graded coalgebras $f\:C\rarrow D$ there are
restriction-of-scalars functors assigning to graded comodules and
contramodules over~$C$ graded $D$\+comodule and $D$\+contramodule
structures on the same graded vector spaces.
 Now let $f=(f,a)\:C\rarrow D$ be a morphism of CDG\+coalgebras and
$M$ be a left CDG\+comodule over~$C$.
 Then the left CDG\+comodule $R_fM$ over $D$ is defined by restricting
scalars in the graded $C$\+comodule $M$ via the morphism of graded
coalgebras~$f\:C\rarrow D$ and changing the differential~$d$ on $M$
by the rule $d'(x)=d(x)+a*x$.
 Analogously, let $N$ be a right CDG\+comodule over~$C$.
 Then the right CDG\+comodule $R_fN$ over $D$ is defined by
restricting scalars in the graded comodule $N$ and changing
the differential~$d$ on $N$ by the rule $d'(y)=d(y)-(-1)^{|y|}y*a$.
 Finally, let $P$ be a graded left CDG\+contramodule over~$C$.
 Then the left CDG\+contramodule $R^fP$ over $D$ is defined by
restricting scalars in the graded contramodule $P$ via the morphism~$f$
and changing the differential~$d$ on $P$ by the rule $d'(p)=d(p)+a*p$.

 Whenever $N$ is a right CDG\+comodule and $M$ is a left CDG\+comodule
over a CDG\+coalgebra $C$, the tensor product $N\oc_CM$ of the graded
comodules $N$ and $M$ over the graded coalgebra $C$ considered as
a subspace of the tensor product $N\ot_k M$ is preserved by
the differential of $N\ot_k M$.
 The restriction of this differential to $N\oc_C M$ has a zero square,
which makes $N\oc_C M$ a complex.
 Whenever $M$ is a left CDG\+comodule and $P$ is a left
CDG\+contramodule over a CDG\+coalgebra $C$, the graded space of
cohomomorphisms $\Cohom_C(M,P)$ is an invariant quotient space of
the graded space $\Hom_k(M,P)$ with respect to the differential on
$\Hom_k(M,P)$.
 The induced differential on $\Cohom_C(M,P)$ has a zero square, which
makes $\Cohom_C(M,P)$ a complex.

 Whenever $N$ is a right CDG\+comodule and $P$ is a left
CDG\+contramodule over a CDG\+coalgebra $C$, the contratensor product
$N\ocn_CP$ of the graded comodule $N$ and the graded contramodule $P$
over the graded coalgebra $C$ is an invariant quotient space of
the tensor product $N\ot_kP$ with respect to the differential on
$N\ot_kP$.
 The induced differential on $N\ocn_CP$ has a zero square, which
makes $N\ocn_CP$ a complex.

 For any left CDG\+comodules $L$ and $M$ over a CDG\+coalgebra $C$,
the graded vector space of homomorphisms between the graded
comodules $L$ and $M$ over the graded coalgebra $C$ considered as
a subspace of the graded space $\Hom_k(L,M)$ is preserved by
the differential on $\Hom_k(L,M)$.
 The induced differential on $\Hom_C(L,M)$ has a zero square, which
makes $\Hom_C(L,M)$ a complex.
 Differentials with zero squares on the graded vector spaces of
homomorphisms $\Hom_C(R,N)$ and $\Hom^C(P,Q)$ for right
CDG\+comodules $R$, $N$ and left CDG\+contramodules $P$, $Q$ over
a CDG\+coalgebra $C$ are constructed in the analogous way.
 These constructions define the DG\+categories $\sDG(C\comodl)$, \
$\sDG(\comodr C)$, and $\sDG(C\contra)$ of left CDG\+comodules,
right CDG\+comodules, and left CDG\+contramodules over $C$.

 For a CDG\+coalgebra $C$, we will sometimes denote by $C^\#$
the graded coalgebra $C$ considered without its differential~$d$
and linear function~$h$ (or with the zero differential and
linear function).
 For left or right CDG\+comodules, or CDG\+contramodules $M$, \ $N$,
or $P$ we will denote by $M^\#$, \ $N^\#$, and $P^\#$
the corresponding graded comodules and contramodules (or
CDG\+comodules and CDG\+contramodules with zero differentials)
over~$C^\#$.
 Notice that for $\sDG$ being the DG\+category of DG\+comodules
or DG\+contramodules over $C$, the corresponding additive category
$Z^0(\sDG^\natural)$ can be identified with the category of graded
comodules or contramodules over $C^\#$; the functor $\natural$ is
identified with the functor of forgetting the differential.

 All shifts, twists, infinite direct sums, and infinite direct products 
exist in the DG\+categories of CDG\+modules and CDG\+contramodules.
 The homotopy category of (the DG\+category of) left CDG\+comodules
over $C$ is denoted by $\Hot(C\comodl)$, the homotopy category of
right CDG\+comodules over $C$ is denoted by $\Hot(\comodr C)$, and
the homotopy category of left CDG\+contramodules over $C$ is
denoted by $\Hot(C\contra)$.
 Notice that there is no obvious way to define derived categories of
CDG\+comodules or CDG\+contramodules, as there is no notion of
cohomology of a CDG\+comodule or a CDG\+contramodule, and hence
no class of acyclic CDG\+comodules or CDG\+contramodules.

\subsection{Coderived and contraderived categories}
 The absolute derived categories of CDG\+comodules and
CDG\+contramodules, the coderived categories of CDG\+comodules,
and the contraderived categories of CDG\+contramodules are
defined in the way analogous to that for CDG\+modules.
 These are all particular cases of the general definition
sketched in Remarks~\ref{cdg-modules-co-contra-derived}
and~\ref{cdg-mod-orthogonality}.
 Let us spell out these definitions in a little more detail.
{\hbadness=1200\par}

 Let $C$ be a CDG\+coalgebra. 
 We will consider exact triples in the abelian categories
$Z^0\sDG(C\comodl)$ and $Z^0\sDG(C\contra)$, i.~e., exact triples
of left CDG\+modules or left CDG\+contramodules over $C$ and
closed morphisms between them.
 An exact triple of CDG\+comodules or CDG\+contramodules can be
viewed as a finite complex of CDG\+comodules or CDG\+contramodules,
so the total CDG\+comodule or CDG\+contramodule is defined for
such an exact triple.

 A left CDG\+comodule or left CDG\+contramodule over $C$ is called
\emph{absolutely acyclic} if it belongs to the minimal thick
subcategory of the homotopy category $\Hot(C\comodl)$ or
$\Hot(C\contra)$ containing the total CDG\+comodules or
CDG\+contramodules of exact triples of left CDG\+comodules or
left CDG\+contra\-modules over~$C$.
 The thick subcategories of absolutely acyclic CDG\+comodules and
CDG\+contramodules are denoted by $\Acycl^\abs(C\comodl)\subset
\Hot(C\comodl)$ and $\Acycl^\abs(C\contra)\subset\Hot(C\contra)$.
 The quotient categories $\sD^\abs(C\comodl)$ and $\sD^\abs(C\contra)$
of the homotopy categories of left CDG\+comodules and left
CDG\+contramodules by these thick subcategories are called
the \emph{absolute derived categories} of left CDG\+comodules and
left CDG\+contramodules over~$C$.

 A left CDG\+comodule over $C$ is called \emph{coacyclic} if it
belongs to the minimal triangulated subcategory of $\Hot(C\comodl)$
containing the total CDG\+comodules of exact triples of left
CDG\+comodules over $C$ and closed under infinite direct sums.
 The thick subcategory formed by all coacyclic CDG\+comodules is
denoted by $\Acycl^\co(C\comodl)\subset\Hot(C\comodl)$.
 The \emph{coderived category} of left CDG\+comodules over $C$ is
defined as the quotient category $\sD^\co(C\comodl)=\Hot(C\comodl)/
\Acycl^\co(C\comodl)$.

 A left CDG\+contramodule over $C$ is called \emph{contraacyclic} if
it belongs to the minimal triangulated subcategory of $\Hot(C\contra)$
containing the total CDG\+contramodules of exact triples of left
CDG\+contramodules over $C$ and closed under infinite products.
 The thick subcategory formed by all contraacyclic CDG\+contramodules
is denoted by $\Acycl^\ctr(C\contra)\subset\Hot(C\contra)$.
 The \emph{contraderived category} of left CDG\+contramodules over $C$
is defined as the quotient category $\sDG^\ctr(C\contra)=\Hot(C\contra)
/\Acycl^\ctr(C\contra)$.

 All the above definitions for left CDG\+comodules can be repeated
verbatim for right CDG\+comodules, so there are thick subcategories
$\Acycl^\co(\comodr C)$ and $\Acycl^\abs(\comodr C)$ in $\Hot(\comodr C)$
with the corresponding quotient categories $\sD^\co(\comodr C)$ and
$\sD^\abs(\comodr C)$.

\subsection{Bounded cases}  \label{dg-coalgebra-bounded-cases}
 Let $C$ be a DG\+coalgebra.
 Denote by $\Hot^+(C\comodl)$ and $\Hot^+(C\contra)$ the homotopy
categories of DG\+comodules and DG\+contramodules over $C$ bounded
from below, and denote by $\Hot^-(C\comodl)$ and $\Hot^-(C\contra)$
the homotopy categories of DG\+comodules and DG\+contramodules
over $C$ bounded from above.
 That is, $M\in\Hot^+(C\comodl)$ iff $M^i=0$ for all $i\ll0$ and
$P\in\Hot^-(C\contra)$ iff $P^i=0$ for all $i\gg0$; similarly for
$\Hot^-(C\comodl)$ and $\Hot^+(C\contra)$.
 Set $\Acycl^\pm(C\comodl)=\Acycl(C\comodl)\cap\Hot^\pm(C\comodl)$
and analogously for $\Acycl^\pm(C\contra)$; also set
$\Acycl^{\co,\pm}(C\comodl)=\Acycl^\co(C\comodl)\cap\Hot^\pm(C\comodl)$
and analogously for $\Acycl^{\ctr,\pm}(C\contra)$.

 Clearly, one has $\Acycl^\co(C\comodl)\subset\Acycl(C\comodl)$ and
$\Acycl^\ctr(C\contra)\subset\Acycl(C\contra)$ for any
DG\+coalgebra~$C$.

{\hbadness=4000
\begin{thm1}
 Assume that $C^i=0$ for $i<0$. Then \par
\textup{(a)} $\Acycl^{\co,+}(C\comodl)=\Acycl^+(C\comodl)$ and\/
$\Acycl^{\ctr,-}(C\contra)=\Acycl^-(C\contra)$; \par
\textup{(b)} the natural functors\/ $\Hot^\pm(C\comodl)/\Acycl^\pm
(C\comodl)\rarrow\sD(C\comodl)$ and\/ $\Hot^\pm(C\contra)/
\Acycl^\pm(C\contra)\rarrow\sD(C\contra)$ are fully faithful; \par
\textup{(c)} the natural functors\/ $\Hot^\pm(C\comodl)/\Acycl^{\co,\pm}
(C\comodl)\rarrow\sD^\co(C\comodl)$ and\/ $\Hot^\pm(C\contra)/
\Acycl^{\ctr,\pm}(C\contra)\rarrow\sD^\ctr(C\contra)$ are
fully faithful.
\end{thm1}}

\begin{proof}
 Analogous to the proof of Theorem~\ref{dg-ring-bounded-cases}.1.
\end{proof}

 For an ungraded coalgebra $E$, the following conditions are
equivalent: (i)~the category of left comodules over $E$ is semisimple;
(ii)~the category of right comodules over $E$ is semisimple;
(iii)~the category of left contramodules over $E$ is semisimple;
(iv)~$E$ is the sum of its cosimple subcoalgebras, where a coalgebra is
called \emph{cosimple} if it contains no nonzero proper subcoalgebras.
 A coalgebra $E$ satisfying these equivalent conditions is called
\emph{cosemisimple}~\cite[Appendix~A]{P}.

{\hbadness=4000
\begin{thm2}
 Assume that $C^i=0$ for $i>0$, the coalgebra $C^0$ is cosemisimple,
and $C^{-1}=0$.  Then \par
\textup{(a)} $\Acycl^{\co,-}(C\comodl)=\Acycl^-(C\comodl)$ and\/
$\Acycl^{\ctr,+}(C\contra)=\Acycl^+(C\contra)$; \par
\textup{(b)} the natural functors\/ $\Hot^\pm(C\comodl)/\Acycl^\pm
(C\comodl)\rarrow\sD(C\comodl)$ and\/ $\Hot^\pm(C\contra)/
\Acycl^\pm(C\contra)\rarrow\sD(C\contra)$ are fully faithful; \par
\textup{(c)} the natural functors\/ $\Hot^\pm(C\comodl)/\Acycl^{\co,\pm}
(C\comodl)\rarrow\sD^\co(C\comodl)$ and\/ $\Hot^\pm(C\contra)/
\Acycl^{\ctr,\pm}(C\contra)\rarrow\sD^\ctr(C\contra)$ are
fully faithful.
\end{thm2}}

\begin{proof}
 Analogous to the proof of Theorem~\ref{dg-ring-bounded-cases}.2.
\end{proof}

\begin{rem}
 The assertion~(b) of Theorem~2 still holds after one replaces
the condition that $C^{-1}=0$ with the weaker condition that
$d(C^{-1})=0$.
 See Remark~\ref{dg-ring-bounded-cases}.
\end{rem}

\subsection{Injective and projective resolutions}
\label{cdg-coalgebra-inj-proj-resolutions}
 Let $C$ be a CDG\+coalgebra.
 Denote by $\Hot(C\comodl_\inj)$ the full triangulated subcategory
of the homotopy category of left CDG\+comodules over $C$ consisting of
all the CDG\+comodules $M$ for which the graded comodule $M^\#$ over 
he graded coalgebra $C^\#$ is injective.
 Analogously, denote by $\Hot(C\contra_\proj)$ the full triangulated
subcategory of the homotopy category of left CDG\+contramodules over $C$
consisting of all the CDG\+contramodules $P$ for which the graded
contramodule $P^\#$ over the graded coalgebra $C^\#$ is projective.
 The next Theorem provides semiorthogonal decompositions of
the homotopy categories $\Hot(C\comodl)$ and $\Hot(C\contra)$, and
describes the coderived category $\sD^\co(C\comodl)$
and the contraderived category $\sD^\ctr(C\contra)$ in terms of
injective and projective resolutions, respectively.

{\hbadness=2000\hfuzz=3pt\emergencystretch=1em
\begin{thm}
\textup{(a)} For any CDG\+comodules $L\in\Acycl^\co(C\comodl)$ and
$M\in\Hot(C\comodl_\inj)$, the complex\/ $\Hom_C(L,M)$ is acyclic. \par
\textup{(b)} For any CDG\+contramodules $P\in\Hot(C\contra_\proj)$ and
$Q\in\Acycl^\ctr(C\contra)$, the complex\/ $\Hom^C(P,Q)$ is acyclic.
\par
\textup{(c)} The composition of functors\/ $\Hot(C\comodl_\inj)\rarrow
\Hot(C\comodl)\rarrow\sD^\co(C\comodl)$ is an equivalence of
triangulated categories. \par
\textup{(d)} The composition of functors\/ $\Hot(C\contra_\proj)
\rarrow\Hot(C\contra)\rarrow\sD^\ctr(C\contra)$ is an equivalence
of triangulated categories.
\end{thm}}

\begin{proof}
 Parts (a) and~(b) are easy; see the proof of
Theorem~\ref{cdg-mod-orthogonality}.
 Parts (c) and~(d) can be proven in the way analogous to
that of Theorems~\ref{noetherian-cdg-ring-case}
and~\ref{coherent-cdg-ring-case}, or alternatively
in the following way.
 Let us first consider the case of a DG\+coalgebra~$C$.
 For any DG\+comodule $M$ over $C$, consider the cobar resolution
$C\ot_k M\rarrow C\ot_kC\ot_kM\rarrow\dsb$.
 This is a complex of DG\+comodules over $C$ and closed morphisms
between them; denote by $J$ the total DG\+comodule of this complex
formed by taking infinite direct sums.
 Then the graded $C^\#$\+comodule $J^\#$ is injective and the cone of
the closed morphism $M\rarrow J$ is coacyclic.
 Analogously, for a DG\+contramodule $P$ over $C$ one considers the bar
resolution $\dsb\rarrow\Hom_k(C\ot_kC\;P)\rarrow\Hom_k(C,P)$ and forms
its total DG\+contramodule by taking infinite products.
 In the case of a CDG\+coalgebra $C$, the construction has to be
modified as follows.
 Let $M$ be a left CDG\+comodule over $C$; consider the graded left
$C^\#$\+comodule $J=\bigoplus_{n=0}^\infty (C^{\ot n+1}\ot_k M)[-n]$,
where the comodule structure on $C^{\ot n+1}\ot_k M$ comes from
the comodule structure on the leftmost factor $C$ and the shift of
the grading introduces the appropriate sign in the graded comodule
structure.
 The differential on $J$ is described as the sum of three components
$\d$, \ $d$, and $\delta$ given by the formulas
$\d(c_0\ot c_1\ot\dsb\ot c_n\ot x) = \mu(c_0)\ot c_1\ot\dsb\ot x
- c_0\ot \mu(c_1)\ot c_2\ot\dsb\ot x + \dsb + (-1)^{n+1}
c_0\ot c_1\ot\dsb\ot c_n\ot\lambda(x)$, where $\mu\:C\rarrow C\ot C$
and $\lambda\:M\rarrow C\ot M$ are the comultiplication and coaction
maps, $(-1)^n d(c_0\ot c_1\ot\dsb\ot c_n\ot x) = d(c_0)\ot c_1\ot
\dsb\ot x + (-1)^{|c_0|}c_0\ot d(c_1)\ot c_2\ot\dsb\ot x + \dsb +
(-1)^{|c_0|+\dsb+|c_n|}c_0\ot\dsb\ot c_n\ot d(x)$, and
$\delta(c_0\ot c_1\ot\dsb\ot c_n\ot x) = h(c_1)c_0\ot c_2\ot c_3
\ot\dsb\ot x - h(c_2)c_0\ot c_1\ot c_3\ot c_4\ot\dsb\ot x + \dsb +
(-1)^{n-1}h(c_n)c_0\ot\dsb c_{n-1}\ot x$.
 The graded $C^\#$\+comodule $J$ with the differential $\d+d+\delta$
is a CDG\+comodule over $C$; it is endowed with a closed morphism
of CDG\+comodules $M\rarrow J$.
 Denote by $\tau_{\le n} J$ the subspaces of canonical filtration
of the vector space $J$ considered as a complex with the grading~$n$
and the differential~$\d$; then $\tau_{\le n}J$ are
CDG\+subcomodules of $J$.
 The quotient CDG\+comodules $\tau_{\le n}J/\tau_{\le n-1}J$
(taken in the abelian category of CDG\+comodules and closed morphisms)
are contractible CDG\+comodules, the contracting homotopy being
given by the map inverse to the differential induced by
the differential~$\d$.
 The only exception is the CDG\+comodule $\tau_{\le 0}J$, which
is isomorphic to $M$.
 It follows that the cone of the closed morphism $M\rarrow J$
is coacyclic.
\end{proof}

\begin{rem}
 The following results, which are particular cases of
Remark~\ref{noetherian-cdg-ring-case},
generalize simultaneously the above Theorem and, to some extent,
Theorems~\ref{noetherian-cdg-ring-case}--\ref{coherent-cdg-ring-case}.
 A topological graded abelian group (with additive topology) is
defined as a graded abelian group endowed with a system of graded
subgroups closed under finite intersections and containing with any
subgroup all the larger graded subgroups; graded subgroups belonging
to the system are called open.
 To any topological graded abelian group one can assign an (ungraded)
topological abelian group by taking the projective limit of the direct
sums of all grading components of the graded quotient groups by
open graded subgroups.
 A topological graded abelian group with a graded ring structure is
called a topological graded ring if its multiplication can be extended
to a topological ring structure on the associated ungraded topological
abelian group.
 Let us restrict ourselves to separated and complete topological
graded rings $B$ where open two-sided graded ideals form a base
of the topology; these are exactly the graded projective limits
of discrete graded rings.
 Let $(B,d,h)$ be a CDG\+ring structure on $B$ such that
the differential~$d$ is continuous; one can easily check that
open two-sided graded differential ideals form a base of
the topology of~$B$ in this case, so $B$ is a projective limit of
discrete CDG\+rings.
 First assume that all discrete graded quotient rings of $B$ are
left Noetherian.
 A graded left $B$\+module is called discrete if the annihilator
of every its homogeneous element is an open left ideal in~$B$.
 Consider the DG\+category $\sDG(B\modl)$ of discrete graded left
$B$\+modules with CDG\+module structures.
 The corresponding coderived category $\sD^\co(B\modl)$ is defined
in the obvious way.
 The graded left $B$\+module of continuous homogeneous
abelian group homomorphisms from $B$ into any (discrete) injective
graded abelian group is an injective object in the category of
discrete graded left $B$\+modules.
 A discrete graded left $B$\+module $M$ is injective if and only if
for any open two-sided graded ideal $J\subset B$ the annihilator
of $J$ in $M$ is an injective graded left $B/J$\+module.
 It follows that there are enough injectives in the abelian
category of discrete graded left $B$\+modules and the class of
injectives is closed under infinite direct sums.
 For any discrete graded left $B$\+module $M$ the graded left
$B$\+module $G^-(M)$ is also discrete.
 So the category $Z^0(\sDG(B\modl)^\natural)$ can be identified
with the category of discrete graded left $B$\+modules and
the result of Remark~\ref{noetherian-cdg-ring-case} applies.
 Thus the coderived category $\sD^\co(B\modl)$ is equivalent
to the homotopy category of discrete left CDG\+modules over~$B$
that are injective as discrete graded modules.
 Now assume that all discrete graded quotient rings of $B$
are right Artinian.
 A pseudo-compact graded right module~\cite{Gab}
(see also~\cite{KelAppx}) over $B$ is a topological graded right
module where open graded submodules form a base of the topology and
all discrete quotient modules have finite length.
 A pseudo-compact right CDG\+module over $B$ is a pseudo-compact
graded right module endowed with a CDG\+module structure such that
the differential is continuous.
 By the result of Remark~\ref{noetherian-cdg-ring-case},
the contraderived category of pseudo-compact right CDG\+modules
over $B$ is equivalent to the homotopy category of pseudo-compact right
CDG\+modules that are projective as pseudo-compact graded modules.
 Furthermore, it is not difficult to define the DG\+category
$\sDG(B\contra)$ of left CDG\+contramodules over~$B$
(cf.~\cite[Remark~A.3]{P}).
 The corresponding contraderived category $\sD^\ctr(B\contra)$ is
equivalent to the homotopy category of left CDG\+contramodules
that are projective as graded contramodules.
 The key step is to show that a graded left contramodule $P$
over $B$ is projective if and only if for any open two-sided graded
ideal $J\subset B$ the maximal quotient contramodule of $P$ whose
$B$\+contramodule structure comes from a $B/J$\+(contra)module
structure is a projective $B/J$\+module.
 More generally, assume that $B$ has a countable base of the topology
and all discrete graded quotient rings of $B$ are right Noetherian.
 Define a contraflat graded left contramodule over $B$ as a graded
contramodule $P$ such that for any $J\subset B$ as above the maximal
quotient contramodule of $P$ whose $B$\+contramodule structure comes
from a $B/J$\+module structure is a flat $B/J$\+module.
 Then the contraderived category $\sD^\ctr(B\contra)$ is equivalent
to the homotopy category of left CDG\+contramodules over $B$
that are projective as graded contramodules provided that all
contraflat graded left contramodules over $B$ have finite projective
dimensions in the abelian category of graded contramodules.
\end{rem}

\subsection{Finite homological dimension case}
\label{finite-homol-dim-cdg-coalgebra}
 Let $E$ be a graded coalgebra.
 Then the homological dimensions of the categories of graded right
$E$\+comodules, graded left $E$\+comodules, and graded left
$E$\+contramodules coincide, as they coincide with the homological
dimensions of the derived functors of cotensor product and
cohomomorphisms on the abelian categories of comodules and
contramodules~\cite[section~0.2.9]{P}.
 The common value of these three homological dimensions we will
call the homological dimension of the graded coalgebra~$E$.

 Let $C$ be a CDG\+coalgebra.
 Assume that the graded coalgebra $C^\#$ has a finite homological
dimension.
 The next Theorem identifies the coderived category of
$C$\+comodules and the contraderived category of $C$\+contramodules
with the corresponding absolute derived categories.

\begin{thm}
\textup{(a)} The two thick subcategories\/ $\Acycl^\co(C\comodl)$ and\/
$\Acycl^\abs(C\comodl)$ in the homotopy category\/ $\Hot(C\comodl)$
coincide. \par
\textup{(b)} The two thick subcategories\/ $\Acycl^\ctr(C\contra)$ and\/
$\Acycl^\abs(C\contra)$ in the homotopy category\/ $\Hot(C\contra)$
coincide.
\end{thm}

\begin{proof}
 The proof is analogous to that of
Theorem~\ref{finite-homol-dim-cdg-ring} and can be based on either
appropriate versions of the constructions from the proof of that
Theorem or the constructions from the proof of
Theorem~\ref{cdg-coalgebra-inj-proj-resolutions}.
\end{proof}

\subsection{Finite-dimensional CDG\+comodules} \label{fin-dim-cdg-comod}
 Let $C$ be a CDG\+coalgebra.
 Denote by $\Hot(C\comodl_\fd)$ the homotopy category of (totally)
finite-dimensional CDG\+comodules over~$C$.
 Let $\Acycl^\abs(C\comodl_\fd)$ denote the minimal thick subcategory
of $\Hot(C\comodl_\fd)$ containing the total CDG\+comodules of
exact triples of finite-dimensional left CDG\+comodules over~$C$.
 The quotient category $\sD^\abs(C\comodl_\fd)=\Hot(C\comodl_\fd)/
\Acycl^\abs(C\comodl_\fd)$ is called the \emph{absolute derived
category of finite-dimensional left CDG\+comodules} over~$C$.

\begin{thm}
 The natural functor\/ $\sD^\abs(C\comodl_\fd)\rarrow\sD^\co(C\comodl)$
is fully faithful.
 In particular, any object of\/ $\Acycl^\co(C\comodl)$ that is homotopy
equivalent to an object of\/ $\Hot(C\comodl_\fd)$ is homotopy equivalent
to an object of\/ $\Acycl^\abs(C\comodl_\fd)$.
\end{thm}

\begin{proof}
 Analogous to the proof of Theorem~\ref{fin-gen-cdg-mod}.1.
\end{proof}

 It is explained in~\ref{dg-coalgebra-resolutions-proof} that
the objects of $\sD^\abs(C\comodl_\fd)$ are compact generators of
the triangulated category $\sD^\co(C\comodl)$.

\subsection{Cotor, Coext, and Ctrtor}  \label{cdg-cotor-coext-ctrtor}
 Let us define the derived functors
\begin{align*}
 \Cotor^C\:\sD^\co(\comodr C)\times\sD^\co(C\comodl)&\lrarrow
k\vect^\sgr \\
 \Coext_C\:\sD^\co(C\comodl)^\op\times\sD^\ctr(C\contra)&\lrarrow
k\vect^\sgr \\
 \Ctrtor^C\:\sD^\co(\comodr C)\times\sD^\ctr(C\contra)&\lrarrow
k\vect^\sgr
\end{align*}
for a CDG\+coalgebra~$C$.
 We denote by $\Hot(\comodrinj C)$ the full subcategory of
$\Hot(\comodr C)$ formed by all the right CDG\+comodules $N$ over $C$
for which the graded $C^\#$\+comodule $N^\#$ is injective.
 To check that the composition of functors $\Hot(\comodrinj C)\rarrow
\Hot(\comodr C)\rarrow\sD^\co(\comodr C)$, one can pass to
the opposite CDG\+coalgebra $C^\rop=(C^\rop,d^\rop,h^\rop)$,
which coincides with $C$ as a graded vector space and has
the comultiplication, differential, and curvature
defined by the formulas $(c^\rop)_{(1)}\ot (c^\rop)_{(2)} = 
(-1)^{|c_{(1)}||c_{(2)}|}c_{(2)}^\rop\ot c_{(1)}^\rop$, \
$d^\rop(c^\rop) = d(c)^\rop$, and $h^\rop(c^\rop) = -h(c)$,
where $c^\rop$ denotes the element of $C^\rop$ corresponding to
an element $c\in C$.

 To define the functor $\Cotor^C$, restrict the functor of cotensor
product $\oc_C\:\allowbreak\Hot(\comodr C)\times\Hot(C\comodl)\rarrow
\Hot(k\vect)$ to either of the full subcategories $\Hot(\comodrinj C)
\times\Hot(C\comodl)$ or $\Hot(\comodr C)\times\Hot(C\comodl_\inj)$.
 The functors so obtained factorize through the localization
$\sD^\co(\comodr C)\times \sD^\co(C\comodl)$ and the two induced
derived functors $\sD^\co(\comodr C)\times \sD^\co(C\comodl)\rarrow
k\vect^\sgr$ are naturally isomorphic to each other.
 Indeed, the cotensor product $N\oc_CM$ is acyclic whenever one
of the CDG\+comodules $N$ and $M$ is coacyclic and the other
is injective as a graded comodule.
 This follows from the fact that the functor of cotensor product
with an injective graded comodule sends exact triples of graded
comodules to exact triples of graded vector spaces. 
 To construct an isomorphism between the two induced derived functors,
it suffices to notice that both of them are isomorphic to the derived
functor obtained by restricting the functor $\oc_C$ to the full
subcategory $\Hot(\comodrinj C)\times\Hot(C\comodl_\inj)$.

 To define the functor $\Coext_C$, restrict the functor of
cohomomorphisms $\Cohom_C\:\allowbreak\Hot(C\comodl)^\op\times
\Hot(C\contra)\rarrow\Hot(k\vect)$ to either of the full
subcategories $\Hot(C\comodl_\inj)^\op\times\Hot(C\contra)$ or
$\Hot(C\comodl)^\op\times\Hot(C\contra_\proj)$.
 The functors so obtained factorize through the localization
$\sD^\co(C\comodl)^\op\times\sD^\ctr(C\contra)$ and the two induced
derived functors $\sD^\co(C\comodl)^\op\times\sD^\ctr(C\contra)\rarrow
k\vect^\sgr$ are naturally isomorphic.
 Indeed, the complex of cohomomorphisms $\Cohom_C(M,P)$ is acyclic
whenever either the CDG\+comodule $M$ is coacyclic and the
CDG\+contra\-module $P$ is projective as a graded contramodule, or
the CDG\+comodule $M$ is injective as a graded comodule and
the CDG\+contramodule $P$ is contraacyclic.

 To define the functor $\Ctrtor^C$, restrict the functor of contratensor
product $\ocn_C\:\allowbreak\Hot(\comodr C)\times\Hot(C\contra)\rarrow
\Hot(k\vect)$ to the full subcategory $\Hot(\comodr C)\times
\Hot(C\contra_\proj)$.
 The functor so obtained factorizes through the localization
$\sD^\co(\comodr C)\times\sD^\ctr(C\contra)$, so one obtains the desired
derived functor.
 Indeed, the contratensor product $N\ocn_C P$ is acyclic whenever
the CDG\+comodule $N$ is coacyclic and the CDG\+contramodule $P$ is
projective as a graded contramodule.
 Notice that one can only obtain the functor $\Ctrtor$ as the derived
functor in its second argument, but apparently not in its first argument,
as comodules adjusted to contratensor product most often do not exist.

{\hbadness=1500
 By Lemma~\ref{semiorthogonal}, one can also compute the functor
$\Ext_C=\Hom_{\sD^\co(C\comodl)}$ of homomorphisms in the coderived
category of left CDG\+comodules in terms of injective resolutions of
the second argument and the functor $\Ext^C = \Hom_{\sD^\ctr(C\contra)}$
of homomorphisms in the contraderived category of left CDG\+contramodules
in terms of projective resolutions of the first argument.
 Namely, one has $\Ext_C(L,M)=H(\Hom_C(L,M))$ whenever the CDG\+comodule
$M$ is injective as a graded $C^\#$\+comodule and $\Ext^C(P,Q) = 
H(\Hom^C(P,Q))$ whenever the CDG\+contramodule $P$ is projective as
a graded $C^\#$\+contramodule. \par}

 For any right CDG\+comodule $N$ over $C$ and any complex of $k$\+vector
spaces $V$ the differential on $\Hom_k(N,V)$ defined by the usual
formula provides a structure of CDG\+contramodule over $C$ on
$\Hom_k(N,V)$.
 The functor $\Hom_k({-},V)$ assigns contraacyclic CDG\+contramodules
to coacyclic CDG\+comodules $N$ and so induces a functor
$\boI\Hom_k({-},V)\:\sD^\co(\comodr C)\rarrow\sD^\ctr(C\contra)$ on
the level of coderived and contraderived categories.
 There are natural isomorphisms of functors of two arguments
$\Hom_k(\Cotor^C(N,M),H(V))\simeq\Coext_C(M\;\boI\Hom_k(N,V))$ and
$\Hom_k(\Ctrtor^C(N,P),H(V))\simeq\Ext^C(P\;\boI\Hom_k(N,V))$,
where $H(V)$ denotes the cohomology of the complex~$V$.

\subsection{Restriction and extension of scalars}
\label{cdg-coalgebra-scalars}
 Let $g\:C\rarrow D$ be a morphism of CDG\+coalgebras.
 Then any CDG\+comodule or CDG\+contramodule over $C$ can be also
considered as a CDG\+comodule or CDG\+contramodule over $D$, as
explained in~\ref{cdg-co-contra-modules}.
 This defines the restriction-of-scalars functors $R_g\:\Hot(C\comodl)
\rarrow\Hot(D\comodl)$ and $R^g\:\Hot(C\contra)\rarrow\Hot(D\contra)$.
 The functor $R_g$ has a right adjoint functor $E_g$ given by
the formula $E_g(N) = C\oc_DN$, while the functor $R^g$ has a left
adjoint functor given by the formula $E^g(Q) = \Cohom_D(C,Q)$;
to define the differentials on $E_g(N)$ and $E^g(Q)$, it is simplest
to decompose $g$ into an isomorphism of CDG\+coalgebras followed by
a morphism of CDG\+coalgebras with a vanishing linear function~$a$.

{\hbadness=1200
 The functors $R_g$ and $R^g$ obviously map coacyclic CDG\+comodules
and contraacyclic CDG\+contramodules to CDG\+comodules and
CDG\+contramodules of the same kind, and so induce functors
$\sD^\co(C\comodl)\rarrow\sD^\co(D\comodl)$ and $\sD^\co(C\contra)
\rarrow\sD^\co(D\contra)$, which we denote by $\boI R_g$ and $\boI R^g$.
 The functor $E_g$ has a right derived functor $\boR E_g$ obtained by
restricting $E_g$ to the full subcategory $\Hot(D\comodl_\inj)\subset
\Hot(D\comodl)$ and composing it with the localization functor
$\Hot(C\comodl)\rarrow\sD^\co(C\comodl)$.
 The functor $E^g$ has a left derived functor $\boL E^g$ obtained by
restricting $E^g$ to the full subcategory $\Hot(D\contra_\proj)\subset
\Hot(D\contra)$ and composing it with the localization functor
$\Hot(C\contra)\rarrow\sD^\ctr(C\contra)$.
 The functor $\boR E_g$ is right adjoint to the functor $\boI R_g$ 
and the functor $\boL E^g$ is left adjoint to the functor $\boI R^g$.
\par}

 For any two CDG\+coalgebras $E$ and $F$, a \emph{CDG\+bicomodule} $K$
over $E$ and $F$ is a graded vector space endowed with commuting
structures of a graded left $E$\+comodule and a graded right
$F$\+comodule and a differential $d$ compatible with both
the differentials in $E$ and $F$ and satisfying the equation
$d^2(x) = h_E*x - x*h_F$ for all $x\in K$.
 Notice that a CDG\+bicomodule over $E$ and $F$ has no natural
structures of a left CDG\+comodule over $E$ or right CDG\+comodule
over $F$, as the equations for $d^2$ are different for CDG\+comodules
and CDG\+bicomodules.

{\hbadness=1250
 CDG\+bicomodules over $E$ and $F$ form a DG\+category with morphisms
of CDG\+bicomodules being homogeneous linear maps satisfying
the compatibility equations for both the graded left $E$\+comodule
and graded right $F$\+comodule structures; the differential on
morphisms of CDG\+bicomodules is defined by the usual formula.
 The class of \emph{coacyclic} CDG\+bicomodules over $E$ and $F$ is
constructed in the same way as the class of coacyclic CDG\+comodules,
i.~e., one considers exact triples of CDG\+bicomodules and closed
morphisms between them, and generates the minimal triangulated
subcategory of the homotopy category of CDG\+bicomodules containing
the total CDG\+bicomodules of exact triples of CDG\+bicomodules and
closed under infinite direct sums.
 A CDG\+bicomodule over $E$ and $E$ is called simply a CDG\+bicomodule
over~$E$.
\par}

 Assume that a CDG\+coalgebra $C$ is endowed with an increasing
filtration by graded vector subspaces $F_0C\subset F_1C\subset\dsb$, \
$C=\bigcup_n F_nC$ that is compatible with the comultiplication and
the differential, that is $\mu(F_nC)\subset\sum_{p+q=n}F_pC\ot F_qC$
and $d(F_nC)\subset F_nC$, where $\mu$ denotes the comultiplication map.
 Then the associated quotient object $\gr_FC = \bigoplus_n
F_nC/F_{n-1}C$ becomes a CDG\+coalgebra with the comultiplication and
differential induced by those in $C$, and the counit and the curvature
linear function~$h$ obtained by restricting the corresponding linear
functions on $C$ to $F_0C$.
 In particular, $F_0C$ is also a CDG\+coalgebra; it is simultaneously
a CDG\+subcoalgebra of both $C$ and $\gr_FC$ and a quotient
CDG\+coalgebra of $\gr_FC$.
 The associated quotient space $\gr_FC$, in addition to a CDG\+coalgebra
structure, has a structure of CDG\+bicomodule over $F_0C$. 

 Now suppose that both CDG\+coalgebras $C$ and $D$ are endowed with
increasing filtrations $F$ as above and the morphism of CDG\+coalgebras
$g\:C\rarrow D$ preserves the filtrations.
 Moreover, let us assume that the morphism of CDG\+coalgebras
$F_0C\rarrow F_0D$ induced by~$g$ is an isomorphism and the cone
of the morphism $\gr_FC\rarrow\gr_FD$ of CDG\+bicomodules over
$F_0D$ is a coacyclic CDG\+bicomodule.

{\hbadness=4000
\begin{thm}
\textup{(a)} Assume that\/ $\gr_FC^\#$ and\/ $\gr_FD^\#$ are injective
graded right $F_0D^\#$\+comodules.
 Then the adjoint functors\/ $\boI R_g$ and\/ $\boR E_g$ are equivalences
of triangulated categories. \par
\textup{(b)} Assume that\/ $\gr_FC^\#$ and\/ $\gr_FD^\#$ are injective
graded left $F_0D^\#$\+comodules.
 Then the adjoint functors\/ $\boI R^g$ and\/ $\boL E^g$ are equivalences
of triangulated categories.
\end{thm}}

\begin{proof}
 We will prove~(a); the proof of~(b) is analogous.
 Let $N$ be a CDG\+comodule over $D$ such that the graded comodule
$N^\#$ over $D^\#$ is injective.
 We have to check that the cone of the morphism of CDG\+comodules
$R_g(E_g(N)) = C\oc_D N\rarrow N$ is a coacyclic CDG\+comodule
over~$D$.
 Introduce an increasing filtration $F$ on $N$ by the rule
$F_nN=\lambda^{-1}(F_nD\ot_kN)$, where $\lambda\:N\rarrow D\ot_kN$
denotes the left coaction map.
 The filtration $F$ is compatible with the graded comodule structure
and the differential on~$N$.
 The induced filtration $F$ on the cotensor product $C\oc_D N$ can be
obtained by the same construction applied to the CDG\+comodule
$C\oc_DN$ over~$C$.
 It suffices to check that the associated quotient object
$\gr_F\cone(C\oc_DN\to N)$ is a coacyclic CDG\+comodule over $F_0D$.
 But this associated quotient object is isomorphic to the cotensor
product $\cone(\gr_FC\to\gr_FD)\oc_{F_0D}F_0N$, so it remains to
notice that the cotensor product of the CDG\+bicomodule
$\cone(\gr_FC\to\gr_FD)$ over $F_0D$ with any left CDG\+comodule $L$
over $F_0D$ is a coacyclic left CDG\+comodule over $F_0D$ in our
assumptions.
 To check the latter, one can choose a morphism from $L$ into 
a CDG\+comodule that is injective as a graded comodule such that
the cone of that morphism of CDG\+comodules is coacyclic.
 Now let $M$ be a CDG\+comodule over~$C$.
 We have to check that the cone of the morphism $M\rarrow\boR E_g
(\boI R_g(M))$ in the coderived category of CDG\+comodules over~$C$
is trivial.
 To do so, we will need an injective resolution of $R_g(M)$ that is
natural enough, so that filtrations of $C$, \ $D$, and $M$ would
induce a filtration of the resolution in a way compatible with
the passage to the associated quotient objects.
 One can use either the construction from the proof of
Theorem~\ref{cdg-coalgebra-inj-proj-resolutions}, or a version of
the construction from the proof of Theorem~\ref{noetherian-cdg-ring-case}
with the coaction map in the role of a natural embedding of a graded
comodule into an injective graded comodule.
 Computing the object $\boR E_g(\boI R_g(M))$ in terms of such
a natural resolution $J$ of the CDG\+comodule $R_g(M)$ over~$D$,
we find out that it suffices to check that the cone of the morphism
$\gr_FM\rarrow \gr_F(C\oc_DJ)$ is a coacyclic CDG\+comodule over $F_0D$.
 But the cones of the morphisms $\gr_FM\rarrow\gr_FJ$ and
$\gr_F(C\oc_DJ)\rarrow\gr_FJ$ are coacyclic CDG\+comodules over $F_0D$,
the latter one in view of the above argument applied to
the CDG\+comodule $N=\gr_FJ$ over the CDG\+coalgebra $\gr_FD$ endowed
with a morphism of CDG\+coalgebras $\gr_FC\rarrow\gr_FD$.
\end{proof}

\Section{Comodule-Contramodule Correspondence}

\subsection{Functors $\Phi_C$ and $\Psi_C$} 
\label{functors-phi-and-psi}
 Let $C$ be a CDG\+coalgebra over a field~$k$.
 For any left CDG\+contramodule $P$ over $C$, let $\Phi_C(P)$ denote
the left CDG\+comodule over $C$ constructed in the following way.
 The underlying graded vector space of $\Phi_C(P)$ is, by
the definition, the contratensor product $C\ocn_CP$, which is a graded
quotient space of the tensor product $C\ot_kP$.
 The graded left $C$\+comodule structure and the differential on
$\Phi_C(P)$ are induced by the graded left $C$\+comodule structure and
the differential on $C\ot_kP$.
 The graded left $C$\+comodule structure on $C\ot_kP$ comes from
the graded left $C$\+comodule structure on $C$, while the differential
on $C\ot_kP$ is given by the standard rule as the differential on
the tensor product of graded vector spaces with differentials.
 It is straightforward to check that $\Phi_C(P)$ is a left
CDG\+comodule.
 
 For any left CDG\+comodule $M$ over $C$, let $\Psi_C(M)$ denote
the left CDG\+contra\-module over $C$ constructed as follows.
 The underlying graded vector space of $\Psi_C(M)$ is, by
the definition, the space of comodule homomorphisms $\Hom_C(C,M)$.
 The graded left contramodule structure and the differential on
$\Psi_C(P)$ are obtained by restricting the graded left contramodule
structure and the differential on $\Hom_k(C,M)$ to this graded vector
subspace.
 The graded left $C$\+contramodule structure on $\Hom_k(C,M)$ is
induced by the graded right $C$\+comodule structure on $C$, while
the differential on $\Hom_k(C,M)$ is given by the standard rule as
the differential on the space of homogeneous linear maps between
graded vector spaces with differentials.

 To a morphism $f\:L\rarrow M$ in the DG\+category of left
CDG\+comodules over $C$, one assigns the morphism $\Phi_C(f)$ given by
the formula $c\ot x\mpsto (-)^{|f||c|}c\ot f(x)$.
 To a morphism $f\:P\rarrow Q$ in the DG\+category of left
CDG\+contramodules over $C$, one assigns the morphism $\Psi_C(f)$ given
by the formula $g\mpsto f\circ g$.
 These rules define DG\+functors $\Phi_C\:\sDG(C\contra)\rarrow
\sDG(C\comodl)$ and $\Psi_C\:\sDG(C\comodl)\rarrow\sDG(C\contra)$.
 The isomorphism between the complexes of morphisms induced by our
standard isomorphism $\Hom_k(C\ot_k P\;M)\simeq\Hom_k(P,\Hom_k(C,M))$
makes the DG\+functor $\Phi_C$ left adjoint to the DG\+functor $\Psi_C$.
 So there are induced adjoint functors $\Hot(C\contra)\rarrow
\Hot(C\comodl)$ and $\Hot(C\comodl)\rarrow\Hot(C\contra)$, which
we also denote by $\Phi_C$ and $\Psi_C$.

\begin{ex}
 Let $C$ be the graded coalgebra over a field~$k$ for which the graded
dual algebra $C^*=k[x]$ is the algebra of polynomials in one
variable~$x$  of degree~$1$.
 Consider the category of vector spaces $V$ over~$k$ endowed with
a complete and cocomplete filtration $F$ indexed by the integers,
that is $\dsb\subset F_{i-1}V\subset F_iV\subset\dsb\subset V$ and
$\varinjlim F_iV\simeq V\simeq \varprojlim V/F_{i-1}V$.
 Then the graded vector space with the components $F_iV$ has a natural
structure of a graded contramodule over $C$ and the graded vector space
with the components $V/F_{i-1}V$ has a natural structure of a graded
comodule over $C$, for any such filtered vector space~$V$.
 These constructions define an equivalence between the categories of
(complete and cocomplete) filtered vector spaces, projective graded
$C$\+contramodules, and injective graded $C$\+comodules.
 The equivalence between the latter two categories assigning
the injective graded comodule $\bigoplus_i V/F_{i-1}V$ to
the projective graded contramodule $\bigoplus_i F_iV$ and vice versa
is given by the functors $P\mpsto C\ocn_CP$ and $M\mpsto\Hom_C(C,M)$.
 The category of filtered vector spaces has a natural exact category
structure; in fact, all exact triples in this exact category split.
 It follows from the graded version of
Theorem~\ref{acyclic-twisting-cochain} below (see also
Appendix~\ref{homogeneous-koszul-appendix}) that the derived
(or homotopy) category of filtered vector spaces is equivalent
to the derived category of the abelian category of graded modules
over the graded ring $k[\eps]/\eps^2$ with $\deg \eps=1$.
 The functor of forgetting the filtration, acting from the derived
category of filtered vector spaces to the derived category of
vector spaces, can be interpreted in terms of the derived category
of graded $k[\eps]/\eps^2$\+modules as the Tate cohomology
functor (cf.~\ref{gorenstein-cdg-ring-case}).
 Besides, the functor of passing to the associated graded space on
the derived category of filtered vector spaces corresponds to
the functor of forgetting the action of~$\eps$ on the derived
category of graded $k[\eps]/\eps^2$\+modules.
\end{ex}

\subsection{Correspondence Theorem}  \label{co-contra-corr-theorem}
 Restricting the functor $\Phi_C$ to the full triangulated subcategory
$\Hot(C\contra_\proj)\subset\Hot(C\contra)$ and composing it with
the localization functor $\Hot(C\comodl)\rarrow\sD^\co(C\comodl)$, we
obtain the left derived functor $\boL\Phi_C\:\sD^\ctr(C\contra)\rarrow
\sD^\co(C\comodl)$.
 Restricting the functor $\Psi_C$ to the full triangulated subcategory
$\Hot(C\comodl_\inj)\subset\Hot(C\comodl)$ and composing it with
the localization functor $\Hot(C\contra)\rarrow\sD^\ctr(C\contra)$,
we obtain the right derived functor $\boR\Psi_C\:\sD^\co(C\comodl)
\rarrow\sD^\ctr(C\contra)$.

\begin{thm}
 The functors $\boL\Phi_C$ and $\boR\Psi_C$ are mutually inverse 
equivalences between the coderived category $\sD^\co(C\comodl)$ and
the contraderived category $\sD^\ctr(C\contra)$.
\end{thm}

\begin{proof}
 One can easily see that the functors $\Phi_C$ and $\Psi_C$ between
the homotopy categories of CDG\+contramodules and CDG\+comodules over
$C$ map the full triangulated subcategories $\Hot(C\contra_\proj)$ and
$\Hot(C\comodl_\inj)$ into each other and their restrictions to
these subcategories are mutually inverse equivalences between them.
 More precisely, the functors $P\mpsto C\ocn_CP$ and $M\mpsto 
\Hom_C(C,M)$ transform the free graded contramodule $\Hom_k(C,V)$
into the cofree graded comodule $C\ot_kV$ and vice versa, for any
graded vector space~$V$.
\end{proof}

 In particular, let $B=(B,d,h)$ be a CDG\+algebra over a field~$k$
such that the graded algebra $B$ is (totally) finite-dimensional.
 Then there is a natural equivalence of triangulated categories
$\sD^\co(B\modl)\simeq\sD^\ctr(B\modl)$.
 Indeed, let $C=B^*$ be the graded dual vector space to~$B$ with
the CDG\+coalgebra structure defined by the formulas
$c(ab)=c_{(2)}(a)c_{(1)}(b)$, \ $d_C(c)(b)=(-1)^{|c|}c(d(b))$, and
$h_C(c)=c(h)$ for $a$, $b\in B$ and $c\in C$.
 Then the DG\+categories of left CDG\+modules over $B$, left
CDG\+comodules over $C$, and left CDG\+contramodules over $C$ are
all isomorphic, so $\sD^\co(B\modl)=\sD^\co(C\comodl)\simeq
\sD^\ctr(C\contra)=\sD^\ctr(B\modl)$.
 This is a particular case of Theorem~\ref{finite-over-gorenstein}.

\subsection{Coext and Ext, Cotor and Ctrtor}
\label{coext-and-ext-cotor-and-ctrtor}
 For any left CDG\+comodules $L$ and $M$ over a CDG\+coalgebra $C$,
there is a natural closed morphism of complexes of vector spaces
$\Cohom_C(L,\Psi_C(M))\rarrow\Hom_C(L,M)$, which is an isomorphism
whenever either of the graded left $C^\#$\+comodules $L^\#$ and $M^\#$
is injective.
 For any left CDG\+contramodules $P$ and $Q$ over $C$, there is 
a natural morphism of complexes of vector spaces $\Cohom_C(\Phi_C(P),Q)
\rarrow\Hom^C(P,Q)$, which is an isomorphism whenever either of
the graded left $C^\#$\+contramodules $P^\#$ and $Q^\#$ is projective.

 For any right CDG\+comodule $N$ and left CDG\+contramodule $P$
over $C$, there is a natural closed morphism of complexes of vector
spaces $N\ocn_C P\rarrow N\oc_C\Phi_C(P)$, which is an isomorphism
whenever either the graded right $C^\#$\+comodule $N$ is injective,
or the graded left $C^\#$\+contramodule $P$ is projective.

 It follows that there are natural isomorphisms of derived functors
of two arguments $\Ext_C(M,\boL\Phi_C(P))\simeq\Coext_C(M,P)\simeq
\Ext^C(\boR\Psi_C(M),P)$ and $\Cotor^C(N,M)\simeq
\Ctrtor^C(N,\boR\Psi_C(M))$ for $M\in\sD^\co(C\comodl)$,
\ $N\in\sD^\co(\comodr C)$, and $P\in\sD^\ctr(C\contra)$.
 In other words, the comodule-contramodule correspondence transforms
the functor $\Coext_C$ into the functors $\Ext_C$ and $\Ext^C$, and
also it transforms the functor $\Cotor^C$ into the functor $\Ctrtor^C$.

\subsection{Relation with extension of scalars}
\label{co-contra-ext-scalars}
 Let $g\:C\rarrow D$ be a CDG\+coalgebra morphism.
 For any left CDG\+comodule $N$ over $D$ such that the graded
$D^\#$\+comodule $N^\#$ is injective, there is a natural closed
isomorphism $\Psi_C(E_g(N))\simeq E^g(\Psi_D(N))$ of 
CDG\+contramodules over $C$ provided by the isomorphisms
$\Hom_C(C\;C\oc_DN)\simeq\Hom_D(C,N)\simeq\Cohom_D(C,\Hom_D(D,N))$
of graded vector spaces.
 Analogously, for any left CDG\+contramodule $Q$ over $D$ such that
the graded $D^\#$\+contramodule $Q^\#$ is projective, there is a natural
closed isomorphism $\Phi_C(E^g(Q))\simeq E_g(\Phi_D(Q))$ of
CDG\+comodules over $C$ provided by the isomorphisms
$C\ocn_C\Cohom_D(C,Q)\simeq C\ocn_DQ\simeq C\oc_D (D\ocn_DQ)$.

 Notice also that the functors $E_g$ and $E^g$ preserve the classes
of CDG\+comodules and CDG\+contramodules that are injective or
projective as graded comodules and contramodules, while the functors
$\Phi$ and $\Psi$ map these classes into each other.
 Thus we have natural isomorphisms of compositions of derived
functors $\boR\Psi_C\circ\boR E_g\simeq\boL E^g\circ\boR\Psi_D$
and $\boL\Phi_C\circ\boL E^g\simeq\boR E_g\circ\boL\Phi_D$;
in other words, the equivalences between the coderived and contraderived
categories of comodules and contramodules transform the derived functor
$\boR E_g$ into the derived functor $\boL E^g$.
 Identifying $\sD^\co(C\comodl)$ with $\sD^\ctr(C\contra)$ and
$\sD^\co(D\comodl)$ with $\sD^\ctr(D\contra)$, one can say that
the functor $\boI R_g$ is left adjoint to the functor
$\boR E_g=\boL E^g$, while the functor $\boI R^g$ is right adjoint
to the same functor $\boR E_g=\boL E^g$.

\subsection{Proof of Theorem~\ref{dg-coalgebra-resolutions}}
\label{dg-coalgebra-resolutions-proof}
 Our first aim is to show that the triangulated category
$\sD^\co(C\comodl)\simeq\sD^\ctr(C\contra)$ is compactly generated
for any CDG\+coalgebra~$C$.
 A triangulated category $\sD$ where all infinite direct sums exist is
said to be compactly generated if it contains a set of compact objects
$\sC$ (see Remark~\ref{dg-module-t-structure}.3 for the definition)
such that $\sD$ coincides with the minimal triangulated subcategory 
of $\sD$ containing $\sC$ and closed under infinite direct sums.

 We will work with the coderived category $\sD^\co(C\comodl)$.
 It follows from Theorem~\ref{cdg-coalgebra-inj-proj-resolutions} that
any finite-dimensional CDG\+comodule over $C$ represents a compact
object in $\sD^\co(C\comodl)$, since the full triangulated subcategory
$\Hot(C\comodl_\inj)\subset\Hot(C\comodl)$ is closed under infinite
direct sums.
 Let us check that any CDG\+co\-module over $C$ up to an isomorphism in
$\sD^\co(C\comodl)$ can be obtained from finite-dimensional
CDG\+comodules by iterated operations of cone and infinite direct sum.
 (Another proof of this result proceeds along the lines of
the proof of Theorem~\ref{fin-gen-cdg-mod}.2.)

 A graded coalgebra $E$ is called \emph{cosemisimple} if its
homological dimension (see~\ref{finite-homol-dim-cdg-coalgebra})
is equal to zero.
 When the abelian group in which the grading of $E$ takes values has
no torsion of the order equal to the characteristic of~$k$, a graded
coalgebra $E$ is cosemisimple if and only if it is cosemisimple as
an ungraded coalgebra.
 For any graded coalgebra $E$ there exists a unique maximal
cosemisimple graded subcoalgebra $E^\ss\subset E$, which coincides with
the maximal cosemisimple subcoalgebra of the ungraded coalgebra~$E$
in the above-mentioned assumption.
 The quotient coalgebra (without counit) $E/E^\ss$ is
\emph{conilpotent}, i.~e., for any element $e\in E/E^\ss$ the image
of $e$ under the iterated comultiplication map
$E/E^\ss\rarrow(E/E^\ss)^{\ot n}$ vanishes for $n$ large enough.
 One can easily prove these results, e.~g., using the fact that any
graded coalgebra is the union of its finite-dimensional graded
subcoalgebras together with the graded version of the structure theory
of finite-dimensional associative algebras.

 Let $E=C\sptilde$ be the graded coalgebra for which the category of
CDG\+comodules over $C$ and closed morphisms between them is equivalent
to the category of graded comodules over~$E$
(see~\cite[sections 0.4.4 and~11.2.2]{P} for an explicit
construction).
 Let $F_nE\subset E$ be the graded subspace formed by all elements
$e\in E$ whose images vanish in $(E/E^\ss)^{\ot n+1}$.
 Then $F_0E=E^\ss$, \ $E=\bigcup_n F_nE$, and the filtration $F_nE$ is
compatible with the coalgebra structure on~$E$.
 For any graded left comodule $M$ over $E$, set $F_nM$ to be the full
preimage of $F_nE\ot_kM$ under the comultiplication map
$M\rarrow E\ot_kM$.
 Then the filtrations $F$ on $E$ and $M$ are compatible with
the coaction map; in paricular, $F_nM$ are $E$\+subcomodules of $M$ and
the quotient comodules $F_nM/F_{n-1}M$ are comodules over $F_0E$.
 Since $F_0E$ is a cosemisimple graded coalgebra, the comodules
$F_nM/F_{n-1}M$ are direct sums of irreducible comodules, which are
finite-dimensional.

 Any left CDG\+comodule $M$ over $C$ can be viewed as a left graded
comodule over $E$; the above construction provides a filtration $F$ on
$M$ such that $F_nM$ are CDG\+co\-modules over $C$, the embeddings
$F_nM\rarrow M$ are closed morphisms, and the quotient CDG\+comodules
$F_nM/F_{n-1}M$ taken in the abelian category $Z^0\sDG(C\comodl)$ are
direct sums of finite-dimensional CDG\+comodules over~$C$.
 It follows that $M$ belongs to the minimal triangulated subcategory of
$\sD^\co(C\comodl)$ containing finite-dimensional CDG\+comodules and
closed under infinite direct sums.
 So we have proven that $\sD^\co(C\comodl)$ is compactly generated.
 The full subcategory formed by finite-dimensional CDG\+comodules
in $\sD^\co(C\comodl)$ is described in~\ref{fin-dim-cdg-comod}.

 Let us point out that in the similar way one can prove that
$\sD^\ctr(C\contra)$ is the minimal triangulated subcategory of
itself containing finite-dimensional CDG\+contramodules and closed
under infinite products.
 Just as for comodules, the category of CDG\+contramodules over $C$
and closed morphisms between them is equivalent to the category of
graded contramodules over~$E$.
 Even though the natural decreasing filtration $F^nP=\im
\Hom_k(E/F_{n-1}E,P)$ on a graded contramodule $P$ over~$E$ associated
with the filtration $F$ of~$E$ is not always separated, it is always
separated and complete for projective graded contramodules and their
graded subcontramodules, which is sufficient for the argument
to work~\cite[Appendix~A]{P}.

{\hbadness=1200
 Now let $C$ be a DG\+coalgebra.
 To prove Theorem~\ref{dg-coalgebra-resolutions}(a), notice that
the class of quasi-isomorphisms of DG\+comodules in the homotopy
category $\Hot(C\comodl)$ is locally small~\cite[10.4.4--5]{Wei},
hence morphisms between any given two objects in $\sD(C\comodl)$
form a set rather than a class.
 The localization functor $\sD^\co(C\comodl)\rarrow\sD(C\comodl)$
preserves infinite direct sums, since the thick subcategory of acyclic
DG\+comodules in $\sD^\co(C\comodl)$ is closed under infinite
direct sums.
 So it follows from the Brown representability theorem for
the compactly generated triangulated category
$\sD^\co(C\comodl)$~\cite{Neem2,Kra0} that the localization functor
$\sD^\co(C\comodl)\rarrow\sD(C\comodl)$ has a right adjoint functor.
 The localization functor $\Hot(C\comodl)\rarrow\sD^\co(C\comodl)$ has
a right adjoint by Theorem~\ref{cdg-coalgebra-inj-proj-resolutions},
thus the localization functor $\Hot(C\comodl)\rarrow\sD(C\comodl)$
also has a right adjoint functor.
 Then it remains to use Lemma~\ref{semiorthogonal}.
\par}

 Alternatively, notice that $\sD(C\comodl)$ is the quotient category of
$\sD^\co(C\comodl)$ by the thick subcategory which can be represented as
the kernel of the forgetful functor $\sD^\co(C\comodl)\rarrow
\sD(k\vect)$ or the kernel of the homological functor
$H\:\sD^\co(C\comodl)\rarrow k\vect^\sgr$.
 Both this forgetful functor and this homological functor preserve
infinite direct sums.
 It follows that this thick subcategory is
well-generated~\cite[section~7]{Kra} and therefore the localization
functor $\sD^\co(C\comodl)\rarrow\sD(C\comodl)$ has a right
adjoint functor.

 To prove Theorem~\ref{dg-coalgebra-resolutions}(b), one can also
notice that the class of quasi-isomorphisms of DG\+contramodules is
locally small in $\Hot(C\contra)$, the localization functor
$\sD^\ctr(C\contra)\rarrow\sD(C\contra)$ preserves infinite products,
and the covariant Brown representability~\cite{Kra0} for the compactly
generated triangulated category $\sD^\ctr(C\contra)$ implies existence
of a left adjoint functor to the localization functor
$\sD^\ctr(C\contra)\rarrow\sD(C\contra)$.
 But the following argument is more illuminating.

 Consider the object $P=\Hom_k(C,k)\in\sD^\ctr(C\contra)$.
 Notice that $\sD(C\contra)$ is the quotient category of
$\sD^\ctr(C\contra)$ by the thick subcategory of all objects $Q$
such that $\Hom_{\sD^\ctr(C\contra)}(P,Q)=0$.
 Consider the minimal triangulated subcategory of $\sD^\ctr(C\contra)$
containing $P$ and closed under infinite direct sums.
 This triangulated category is well-generated and therefore the functor
of its embedding into $\sD^\ctr(C\contra)$ has a right adjoint functor.
 It follows that the localization functor $\sD^\ctr(C\contra)\rarrow
\sD(C\contra)$ has a left adjoint functor whose image coincides with
the minimal triangulated subcategory of $\sD^\ctr(C\contra)$ containing
$P$ and closed under infinite direct sums.
 The localization functor $\Hot(C\comodl)\rarrow\sD^\ctr(C\contra)$ has
a left adjoint functor by
Theorem~\ref{cdg-coalgebra-inj-proj-resolutions}, thus the localization
functor $\Hot(C\contra)\rarrow\sD(C\contra)$ also has a left adjoint.

 In addition to the assertions of Theorem, we have proven that
the triangulated subcategory $\Hot(C\contra)_\proj$ coincides with
the minimal triangulated subcategory of $\Hot(C\contra)$ containing
the DG\+contramodule $\Hom_k(C,k)$ and closed under infinite direct
sums.
 Indeed, this is so in the triangulated category $\sD^\ctr(C\contra)$
and the triangulated subcategory $\Hot(C\contra_\proj)\subset
\Hot(C\contra)$, which is equivalent to $\sD^\ctr(C\contra)$,
contains $\Hom_k(C,k)$ and is closed under infinite direct sums.
 We do not know whether the triangulated subcategory
$\Hot(C\comodl)_\inj$ coincides with the minimal triangulated
subcategory of $\Hot(C\comodl)$ containing the DG\+comodule $C$ and
closed under infinite products; the former subcategory certainly
contains the latter one.
 Notice that the DG\+comodule $C$ and the DG\+contramodule
$\Hom_k(C,k)$ correspond to each other under the equivalence of
categories $\sD^\co(C\comodl)\simeq\sD^\ctr(C\contra)$. \qed

\medskip

\textit{Note added three years later}:
 Assuming some large cardinal axioms, one can indeed show that
$\Hot(C\comodl)_\inj$ is the minimal triangulated subcategory in
$\Hot(C\comodl)$ containing the DG\+comodule $C$ and closed under
infinite products.
 Equivalently, the coderived category $\sD^\co(C\comodl)$ has
a semiorthogonal decomposition formed by the full triangulated
subcategory of acyclic DG\+comodules and the minimal triangulated
subcategory containing the DG\+comodule $C$ and closed under
infinite products.

 It follows that the derived category $\sD(C\comodl)$ can be
described as the minimal triangulated subcategory in
$\sD^\co(C\comodl)\simeq\sD^\ctr(C\contra)$ containing
the canonical object $C\in\sD^\co(C\comodl)$, while, as we have
seen, the derived category $\sD(C\contra)$ is described as
the minimal triangulated subcategory in the same triangulated
category $\sD^\co(C\comodl)\simeq\sD^\ctr(C\contra)$ containing
the same canonical object $\Hom_k(C,k)\in\sD^\ctr(C\contra)$
and closed under infinite direct sums.
 Notice again that the coderived~$\simeq$~contraderived category is
compactly generated, but the canonical object
$C\longleftrightarrow\Hom_k(C,k)$ is not compact in it, in general.

 Indeed, a DG\+comodule $M$ over a DG\+coalgebra $C$ is acyclic
if and only if one has $\Hom_{\sD^\co(C\comodl)}(M,C)=0$.
 Hence it suffices to show that the minimal triangulated
subcategory in $\sD^\co(C\comodl)$ containing the object $C$
and closed under infinite products together with its left
orthogonal complement form a semiorthogonal decomposition
of $\sD^\co(C\comodl)$.
 By~\cite[Theorem~2.4]{CGR}, in the assumption of Vop\v enka's
principle~\cite{AR}, every triangulated subcategory closed under
infinite products in the homotopy category of a locally presentable
category with a stable model category structure has this property.
 For any CDG\+coalgebra $C$ over~$k$, a model category structure
on a locally Noetherian Grothendieck abelian category with
the homotopy category equivalent to $\sD^\co(C\comodl)$ is
provided by Theorem~\ref{cdg-co-contra-model-categories}(a).

 Moreover, this model category is cofibrantly generated with
injective morphisms of finite-dimensional CDG\+comodules being
the generating cofibrations and injective morphisms of
finite-dimensional CDG\+comodules with the cokernels belonging
to $\Acycl^\abs(C\comodl_\fd)$ in the role of generating
trivial cofibrations (see~\ref{fin-dim-cdg-comod} for
the notation, and the proof of Theorem~\ref{fin-gen-cdg-mod}.1
for a relevant argument).
 In fact, the model category of CDG\+comodules over $C$ is even
finitely generated by these generating sets of morphisms
(see~\cite{Hov} for the background definitions).
 Indeed, any injective map of CDG\+co\-modules with
a finite-dimensional cokernel is a pushout of an injective map
of finite-dimensional CDG\+comodules with the same cokernel.
 So any cofibration of CDG\+co\-modules is a transfinite composition
of our generating cofibrations and any trivial cofibration is
a transfinite composition of our generating trivial cofibrations.
 A similar and, essentially, much more general result about
the coderived (and also the contraderived) model structure on
CDG\+modules over a CDG\+ring being cofibrantly generated
can be found in~\cite[Proposition~1.3.6]{Bec}
(the proof is based on~\cite[Theorem~2.13]{SS};
cf.\ Remark~\ref{cdg-module-model-categories}).

 It would be interesting to know how much can one weaken
the set-theoretical axioms used in the above argument.
 The assertion of~\cite[Theorem~9.5]{BCMR} seems to suggest
that supercompact cardinals might suffice.
 On the other hand, some large cardinals may be
necessary~\cite[2.6\+-2.7]{CGR}.

 I am grateful to Greg Stevenson for directing my attention
to the paper~\cite{CGR}, and to Carles Casacuberta for his
explanations and comments on the results contained in this
and several other recent papers of his, including~\cite{BCMR}.
 I also wish to thank Hanno Becker for very interesting
correspondence.

\Section{Koszul Duality: Conilpotent and Nonconilpotent Cases}

\subsection{Bar and cobar constructions}  \label{bar-cobar-constr}
 Let $B=(B,d,h)$ be a CDG\+algebra over a field~$k$.
 We assume that $B$ is nonzero, i.~e., the unit element $1\in B$ is not
equal to~$0$, and consider $k=k\cdot 1$ as a graded vector subspace
in~$B$.
 Let $v\:B\rarrow k$ be a homogeneous $k$\+linear retraction of
the graded vector space $B$ to its subspace $k$;
set $V=\ker v\subset B$.
 The direct sum decomposition $B=V\oplus k$ allows one to split
the multiplication map $m\:V\ot_k V\rarrow B$, the differential
$d\:V\rarrow B$, and the curvature element $h\in B$ into
the components $m=(m_V,m_k)$, \ $d=(d_V,d_k)$, and $h=(h_V,h_k)$,
where $m_V\:V\ot_kV\rarrow V$, \ $m_k\:V\ot_kV\rarrow k$, \ 
$d_V\:V\rarrow V$, \ $d_k\:V\rarrow k$, \ $h_V\in V$, and $h_k\in k$.
 Notice that the restrictions of the multiplication map and
the differential to $k\ot_kV$, \ $V\ot_kk$, \ $k\ot_kk$, and $k$
are uniquely determined by the axioms of a graded algebra and its
derivation.
 One has $h_k=0$ for the dimension reasons when $B$ is
$\boZ$\+graded, but $h_k$ may be nonzero when $B$ is $\boZ/2$\+graded
(see Remark~\ref{dg-rings-modules}).

 Set $B_+=B/k$. 
 Let $\Br(B)=\bigoplus_{n=0}^\infty B_+^{\ot n}$ be the tensor
coalgebra generated by the graded vector space~$B_+$.
 The comultiplication in $\Br(B)$ is given by the rule $b_1\ot\dsb\ot
b_n \mpsto \sum_{j=0}^n (b_1\ot\dsb\ot b_j)\ot (b_{j+1}\ot\dsb\ot b_n)$
and the counit is the projection to the component $B_+^{\ot 0}
\simeq k$.
 The coalgebra $\Br(B)$ is a graded coalgebra with the grading given
by the rule $\deg(b_1\ot\dsb\ot b_n)=\deg(b_1)+\dsb+\deg(b_n)-n$.

 Odd coderivations of degree~$1$ on $\Br(B)$ are determined by their
compositions with the projection of $\Br(B)$ to the component
$B_+^{\ot 1}\simeq B_+$; conversely, any linear map $\Br(B)\to B_+$ of
degree~$2$ gives rise to an odd coderivation of degree~$1$ on $\Br(B)$.
 Let $d_{\Br}$ be odd coderivation of degree~$1$ on $\Br(B)$ whose
compositions with the projection $\Br(B)\rarrow B_+$ are given by
the rules $b_1\ot\dsb\ot b_n\mpsto 0$ for $n\ge3$, \
$b_1\ot b_2\mpsto (-1)^{|b_1|+1}m_V(b_1\ot b_2)$, \ $b\mpsto -d_V(b)$,
and $1\mpsto h_V$, where $B_+$ is identified with $V$ and 
$1\in B_+^{\ot0}$.
 Let $h_{\Br}\:\Br(B)\rarrow k$ be the linear function given by
the formulas $h_{\Br}(b_1\ot\dsb\ot b_n)=0$ for $n\ge 3$, \
$h_{\Br}(b_1\ot b_2)=(-1)^{|b_1|+1}h_k(b_1\ot b_2)$, \
$h_{\Br}(b)=-d_k(b)$, and $h_{\Br}(1)=h_k$.
 Then $\Br_v(B)=(\Br(B),d_{\Br},h_{\Br})$ is a CDG\+coalgebra over~$k$.
 The CDG\+coalgebra $\Br_v(B)$ is called the \emph{bar-construction}
of a CDG\+algebra $B$ endowed with a homogeneous $k$\+linear retraction
$v\:B\rarrow k$.

 A retraction $v\:B\rarrow k$ is called an \emph{augmentation}
of a CDG\+algebra $B$ if $(v,0)\: (B,d,h)\rarrow (k,0,0)$ is a morphism
of CDG\+algebras; equivalently, $v$ is an augmentation if it is
a morphism of graded algebras satisfying the equations $v(d(b))=0$ and
$v(h)=0$.
 A $k$\+linear retraction $v$ is an augmentation if and only if
the CDG\+coalgebra $\Br_v(B)$ is actually a DG\+coalgebra, i.~e.,
$h_{\Br}=0$.

 Let $C=(C,d,h)$ be a CDG\+coalgebra over~$k$.
 We assume that $C$ is nonzero, i.~e., the counit map $\eps\:C\rarrow k$
is a nonzero linear function.
 Let $w\:k\rarrow C$ be a homogeneous $k$\+linear section of
the surjective map of graded vector spaces~$\eps$; set $W=\coker w$.
 The direct sum decomposition $C=W\oplus k$ allows one to split
the comultiplication map $\mu\:C\rarrow W\ot_kW$, the differential
$d\:C\rarrow W$, and the curvature linear function $h\:C\rarrow k$
into the components $\mu=(\mu_W,\mu_k)$, \ $d=(d_W,d_k)$, and
$h=(h_W,h_k)$, where $\mu_W\:W\rarrow W\ot_kW$, \ $\mu_k\in W\ot_kW$,
\ $d_W\:W\rarrow W$, \ $d_k\in W$, \ $h_W\:W\rarrow k$, and $h_k\in k$.
 Notice that the compositions of the comultiplication map with
the projections $C\ot_kC\rarrow k\ot_kW$, $W\ot_kk$, $k\ot_kk$ and
the composition of the differential with the projection (counit)
$C\rarrow k$ are uniquely detemined by the axioms of a graded coalgebra
and a differential compatible with the coalgebra structure.
 One has $h_k=0$ for dimension reasons when $B$ is $\boZ$\+graded, but
$h_k$ may be nonzero when $B$ is $\boZ/2$\+graded.

 Set $C_+=\ker \eps$.
 Let $\Cb(C)=\bigoplus_{n=0}^\infty C_+^{\ot n}$ be the tensor (free
associative) algebra, generated by the graded vector space~$C_+$.
 The multiplication in $\Cb(C)$ is given by the rule $(c_1\ot\dsb\ot c_j)
(c_{j+1}\ot\dsb\ot c_n)=c_1\ot\dsb\ot c_n$ and the unit element is
$1\in k\simeq C_+^{\ot 0}$.
 The algebra $\Cb(C)$ is a graded algebra with the grading given by
the rule $\deg(c_1\ot\dsb\ot c_n)=\deg c_1+\dsb+\deg c_n+n$.

 Odd derivations of degree~$1$ on $\Cb(C)$ are determined by their
restrictions to the component $C_+\simeq C_+^{\ot1}\subset\Cb(C)$; 
conversely, any linear map $C_+\rarrow\Cb(C)$ of degree~$2$ gives
rise to an odd derivation of degree~$1$ on $\Cb(C)$.
 Let $d_{\Cb}$ be the odd derivation on $\Cb(C)$ whose restriction
to $C_+$ is given by the formula $d(c) = (-1)^{|c_{(1,W)}|+1}c_{(1,W)}
\ot c_{(2,W)} - d_W(c) + h_W(c)$, where $C_+$ is identified with $W$
and $\mu_W(c)=c_{(1,W)}\ot c_{(2,W)}$.
 Let $h_{\Cb}\in\Cb(C)$ be the element given by the formula
$h_{\Cb} = (-1)^{|\mu_{(1,k)}|+1}\mu_{(1,k)}\ot\mu_{(2,k)} - d_k + h_k$,
where $\mu_k=\mu_{(1,k)}\ot\mu_{(2,k)}$.
 Then $\Cb_w(C)=(\Cb(C),d_{\Cb},h_{\Cb})$ is a CDG\+algebra over~$k$.
 The CDG\+algebra $\Cb_w(C)$ is called the \emph{cobar-construction}
of a CDG\+coalgebra $C$ endowed with a homogeneous $k$\+linear
section $w\:k\rarrow C$ of the counit map $\eps\:C\rarrow k$.

 A section $w\:k\rarrow C$ is called a \emph{coaugmentation} of
a CDG\+coalgebra $C$ if $(w,0)\:(k,0,0)\rarrow(C,d,h)$ is a morphism
of CDG\+coalgebras; equivalently, $w$ is a coaugmentation if it is
a morphism of graded coalgebras satisfying the equations $d\circ w = 0$
and $h\circ w = 0$.
 A $k$\+linear section $w$ is a coaugmentation if and only if
the CDG\+algebra $\Cb_w(C)$ is actually a DG\+algebra, i.~e., 
$h_{\Cb}=0$.

 For any CDG\+algebra $B$ with a $k$\+linear retraction~$v$,
the $k$\+linear section $w\:k\rarrow\Br_v(B)$ given by the embedding
of $k\simeq B_+^{\ot0}$ into $\Br(B)$ is a coaugmentation of
the CDG\+coalgebra $\Br_v(B)$ if and only if $h=0$ in $B$, i.~e.,
$B$ is a DG\+algebra.
 For any CDG\+coalgebra $C$ with a $k$\+linear section~$w$,
the $k$\+linear retraction $v\:\Cb_w(C)\rarrow k$ given by
the projection of $\Cb(C)$ onto $C_+^{\ot0}\simeq k$ is
an augmentation of the CDG\+algebra $\Cb_w(C)$ if and only if
$h=0$ on $C$, i.~e., $C$ is a DG\+coalgebra.
 So a (co)augmentation on one side of the (co)bar-construction
corresponds to the vanishing of the curvature element on the other side.

 Given a CDG\+algebra $B$, changing a retraction $v\:B\rarrow k$
to another retraction $v'\:B\rarrow k$ given by the formula 
$v'(b)=v(b)+a(b)$ leads to an isomorphism of CDG\+coalgebras
$(\id,a)\:\Br_v(B)\rarrow\Br_{v'}(B)$, where $a\:B_+\rarrow k$ is
a linear function of degree~$0$ identified with the corresponding
linear function $\Br(B)\rarrow B_+\rarrow k$ of degree~$1$.
 Given a CDG\+coalgebra $C$, changing a section $w\:k\rarrow C$
to another section $w'\:k\rarrow C$ given by the rule $w'(1)=w(1)+a$
leads to an isomorphism of CDG\+algebras
$(\id,a)\:\Cb_{w'}(C)\rarrow\Cb_w(C)$, where $a\in C_+$ is an element
of degree~$0$ identified with the corresponding element of
$\Cb(C)\supset C_+$ of degree~$1$.
 To an isomorphism of CDG\+coalgebras of the form $(\id,a)\:
(C,d,h)\rarrow (C,d',h')$ one can assign an isomorphism of
the corresponding cobar-constructions of the form $(f_a,0)\:
\Cb_w(C,d,h)\rarrow\Cb_w(C,d',h')$ with the automorphism~$f_a$ of
the graded algebra $\Cb(C)$ given by the rule $c\mpsto c-a(c)$
for $c\in C_+$.
 Here $a\:C\rarrow k$ is a linear function of degree~$1$.

 Consequently, there is a functor from the category of CDG\+coalgebras
to the category of CDG\+algebras assigning to a CDG\+coalgebra $C$
its cobar-construction $\Cb_w(C)$.
 The cobar-construction is also a functor from the category of
coaugmented CDG\+coalgebras to the category of DG\+algebras, from
the category of DG\+coalgebras to the category of augmented
CDG\+algebras, and from the category of coaugmented DG\+coalgebras
to the category of augmented DG\+algebras.

 Furthermore, let us call a morphism of CDG\+algebras
$(f,a)\:B\rarrow A$ \emph{strict} if one has $a=0$.
 Then there is a functor from the category of CDG\+algebras and strict
morphisms between them to the category of CDG\+coalgebras assigning to
a CDG\+algebra $B$ its bar-construction $\Br_v(B)$.
 The bar-construction is also a functor from the category of
DG\+algebras to the category of coaugmented CDG\+coalgebras, from
the category of augmented CDG\+algebras and strict morphisms between
them to the category of DG\+coalgebras, and from the category of
augmented DG\+algebras to the category of coaugmented DG\+coalgebras.

\begin{rem}
 There is \emph{no} isomorphism of bar-constructions corresponding to
an isomorphism of CDG\+algebras that is not strict.
 The reason is, essentially, that there exist no morphisms of tensor
coalgebras $\Br(B)$ that do not preserve their coaumentations
$k\simeq B_+^{\ot 0}\rarrow\Br(B)$, while there do exist coderivations
of $\Br(B)$ not compatible with the coaugmentation.
 Moreover, for any augmented CDG\+algebra $B=(B,d,h)$ with $h\ne0$
the DG\+coalgebra $\Br_v(B)$ is acyclic, i.~e., its cohomology is
the zero coalgebra.
 Indeed, consider the dual DG\+algebra $\Br_v(B)^*$.
 Its subalgebra of cocycles of degree zero $Z^0(\Br_v(B)^*)$ is complete
in the adic topology of its augmentation ideal $\ker(Z^0(\Br_v(B)^*)
\to k)$, while the ideal of coboundaries $\im d^{-1}\subset Z^0(\Br_v
(B)^*)$ contains elements not belonging to the augmentation ideal.
 Thus $\im d^{-1}=Z^0(\Br_v(B)^*)$ and the unit element $1\in\Br_v(B)^*$
is a coboundary.
  One can show that the DG\+coalgebra $\Br_v(B)$ considered up to
DG\+coalgebra isomorphisms carries no information about a coaugmented
CDG\+algebra $B$ except for the dimensions of its graded components.
 Furthermore, for any CDG\+algebra $(B,d,h)$ with $h\ne0$ and any
left CDG\+module $M$ over $B$, the CDG\+comodule
$\Br_v(B)\ot^{\tau_{B,v}}M$ and the CDG\+contramodule
$\Hom^{\tau_{B,v}}(\Br_v(B),M)$ over the CDG\+coalgebra $\Br_v(B)$ are
contractible, the notation being introduced
in~\ref{twisting-cochains-subsect}.
 (See~Remark~\ref{derived-category-ainfty-modules}.)
\end{rem}

\subsection{Twisting cochains}  \label{twisting-cochains-subsect}
 Let $C=(C,d_C,h_C)$ be a CDG\+coalgebra and $B=(B,d_B,h_B)$ be
a CDG\+algebra over the same field~$k$.
 We introduce a CDG\+algebra structure on the graded vector space
of homogeneous homomorphisms $\Hom_k(C,B)$ in the following way.
 The multiplication in $\Hom_k(C,B)$ is given by the formula
$(fg)(c)=(-1)^{|g||c_{(1)}|}f(c_{(1)})g(c_{(2)})$.
 The differential is given by the standard rule $d(f)(c)=
d_B(f(c))-(-1)^{|f|}f(d_C(c))$.
 The curvature element is defined by the formula $h(c)=\eps(c)h_B
-h_C(c)\cdot 1$, where $1$ is the unit element of $B$ and $\eps$
is the counit map of~$C$.
 A homogeneous linear map $\tau\:C\rarrow B$ of degree~$1$ is
called a \emph{twisting cochain}~\cite{Lef,Hue,Nic} if it satisfies
the equation $\tau^2+d\tau+h=0$ with respect to the above-defined
CDG\+algebra structure on $\Hom_k(C,B)$.

 Let $C$ be a CDG\+coalgebra and $w\:k\rarrow C$ be a homogeneous
$k$\+linear section of the counit map~$\eps$.
 Then the composition $\tau=\tau_{C,w}\:C\rarrow\Cb(C)$ of
the homogeneous linear maps $C\rarrow W\simeq C_+\simeq C_+^{\ot 1}
\rarrow\Cb(C)$ is a twisting cochain for $C$ and $\Cb_w(C)$.
 Let $B$ be a CDG\+algebra and $v\:C\rarrow k$ be a homogeneous
$k$\+linear retraction.
 Then minus the composition $\Br_v(B)\rarrow B_+^{\ot1}\simeq B_+
\simeq V\rarrow B$ is a twisting cochain $\tau=\tau_{B,v}\:\Br_v(B)
\rarrow B$ for $\Br_v(B)$ and~$B$.

 Let $\tau\:C\rarrow B$ be a twisting cochain for a CDG\+coalgebra $C$
and a CDG\+alge\-bra~$B$.
 Then for any left CDG\+module $M$ over $B$ there is a natural structure
of left CDG\+comodule over $C$ on the tensor product $C\ot_kM$.
 Namely, the coaction of $C$ in $C\ot_kM$ is induced by the left
coaction of $C$ in itself, while the differential on $C\ot_kM$ is given
by the formula $d(c\ot x) = d(c)\ot x + (-1)^{|c|}c\ot d(x) +
(-1)^{|c_{(1)}|}c_{(1)}\ot \tau(c_{(2)})x$.
 We will denote the tensor product $C\ot_kM$ with this CDG\+comodule
structure by $C\ot^\tau\!\. M$.
 Furthermore, for any left CDG\+comodule $N$ over $C$ there is a natural
structure of left CDG\+module over $B$ on the tensor product $B\ot_kN$.
 Namely, the action of $B$ in $B\ot_kN$ is induced by the left action
of $B$ in itself, while the differential on $B\ot_kN$ is given by
the formula $d(b\ot y) = d(b)\ot y + (-1)^{|b|}b\ot d(y) +
(-1)^{|b|+1}b\tau(y_{(-1)})\ot y_{(0)}$.
 We will denote the tensor product $B\ot_kN$ with this CDG\+module
structure by $B\ot^\tau\!\.N$.

 The correspondences assigning to a CDG\+module $M$ over $B$
the CDG\+comodule $C\ot^\tau\!\. M$ over $C$ and to a CDG\+comodule $N$
over $C$ the CDG\+module $B\ot^\tau\!\.N$ over $B$ can be extended to
DG\+functors whose action on morphisms is given by the standard
formulas $f_*(c\ot x) = (-1)^{|f||c|}c\ot f_*(x)$ and
$g_*(b\ot y) = (-1)^{|g||b|} b\ot g_*(y)$.
 The DG\+functor $C\ot^\tau\!\.{-}\:\sDG(B\modl)\rarrow\sDG(C\comodl)$
is right adjoint to the DG\+functor $B\ot^\tau{-}\!\.\:\sDG(C\comodl)
\rarrow\sDG(B\modl)$.

 Analogously, for any right CDG\+module $M$ over $B$ there is a natural
structure of right CDG\+comodule over $C$ on the tensor product
$M\ot_kC$.
 The coaction of $C$ in $M\ot_kC$ is induced by the right coaction of
$C$ in itself and the differential on $M\ot_kC$ is given by the formula
$d(x\ot c) = d(x)\ot c + (-1)^{|x|}x\ot d(c) + (-1)^{|x|+1}
x\tau(c_{(1)})\ot c_{(2)}$.
 We will denote the tensor product $M\ot_kC$ with this CDG\+comodule
structure by $M\ot^\tau\!\.C$.
 For any right CDG\+comodule $N$ over $C$ there is a natural structure
of right CDG\+module over $B$ on the tensor product $N\ot_kB$.
 Namely, the action of $B$ in $N\ot_kB$ is induced by the right action
of $B$ in itself and the differential on $N\ot_kB$ is given by
the formula $d(y\ot b) = d(y)\ot b + (-1)^{|y|}y\ot d(b) +
(-1)^{|y_{(0)}|} y_{(0)}\ot\tau(y_{(1)})b$.
 We will denote the tensor product $N\ot_kB$ with this CDG\+module
structure by $N\ot^\tau\!\.B$.

 For any left CDG\+module $P$ over $B$ there is a natural structure of
left CDG\+contra\-module over $C$ on the graded vector space of
homogeneous linear maps $\Hom_k(C,P)$.
 The contraaction of $C$ in $\Hom_k(C,P)$ is induced by the right
coaction of $C$ in itself as explained in~\ref{graded-contramodules}.
 The differential on $\Hom_k(C,P)$ is given by the formula
$d(f)(c) = d(f(c)) - (-1)^{|f|}f(d(c)) + (-1)^{|f||c_{(1)}|}
\tau(c_{(1)})f(c_{(2)})$ for $f\in\Hom_k(C,P)$.
 We will denote the graded vector space $\Hom_k(C,P)$ with this
CDG\+contramodule structure by $\Hom^\tau(C,P)$.
 For any left CDG\+contramodule $Q$ over $C$ there is a natural
structure of left CDG\+module over $B$ on the graded vector space of
homogeneous linear maps $\Hom_k(B,Q)$.
 The action of $B$ in $\Hom_k(B,Q)$ is induced by the right action
of $B$ in itself as explained in \ref{injective-dg-modules}
and~\ref{dg-mod-scalars}.
 The differential on $\Hom_k(B,Q)$ is given by the formula $d(f)(b) =
d(f(b)) - (-1)^{|f|}f(d(b)) + \pi(c\mapsto (-1)^{|f|+1+|c||b|}
f(\tau(c)b))$, where $\pi$ denotes the contraaction map
$\Hom_k(C,Q)\rarrow Q$.
 We will denote the graded vector space $\Hom_k(B,Q)$ with this
CDG\+module structure by $\Hom^\tau(B,Q)$.

 The correspondences assigning to a CDG\+module $P$ over $B$
the CDG\+contramodule $\Hom^\tau(C,P)$ over $C$ and to
a CDG\+contramodule $Q$ over $C$ the CDG\+module $\Hom^\tau(B,Q)$
over $B$ can be extended to DG\+functors whose action on morphisms
is given by the standard formula $g_*(f)=g\circ f$ for $f\:C\rarrow P$
or $f\:B\rarrow Q$.
 The DG\+functor $\Hom^\tau(C,{-})\:\sDG(B\modl)\rarrow\sDG(C\contra)$
is left adjoint to the DG\+functor $\Hom^\tau(B,{-})\:\sDG(C\contra)
\rarrow\sDG(B\modl)$.

\subsection{Duality for bar-construction}
\label{bar-construction-duality}
 Let $A=(A,d)$ be a DG\+algebra over a field~$k$.
 Choose a homogeneous $k$\+linear retraction $v\:A\rarrow k$ and
consider the bar-construction $C=\Br_v(A)$; then $C$ is
a coaugmented CDG\+coalgebra.
 Let $\tau=\tau_{A,v}\:C\rarrow A$ be the natural twisting cochain.

\begin{thm}
\textup{(a)} The functors $C\ot^\tau\!\.{-}\:\Hot(A\modl)\rarrow
\Hot(C\comodl)$ and $A\ot^\tau\!\.{-}\:\allowbreak\Hot(C\comodl)\rarrow
\Hot(A\modl)$ induce functors\/ $\sD(A\modl)\rarrow\sD^\co(C\comodl)$
and\/ $\sD^\co(C\comodl)\rarrow\sD(A\modl)$, which are mutually
inverse equivalences of triangulated categories. \par
\textup{(b)} The functors\/ $\Hom^\tau(C,{-})\:\Hot(A\modl)\rarrow
\Hot(C\contra)$ and\/ $\Hom^\tau(A,{-})\:\allowbreak\Hot(C\contra)
\rarrow\Hot(A\modl)$ induce functors\/ $\sD(A\modl)\rarrow
\sD^\ctr(C\contra)$ and\/ $\sD^\ctr(C\contra)\rarrow\sD(A\modl)$,
which are mutually inverse equivalences of triangulated categories. \par
\textup{(c)} The above equivalences of triangulated categories\/
$\sD(A\modl)\simeq\sD^\co(C\comodl)$ and\/ $\sD(A\modl)\simeq
\sD^\ctr(C\contra)$ form a commutative diagram with the equivalence
of triangulated categories\/ $\sD^\co(C\comodl)\simeq\sD^\ctr(C\contra)$
provided by the derived functors\/ $\boL\Phi_C$ and\/ $\boR\Psi_C$
of Theorem~\ref{co-contra-corr-theorem}.
\end{thm}

\begin{proof}
 Part~(a): first notice that for any coacyclic CDG\+comodule $N$
over $C$ the DG\+module $A\ot^\tau\!\.N$ over $A$ is contractible.
 Indeed, whenever $N$ is the total CDG\+module of an exact triple
of CDG\+modules $A\ot^\tau\!\.N$ is the total DG\+module of
an exact triple of DG\+modules that splits as an exact triple of
graded $A$\+modules.
 Secondly, let us check that for any acyclic DG\+module $M$ over $A$
the CDG\+comodule $C\ot^\tau\!\.M$ over $C$ is coacyclic.
 Introduce an increasing filtration $F$ on the coalgebra $C=\Cb(A)$
by the rule $F_n\Cb(A)=\bigoplus_{j\le n}A_+^{\ot j}$.
 There is an induced filtration on $C\ot^\tau\!\.M$ given by
the formula $F_n(C\ot^\tau\!\.M) = F_nC\ot_k M$.
 This is a filtration by CDG\+subcomodules and the quotient
CDG\+comodules $F_n(C\ot^\tau\!\.M)/F_{n-1}(C\ot^\tau\!\.M)$
have trivial $C$\+comodule structures.
 So these quotient CDG\+comodules can be considered simply as
complexes of vector spaces, and as such they are isomorphic to
the complexes $A_+^{\ot n}\ot_k M$.
 These complexes are acyclic, and hence coacyclic, whenever $M$
is acyclic.
 So the CDG\+comodule $C\ot^\tau\!\.M$ is coacyclic.
 Since it is cofree as a graded $C$\+comodule, it is even
contractible.
 We have shown that there are induced functors $\sD(A\modl)\rarrow
\sD^\co(C\comodl)$ and $\sD^\co(C\comodl)\rarrow\sD(A\modl)$;
it remains to check that they are mutually inverse equivalences.
 For any DG\+module $M$ over $A$, the DG\+module $A\ot^\tau\!\.C
\ot^\tau\!\.M$ is isomorphic to the total DG\+module of
the reduced bar-resolution $\dsb\rarrow A\ot A_+\ot A_+\ot M
\rarrow\ A\ot A_+\ot M\rarrow A\ot M$.
 So the cone of the adjunction morphism $A\ot^\tau\!\.C\ot^\tau\!\.M
\rarrow M$ is acyclic, since the reduced bar-resolution remains exact
after passing to the cohomology $\dsb\rarrow H(A)\ot H(A_+)\ot H(M)
\rarrow H(A)\ot H(M)\rarrow H(M)\rarrow 0$, as explained in the proof
of Theorem~\ref{projective-dg-modules}.
 For a CDG\+comodule $N$ over $C$, let $K$ denote the cone of
the adjunction morphism $N\rarrow C\ot^\tau\!\.A\ot^\tau\!\.N$.
 Let us show that the CDG\+comodule $K$ is absolutely acyclic.
 Introduce a finite increasing filtration on the graded $C$\+comodule
$K$ by the rules $G_{-2}K=0$, \ $G_{-1}K=N[1]$, \ $G_0K=N[1]\oplus
C\ot_k k\ot_k N \subset N[1]\oplus C\ot_k A\ot_k N$, and $G_1K=K$,
where $C\ot_kk\ot_kN$ is embedded into $C\ot_kA\ot_kN$ by the map
induced by the unit element of~$A$.
 The differential $d$ on $K$ does not preserve this filtration; still
one has $d(G_iK)\subset G_{i+1}K$.
 Let $\d$ denote the differential induced by $d$ on the associated
quotient $C$\+comodule $\gr_GK$.
 Then $(\gr_GK,\d)$ is an exact complex of graded $C$\+comodules;
indeed, it is isomorphic to the standard resolution of the graded
comodule $N$ over the graded tensor coalgebra~$C$.
 Set $L=G_{-1}K+d(G_{-1}K)\subset K$; it follows that both $L$ and
$K/L$ are contractible CDG\+comodules over~$C$.
 Part~(a) is proven; the proof of part~(b) is completely analogous
(up to duality).
 To prove~(c), it suffices to notice the natural isomorphisms
$\Psi_C(C\ot^\tau\!\.M)\simeq\Hom^\tau(C,M)$ and $\Phi_C(\Hom^\tau
(C,M))\simeq C\ot^\tau\!\.M$ for a DG\+module $M$ over~$A$.
\end{proof}

\subsection{Conilpotent duality for cobar-construction}
\label{conilpotent-cobar-duality}
 A graded coalgebra $E$ without counit is called \emph{conilpotent}
(cf.~\ref{dg-coalgebra-resolutions-proof}) if it is the union
of the kernels of iterated comultiplication maps $E\rarrow E^{\ot n}$.
 A graded coalgebra $D$ endowed with a coaugmentation (morphism of
coalgebras) $w\:k\rarrow D$ is called \emph{conilpotent} if
the graded coalgebra without counit $D/w(k)$ is conilpotent.
 One can easily see that a conilpotent graded coalgebra has a unique
coaugmentation.

 For a conilpotent graded coalgebra $D$ set $F_nD$ to be the kernel
of the composition $D\rarrow D^{\ot n+1}\rarrow (D/w(k))^{\ot n+1}$;
then the increasing filtration~$F$ on $D$ is compatible with
the coalgebra structure.
 We will call a CDG\+coalgebra $C=(C,d,h)$ \emph{conilpotent} if it is
conilpotent as a graded coalgebra and coaugmented as a CDG\+coalgebra.
 A DG\+coalgebra is \emph{conilpotent} if it is conilpotent as
a CDG\+coalgebra.
 For a conilpotent CDG\+coalgebra $C$, the filtration $F$ defined
above is compatible with the CDG\+coal\-gebra structure, i.~e., one
has $d(F_nC)\subset F_nC$, and in addition, $h(F_0C)=0$.

 Let $C$ be a conilpotent CDG\+coalgebra and $w\:k\rarrow C$ be its
coaugmentation map.
 Consider the cobar-construction $A=\Cb_w(C)$; then $A$ is 
a DG\+algebra.
 Let $\tau=\tau_{C,w}\:C\rarrow A$ be the natural twisting cochain.

\begin{thm}
 All the assertions of Theorem~\ref{bar-construction-duality} hold
for the DG\+algebra $A$, CDG\+co\-algebra $C$, and twisting cochain\/
$\tau$ as above in place of $A$, \ $C$, and\/ $\tau$
from~\ref{bar-construction-duality}.
\end{thm}

\begin{proof}
 Just as in the proof of Theorem~\ref{bar-construction-duality} one
shows that for any coacyclic CDG\+comodule $N$ over $C$ the DG\+module
$A\ot^\tau\!\.N$ over $A$ is contractible.
 To check that the CDG\+comodule $C\ot^\tau\!\.M$ is coacyclic (and
even contractible) for any acyclic DG\+module $M$ over~$A$, one uses
the filtration $F$ on the coalgebra $C$ that was constructed above
and the induced filtration of the CDG\+comodule $C\ot^\tau\!\.M$
by its CDG\+subcomodules $F_n(C\ot^\tau\!\.M)=F_nC\ot_kM$.
 The quotient CDG\+comodules $F_n(C\ot^\tau\!\.M)/F_{n-1}(C\ot^\tau
\!\.M)$ are simply the complexes $F_nC/F_{n-1}C\ot_kM$ with
the trivial $C$\+comodule structures, so they are coacyclic
whenever $M$ is acyclic.
 For any CDG\+comodule $N$ over $C$, the CDG\+comodule $C\ot^\tau\!\.
A\ot^\tau\!\.N$ is isomorphic to the reduced version of the curved
cobar-resolution introduced in the proof of
Theorem~\ref{cdg-coalgebra-inj-proj-resolutions}.
 So the same argument with the canonical filtration with respect to
the cobar differential~$\d$ proves that the cone of the adjunction
morphism $N\rarrow C\ot^\tau\!\.A\ot^\tau\!\.N$ is coacyclic.
 For a DG\+module $M$ over $A$, denote by $K$ the cocone of
the adjunction morphism $A\ot^\tau\!\.C\ot^\tau\!\.M\rarrow M$.
 We will show that the DG\+module $K$ is absolutely acyclic.
 Introduce a finite decreasing filtration $G$ on the graded
$A$\+module $K$ by the rules $K/G^{-1}K=0$, \ $K/G^0K=M[-1]$, \
$K/G^1K=A\ot_kk\ot_kM\oplus M[-1]$, and $K/G^2K=K$, where
$A\ot_kC\ot_kM$ maps onto $A\ot_kk\ot_kM$ by the map induced by
the counit of~$C$.
 The differential $d$ on $K$ does not preserve this filtration; still
one has $d(G^iK)\subset G^{i-1}K$.
 Let $\d$ denote the differential induced by $d$ on the associated
quotient $A$\+module $\gr_GK$.
 Then $(\gr_GK,\d)$ is an exact complex of graded $A$\+modules;
indeed, it is isomorphic to the standard resolution of the graded
module $M$ over the graded tensor algebra~$A$.
 Set $L=G^1K+d(G^1K)\subset L$; it follows that both $L$ and $K/L$
are contractible DG\+modules over~$A$.
 The proof of part~(b) is similar (up to duality), and the proof
of~(c) is analogous to the proof of
Theorem~\ref{bar-construction-duality}(c).
\end{proof}

\subsection{Acyclic twisting cochains}  \label{acyclic-twisting-cochain}
 Let $C$ be a coaugmented CDG\+coalgebra with a coaugmentation~$w$
and $A$ be a DG\+algebra.
 Then there is a natural bijective correspondence between
morphisms of DG\+algebras $\Cb_w(C)\rarrow A$ and twisting cochains
$\tau\:C\rarrow A$ such that $\tau\circ w=0$.
 Whenever $C$ is a DG\+coalgebra, so that $\Cb_w(C)$ is an augmented
DG\+algebra, and $A$ is also an augmented DG\+algebra with
an augmentation~$v$, a morphism of DG\+algebras $\Cb_w(C)\rarrow A$
preserves the augmentations if and only if one has $v\circ\tau=0$
for the corresponding twisting cochain~$\tau$.

 Let us assume from now on that $C$ is a conilpotent CDG\+coalgebra.
 Then a twisting cochain $\tau\:C\rarrow A$ with $\tau\circ w=0$ is
said to be \emph{acyclic} if the corresponding morphism of
DG\+algebras $\Cb(C)\rarrow A$ is a quasi-isomorphism.

\begin{thm}
\textup{(a)} The functors $C\ot^\tau\!\.{-}\:\Hot(A\modl)\rarrow
\Hot(C\comodl)$ and $A\ot^\tau\!\.{-}\:\allowbreak\Hot(C\comodl)\rarrow
\Hot(A\modl)$ induce functors\/ $\sD(A\modl)\rarrow\sD^\co(C\comodl)$
and\/ $\sD^\co(C\comodl)\rarrow\sD(A\modl)$, the former of which is
right adjoint to the latter.
 These functors are mutually inverse equivalences of triangulated
categories if and only if the twisting cochain $\tau$ is acyclic. \par
\textup{(b)} The functors\/ $\Hom^\tau(C,{-})\:\Hot(A\modl)\rarrow
\Hot(C\contra)$ and\/ $\Hom^\tau(A,{-})\:\allowbreak\Hot(C\contra)
\rarrow\Hot(A\modl)$ induce functors\/ $\sD(A\modl)\rarrow
\sD^\ctr(C\contra)$ and\/ $\sD^\ctr(C\contra)\rarrow\sD(A\modl)$,
the former of which is left adjoint to the latter.
 These functors are mutually inverse equivalences of triangulated
categories if and only if the twisting cochain $\tau$ is acyclic. \par
\textup{(c)} Whenever the twisting cochain $\tau$ is acyclic,
the above equivalences of triangulated categories\/
$\sD(A\modl)\simeq\sD^\co(C\comodl)$ and\/ $\sD(A\modl)\simeq
\sD^\ctr(C\contra)$ form a commutative diagram with the equivalence
of triangulated categories\/ $\sD^\co(C\comodl)\simeq\sD^\ctr(C\contra)$
provided by the derived functors\/ $\boL\Phi_C$ and\/ $\boR\Psi_C$
of Theorem~\ref{co-contra-corr-theorem}.
\end{thm}

 So in particular the twisting cochain $\tau=\tau_{A,v}$
of~\ref{bar-construction-duality} is acyclic; the twisting cochain
$\tau=\tau_{C,w}$ of~\ref{conilpotent-cobar-duality} is acyclic
by the definition.

 Notice that for any acyclic twisting cochain~$\tau$ the above
equivalences of derived categories (of the first and the second kind)
transform the trivial CDG\+comodule $k$ over $C$ into the free
DG\+module $A$ over $A$ and the trivial CDG\+contramodule $k$
over $C$ into the cofree DG\+module $\Hom_k(A,k)$ over~$A$.
 When $C$ is a DG\+coalgebra, $A$ is an augmented DG\+algebra with
an augmentation $v$, and one has $v\circ\tau=0$, these equivalences of
exotic derived categories also transform the trivial DG\+module $k$
over $A$ into the cofree DG\+comodule $C$ over $C$ and into
the free DG\+contramodule $\Hom_k(C,k)$ over~$C$.
 Here the trivial comodule, contramodule, and module structures on~$k$
are defined in terms of the coaugmentation~$w$ and augmentation~$v$.

\begin{proof}
 Part~(a): Just as in the proofs of
Theorems~\ref{bar-construction-duality}
and~\ref{conilpotent-cobar-duality} one shows that the functor
$N\mpsto A\ot^\tau\!\.M$ sends coacyclic CDG\+comodules to contractible
DG\+mod\-ules and the functor $M\mpsto C\ot^\tau\!\.M$ sends acyclic
DG\+modules to contractible CDG\+comodules.
 In order to see that the induced functors are adjoint it suffices to
recall that adjointness of functors can be expressed in terms of
adjunction morphisms and equations they satisfy; these morphisms
obviously continue to exist and the equations continue to hold after
passing to the induced functors between the quotient categories.
 To prove that these functors are equivalences of triangulated
categories if and only if $\tau$ is an acyclic twisting cochain,
it suffices to apply Theorem~\ref{conilpotent-cobar-duality}
and Theorem~\ref{dg-mod-scalars} for the morphism of DG\+algebras
$f\:\Cb_w(C)\rarrow A$.
 Indeed, there are obvious isomorphisms of functors
$C\ot^\tau\!\.M\simeq C\ot^{\tau_{C,w}}\!\.R_f(M)$ for
a DG\+module $M$ over $A$ and $A\ot^\tau\!\.N\simeq
E_f(\Cb_w(C)\ot^{\tau_{C,w}}\!\.N)$ for a CDG\+comodule $N$ over~$C$.
 The proof of part~(b) is completely analogous and uses
the functor $E^f$ instead of $E_f$.
 Notice that for any twisting cochain~$\tau$, for any CDG\+comodule $N$
over $C$ the DG\+module $A\ot^\tau\!\.N$ over $A$ is projective and
for any CDG\+contramodule $Q$ over $C$ the DG\+module $\Hom^\tau(A,Q)$
over $A$ is injective, as one can prove using either the adjointness of
the $\tau$\+related functors between the homotopy categories, or
the facts that $\sD^\co(C\comodl)$ is generated by the trivial
CDG\+comodule $k$ as a triangulated category with infinite direct sums
and $\sD^\ctr(C\contra)$ is generated by the trivial CDG\+contramodule
$k$ as a triangulated category with infinite products
(see~\ref{dg-coalgebra-resolutions-proof} for some details).
 The proof of part~(c) is identical to the proof of
Theorem~\ref{bar-construction-duality}(c).
\end{proof}

 Notice that for any acyclic twisting cochain~$\tau\:C\rarrow A$ and
any left CDG\+comod\-ule $N$ over $C$ the complex $A\ot^\tau\!\.N$
computes $\Cotor^C(k,N)\simeq\Ext_C(k,N)$.
 Indeed, the CDG\+comodule $C\ot^\tau\!\. A\ot^\tau\!\. N$ is
isomorphic to $N$ in the coderived category of CDG\+comodules over~$C$.
 This CDG\+comodule is also cofree as a graded $C$\+comodule and
one has $A\ot^\tau\!\.N\simeq k\oc_C (C\ot^\tau\!\.A\ot^\tau\!\.N)
\simeq \Hom_C(k\;C\ot^\tau\!\.A\ot^\tau\!\.N)$.
 Analogously, for any acyclic twisting cochain $\tau$ and any left
CDG\+contramodule $Q$ over $C$ the complex $\Hom^\tau(A,Q)\simeq
\Cohom_C(k\;\Hom^\tau(C,\Hom^\tau(A,Q)))\simeq k\ocn_C\Hom^\tau(C,
\Hom^\tau(A,Q))$ computes $\Coext_C(k,Q)\simeq\Ctrtor^C(k,Q)$.
 The DG\+algebra $A$ itself, considered as a complex, computes
$\Cotor^C(k,k)\simeq\Ext_C(k,k)$.

 Now let $C$ be a conilpotent DG\+coalgebra, $A$ be an augmented
DG\+algebra with an augmentation~$v$, and $\tau\:C\rarrow A$
be an acyclic twisting cochain for which $v\circ\tau=0$.
 Then for any left DG\+module $M$ over $A$ the complex
$C\ot^\tau\!\.M$ computes $\Tor^A(k,M)$.
 Indeed, the DG\+module $A\ot^\tau\!\.C\ot^\tau\!\.M$ is isomorphic
to $M$ in the derived category of DG\+modules over~$A$.
 This DG\+module is also projective, as mentioned in the above proof,
and one has $C\ot^\tau\!\.M\simeq k\ot_A(A\ot^\tau\!\.C\ot^\tau\!\.M)$.
 Analogously, for any left DG\+module $P$ over $A$ the complex
$\Hom^\tau(C,P)\simeq\Hom_A(k\;\Hom^\tau(A,\Hom^\tau(C,M)))$ computes
$\Ext_A(k,P)$.
 The DG\+coalgebra $C$ itself, considered as a complex, computes
$\Tor^A(k,k)$.

 It follows from the above Theorem that our definion of the coderived
category of CDG\+comodules is equivalent to the definition of
Lef\`evre-Hasegawa and Keller~\cite{Lef,Kel2} for a conilpotent
CDG\+coalgebra $C$.

\subsection{Koszul generators}  \label{koszul-generators}
 Let $A$ be a DG\+algebra over a field~$k$.
 Suppose that $A$ is endowed with an increasing filtration by graded
subspaces $k=F_0A\subset F_1A\subset F_2A\subset\dsb\subset A$ which
is compatible with the multiplication, preserved by the differential,
and cocomplete, i.~e., $A=\bigcup_n F_nA$.
 Let $\gr_FA=\bigoplus_n F_nA/F_{n-1}A$ be the associated quotient
algebra; it is a bigraded algebra with a grading~$i$ induced by
the grading of $A$ and a nonnegative grading~$n$ coming from
the filtration~$F$.
 Assume that the algebra $\gr_FA$ is Koszul~\cite{PP,P0,P} in its
nonnegative grading~$n$.

 Choose a graded subspace $\overline{V}\subset F_1A$ complementary
to $k=F_0A$ in $F_1A$.
 Notice that the filtration $F$ on~$A$ is determined by the subspace
$\overline{V}\subset A$, as a Koszul algebra is generated by its
component of degree~$1$.
 We will call $F$ a \emph{Koszul filtration} and $\overline{V}$
a \emph{Koszul generating subspace} of~$A$.
 Extend $\overline{V}$ to a subspace $V\subset A$ complementary
to~$k$ in $A$ and denote by $v\:A\rarrow k$ the projection of $A$
to~$k$ along~$V$.

 Let $C\subset\bigoplus_n (F_1A/k)^{\ot n}$ be the Koszul coalgebra
quadratic dual to $\gr_FA$.
 Recall that $C$ is constructed as the direct sum of intersections
of the form $C=\bigoplus_{n=0}^\infty\bigcap_{j=1}^{n-1}
(F_1A/k)^{\ot j-1}\ot_k R\ot_k(F_1A/k)^{\ot n-j-1}$, where
$R\subset (F_1A/k)\ot_k (F_1A/k)$ is the kernel of
the multiplication map $(F_1A/k)^{\ot2}\rarrow F_2A/F_1A$.
 In particular, $C_0=k$, \ $C_1=F_1A/k$, and $C_2=R$ are
the low-degree components of $C$ in the grading~$n$.
 We will consider $C$ as a subcoalgebra of the tensor coalgebra
$\Br(A)=\bigoplus_n (A/k)^{\ot n}$ and endow $C$ with the total
grading inherited from the grading of $\Br(A)$.

 One can easily see that the graded subcoalgebra $C\subset\Br_v(A)$
is preserved by the differential of $\Br_v(A)$, which makes it
a CDG\+algebra and a CDG\+subcoalgebra of $\Br_v(A)$.
 The CDG\+algebra structure on $C$ does not depend on the choice
of a subspace $V\subset A$, but only on the subspace $\overline{V}
\subset F_1A$.
 Define the homogeneneous linear map $\tau\:C\rarrow A$ of degree~$1$
as minus the composition $C\rarrow C_1=F_1A/k\simeq\overline{V}\rarrow
F_1A\subset A$.
 Clearly, $C$ is a conilpotent CDG\+coalgebra with the coaugmentation
$w\:k\simeq C_0\rarrow C$ and $\tau\circ w=0$.

\begin{thm}
 The map $\tau$ is an acyclic twisting cochain.
\end{thm}

\begin{proof}
 The element $\tau\in\Hom_k(C,A)$ satisfies the equation
$\tau^2+d\tau+h=0$, since it is the image of the twisting cochain
$\tau_{A,v}\in\Hom_k(\Br_v(A),A)$ under the natural strict surjective
morphism of CDG\+algebras $\Hom_k(\Br_v(A),A)\rarrow\Hom_k(C,A)$
induced by the embedding $C\rarrow\Br_v(A)$.
 To check that the morphism of DG\+algebras $\Cb_w(C)\rarrow A$
is a quasi-isomorphism, it suffices to pass to the associated
quotient objects with respect to the increasing filtration~$F$ on
$A$ and the increasing filtration $F$ on $\Cb_w(C)$ induced by
the increasing filtration $F$ on $C$ associated with the grading~$n$.
 Then it remains to use the fact that the coalgebra $C$ is
Koszul~\cite{P0}.
\end{proof}

 Let $A$ be a DG\+algebra with a Koszul generating
subspace~$\overline{V}$ and the corresponding Koszul filtration~$F$.
 Consider the CDG\+coalgebra $C$ and the twisting cochain
$\tau\:C\rarrow A$ constructed above.
 Let $M$ be a left DG\+module over $A$; suppose that $M$ is endowed
with an increasing filtration by graded subspaces $F_0M\subset F_1M
\subset\dsb\subset M$ that is compatible with the filtration on $A$
and the action of $A$ on $M$, preserved by the differential on~$M$,
and cocomplete, i.~e., $M=\bigcup_n F_nM$.
 Assume that the associated quotient module $\gr_FM$ over
the associated quotient algebra $\gr_FA$ is Koszul in its
nonnegative grading~$n$.
 Define the graded subcomodule $N\subset C\ot^\tau\!\.M$ as
the intersection $C\ot_k F_0M\cap C\ot_k S\subset C\ot_k M$,
where $S\subset (F_1A/k)\ot F_0M$ is the kernel of the action map
$(F_1A/k)\ot_k F_0M\rarrow F_1M/F_0M$.
 This is the Koszul comodule quadratic dual to the Koszul module
$\gr_FM$ over $\gr_FA$.
 The subcomodule $N$ is preserved by the differential on
$C\ot^\tau\!\.M$, so it is a CDG\+comodule over~$C$.
 The natural morphism of DG\+modules $A\ot^\tau\!\.N\rarrow M$
over $A$ is a quasi-isomorphism, as one can show in the way
analogous to the proof of the above Theorem.
 A dual result holds for DG\+modules $P$ over $A$ endowed with
a complete decreasing filtration satisfying the Koszulity condition
and the CDG\+contramodules $P$ quadratic dual to them.

\begin{ex}
 Let $\g$ be a Lie algebra and $A=U\g$ be its universal enveloping
algebra considered as a DG\+algebra concentrated in degree~$0$.
 Let $F$ be the standard filtration on $U\g$ and $\overline V = \g
\subset U\g$ be the standard generating subspace; they are well-known
to be Koszul.
 Then $A$ is augmented, so $C$ is a DG\+coalgebra; it can be
identified with the standard homological complex $C_*(\g)$.
 The functors $M\mpsto C\ot^\tau\!\.M$ and $P\mpsto\Hom^\tau(C,P)$
are isomorphic to the functors of standard homological and
cohomological complexes $M\mpsto C_*(\g,M)$ and $P\mpsto C^*(\g,P)$
with coefficients in complexes of $\g$\+modules $M$ and~$P$.
 Hence we see that these functors induce equivalences between
the derived category of $\g$\+modules, the coderived category of
DG\+comodules over $C_*(\g)$, and the contraderived category of
DG\+contramodules over $C_*(\g)$.
 When $\g$ and consequently $C_*(\g)$ are finite-dimensional,
DG\+comodules and DG\+contramodules over $C_*(\g)$ can be identified
with DG\+modules over the standard cohomological complex $C^*(\g)$,
so the functors $M\mpsto C_*(\g,M)$ and $P\mpsto C^*(\g,P)$ induce
equivalences between the derived category of $\g$\+modules,
the coderived category of DG\+modules over $C^*(\g)$, and
the contraderived category of DG\+modules over $C^*(\g)$.
 These results can be extended to the case of a central extension
of Lie algebras $0\rarrow k\rarrow\g'\rarrow\g\rarrow 0$ with
the kernel~$k$ and the enveloping algebra $U'\g = U\g'/(1_{U\g'}-
1_{\g'})$ governing representations of $\g'$ where the central
element $1\in k\subset\g'$ acts by the identity.
 Choose a section $\g\rarrow\g'$ of our central extension and define
the generating subspace $\overline{V}\subset\g'\subset U'\g$
accordingly; then the corresponding CDG\+coalgebra $C$ coincides
with the DG\+coalgebra $C_*(\g)$ as a graded coalgebra with
a coderivation; the $2$\+cochain $C_2(\g)\rarrow k$ of the central
extension $\g'\rarrow\g$ is the curvature linear function of~$C$.
 The derived category of $U'\g$\+modules is equivalent to
the coderived category of CDG\+comodules and the contraderived
category of CDG\+contramodules over this CDG\+coalgebra~$C$.
 Furthermore, for a finite-dimensional Lie algebra $\g$ the bounded
derived category of finitely generated $U'\g$\+modules is
equivalent to the absolute derived category of finite-dimensional
CDG\+comodules over the (finite-dimensional) CDG\+coalgebra $C$,
since the algebra $U'\g$ is Noetherian and has a finite homological
dimension (cf.~\ref{coherent-d-modules}).
\end{ex}

\subsection{Nonconilpotent duality for cobar-construction}
\label{nonconilpotent-duality}
 Let $C$ be a CDG\+coalge\-bra endowed with a homogeneous $k$\+linear
section $w\:k\rarrow C$ of the counit map~$\eps$, and let $B$ be
a CDG\+algebra.
 Then there is a natural bijective correspondence between morphisms
of CDG\+algebras $\Cb_w(C)\rarrow B$ and twisting cochains
$\tau\:C\rarrow B$.
 A morphism of CDG\+algebras $\Cb_w(C)\rarrow B$ is strict if and only
if one has $\tau\circ w=0$ for the corresponding twisting
cochain~$\tau$.
 Whenever $C$ is a DG\+coalgebra, so that $\Cb_w(C)$ is an augmented
CDG\+algebra, and $B$ is also an augmented CDG\+algebra with
an augmentation~$v$, a morphism of CDG\+algebras $\Cb_w(C)\rarrow B$
preserves the augmentations if and only if one has $v\circ\tau=0$
for the corresponding twisting cochain~$\tau$.

 Given a CDG\+coalgebra $C$ with a $k$\+linear section~$w\:k\rarrow C$,
set $B=\Cb_w(C)$ and $\tau=\tau_{C,w}\:C\rarrow B$.

{\emergencystretch=0em \hfuzz=2.5pt
\begin{thm}
\textup{(a)} The functors $C\ot^\tau\!\.{-}\:\Hot(B\modl)\rarrow
\Hot(C\comodl)$ and $B\ot^\tau\!\.{-}\:\allowbreak\Hot(C\comodl)\rarrow
\Hot(B\modl)$ induce functors\/ $\sD^\co(B\modl)\rarrow
\sD^\co(C\comodl)$ and\/ $\sD^\co(C\comodl)\rarrow\sD^\co(B\modl)$,
which are mutually inverse equivalences of triangulated categories. \par
\textup{(b)} The functors\/ $\Hom^\tau(C,{-})\:\Hot(B\modl)\rarrow
\Hot(C\contra)$ and\/ $\Hom^\tau(B,{-})\:\allowbreak\Hot(C\contra)
\rarrow\Hot(B\modl)$ induce functors\/ $\sD^\ctr(B\modl)\rarrow
\sD^\ctr(C\contra)$ and\/ $\sD^\ctr(C\contra)\rarrow\sD^\ctr(B\modl)$,
which are mutually inverse equivalences of triangulated categories. \par
\textup{(c)} The above equivalences of triangulated categories\/
$\sD^\abs(B\modl)\simeq\sD^\co(C\comodl)$ and\/ $\sD^\abs(B\modl)\simeq
\sD^\ctr(C\contra)$ form a commutative diagram with the equivalence
of triangulated categories\/ $\sD^\co(C\comodl)\simeq\sD^\ctr(C\contra)$
provided by the derived functors\/ $\boL\Phi_C$ and\/ $\boR\Psi_C$
of Theorem~\ref{co-contra-corr-theorem}.
\end{thm}}

 Whenever $C$ is an coaugmented CDG\+coalgebra and accordingly $B$ is
a DG\+algebra, the above equivalences of triangulated categories
transform the trivial CDG\+comodule $k$ over $C$ into the free
DG\+module $B$ over $B$ and the trivial CDG\+contramodule $k$ over $C$
into the cofree DG\+module $\Hom_k(B,k)$ over~$B$.
 Whenever $C$ is a DG\+coalgebra and accordingly $B$ is an augmented
CDG\+algebra, the same equivalences of triangulated categories
transform the trivial CDG\+module $k$ over $B$ into the cofree
DG\+comodule $C$ over $C$ and into the free DG\+contramodule
$\Hom_k(C,k)$ over~$C$.

\begin{proof}
 The assertions about existence of induced functors in (a) and~(b)
hold for any CDG\+coalgebra $C$, \ CDG\+algebra $B$, and twisting
cochain $\tau\:C\rarrow B$.
 Indeed, the functors $M\mpsto C\ot^\tau\!\.M$ and $N\mpsto
B\ot^\tau\!\.N$ send coacyclic objects to contractible ones, while
the functors $P\mpsto\Hom^\tau(C,P)$ and $Q\mpsto\Hom^\tau(B,Q)$
send contraacyclic objects to contractible ones, for the reasons
explained in the proof of Theorem~\ref{bar-construction-duality}.
 One also finds that the induced functors are adjoint to each other,
as explained in the proof of Theorem~\ref{acyclic-twisting-cochain}.
 Now when $B=\Cb_w(C)$ and $\tau=\tau_{C,w}$, the adjunction
morphisms are isomorphisms, as it was shown in the proof of
Theorem~\ref{conilpotent-cobar-duality}.
 Notice that $\sD^\co(B\modl)=\sD^\abs(B\modl)=\sD^\ctr(B\modl)$
in this case by Theorem~\ref{finite-homol-dim-cdg-ring}.
 The proof of part~(c) is identical to that of
Theorem~\ref{bar-construction-duality}(c).
\end{proof}

 We continue to assume that $B=\Cb_w(C)$ and $\tau=\tau_{C,w}$.
 Whenever $C$ is a coaug\-mented CDG\+coalgebra, for any left
CDG\+comodule $N$ over $C$ the complex $B\ot^\tau\!\.N$ computes
$\Cotor^C(k,N)\simeq\Ext_C(k,N)$ and for any left CDG\+contramodule
$Q$ over $C$ the complex $\Hom^\tau(B,Q)$ computes $\Coext_C(k,Q)
\simeq\Ctrtor^C(k,Q)$, as explained in~\ref{acyclic-twisting-cochain}.
 The DG\+algebra $B$ itself, considered as a complex, computes
$\Cotor^C(k,k)\simeq\Ext_C(k,k)$.
 Whenever $C$ is a DG\+coalgebra, for any left CDG\+module $M$
over $B$ the complex $C\ot^\tau\!\.M$ computes $\Tor^{B,I\!I}(k,M)$
and for any left CDG\+module $P$ over $B$ the complex $\Hom^\tau(C,P)$
computes $\Ext_B^{I\!I}(k,P)$.
 The DG\+coalgebra $C$ itself, considered as a complex, computes
$\Tor^{B,I\!I}(k,k)$.

\begin{cor}
 Let $C$ be a conilpotent CDG\+coalgebra, $w\:k\rarrow C$ be its
coaugmentation, and $A=\Cb_w(C)$ be its cobar-construction.
 Then the derived category\/ $\sD(A\modl)$ and the absolute derived
category\/ $\sD^\abs(A\modl)$ coincide; in other words, any acyclic
DG\+module over~$A$ is absolutely acyclic.
\end{cor}

\begin{proof}
 Compare Theorem~\ref{conilpotent-cobar-duality} and
Theorem~\ref{nonconilpotent-duality}.
\end{proof}

 For more general results of this kind, see~\ref{koszul-cogenerators}
and~\ref{cofibrant-dg-alg}.
 For a counterexample showing that the conilpotency condition is
necessary in this Corollary, see Remark~\ref{co-algebra-bar-duality}.

 Now let $\tau\:C\rarrow B$ be any twisting cochain between
a CDG\+coalgebra $C$ and a CDG\+algebra $B$.
 Let us discuss the adjunction properties of our functors between
the homotopy categories in some more detail.
 Notice that the functor $M\mpsto C\ot^\tau\!\. M$ is the composition
of left adjoint functors $M\mpsto \Hom^\tau(C,M)$ and $Q\mpsto
\Phi_C(Q)$.
 So the corresponding composition of right adjoint functors
$N\mpsto \Hom^\tau(B,\Psi_C(N))$ is right adjoint to the functor
$M\mpsto C\ot^\tau\!\.M$.
 At the same time, the functor $N\mpsto B\ot^\tau\!\.N$ is left
adjoint to the functor $M\mpsto C\ot^\tau\!\.M$.
 Analogously, the functor $M\mpsto\Hom^\tau(C,M)$ is the composition
of right adjoint functors $M\mpsto C\ot^\tau\!\.M$ and $N\mpsto
\Psi_C(N)$.
 So the corresponding composition of left adjoint functors
$Q\mpsto B\ot^\tau\!\.\Phi_C(Q)$ is left adjoint to the functor
$M\mpsto\Hom^\tau(C,M)$.
 At the same time, the functor $Q\mpsto\Hom^\tau(B,Q)$ is right
adjoint to the functor $M\mpsto\Hom^\tau(C,M)$.
 
\subsection{Koszul cogenerators}  \label{koszul-cogenerators}
 Let $C$ be a CDG\+coalgebra over a field~$k$.
 Suppose $C$ is endowed with a decreasing filtration by graded
subspaces $C=G^0C\supset G^1C\subset G^2C\supset\dsb\supset G^NC
\supset G^{N+1}C=0$ such that $G$ is compatible with
the comultiplication and the counit, $d_C(G^nC)\subset G^{n-1}C$,
and $h_C(G^3C)=0$.
 Let $\gr_GC=\bigoplus_n G^nC/G^{n+1}C$ be the associated quotient
coalgebra; it is a bigraded coalgebra with a grading~$i$ induced by
the grading of $C$ and a grading $0\le n\le N$ coming from
the filtration~$G$.
 Assume that the coalgebra $\gr_GC$ is Koszul~\cite{P0,P} in
its nonnegative grading~$n$.

 The differential $d_C$ on $C$ induces a differential $d_{\gr C}$
on $\gr_GC$ having the degrees $i=1$ and $n=-1$.
 The linear function $h_C$ on $C$ induces a linear function
$h_{\gr C}$ on $\gr_GC$ having the degrees $i=2$ and $n=-2$.
 This defines a structure of CDG\+coalgebra (with respect to
the grading~$i$) on $\gr_GC$.
 The natural map $k\simeq C/G^1G\rarrow\gr_GC$ provides
a coaugmentation~$w$ for $\gr_GC$.
 Consider the DG\+algebra $\Cb_w(\gr_GC)$.
 The increasing filtration $F$ on $\gr_GC$ coming from
the grading~$n$ is compatible with the CDG\+coalgebra structure
and induces an increasing filtration $F$ on $\Cb_w(\gr_GC)$.
 Considering the associated graded DG\+algebra $\gr_F\Cb_w(\gr_GC)$,
one can show~\cite[Section~7 of Chapter~5]{PP} that the cohomology
algebra $B$ of the DG\+algebra $\Cb_w(\gr_GC)$ is concentrated in
degree~$0$ with respect to the difference of the grading~$n$
induced by the grading~$n$ of $\gr_GC$ and the grading~$m$ of
$\Cb_w(\gr_GC)$ by the tensor degrees.
 Moreover, the associated graded algebra $\gr_FB$ is quadratic dual
to the Koszul coalgebra $\gr_GC$.
 The component $F_1B$ is identified with the direct sum $k\oplus V$
of the subspace $k=F_0B$ and a graded subspace $V\subset B$ naturally
isomorphic to $(G^1C/G^2C)[-1]$.

 Now choose a homogeneous $k$\+linear section $w\:k\rarrow C$ of
the counit map of~$C$.
 Define a graded algebra morphism $\Cb_w(C)\rarrow B$ by the rule
that the map $C_+[-1]\rarrow B$ is the composition $C_+[-1] = G^1C[-1]
\rarrow (G^1C/G^2C)[-1]\simeq V\rarrow B$.
 We claim that the CDG\+algebra structure on $\Cb_w(C)$ induces such
a structure on~$B$.
 Indeed, the graded algebra $B$ is the quotient algebra of
the graded algebra $\Cb_w(C)$ by the ideal generated by
$G^2C[-1]+d_{\Cb}(G^2C[-1])$.
 Furthermore, one has $d(F_nB)\subset F_{n+1}B$ and $h\in F_2B$.
 The associated quotient algebra $\gr_FB$ is Koszul with respect
to the grading~$n$ coming from the filtration~$F$.
 In addition, $\gr_FB$ has a finite homological dimension as
a graded algebra with respect to the grading~$n$, or equivalently,
as a bigraded algebra.
 Conversely, given a CDG\+algebra $B$ with an increasing filtration
$F$ with the above properties one can recover a CDG\+coalgebra~$C$.
 One proves this using the dual version of the spectral sequence
argument from~\cite{PP}, which works in the assumption of finite
homological dimension.

 The graded algebra $B$ has a finite (left and right) homological
dimension, since one can compute the spaces $\Ext$ over it in
terms of nonhomogeneous Koszul complexes.
 Let $\tau\:C\rarrow B$ be the twisting cochain corresponding to
the morphism of CDG\+algebras $\Cb_w(C)\rarrow B$.

\begin{thm}
 The assertions (a-c) of Theorem~\ref{nonconilpotent-duality} hold
for the CDG\+coalgebra $C$, the DG\+algebra $B$, and the twisting
cochain~$\tau$.
\end{thm}

\begin{proof}
 Let us check that the adjunction morphism $N\rarrow C\ot^\tau\!\.
B\ot^\tau\!\.N$ has a coacyclic cone for any CDG\+comodule $N$
over~$C$.
 Introduce an increasing filtration $F$ on this cone $K$ by the rules
$F_{-2}K=0$, \ $F_{-1}K=N[1]$, and $F_nK=\cone(N\to C\ot_k F_nB
\ot_k N)$ for $n\ge0$.
 The differential $d$ on $K$ does not preserve this filtration;
still one has $d(F_nK)\subset F_{n+1}K$.
 Let $\d$ denote the induced differential on the associated quotient
$C$\+comodule $\gr_FK$.
 Then $(\gr_FK,\d)$ is an exact complex of graded $C$\+comodules;
indeed, it is isomorphic to the nonhomogeneous Koszul resolution of
the graded $C$\+comodule $N$.
 Consider the filtration of the CDG\+comodule $K$ over $C$ by
the CDG\+subcomodules $F_nK+d(F_nK)$.
 Then the associated quotient CDG\+comodules of this filtration are
contractible, hence $K$ is coacyclic.
 To check that the adjunction morphism $B\ot^\tau\!\.C\ot^\tau\!\.M
\rarrow M$ has an absolutely acyclic cone for any CDG\+module $M$
over $B$, one argues in the analogous way.
 Introduce a decreasing filtration $G$ on this cone $K$ by
the rules $K/G^{-1}K=0$, \ $K/G^0K=M$, and $K/G^{n+1}K=
\cone(B\ot^\tau\!\.C/G^{n+1}C\ot^\tau\!\.M\to M)$ for $n\ge0$.
 Then the associated quotient CDG\+modules of the finite filtration
of $K$ by the CDG\+submodules $G^nK+d(G^nK)$ are contractible.
 The rest is explained in the proof of
Theorem~\ref{nonconilpotent-duality} and the previous theorems.
\end{proof}

 Whenever $C$ is a coaugmented coalgebra and $w$~is its coaugmentation,
the CDG\+algebra $B$ is, in fact, a DG\+algebra.
 In this case, for any left CDG\+comodule $N$ over $C$ the complex
$B\ot^\tau\!\.N$ computes $\Cotor^C(k,N)\simeq\Ext_C(k,N)$ and for any
left CDG\+contramodule $Q$ over $C$ the complex $\Hom^\tau(B,Q)$
computes $\Coext_C(k,Q)\simeq\Ctrtor^C(k,Q)$, as explained
in~\ref{acyclic-twisting-cochain}.
 The DG\+algebra $B$ itself, considered as a complex, computes
$\Cotor^C(k,k)\simeq\Ext_C(k,k)$.
 Whenever $C$ is, actually, a DG\+coalgebra, the CDG\+algebra $B$ is
augmented.
 In this case, for any left CDG\+module $M$ over $B$ the complex
$C\ot^\tau\!\.M$ computes $\Tor^{B,I\!I}(k,M)$ and for any left
CDG\+module $P$ over $B$ the complex $\Hom^\tau(C,P)$ computes
$\Ext_B^{I\!I}(k,P)$.
 The DG\+coalgebra $C$ itself, considered as a complex, computes
$\Tor^{B,I\!I}(k,k)$.

 Let $\Hot(B\modl_\fgp)$ denote the full triangulated subcategory
of the homotopy category $\Hot(B\modl)$ consisting of all CDG\+modules
that are projective and finitely generated as graded $B$\+modules.
 The category $\Hot(B\modl_\fgp)$ can be also considered as a full 
triangulated subcategory of $\sD^\co(B\modl)=\sD^\ctr(B\modl)=
\sD^\abs(B\modl)$.

\begin{cor1}
 For a CDG\+algebra $B$ as above, the objects of the triangulated
subcategory\/ $\Hot(B\modl_\fgp)\subset\sD^\abs(B\modl)$ are compact
generators of the triangulated category\/ $\sD^\abs(B\modl)$.
\end{cor1}

\begin{proof}
 It was shown in~\ref{dg-coalgebra-resolutions-proof} (see
also~\ref{fin-dim-cdg-comod}) that the objects of the full triangulated
subcategory $\sD^\abs(C\comodl_\fd)\subset\sD^\co(C\comodl)$ are
compact generators of the coderived category $\sD^\co(C\comodl)$.
 Obviously, the Koszul duality functor $N\mpsto B\ot^\tau\!\.M$ from
the above Theorem maps $\sD^\abs(C\comodl_\fd)$ into $\Hot(B\modl_\fgp)
\subset\sD^\abs(B\modl)$.
 On the other hand, by Theorem~\ref{finite-homol-dim-cdg-ring}
the absolute derived category $\sD^\abs(B\modl)$ is equivalent to
the homotopy category $\Hot(B\modl_\proj)$ of CDG\+modules that
are projective as graded $B$\+modules.
 It is clear that the objects of $\Hot(B\modl_\fgp)$ are compact in
$\Hot(B\modl_\proj)$.
 The assertion of Corollary follows from that of Theorem, part~(a)
together with these observations.
\end{proof}

 For a discussion of questions and results related to Corollary~1,
see Question~\ref{coherent-cdg-ring-case}.

\begin{cor2}
 In the above notation, suppose that $C$ is a conilpotent
CDG\+coalgebra and $w\:k\rarrow C$ is its coaugmentation.
 Then the derived category\/ $\sD(B\modl)$ and the absolute derived
category\/ $\sD^\abs(B\modl)$ coincide; in other words, any acyclic
DG\+module over $B$ is absolutely acyclic.
\end{cor2}

\begin{proof}
 It suffices to show that the CDG\+comodule $C\ot^\tau\!\.M$ over~$C$
is coacyclic whenever a DG\+module $M$ over $B$ is acyclic.
 It was explained in the proof of
Theorem~\ref{conilpotent-cobar-duality} how to do so.
\end{proof}

 For other results similar to that of Corollary~2,
see~\ref{cofibrant-dg-alg}.

\subsection{Cotor and Tor, Coext and Ext, Ctrtor and Tor; restriction
and extension of scalars} \label{cotor-and-tor-coext-and-ext}
 Let $C$ be a conilpotent CDG\+coalgebra, $A$ be a DG\+algebra, and
$\tau\:C\rarrow A$ be an acyclic twisting cochain.
 Notice that by the right version of
Theorem~\ref{acyclic-twisting-cochain} the functors $M\mpsto
M\ot^\tau\!\.C$ and $N\mpsto N\ot^\tau\!\.A$ induce an equivalence
of triangulated categories $\sD(\modr A)\simeq\sD^\co(\comodr C)$.

\begin{thm1}
\textup{(a)} The equivalences of triangulated categories\/
$\sD^\co(\comodr C)\simeq\sD(\modr A)$ and\/ $\sD^\co(C\comodl)
\simeq\sD(A\modl)$ transform the functor\/ $\Cotor^C$ into
the functor\/ $\Tor^A$. \par
\textup{(b)} The equivalences of triangulated categories\/
$\sD^\co(C\comodl)\simeq\sD(A\modl)$ and\/ $\sD^\ctr(C\contra)
\simeq\sD(A\modl)$ transform the functor\/ $\Coext_C$ into
the functor\/ $\Ext_A$. \par
\textup{(c)} The equivalences of triangulated categories\/
$\sD^\co(\comodr C)\simeq\sD(\modr A)$ and\/ $\sD^\ctr(C\contra)
\simeq\sD(A\modl)$ transform the functor\/ $\Ctrtor^C$ into
the functor\/ $\Tor^A$.
\end{thm1}

\begin{proof}
 To prove part~(a), it suffices to use either of the natural
isomorphisms of complexes $N'\oc_C(C\ot^\tau\!\.M'')\simeq
(N'\ot^\tau\!\.A)\ot_A M''$ or $(M'\ot^\tau\!\.C)\oc_C N''\simeq
M'\ot_A(A\ot^\tau\!\.N'')$.
 To check that one obtains the same isomorphism of functors in
these two ways, notice that the two compositions $N'\oc_C N''\rarrow
N'\oc_C(C\ot^\tau\!\.A\ot^\tau\!\.N'')\simeq(N'\ot^\tau\!\.A)\ot_A
(A\ot^\tau\!\.N'')$ and $N'\oc_C N''\rarrow (N'\ot^\tau\!\.A\ot^\tau
\!\.C)\oc_C N''\rarrow (N'\ot^\tau\!\.A)\ot_A(A\ot^\tau\!\.N'')$
coincide.
 To prove~(b), use either of the isomorphisms $\Cohom_C(C\ot^\tau\!\.M\;
Q)\simeq\Hom_A(M,\Hom^\tau(A,Q))$ or $\Cohom_C(N,\Hom^\tau(C,P))\simeq
\Hom_A(A\ot^\tau\!\.N,P)$.
 Alternatively, use the result of~\ref{coext-and-ext-cotor-and-ctrtor}
and Theorem~\ref{acyclic-twisting-cochain}.
 To check~(c), use the natural isomorphism $(N\ot^\tau\!\.A)\ot_AP
\simeq N\ocn_C\Hom^\tau(C,P)$.
\end{proof}

 Let $C$ and $D$ be conilpotent CDG\+coalgebras, $A$ and $B$ be
DG\+algebras, $\tau\:C\rarrow A$ and $\sigma\:D\rarrow B$ be
acyclic twisting cochains.
 Let $f\:A\rarrow B$ be a morphism of DG\+algebras and $g=(g,a)\:
C\rarrow D$ be a morphism of CDG\+coalgebras.
 Assume that the following commutativity equation holds:
the difference $\sigma g - f\tau$ is equal to the composition
of the map $a\:C\rarrow k$ and the unit map $k\rarrow B$.

\begin{prop} {\hbadness=3000
\textup{(a)} The equivalences of triangulated categories\/
$\sD(A\modl)\simeq\sD^\co(C\comodl)$ and\/ $\sD(B\modl)\simeq
\sD^\co(D\comodl)$ transform the functor\/ $\boI R_f$ into
the functor\/ $\boR E_g$ and the functor\/ $\boL E_f$ into
the functor\/ $\boI R_g$. \par}
\textup{(b)} The equivalences of triangulated categories\/
$\sD(A\modl)\simeq\sD^\ctr(C\contra)$ and\/ $\sD(B\modl)\simeq
\sD^\ctr(D\contra)$ transform the functor\/ $\boI R_f$ into
the functor\/ $\boL E^g$ and the functor $\boR E^f$ into
the functor\/ $\boI E^g$.
\end{prop}

\begin{proof}
 To prove part~(a), use the natural isomorphisms
$E_g(D\ot^\sigma M)\simeq C\ot^\tau R_f(M)$ and
$E_f(A\ot^\tau N)\simeq B\ot^\sigma R_g(N)$ for a DG\+module $M$
over $B$ and a CDG\+comodule $N$ over~$C$.
 The proof of part~(b) is similar.
\end{proof}

 Now let $C$ be a CDG\+algebra endowed with a $k$\+linear section
$w\:k\rarrow C$ and $B=\Cb_w(C)$ be its cobar-construction.
 More generally, one can assume that $C$ and $B$ are a CDG\+coalgebra
and a CDG\+algebra related by the construction
of~\ref{koszul-cogenerators}.
 By the right version of Theorem~\ref{nonconilpotent-duality}
or~\ref{koszul-cogenerators}, the functors $M\mpsto M\ot^\tau\!\.C$
and $N\mpsto N\ot^\tau\!\.B$ induce an equivalence of triangulated
categories $\sD^\abs(\modr B)\simeq \sD^\co(\comodr C)$.

\begin{thm2}
\textup{(a)} The equivalences of triangulated categories\/
$\sD^\co(\comodr C)\simeq\sD^\abs(\modr B)$ and\/ $\sD^\co(C\comodl)
\simeq\sD^\abs(B\modl)$ transform the functor\/ $\Cotor^C$ into
the functor\/ $\Tor^{B,I\!I}$. \par
\textup{(b)} The equivalences of triangulated categories\/
$\sD^\co(C\comodl)\simeq\sD^\abs(B\modl)$ and\/ $\sD^\ctr(C\contra)
\simeq\sD^\abs(B\modl)$ transform the functor\/ $\Coext_C$ into
the functor\/ $\Ext_B^{I\!I}$. \par
\textup{(c)} The equivalences of triangulated categories\/
$\sD^\co(\comodr C)\simeq\sD(\modr B)$ and\/ $\sD^\ctr(C\contra)
\simeq\sD(B\modl)$ transform the functor\/ $\Ctrtor^C$ into
the functor\/ $\Tor^{B,I\!I}$.
\end{thm2}

\begin{proof}
 See the proof of Theorem~1.
\end{proof}

\subsection{Bar duality between algebras and coalgebras}
\label{co-algebra-bar-duality}
 Graded tensor coalgebras are cofree objects in the category of
conilpotent graded coalgebras.
 More precisely, for any conilpotent graded coalgebra $C$ with
the coaugmentation $w\:k\rarrow C$ and any graded vector space $U$ 
there is a bijective correspondence between graded coalgebra
morphisms $C\rarrow\bigoplus_{n=0}^\infty U^{\ot n}$ and
homogeneous $k$\+linear maps $C/w(k)\rarrow U$ of degree zero.
 Notice that the graded tensor coalgebra $\bigoplus_{n=0}^\infty
U^{\ot n}$ is conilpotent and any morphism of conilpotent graded
coalgebras preserves the coaugmentations.
 Let us emphasize that the above assertion in \emph{not} true
when the graded coalgebra $C$ is not conilpotent.

 Let $B$ be a CDG\+algebra and $v\:B\rarrow k$ be a homogeneous
$k$\+linear retraction.
 Let $C$ be a CDG\+coalgebra that is conilpotent as a graded
coalgebra; denote by $w\:k\rarrow C$ the coaugmentation map (which
does not have to be a coaugmentation of $C$ as a CDG\+algebra).
 Then there is a natural bijective correspondence between morphisms
of CDG\+coalgebras $C\rarrow \Br_v(B)$ and twisting cochains
$\tau\:C\rarrow B$ such that $\tau\circ w=0$.
 Whenever $C$ is a DG\+coalgebra and $v$ is an augmentation of~$B$,
a CDG\+coalgebra morphism $C\rarrow\Br_v(B)$ is actually a morphism
of DG\+coalgebras if and only if one has $v\circ\tau=0$ for
the corresponding twisting cochain~$\tau$.

 Let $k\coalg_\cdg^\conilp$ denote the category of conilpotent
CDG\+coalgebras, $k\coalg_\dg^\conilp$ denote the category of
conilpotent DG\+coalgebras, $k\alg_\dg^+$ denote the category of
DG\+algebras with nonzero units, and $k\alg_\dg^\aug$ denote
the category of augmented DG\+algebras (over the ground field~$k$).
 It follows from the above that the functor of conilpotent
cobar-construction $\Cb_w\:k\coalg_\cdg^\conilp\rarrow k\alg_\dg^+$
is left adjoint to the functor of DG\+algebra bar-construction
$\Br_v\:k\alg_\dg^+\rarrow k\coalg_\cdg^\conilp$.
 Analogously, the functor of conilpotent DG\+coalgebra
cobar-construction $\Cb_w\:k\coalg_\dg^\conilp\rarrow k\alg_\dg^\aug$
is right adjoint to the functor of augmented DG\+algebra
bar-construction $\Br_v\:k\alg_\dg^\aug\rarrow k\coalg_\dg^\conilp$.

 A morphism of conilpotent CDG\+coalgebras $(f,a)\:C\rarrow D$
is called a \emph{filtered quasi-isomorphism} if there exist
increasing filtrations $F$ on $C$ and $D$ satisfying the following
conditions.
 The filtrations $F$ must be compatible with the comultiplications
and differentials on $C$ and $D$; one must have $F_0C = w_C(k)$ and
$F_0D = w_D(k)$, so that, in particular, the associated quotient
objects $\gr_FC$ and $\gr_FD$ are DG\+coalgebras; and the induced
morphism $\gr_Ff\:\gr_FC\rarrow\gr_FD$ must be a quasi-isomorphism of
graded complexes of vector spaces.
 A morphism of conilpotent DG\+coalgebras is a filtered
quasi-isomorphism if it is a filtered quasi-isomorphism as
a morphism of conilpotent CDG\+coalgebras.
 The classes of filtered quasi-isomorphisms will be denoted by
$\FQuis\subset k\coalg_\cdg^\conilp$ and $\FQuis\subset
k\coalg_\dg^\conilp$.
 Let us emphasize that there is \emph{no} claim that the classes
of filtered quasi-isomorphisms are closed under composition of
morphisms.
 The classes of quasi-isomorphisms of DG\+algebras and augmented
DG\+algebras will be denoted by $\Quis\subset k\alg_\dg^+$ and
$\Quis\subset k\alg_\dg^\aug$.

\begin{thm}
 \textup{(a)} The functors\/ $\Cb_w\:k\coalg_\cdg^\conilp\rarrow
k\alg_\dg^+$ and\/ $\Br_v\:k\alg_\dg^+\rarrow k\coalg_\cdg^\conilp$
induce functors between the localized categories\/
$k\coalg_\cdg^\conilp[\FQuis^{-1}]\rarrow k\alg_\dg^+[\Quis^{-1}]$ and\/
$k\alg_\dg^+[\Quis^{-1}]\rarrow k\coalg_\cdg^\conilp[\FQuis^{-1}]$,
which are mutually inverse equivalences of categories. \par
 \textup{(b)} The functors\/ $\Cb_w\:k\coalg_\dg^\conilp\rarrow
k\alg_\dg^\aug$ and\/ $\Br_v\:k\alg_\dg^\aug\rarrow k\coalg_\dg^\conilp$
induce functors between the localized categories\/
$k\coalg_\dg^\conilp[\FQuis^{-1}]\rarrow k\alg_\dg^\aug[\Quis^{-1}]$
and\/ $k\alg_\dg^\aug[\Quis^{-1}]\rarrow k\coalg_\dg^\conilp
[\FQuis^{-1}]$, which are mutually inverse equivalences of categories.
\end{thm}

\begin{proof}
 We will prove part~(a); the proof of part~(b) is similar.
 For any filtered quasi-isomorphism of conilpotent CDG\+coalgebras
$(f,a)\:C\rarrow D$ the induced morphism of cobar-constructions
$\Cb_w(f,a)\:\Cb_w(C)\rarrow\Cb_w(D)$ is a quasi-isomorphism of
DG\+algebras.
 Indeed, let $F$ denote the increasing filtrations on
the cobar-constructions induced by the filtrations $F$ of $C$ and $D$;
then the morphism of associated graded DG\+algebras
$\gr_F\Cb_w(f,a)$ is a quasi-isomorphism, since the tensor products
and the cones of morphisms of complexes preserve quasi-isomorphisms.
 Conversely, for any quasi-isomorphism of DG\+algebras $g\:A\rarrow B$
the induced morphism of bar-constructions $\Br_v(g)\:\Br_v(A)\rarrow
\Br_v(B)$ is a filtered quasi-isomorphism.
 Indeed, it suffices to consider the increasing filtrations of
bar-constructions associated with their nonnegative gradings~$n$
by the number of factors in tensor powers.
 So the induced functors exist; it remains to check that they are
mutually inverse equivalences.
 For any DG\+algebra $A$, the adjunction morphism $\Cb_w(\Br_v(A))
\rarrow A$ is a quasi-isomorphism.
 One can prove this by passing to the associated quotients with
respect to the increasing filtration $F$ on $A$ defined by
the rules $F_0A=k$ and $F_1A=A$, and the induced filtration on
$\Cb_w(\Br_v(A))$.
 Finally, for any conilpotent CDG\+coalgebra $C$, the adjunction
morphism $C\rarrow\Br_v(\Cb_w(C))$ is a filtered quasi-isomorphism.
 Indeed, consider the natural increasing filtration $F$ on $C$
defined in~\ref{conilpotent-cobar-duality} and the induced
filtration $F$ on $\Br_v(\Cb_w(C))$.
 We have to prove that our adjunction morphism becomes
a quasi-isomorphism after passing to the associated quotient
objects, i.~e., the morphism of graded DG\+coalgebras
$\gr_FC\rarrow\Br_v(\Cb_w(\gr_FC))$ is a quasi-isomorphism.
 Here it suffices to consider the decreasing filtration $G$ on
$\gr_FC$ defined by the rules $G^0\gr_FC=\gr_FC$, \ 
$G^1\gr_FC=\ker(\eps\:\gr_FC\to k)$, and $G^2\gr_FC=0$.
 The induced filtration on $\Br_v(\Cb_w(\gr_FC))$ stabilizes at
every degree of the nonnegative grading coming from
the filtration~$F$ and the morphism
$\gr_FC\rarrow\Br_v(\Cb_w(\gr_FC))$ can be easily seen to become
a quasi-isomorphism after passing to the associated quotient
objects with respect to the filtration~$G$.
\end{proof}

\begin{rem}
 Notice that the notion of a filtered quasi-isomorphism makes sense
for conilpotent CDG\+coalgebras only, as any CDG\+coalgebra
admitting an increasing filtration $F$ satisfying the conditions in
the definition of a filtered quasi-isomorphism is conilpotent.
 And one cannot even speak about conventional (nonfiltered)
quasi-isomorphisms of CDG\+coalgebras, as the latter are not
complexes.
 Furthermore, the assertions of Theorem do not hold with
the filtered quasi-isomorphisms replaced with conventional
quasi-isomorphisms of DG\+coalgebras, or with the conilpotency
condition dropped.
 Indeed, let $A$ be any DG\+algebra; consider it as a DG\+algebra
without unit and add a unit formally to it, obtaining an augmented
DG\+algebra $k\oplus A$ with the augmentation~$v$.
 Then the morphisms of augmented DG\+algebras $k\rarrow k\oplus A
\rarrow k$ induce quasi-isomorphisms of bar-constructions
$\Br_v(k)\rarrow\Br_v(k\oplus A)\rarrow\Br_v(k)$; applying
the cobar-construction, we find that the morphisms
$\Cb_w(\Br_v(k))\rarrow\Cb_w(\Br_v(k\oplus A))\rarrow\Cb_w(\Br_v(k))$
are \emph{not} quasi-isomorphisms, since the middle term is
quasi-isomorphic to $k\oplus A$.
 Analogously, let $D$ be any DG\+coalgebra; consider it as
a DG\+coalgebra without counit and add a counit formally to it,
obtaining a coaugmented DG\+coalgebra $k\oplus D$ with
the coaugmentation~$w$.
 Then the morphisms of augmented DG\+algebras $k\rarrow\Cb_w(k\oplus D)
\rarrow k$ are quasi-isomorphisms and it follows that the induced
morphisms of bar-constructions $\Br_v(k)\rarrow\Br_v(\Cb_w(k\oplus D))
\rarrow\Br_v(k)$ are also quasi-isomorphisms, hence the cohomology of
the DG\+coalgebra $\Br_v(\Cb_w(k\oplus D))$ is different from that
of~$k\oplus D$.
 And there is even no natural morphism between $C$ and
$\Br_v(\Cb_w(C))$ for a nonconilpotent DG\+coalgebra $C$.
 Finally, let $D$ be a DG\+coalgebra and $N$ be a left DG\+comodule
over $D$.
 Consider $N$ as a DG\+comodule over the above DG\+coalgebra
$k\oplus D$; then the DG\+module $\Cb_w(k\oplus D)
\ot^{\tau_{k\oplus D,w}}\!\.N$ over $\Cb_w(k\oplus D)$ is acyclic.
 It follows that the assertions of
Theorem~\ref{conilpotent-cobar-duality} and
Corollary~\ref{nonconilpotent-duality} do not hold without
the conilpotency assumption on the coaugmented CDG\+coalgebra~$C$.
\end{rem}

\Section{$\Ainfty$-Algebras and Curved $\Ainfty$-Coalgebras}

\subsection{Nonunital $\Ainfty$\+algebras} \label{non-1-ainfty-algebras}
 Let $A$ be a graded vector space over a field~$k$.
 Consider the graded tensor coalgebra (cofree conilpotent coassociative
graded coalgebra) $\bigoplus_{n=0}^\infty A[1]^{\ot n}$ with its
coaugmentation $w\:k\simeq A[1]^{\ot 0}\rarrow \bigoplus_n A[1]^{\ot n}$.
 A \emph{nonunital $\Ainfty$\+algebra} structure on $A$ is, by
the definition, a coaugmented DG\+coalgebra structure on
$\bigoplus_n A[1]^{\ot n}$, i.~e., an odd coderivation~$d$ of degree~$1$
on $\bigoplus_n A[1]^{\ot n}$ such that $d^2=0$ and $d\circ w=0$.
 Since a coderivation of $\bigoplus_n A[1]^{\ot n}$ is uniquely
determined by its composition with the projection $\bigoplus_n A[1]^
{\ot n}\rarrow A[1]^{\ot 1}\simeq A[1]$, a nonunital $\Ainfty$\+algebra
structure on $A$ can be considered as a sequence of linear maps
$m_n\:A^{\ot n}\rarrow A$, \ $n=1$, $2$,~\dots\ of degree $2-n$.
 More precisely, define the maps $m_n$ by the rule that the image of
the element $d(a_1\ot\dsb\ot a_n)$ under the projection to $A$ equals
$(-1)^{n+\sum_{j=1}^n(n-j)(|a_j|+1)}m_n(a_1\ot\dsb\ot a_n)$ for
$a_j\in A$.
 The sequence of maps $m_n$ must satisfy a sequence of quadratic
equations corresponding to the equation $d^2=0$ on the coderivation~$d$.
 We will not write down these equations explicitly.

 A morphism of nonunital $\Ainfty$\+algebras $f\:A\rarrow B$ over~$k$
is, by the definition, a morphism of (coaugmented) DG\+coalgebras
$\bigoplus_n A[1]^{\ot n}\rarrow \bigoplus_n B[1]^{\ot n}$.
 Since a graded coalgebra morphism into a graded tensor coalgebra
$\bigoplus_n B[1]^{\ot n}$ is determined by its composition
with the projection $\bigoplus_n B[1]^{\ot n}\rarrow B[1]^{\ot 1}
\simeq B[1]$ and any morphism of conilpotent graded coalgebras preserves
coaugmentations, a morphism of nonunital $\Ainfty$\+algebras
$f\:A\rarrow B$ can be considered as a sequence of linear maps
$f_n\:A^{\ot n}\rarrow B$, \ $n=1$, $2$,~\dots\ of degree $1-n$. 
 More precisely, define the maps $f_n$ by the rule that the image of
the element $f(a_1\ot\dsb\ot a_n)$ under the projection to $B$ equals
$(-1)^{n-1+\sum_{j=1}^n(n-j)(|a_j|+1)}f_n(a_1\ot\dsb\ot a_n)$ for
$a_j\in A$.
 The sequence of maps $f_n$ must satisfy a sequence of polynomial
equations corresponding to the equation $d\circ f=f\circ d$ on
the morphism~$f$.

 Let $A$ be a nonunital $\Ainfty$\+algebra over a field~$k$ and $M$
be a graded vector space over~$k$.
 A structure of \emph{nonunital left $\Ainfty$\+module} over $A$
on $M$ is, by the definition, a structure of DG\+comodule over
the DG\+coalgebra $\bigoplus_n A[1]^{\ot n}$ on the cofree graded
left comodule $\bigoplus_n A[1]^{\ot n}\ot_kM$ over the graded
coalgebra $\bigoplus_n A[1]^{\ot n}$.
 Analogously, a structure of \emph{nonunital right $\Ainfty$\+module}
over $A$ on a graded vector space $N$ is defined as a structure of
DG\+comodule over $\bigoplus_n A[1]^{\ot n}$ on the cofree graded
right comodule $N\ot_k\bigoplus_n A[1]^{\ot n}$.
 Since a coderivation of a cofree graded comodule $\bigoplus_n
A[1]^{\ot n}\ot_kM$ compatible with a given coderivation of
the graded coalgebra $\bigoplus_n A[1]^{\ot n}$ is determined by
its composition with the projection $\bigoplus_n A[1]^{\ot n}\ot_kM
\rarrow M$ induced by the counit map $\bigoplus_n A[1]^{\ot n}
\rarrow k$, a nonunital left $\Ainfty$\+module structure on $M$ can be
considered as a sequence of linear maps $l_n\:A^{\ot n}\ot M\rarrow M$,
\ $n=0$, $1$,~\dots\ of degree $1-n$.
 More precisely, define the maps $l_n$ by the rule that the image of
the element $d(a_1\ot\dsb\ot a_n\ot x)$ under the projection to $M$
equals $(-1)^{n+\sum_{j=1}^n(n-j)(|a_j|+1)}l_n(a_1\ot\dsb\ot a_n\ot x)$
for $a_j\in A$ and $x\in M$.
 The sequence of maps $l_n$ must satisfy a system of nonhomogeneous
quadratic equations corresponding to the equation $d^2=0$ on
the coderivation~$d$ on $\bigoplus_n A[1]^{\ot n}\ot_kM$.
 Analogously, a nonunital right $\Ainfty$\+module structure on $N$ can
be considered as a sequence of linear maps $r_n\:N\ot A^{\ot n}\rarrow
N$ defined by the rule that the image of the element
$d(y\ot a_1\ot\dsb\ot a_n)$ under the projection $N\ot_k\bigoplus_n
A[1]^{\ot n}\rarrow N$ equals $(-1)^{n|y|+\sum_{j=1}^n(n-j)(|a_j|+1)}
r_n(y\ot a_1\ot\dsb\ot a_n)$ for $a_j\in A$ and $y\in N$.

 The complex of morphisms between nonunital left $\Ainfty$\+modules
$L$ and $M$ over a nonunital $\Ainfty$\+algebra $A$ is, by
the definition, the complex of morphisms between left DG\+comodules
$\bigoplus_n A[1]^{\ot n}\ot_kL$ and $\bigoplus_n A[1]^{\ot n}\ot_kM$
over the DG\+coalgebra $\bigoplus_n A[1]^{\ot n}$.
 Analogously, the complex of morphisms between nonunital right
$\Ainfty$\+modules $R$ and $N$ over $A$ is, by the definition,
the complex of morphisms between right DG\+comodules
$R\ot_k\bigoplus_n A[1]^{\ot n}$ and $N\ot_k\bigoplus_n A[1]^{\ot n}$
over the DG\+coalgebra $\bigoplus_n A[1]^{\ot n}$.
 A morphism of nonunital left $\Ainfty$\+modules $f\:L\rarrow M$ of
degree~$i$ is the same that a sequence of linear maps $f_n\:A^{\ot n}
\ot_k L \rarrow M$, \ $n=0$, $1$,~\dots\ of degree $i-n$.
 More precisely, define the maps $f_n$ by the rule that the image of
the element $f(a_1\ot\dsb\ot a_n\ot x)$ under the projection
$\bigoplus_n A[1]^{\ot n}\ot_kM\rarrow M$ equals $(-1)^{n+\sum_{j=1}^n
(n-j)(|a_j|+1)}f_n(a_1\ot\dsb\ot a_n\ot x)$ for $a_j\in A$ and $x\in L$.
 Any sequence of linear maps $f_n$ corresponds to a (not necessarily
closed) morphism of nonunital $\Ainfty$\+modules~$f$.
 Analogously, a morphism of nonunital right $\Ainfty$\+modules
$g\:R\rarrow N$ of degree~$i$ is the same that a sequence of linear
maps $g_n\:A^{\ot n}\ot_k R\rarrow N$ of degree~$i-n$.
 More precisely, define the maps $g_n$ by the rule that the image of
the element $g(y\ot a_1\ot\dsb\ot a_n)$ under the projection
$N\ot_k\bigoplus_nA[1]^{\ot n}\rarrow N$ equals
$(-1)^{n|y|+\sum_{j=1}^n(n-j)(|a_j|+1)}g_n(y\ot a_1\ot\dsb \ot a_n)$
for $a_j\in A$ and $y\in R$.

 For any CDG\+coalgebra $C$, the functors $\Phi_C$ and $\Psi_C$
of~\ref{functors-phi-and-psi}--\ref{co-contra-corr-theorem} provide
an equivalence between the DG\+category of left CDG\+comodules over $C$
that are cofree as graded $C$\+comodules and the DG\+category of left
CDG\+contramodules over $C$ that are free as graded $C$\+contramodules.
 So one can alternatively define a nonunital left
$\Ainfty$\+module $M$ over a nonunital $\Ainfty$\+algebra $A$
as a graded vector space for which a structure of DG\+contramodule
over the DG\+coalgebra $\bigoplus_n A[1]^{\ot n}$ is given on
the free graded contramodule $\Hom_k(\bigoplus_n A[1]^{\ot n},M)$
over the graded coalgebra $\bigoplus_n A[1]^{\ot n}$.
 Since a contraderivation of a free graded contramodule
$\Hom_k(\bigoplus_nA[1]^{\ot n},M)$ compatible with a given
coderivation of the graded coalgebra $\bigoplus_n A[1]^{\ot n}$
is determined by its restriction to the graded subspace
$M\subset \Hom_k(\bigoplus_n A[1]^{\ot n},M)$, a nonunital left
$\Ainfty$\+module structure on $M$ can be considered as a sequence
of linear maps $p_n\:M\rarrow\Hom_k(A^{\ot n},M)$, \ $n=0$,
$1$,~\dots\ of degree~$1-n$.
 More precisely, define the maps $p_n$ by the formula
$p_n(x)(a_1\ot\dsb\ot a_n) = (-1)^{n+n|x|+\sum_{j=1}^n(n-j)(|a_j|+1)}
d(x)(a_1\ot\dsb\ot a_n)$ for $a_j\in A$ and $x\in M$.
 Then the maps $p_n$ are related to the above maps $l_n\:A^{\ot n}\ot_kM
\rarrow M$ by the rule $p_n(x)(a_1\ot\dsb\ot a_n) =
(-1)^{|x|\sum_{j=1}^n|a_j|}l_n(a_1\ot\dsb \ot a_n\ot x)$.

 Furthermore, one can alternatively define the complex of morphisms
between nonunital left $\Ainfty$\+modules $L$ and $M$ over a nonunital
$\Ainfty$\+algebra $A$ as the complex of morphisms between left
DG\+contramodules $\Hom_k(\bigoplus_n A[1]^{\ot n},L)$ and
$\Hom_k(\bigoplus_n A[1]^{\ot n},M)$ over the DG\+coalgebra
$\bigoplus_n A[1]^{\ot n}$.
 Thus a (not necessarily closed) morphism of nonunital left
$\Ainfty$\+modules $f\:L\rarrow M$ of degree $i$ is the same that
a sequence of linear maps $f^n\:L \rarrow \Hom_k(A^{\ot n},M)$, \ 
$n=0$, $1$,~\dots\ of degree $i-n$.
 More precisely, define the maps $f^n$ by the formula 
$f^n(x)(a_1\ot\dsb\ot a_n) = (-1)^{n+n|x|+\sum_{j=1}^n(n-j)(|a_j|+1)}
f(x)(a_1\ot\dsb\ot a_n)$ for $a_j\in A$ and $x\in L$, where $L$ is
considered as a graded subspace in $\Hom_k(\bigoplus_n A[1]^{\ot n},L)$.
 Then the maps $f^n$ are related to the above maps $f_n\:A^{\ot n}
\ot_k L\rarrow M$ by the rule $f^n(x)(a_1\ot\dsb\ot a_n) =
(-1)^{|x|\sum_{j=1}^n|a_j|}f_n(a_1\ot\dsb \ot a_n\ot x)$.

\subsection{Strictly unital $\Ainfty$\+algebras}
\label{strictly-1-ainfty-algebras}
 Let $A$ be a nonunital $\Ainfty$\+algebra over a field~$k$.
 An element $1\in A$ of degree~$0$ is called a \emph{strict unit}
if one has $m_2(1\ot a)=a=m_2(a\ot 1)$ for all $a\in A$ and
$m_n(a_1\ot\dsb\ot a_{j-1}\ot 1\ot a_{j+1}\ot\dsb\ot a_n)=0$
for all $n\ne2$, \ $1\le j\le n$, and $a_t\in A$.
 Obviously, a strict unit is unique if it exists.
 A \emph{strictly unital $\Ainfty$\+algebra} is a nonunital
$\Ainfty$\+algebra that has a strict unit.
 A morphism of strictly unital $\Ainfty$\+algebras $f\:A\rarrow B$
is a morphism of nonunital $\Ainfty$\+algebras
such that $f_1(1_A)=1_B$ and $f_n(a_1\ot\dsb\ot a_{j-1}\ot 1_A
\ot a_{j+1}\ot\dsb\ot a_n)=0$ for all $n>1$ and $a_t\in A$.
 Notice that for a strictly unital $\Ainfty$\+algebra $A$ with
the unit $1_A$ one has $1_A=0$ if and only if $A=0$.
 We will assume our strictly unital $\Ainfty$\+algebras to have
nonzero units.

 A \emph{strictly unital left $\Ainfty$\+module} $M$ over
a strictly unital $\Ainfty$\+algebra $A$ is a nonunital left
$\Ainfty$\+module such that $l_1(1\ot x)=x$ and 
$l_n(a_1\ot\dsb\ot a_{j-1}\ot 1\ot a_{j+1}\ot\dsb\ot a_n\ot x)=0$
for all $n>1$, \ $1\le j\le n$, \ $a_t\in A$, and $x\in M$.
 Equivalently, one must have $p_1(x)(1)=x$ and $p_n(x)
(a_1\ot\dsb\ot a_{j-1}\ot 1\ot a_{j+1}\ot\dsb\ot a_n)=0$.
 Analogously, a \emph{strictly unital right $\Ainfty$\+module} $N$
over $A$ is a nonunital right $\Ainfty$\+module such that
$r_1(y\ot 1)=y$ and $r_n(y\ot a_1\ot\dsb\ot a_{j-1}\ot 1\ot a_{j+1}
\ot\dsb\ot a_j)=0$ for all $n>1$, \ $a_t\in A$, and $y\in N$.

 The complex of morphisms between strictly unital left
$\Ainfty$\+modules $L$ and $M$ over a strictly unital
$\Ainfty$\+algebra $A$ is the subcomplex of the complex of
morphisms between $L$ and $M$ as nonunital $\Ainfty$\+modules
consisting of all morphisms $f\:L\rarrow M$ such that
$f_n(a_1\ot\dsb\ot a_{j-1}\ot 1\ot a_{j+1}\ot\dsb\ot a_n\ot x)=0$
for all $n>0$,  \ $1\le j\le n$, \ $a_t\in A$, and $x\in L$.
 Equivalently, one must have $f^n(x)(a_1\ot\dsb\ot a_{j-1}\ot 1\ot
a_{j+1}\ot\dsb\ot a_n)=0$.
 Analogously, the complex of morphisms between strictly unital right
$\Ainfty$\+modules $R$ and $N$ over $A$ is the subcomplex of
the complex of morphisms between $R$ and $N$ as nonunital
$\Ainfty$\+modules consisting of all morphisms $g\:R\rarrow N$ such
that $g_n(y\ot a_1\ot\dsb\ot a_{j-1}\ot 1\ot a_{j+1}\ot\dsb\ot a_n)=0$
for all $n>0$, \ $a_t\in A$, and $y\in R$.

 Let $A$ be a nonunital $\Ainfty$\+algebra and $1_A\in A$ be a nonzero
element of degree~$0$.
 Set $A_+=A/k\cdot 1_A$.
 Then the graded tensor coalgebra $\bigoplus_n A_+[1]^{\ot n}$ is 
a quotient coalgebra of the tensor coalgebra $\bigoplus_n A[1]^{\ot n}$.
 Denote by $K_A$ the kernel of the natural surjection
$\bigoplus_n A[1]^{\ot n}\rarrow\bigoplus_n A_+[1]^{\ot n}$ and
by $\kappa_A\:K_A\rarrow k$ the homogeneous linear function of
degree~$1$ sending $1_A\in K_A\cap A[1]$ to $1\in k$ and annihilating
$K_A\cap A[1]^{\ot n}$ for all $n>1$.
 Let $\theta_A\:\bigoplus_n A[1]^{\ot n}\rarrow k$ be any homogeneous
linear function of degree~$1$ extending the linear function $\kappa_A$
on~$K_A$.
 Then the element $1_A\in A$ is a strict unit if and only if
the odd coderivation $d'(c) = d(c)+\theta_A*c-(-1)^{|c|}c*\theta_A$
of degree~$1$ on the tensor coalgebra $\bigoplus_n A[1]^{\ot n}$
preserves the subspace $K_A$ and the linear function $h'(c) = 
\theta_A(d(c))+\theta_A^2(c)$ of degree~$2$ on
$\bigoplus_n A[1]^{\ot n}$ annihilates $K_A$.
 This condition does not depend on the choice of~$\theta_A$.
 For strictly unital $\Ainfty$\+algebras $A$ and $B$, a morphism of
nonunital $\Ainfty$\+algebras $f\:A\rarrow B$ is a morphism of
strictly unital $\Ainfty$\+algebras if and only if $f(K_A)\subset K_B$
and $\kappa_B\circ f|_{K_A} = \kappa_A$.

 Let $A$ be a strictly unital $\Ainfty$\+algebra and $M$ be
a nonunital left $\Ainfty$\+module over~$A$.
 Then $M$ is a strictly unital $\Ainfty$\+module if and only if
the odd coderivation $d'(z)=d(z)+\theta_A*z$ of degree~$1$ on
the cofree comodule $\bigoplus_n A[1]^{\ot n}\ot_k M$ compatible with
the coderivation $d'$ of the coalgebra $\bigoplus_n A[1]^{\ot n}$
preserves the subspace $K_A\ot_k M\subset\bigoplus_n A[1]^{\ot n}
\ot_k M$.
 Equivalently, the odd contraderivation $d'(q)=d(q)+\theta_A*q$ of
degree~$1$ on the free contramodule $\Hom_k(\bigoplus_n A[1]^{\ot n},
M)$ compatible with the coderivation $d'$ of the coalgebra
$\bigoplus_n A[1]^{\ot n}$ must preserve the subspace
$\Hom_k(\bigoplus_n A_+[1]^{\ot n},M)\subset
\Hom_k(\bigoplus_n A[1]^{\ot n},M)$.
 Analogously, a nonunital right $\Ainfty$\+module $N$ over~$A$ is
a strictly unital $\Ainfty$\+module if and only if the odd
coderivation $d'(z)=d(z)-(-1)^{|z|}z*\theta_A$ of degree~$1$ on
the cofree comodule $N\ot_k \bigoplus_n A[1]^{\ot n}$ compatible with
the coderivation $d'$ of the coalgebra $\bigoplus_n A[1]^{\ot n}$
preserves the subspace $N\ot_k K_A\subset N\ot_k \bigoplus_n
A[1]^{\ot n}$.
 For strictly unital left $\Ainfty$\+modules $L$ and $M$ over $A$,
a (not necessarily closed) morphism of nonunital $\Ainfty$\+modules
$f\:L\rarrow M$ is a morphism of strictly unital $\Ainfty$\+modules
if and only if one has $f(K_A\ot_k L)\subset K_A\ot_k M$, or
equivalently, $f(\Hom_k(\bigoplus_n A_+[1]^{\ot n},L))\subset
\Hom_k(\bigoplus_n A_+[1]^{\ot n},M)$.
 Analogously, for strictly unital right $\Ainfty$\+modules $R$ and $N$
over $A$, a (not necessarily closed) morphism of nonunital
$\Ainfty$\+modules $g\:R\rarrow N$ is a morphism of strictly unital
$\Ainfty$\+modules if and only if one has $g(R\ot_k K_A)\subset
N\ot_k K_A$.

 Let $A$ be a strictly unital $\Ainfty$\+algebra.
 Identify $k$ with the subspace $k\cdot 1_A\subset A$ and choose
a homogeneous $k$\+linear retraction $v\:A\rarrow k$.
 Define the homogeneous linear function $\theta_A\:\bigoplus_n A[1]^
{\ot n}\rarrow k$ of degree~$1$ by the rules $\theta_A(a)=v(a)$ and
$\theta_A(a_1\ot\dsb\ot a_n)=0$ for $n\ne 1$.
 Then the linear function $\theta_A$ is an extension of
the linear function $\kappa_A\:K_A\rarrow k$.
 Let $d\:\bigoplus_n A_+[1]^{\ot n}\rarrow \bigoplus_n A_+[1]^{\ot n}$
be the map induced by the odd coderivation $d'$ of
$\bigoplus_n A[1]^{\ot n}$ defined by the above formula, and let
$h\:\bigoplus_n A_+[1]^{\ot n}\rarrow k$ be the linear function 
induced by the above linear function~$h'$.
 Then $\Br_v(A)=(\bigoplus_n A_+[1]^{\ot n}\;d\;h)$ is a coaugmented
(and consequenly, conilpotent) CDG\+coalgebra with the coaugmentation
$k\simeq A_+[1]^{\ot 0}\rarrow \bigoplus_n A_+[1]^{\ot n}$.
 The CDG\+coalgebra $\Br_v(A)$ is called the \emph{bar-construction}
of a strictly unital $\Ainfty$\+algebra~$A$.

 Let $f\:A\rarrow B$ be a morphism of strictly unital
$\Ainfty$\+algebras.
 Let $v\:A\rarrow k$ and $v\:B\rarrow k$ be homogeneous $k$\+linear
retractions, and let $\theta_A$ and $\theta_B$ be the corresponding
homogeneous linear functions of degree~$1$ on the graded tensor
coalgebras.
 The morphism of tensor coalgebras $f\:\bigoplus_n A[1]^{\ot n}
\rarrow\bigoplus_n B[1]^{\ot n}$ maps $K_A$ into $K_B$, so it induces
a morphism of graded tensor coalgebras $\bigoplus_n A_+[1]^{\ot n}
\rarrow\bigoplus_n B_+[1]^{\ot n}$, which we will denote also by~$f$.
 The linear function $\theta_B\circ f-\theta_A\:\bigoplus_n A[1]^
{\ot n}\rarrow k$ annihilates $K_A$, so it induces a linear function
$\bigoplus_n A_+[1]^{\ot n}\rarrow k$, which we will denote by~$\eta_f$.
 Then the pair $(f,\eta_f)$ is a morphism of CDG\+coalgebras
$\Br_v(A)\rarrow\Br_v(B)$.
 Thus the bar-construction $A\mpsto \Br_v(A)$ is a functor from
the category of strictly unital $\Ainfty$\+algebras with nonzero
units to the category of coaugmented CDG\+coalgebras whose
underlying graded coalgebras are graded tensor coalgebras.
 One can easily see that this functor is an equivalence of categories.
 Alternatively, one can use any linear function~$\theta_A$ of
degree~$1$ extending the linear function~$\kappa_A$ in the construction
of this equivalence of categories.

 To obtain the inverse functor, assign to a conilpotent CDG\+coalgebra
$(D,d_D,h_D)$ the conilpotent DG\+coalgebra $(C,d_C)$ constructed as
follows.
 First, adjoin to $D$ a single cofree cogenerator of degree~$-1$,
obtaining a conilpotent graded coalgebra $C$ endowed with a graded
coalgebra morphism $C\rarrow D$ and a homogeneous linear function
$\theta\:C\rarrow k$ of degree~$1$.
 Second, define the odd coderivation $d'_C$ of degree~$1$ on
the graded coalgebra $C$ by the conditions that $d'_C$ must preserve
the kernel of the graded coalgebra morphism $C\rarrow D$ and induce
the differential $d_D$ on~$D$, and that the equation $\theta(d'_C(c))
= \theta^2(c) + h_D(c)$ must hold for all $c\in C$, where $h_D$ is
considered as a linear function on~$C$.
 Finally, set $d_C(c)=d'_C(c)-\theta*c+(-1)^{|c|}c*\theta$ for all
$c\in C$.

 Let $A$ be a strictly unital $\Ainfty$\+algebra, $v\:A\rarrow k$ be
a homogeneous $k$\+linear retraction, and $\theta_A\:\bigoplus_n A[1]^
{\ot n}\rarrow k$ be the corresponding homogeneous linear function
of degree~$1$.
 Let $M$ be a strictly unital left $\Ainfty$\+module over~$A$.
 Set $d\:\bigoplus_n A_+[1]^{\ot n}\ot_k M\rarrow\bigoplus_n
A_+[1]^{\ot n}\ot_k M$ to be the map induced by the differential~$d'$
on $\bigoplus_n A[1]^{\ot n}\ot_k M$ defined by the above formula.
 Then $\Br_v(A,M)=(\bigoplus_n A_+[1]^{\ot n}\ot_k M\;d)$ is a left
CDG\+comodule over the CDG\+coalgebra $\Br_v(A)$.
 Furthermore, set $d\:\Hom_k(\bigoplus_n A_+[1]^{\ot n},M)\rarrow
\Hom_k(\bigoplus_n A_+[1]^{\ot n},M)$ to be the restriction of
the differential $d'$ on $\Hom_k(\bigoplus_n A[1]^{\ot n},M)$
defined above.
 Then $\Cb^v(A,M)=(\Hom_k(\bigoplus_n A_+[1]^{\ot n},M)\;d)$ is
a left CDG\+contramodule over the CDG\+coalgebra $\Br_v(A)$.
 Analogously, for a strictly unital right $\Ainfty$\+module $N$
over~$A$ set $d\:N\ot_k\bigoplus_n A_+[1]^{\ot n}\rarrow N\ot_k
\bigoplus_n A_+[1]^{\ot n}$ to be the map induced by
the differential~$d'$ on $N\ot_k\bigoplus_n A[1]^{\ot n}$ defined above.
 Then $\Br_v(N,A)=(N\ot_k\bigoplus_n A_+[1]^{\ot n}\;d)$ is a right
CDG\+comodule over the CDG\+coalgebra $\Br_v(A)$.

 To a (not necessarily closed) morphism of strictly unital left
$\Ainfty$\+modules $f\:L\rarrow M$ over $A$ one can assign the induced
maps $\bigoplus_n A_+[1]^{\ot n}\ot_k L\rarrow \bigoplus_n A_+[1]^
{\ot n}\ot_k M$ and $\Hom_k(\bigoplus_n A_+[1]^{\ot n},L)\rarrow
\Hom_k(\bigoplus_n A_+[1]^{\ot n},M)$.
 These are a (not necessarily closed) morphism of CDG\+comodules
$\Br_v(A,L)\rarrow\Br_v(A,M)$ and a (not necessarily closed) morphism
of CDG\+contramodules $\Cb^v(A,L)\rarrow\Cb^v(A,M)$ over
the CDG\+coalgebra $\Br_v(A)$.
 So we obtain the DG\+functor $M\mpsto \Br_v(A,M)$, which is
an equivalence between the DG\+category of strictly unital left
$\Ainfty$\+modules over $A$ and the DG\+category of left CDG\+comodules
over $\Br_v(A)$ that are cofree as graded comodules, and the DG\+functor
$M\mpsto\Cb^v(A,M)$, which is an equivalence between the DG\+category
of strictly unital left $\Ainfty$\+modules over $A$ and the DG\+category
of left CDG\+contramodules over $\Br_v(A)$ that are free as graded
contramodules.
 These two equivalences of DG\+categories form a commutative diagram
with the equivalence between the DG\+category of CDG\+comodules
that are cofree as graded comodules and the DG\+category of
CDG\+contramodules that are free as graded contramodules provided
by the functors $\Psi_{\Br_v(A)}$ and $\Phi_{\Br_v(A)}$.
 Analogously, to a (not necessarily closed) morphism of strictly
unital right $\Ainfty$\+modules $g\:R\rarrow N$ over $A$ one can assign
the induced map $R\ot_k\bigoplus_n A_+[1]^{\ot n}\rarrow 
N\ot_k\bigoplus_n A_+[1]^{\ot n}$.
 This is a (not necessarily closed) morphism of CDG\+comodules
$\Br_v(R,A)\rarrow\Br_v(N,A)$.
 The DG\+functor $N\mpsto\Br_v(N,A)$ is an equivalence between
the DG\+category of strictly unital right $\Ainfty$\+modules over $A$
and the DG\+category of right CDG\+comodules over $\Br_v(A)$ that are
cofree as graded comodules.

 Now let $A$ be a DG\+algebra with nonzero unit, $M$ be a left 
DG\+module over $A$, and $N$ be a right DG\+module over~$A$.
 Let $v\:A\rarrow k$ be a homogeneous $k$\+linear retraction.
 Define a strictly unital $\Ainfty$\+algebra structure on $A$ by
the rules $m_1(a)=d(a)$, \ $m_2(a_1\ot a_2) = a_1a_2$, and $m_n=0$
for $n>2$.
 Define a structure of a strictly unital left $\Ainfty$\+module
over $A$ on $M$ by the rules $l_0(x)=d(x)$, \ $l_1(a\ot x)=ax$, and
$l_n=0$ for $i>1$, where $a\in A$ and $x\in M$.
 Analogously, define a structure of a strictly unital right
$\Ainfty$\+module over $A$ on $N$ by the rules $r_0(y)=d(y)$, \ 
$r_1(y\ot a)=ya$, and $r_n=0$ for $n>1$, where $a\in A$ and $y\in N$.
 Then the CDG\+coalgebra structure $\Br_v(A)$ on the graded tensor
coalgebra $\bigoplus_n A_+[1]^{\ot n}$ that was defined
in~\ref{bar-cobar-constr} coincides with the CDG\+coalgebra
structure $\Br_v(A)$ constructed above, so our notation is consistent.
 The left CDG\+comodule structure $\Br_v(A)\ot^{\tau_{A,v}}\!\.M$ on
the cofree graded comodule $\bigoplus_n A_+[1]^{\ot n}\ot_k M$ that
was defined in~\ref{twisting-cochains-subsect} coincides with
the left CDG\+comodule structure $\Br_v(A,M)$.
 The left CDG\+contramodule structure $\Hom^{\tau_{A,v}}(\Br_v(A),M)$
on the free graded contramodule $\Hom_k(\bigoplus_n A_+[1]^{\ot n},M)$
coincides with the CDG\+contramodule structure $\Cb^v(A,M)$.
 The right CDG\+comodule structure $N\ot^{\tau_{A,v}}\!\.\Br_v(A)$ on
the cofree graded comodule $N\ot_k\bigoplus_n A_+[1]^{\ot n}$ coincides
with the right CDG\+comodule structure $\Br_v(N,A)$.

 A morphism of strictly unital $\Ainfty$\+algebras $f\:A\rarrow B$
is called \emph{strict} if $f_n=0$ for all $n>1$.
 An \emph{augmented} strictly unital $\Ainfty$\+algebra $A$ is
a strictly unital $\Ainfty$\+algebra endowed with a morphism of
strictly unital $\Ainfty$\+algebras $A\rarrow k$, where the strictly
unital $\Ainfty$\+algebra structure on~$k$ comes from its structure
of DG\+algebra with zero differential.
 An augmented strictly unital $\Ainfty$\+algebra is \emph{strictly
augmented} if the augmentation morphism is strict.
 A morphism of augmented or strictly augmented strictly unital
$\Ainfty$\+algebras is a morphism of strictly unital
$\Ainfty$\+algebras forming a commutative diagram with
the augmentation morphisms.
 The categories of augmented strictly unital $\Ainfty$\+algebras,
strictly augmented strictly unital $\Ainfty$\+algebras, and
nonunital $\Ainfty$\+algebras are equivalent.
 The equivalence of the latter two categories is provided by
the functor of formal adjoining of the strict unit, and
the equivalence of the former two categories can be deduced from
the equivalence between the categories of DG\+coalgebras $C$ and
CDG\+coalgebras $C$ endowed with a CDG\+coalgebra morphism
$C\rarrow k$.
 The DG\+category of strictly unital $\Ainfty$\+modules over
an augmented strictly unital $\Ainfty$\+algebra $A$ is equivalent
to the DG\+category of nonunital $\Ainfty$\+modules over
the corresponding nonunital $\Ainfty$\+algebra.

\subsection{Derived category of $\Ainfty$\+modules}
\label{derived-category-ainfty-modules}
 Let $A$ be a strictly unital $\Ainfty$\+algebra over a field~$k$.
 A (not necessarily closed) morphism of strictly unital left
$\Ainfty$\+modules $f\:L\rarrow M$ over $A$ is called \emph{strict}
if one has $f_n=0$ and $(df)_n=0$ for all $n>0$, or equivalently,
$f^n=0$ and $(df)^n=0$ for all $n>0$.
 Strictly unital left $\Ainfty$\+modules and strict morphisms between
them form a DG\+subcategory of the DG\+category of strictly unital
left $\Ainfty$\+modules and their morphisms.

 A closed strict morphism of strictly unital $\Ainfty$\+modules is
called a \emph{strict homotopy equivalence} if it is a homotopy
equivalence in the DG\+category of strictly unital $\Ainfty$\+modules
and strict morphisms between them.
 A triple $K\rarrow L\rarrow M$ of strictly unital $\Ainfty$\+modules
with closed strict morphisms between them is said to be \emph{exact}
if $K\rarrow L\rarrow M$ is an exact triple of graded vector spaces.
 The total strictly unital $\Ainfty$\+module of such an exact triple
is defined in the obvious way.

 Any strictly unital left $\Ainfty$\+module $M$ over $A$ can be
considered as a complex with the differential $l_0=p_0\:M\rarrow M$,
since one has $l_0^2=0$.
 A strictly unital left $\Ainfty$\+module $M$ is called \emph{acyclic}
if it is acyclic as a complex with the differential~$l_0$.
 For any closed morphism of strictly unital left $\Ainfty$\+modules
$f\:L\rarrow M$ the map $f_0=f^0\:L\rarrow M$ is a morphism of
complexes with respect to~$l_0$.
 The morphism $f$ is called a \emph{quasi-isomorphism} if $f_0$ is
a quasi-isomorphism of complexes.

 Let $v\:A\rarrow k$ be a homogeneous $k$\+linear retraction and
$C=\Br_v(A)$ be the corresponding CDG\+coalgebra structure on
the graded tensor coalgebra $\bigoplus_n A_+[1]^{\ot n}$.

\begin{thm1}
 The following six definitions of the \emph{derived category\/
$\sD(A\modl)$ of strictly unital left $\Ainfty$\+modules} over $A$
are equivalent, i.~e., lead to naturally isomorphic (triangulated)
categories: \par
\textup{(a)} the homotopy category of the DG\+category of strictly
unital left $\Ainfty$\+modules over $A$ and their morphisms; \par
\textup{(b)} the localization of the category of strictly unital
left $\Ainfty$\+modules over $A$ and their closed morphisms
by the class of quasi-isomorphisms; \par
\textup{(c)} the localization of the category of strictly unital
left $\Ainfty$\+modules over $A$ and their closed morphisms by
the class of strict homotopy equivalences; \par
\textup{(d)} the quotient category of the homotopy category of
the DG\+category of strictly unital left $\Ainfty$\+modules over $A$
and strict morphisms between them by the thick subcategory of
acyclic $\Ainfty$\+modules; \par
\textup{(e)} the localization of the category of strictly unital left
$\Ainfty$\+modules over $A$ and their closed strict morphisms by
the class of strict quasi-isomorphisms; \par
\textup{(f)} the quotient category of the homotopy category of
the DG\+category of strictly unital left $\Ainfty$\+modules over $A$
and strict morphisms between them by its minimal triangulated
subcategory containing all the total strictly unital $\Ainfty$\+modules
of exact triples of strictly unital $\Ainfty$\+modules with closed
strict morphisms between them. \par
 The derived category\/ $\sD(A\modl$) is also naturally equivalent
to the following triangulated categories: \par
\textup{(g)} the coderived category\/ $\sD^\co(C\comodl)$ of left
CDG\+comodules over~$C$; \par
\textup{(h)} the contraderived category\/ $\sD^\ctr(C\contra)$ of left
CDG\+contramodules over~$C$; \par
\textup{(i)} the absolute derived category\/ $\sD^\abs(C\comodl)$ of
left CDG\+comodules over~$C$; \par
\textup{(j)} the absolute derived category\/ $\sD^\abs(C\contra)$ of
left CDG\+contramodules over~$C$.
\end{thm1}

\begin{proof}
 The equivalence of (a-h) holds in the generality of CDG\+comodules
and CDG\+contramodules over an arbitrary conilpotent CDG\+coalgebra~$C$.
 More precisely, let us consider CDG\+comodules over $C$ that are
cofree as graded comodules, or equivalently, CDG\+contramodules over
$C$ that are free as graded contramodules, in place of strictly
unital $\Ainfty$\+modules.
 The equivalence of (a), (g), and~(h) follows from (the proof of)
Theorem~\ref{cdg-coalgebra-inj-proj-resolutions}.

 There is a natural increasing filtration $F$ on a conilpotent
CDG\+coalgebra $C$ that was defined in~\ref{conilpotent-cobar-duality},
and there are induced increasing filtrations
$F_nK=\lambda^{-1}(F_nC\ot_kK)$ on all CDG\+comodules $K$ over $C$
and decreasing filtrations $F^nQ=\pi(\Hom_k(C/F_{n-1}C,Q))$
on all CDG\+contramodules $Q$ over $C$.
 In particular, from any CDG\+comodule $C\ot_kM$ that is cofree as
a graded comodule and the corresponding CDG\+contramodule $\Hom_k(C,M)$
that is free as a graded contramodule one can recover the complex $M$
as $M \simeq F_0(C\ot_kM) \simeq \Hom_k(C,M)/F^1\Hom_k(C,M)$.
 A closed morphism of CDG\+comodules $C\ot_kL\rarrow C\ot_kM$ and
the corresponding closed morphism of CDG\+contramodules
$\Hom_k(C,L)\rarrow\Hom_k(C,M)$ are homotopy equivalences if and
only if the corresponding morphism of complexes $L\rarrow M$ is
a quasi-isomorphism.
 Indeed, let us pass to the cones and check that a cofree CDG\+comodule
$C\ot_kM$ is contractible if and only if the complex $M$ is acyclic.
 The ``only if'' is clear, and ``if'' follows from the fact that
$C\ot_kM$ is coacyclic whenever $M$ is acyclic.
 To check the latter, notice that the quotient CDG\+comodules
$F_n(C\ot_kM)/F_{n-1}(C\ot_kM)$ are just the tensor products of
complexes of vector spaces $F_nC/F_{n-1}C\ot_kM$ with the trivial
CDG\+comodule structures.

 This proves the equivalence of (a) and~(b), since for any DG\+category
$\sDG$ with shifts and cones the homotopy category $H^0(\sDG)$ can
be also obtained by inverting homotopy equivalences in the category
of closed morphisms $Z^0(\sDG)$.
 The equivalence of (d) and~(e) also follows from the latter result
about DG\+categories; and to prove the equivalence of (a) and~(c)
the following slightly stronger formulation of that result is
sufficient. 
 For any DG\+category $\sDG$ with shifts and cones consider the class
of morphisms of the form $(\id_X,0)\:X\oplus\cone(\id_X)\rarrow X$.
 Then by formally inverting all the morphisms in this class one
obtains the homotopy category $H^0(\sDG)$.

 A morphism of CDG\+comodules $f'\:C\ot_kL\rarrow C\ot_kM$ and
the corresponding morphism of CDG\+contramodules $f''\:\Hom_k(C,L)
\rarrow\Hom_k(C,M)$ can be called strict if both $f'$ and $df'$
as maps of graded vector spaces can be obtained by applying the functor
$C\ot{-}$ to certain maps $L\rarrow M$, or equivalently, both $f''$
and $df''$ as maps of graded vector spaces can be obtained by applying
the functor $\Hom_k(C,{-})$ to (the same) maps $L\rarrow M$.
 Let $w\:k\rarrow C$ be the coaugmentation map; consider
the DG\+algebra $U=\Cb_w(C)$.
 When $C=\Br_v(A)$ is the bar-construction of a strictly unital
$\Ainfty$\+algebra $A$, the DG\+algebra $U$ is called
the \emph{enveloping DG\+algebra} of~$A$.
 For any conilpotent CDG\+coalgebra $C$, consider the DG\+functors
$C\ot^{\tau_{C,w}}{-}$ and $\Hom^{\tau_{C,w}}(C,-)$ assigning
CDG\+comodules and CDG\+contramodules over $C$ to DG\+modules over~$U$.
 These two DG\+functors are equivalences between the DG\+categories
of left DG\+modules over $U$, left CDG\+comodules over $C$ that are
cofree as graded comodules with strict morphisms between them, and
left CDG\+contramodules over $C$ that are free as graded contramodules
with strict morphisms between them.
 So the equivalence of (a) and~(d) follows from
Theorem~\ref{conilpotent-cobar-duality}, and the equivalence of
(d) and~(f) follows from Corollary~\ref{nonconilpotent-duality}.

 Finally, the equivalences (g)$\Longleftrightarrow$(i) and
(h)$\Longleftrightarrow$(j) for $C=\Br_v(A)$ are provided by
Theorem~\ref{finite-homol-dim-cdg-coalgebra}.
\end{proof}

 Let $A$ be a DG\+algebra over~$k$; it can be considered as a strictly
unital $\Ainfty$\+algebra and left DG\+modules over it can be
considered as strictly unital $\Ainfty$\+modules as expained
in~\ref{strictly-1-ainfty-algebras}.
 It follows from Theorem~\ref{bar-construction-duality} that
the derived category of left DG\+modules over $A$ is equivalent to
the derived category of left $\Ainfty$\+modules, so our notation
$\sD(A\modl)$ is consistent.

 Any strictly unital $\Ainfty$\+algebra $A$ can be considered as
a complex with the differential $m_1\:A\rarrow A$, since $m_1^2=0$.
 For any morphism of strictly unital $\Ainfty$\+algebras
$f\:A\rarrow B$ the map $f_1\:A\rarrow B$ is a morphism of complexes
with respect to~$m_1$.
 A morphism~$f$ of strictly unital $\Ainfty$\+algebras is called
a \emph{quasi-isomorphism} if $f_1\:A\rarrow B$ is a quasi-isomorphism
of complexes, or equivalently, $f_{1,+}\:A_+\rarrow B_+$ is
a quasi-isomorphism of complexes.
 
 Let $f\:A\rarrow B$ be a morphism of strictly unital
$\Ainfty$\+algebras and $g\:\Br_v(A)\rarrow\Br_v(B)$ be
the corresponding morphism of CDG\+coalgebras.
 Any strictly unital left $\Ainfty$\+module $M$ over $B$ can be
considered as a strictly unital left $\Ainfty$\+module over $A$;
this corresponds to the extension-of-scalarars functors $E_g$ on
the level of CDG\+comodules that are cofree as graded comodules
and $E^g$ on the level of CDG\+contramodules that are free as
graded contramodules.
 Denote the induced functor on derived categories by $\boI R_f\:
\sD(B\modl)\rarrow\sD(A\modl)$.
 The functor $\boI R_f$ has left and right adjoint functors
$\boL E_f$ and $\boR E^f\:\sD(A\modl)\rarrow\sD(B\modl)$ that can
be constructed as the functors $\boI R_g$ and $\boI R^g$ on
the level of coderived categories of CDG\+comodules and
contraderived categories of CDG\+contramodules
(see~\ref{co-contra-ext-scalars}).

\begin{thm2}
 The functor $R_f$ is an equivalence of triangulated categories if
and only if a morphism $f$ of strictly unital $\Ainfty$\+algebras
is a quasi-isomorphism.
\end{thm2}

\begin{proof}
 The ``if'' part follows easily from
Theorem~\ref{cdg-coalgebra-scalars}.
 Both ``if'' and ``only if'' can be deduced from
Theorem~\ref{dg-mod-scalars} in the following way.
 For any strictly unital $\Ainfty$\+algebra $A$ and the corresponding
CDG\+coalgebra $C=\Br_v(A)$ with its coaugmentation~$w$,
the adjunction morphism $C\rarrow\Br_v(\Cb_w(C))$ corresponds to
a morphism of strictly unital $\Ainfty$\+algebras $u\:A\rarrow U(A)$
from $A$ to its the enveloping DG\+algebra $U(A)$ (see the proof of
Theorem~1).
 The morphism $u$ is a quasi-isomorphism, as one can see by
considering the increasing filtration $F$ on $A$ defined by the rules
$F_0A=k$ and $F_1A=A$, and the induced filtration on $U(A)$.
 The functor $\boI R_u$ is an equivalence of triangulated categories,
as it follows from Theorems~\ref{bar-construction-duality}
and~\ref{conilpotent-cobar-duality}, or as we have just proved.
 It remains to apply Theorem~\ref{dg-mod-scalars} to the morphism
of DG\+algebras $U(f)\:U(A)\rarrow U(B)$.
\end{proof}

 Let $A$ be a strictly unital $\Ainfty$\+algebra and $C=\Br_v(A)$
be the corresponding CDG\+coalgebra.
 All the above results about strictly unital left $\Ainfty$\+modules
over $A$ apply to strictly unital right $\Ainfty$\+modules as well,
since one can pass to the opposite CDG\+coalgebra $C^\rop$ as defined
in~\ref{cdg-cotor-coext-ctrtor}.
 In particular, the derived category of strictly unital right
$\Ainfty$\+modules $\sD(\modr A)$ is defined and naturally equivalent
to the coderived category $\sD^\co(\comodr C)$.

{\hbadness=2500
 The functor $\Tor^A\:\sD(\modr A)\times\sD(A\modl)\rarrow k\vect^\sgr$
can be constructed either by restricting the functor of cotensor
product $\oc_C\:\Hot(\comodr C)\times\Hot(C\comodl)\rarrow \Hot(k\vect)$
to the Cartesian product of the homotopy categories of CDG\+comodules
that are cofree as graded comodules, or by restricting the functor of
contratensor product $\ocn_C\:\Hot(\comodr C)\times\Hot(C\contra)
\rarrow\Hot(k\vect)$ to the Cartesian product of the homotopy
categories of CDG\+comodules that are cofree as graded comodules and
CDG\+contramodules that are free as graded contramodules.
 The functors one obtains in these two ways are naturally isomorphic
by the result of~\ref{coext-and-ext-cotor-and-ctrtor}.
 This definition of the functor $\Tor^A$ agrees with the definition
of functor $\Tor^A$ for DG\+algebras $A$ by
Theorem~\ref{cotor-and-tor-coext-and-ext}.1.
\par}

{\hfuzz=2.5pt\emergencystretch=0em
 The functor $\Ext_A=\Hom_{\sD(A\modl)}\:\sD(A\modl)^\op\times
\sD(A\modl)\rarrow k\vect^\sgr$ can be computed in three ways.
 One can either restrict the functor $\Hom_C\:\Hot(C\comodl)^\op\times
\Hot(C\comodl)\rarrow\Hot(k\vect)$ to the Cartesian product of
the homotopy categories of CDG\+comodules that are cofree as graded
comodules, or restrict the functor $\Hom^C\:\Hot(C\contra)^\op\times
\Hot(C\contra)\rarrow\Hot(k\vect)$ to the Cartesian product of
the homotopy categories of CDG\+contramodules that are free as graded
contramodules, or restrict the functor $\Cohom_C\:\Hot(C\comodl)^\op
\times\Hot(C\contra)\rarrow\Hot(k\vect)$ to the Cartesian product
of the homotopy categories of CDG\+comodules that are cofree as
graded comodules and CDG\+contramodules that are free as graded
contramodules.
 The functors one obtains in these three ways are naturally isomorphic
by the result of~\ref{coext-and-ext-cotor-and-ctrtor}, are
isomorphic to the functor $\Hom_{\sD(A\modl)}$ by Theorem~1 above,
and agree with the functor $\Ext_A$ for DG\+algebras $A$ by
Theorem~\ref{bar-construction-duality} or
Theorem~\ref{cotor-and-tor-coext-and-ext}.1.
\par}
 
\begin{rem}
 One can define a nonunital curved $\Ainfty$\+algebra $A$ as
a structure of not necessarily coaugmented DG\+coalgebra on
$\bigoplus_n A[1]^{\ot n}$; such a structure is given by a sequence
of linear maps $m_n\:A^{\ot n}\rarrow A$, \ $n=0$, $1$,~\dots, where
$m_0\:k\rarrow A$ may be a nonzero map (corresponding to
the curvature element of~$A$).
 Any morphism of DG\+coalgebras $f\:\bigoplus_n A[1]^{\ot n}
\rarrow B[1]^{\ot n}$ preserves the coaugmentations of the graded
tensor coalgebras, though, so a morphism of nonunital curved 
$\Ainfty$\+algebras $f\:A\rarrow B$ is given by a sequence of 
linear maps $f_n\:A^{\ot n}\rarrow B$, \ $n=1$, $2$,~\dots{}
 All the definitions of~\ref{non-1-ainfty-algebras}--%
\ref{strictly-1-ainfty-algebras} can be generalized straightforwardly
to the curved situation, and all the results
of~\ref{strictly-1-ainfty-algebras} hold in this case.
 However, this theory is largely trivial.
 For any strictly unital curved $\Ainfty$\+algebra $A$ with $m_0\ne0$,
every object of the DG\+category of strictly unital curved
$\Ainfty$\+modules over~$A$ is contractible.
 In particular, the same applies to nonunital curved
$\Ainfty$\+modules over a nonunital curved $\Ainfty$\+algebra.
 Two cases have to be considered separately, the case when the image
of $m_0$ coincides with $k\cdot 1_A\subset A$, and the case when
$m_0(1)$ and $1_A$ are linearly independent.
 The former case cannot occur in the $\boZ$\+graded situation for
dimension reasons, but in the $\boZ/2$\+graded situation it is possible.
 In this case the differential on the CDG\+coalgebra $C=\bigoplus_n
A_+[1]^{\ot n}$ is compatible with the coaugmentation $w\:k\simeq
A_+[1]^{\ot 0}\rarrow C$, but the curvature linear function
$h\:C\rarrow k$ is not, i.~e., $h\circ w\ne 0$.
 Let $C\ot_kM$ be a left CDG\+comodule over $C$ that is cofree
as a graded comodule and $F_nC\ot_kM$ be its natural increasing
filtration induced by the natural increasing filtration $F$ of
the conilpotent graded coalgebra~$C$.
 Then the filtrations $F$ on both $C$ and $C\ot_kM$ are preserved
by the differentials, since the differential on $C$ is compatible
with the coaugmentation.
 The induced differential $l_0\:M\rarrow M$ on $M=F_0C\ot_kM$
has the square equal to a nonzero constant from~$k$ times
the identity endomorphism of~$M$.
 The CDG\+comodule $(M,l_0)$ over the CDG\+coalgebra $F_0C$ is clearly
contractible; let $t_0$ be its contracting homotopy.
 Set $t=\id\ot t_0\:C\ot_kM\rarrow C\ot_kM$; then $t$ is a nonclosed
endomorphism of $M$ of degree~$-1$ and the endomorphism $d(t)=
dt+td$ is invertible, hence $C\ot_kM$ is contractible.
 This argument is applicable to any conilpotent graded coalgebra~$C$.
 In the case when $m_0(1)$ and $1_A$ are linearly independent,
the theory trivializes even further.
 The author learned the idea of the following arguments from
M.~Kontsevich\-.
 All strictly unital curved $\Ainfty$\+algebra structures with
$m_0(1)$ and $1_A$ linearly independent on a given graded vector
space $A$ are isomorphic, and all structures of a strictly unital
curved $\Ainfty$\+module over $A$ on a given graded vector space $M$
are isomorphic.
 Indeed, consider the component of the tensor degree $n=1$ of
the differential on $C=\bigoplus_n A_+[1]^{\ot n}$; it is determined
by~$m_0$.
 This differential makes $\bigoplus_n A_+[1]^{\ot n}$ into a complex,
and this complex is acyclic.
 Taking this fact into account, one can first find a CDG\+coalgebra
isomorphism of the form $(\id,a)$, \ $a\:C\rarrow k$ between 
a given CDG\+coalgebra structure on $C$ and a certain CDG\+coalgebra
structure with $h=0$, i.~e., a DG\+coalgebra structure.
 One proceeds step by step, killing the component $h_n\:A_+^{\ot n}
\rarrow k$ of the linear function~$h$ using a linear function~$a$
with the only component $a_{n+1}\:A_+^{\ot n+1}\rarrow k$.
 Having obtained a DG\+coalgebra structure on $C$, one subsequently
kills all the components $m_n\:A_+^{\ot n}\rarrow A$ of the
differential~$d$ with $n>0$ using graded tensor coalgebra
automorphisms~$f$ of $C$ with the only nonzero components
$f_1=\id_{A_+}$ and $f_{n+1}\:A_+^{\ot n+1}\rarrow A_+$.
 Analogously one shows that any DG\+comodule over $C$ that is cofree
as a graded comodule is isomorphic to a direct sum of shifts of
the DG\+comodule~$C$.
 Since $C$ is acyclic, such DG\+comodules are clearly contractible.
 Consequently, the coderived category $\sD^\co(C\comodl)$ vanishes.
 Alternatively, one could consider the DG\+category of strictly
unital curved $\Ainfty$\+modules over a strictly unital curved
$\Ainfty$\+algebra $A$ with strict morphisms between the curved
$\Ainfty$\+modules.
 This DG\+category is equivalent to the DG\+category of
CDG\+modules over the CDG\+coalgebra $U=\Cb_w(C)$.
 Its homotopy category $\Hot(U\modl)$ can well be nonzero, but
the corresponding absolute derived category $\sD^\abs(U\modl)$
is zero by Theorem~\ref{nonconilpotent-duality}.
 So in the category of (strictly unital or nonunital) curved
$\Ainfty$\+algebras over a field there are too few and too many
morphisms at the same time: there are no ``change-of-connection''
morphisms, and in particular no morphisms corresponding to
nonstrict morphisms of CDG\+algebras, and still there are enough
morphisms to trivialize the theory in almost all cases.
 One way out of this predicament is to restrict oneself to the curved
$\Ainfty$\+algebras (modules, morphisms) for each of which
the structure maps $m_n$, $l_n$, $f_n$ vanish for $n$ large
enough.
 For such curved $\Ainfty$\+algebras $A$, it actually makes sense
to consider morphisms $f\:A\rarrow B$ with nonvanishing
change-of-connection components $f_0\:k\rarrow B$.
 This theory can be interpreted in terms of the topological
tensor coalgebras $\prod_n A[1]^{\ot n}$ (cf.\
Remark~\ref{curved-ainfty-co-contra-derived}).
 Another solution would be to consider curved $\Ainfty$\+algebras
over the ring of formal power series $k[[t]]$ and require
the curvature and change-of-connection elements to be divisible
by~$t$.
\end{rem}

\subsection{Noncounital curved $\Ainfty$\+coalgebras}
\label{non-co1-curved-ainfty-coalgebras}
 Let $C$ be a graded vector space over a field~$k$.
 Consider the graded tensor algebra (free associative graded algebra)
$\bigoplus_{n=0}^\infty C[-1]^{\ot n}$ generated by the graded vector
space $C[-1]$.
 A \emph{noncounital curved $\Ainfty$\+coalgebra} structure on $C$ is,
by the definition, a DG\+algebra structure on $\bigoplus_n
C[-1]^{\ot n}$, i.~e., an odd derivation~$d$ of degree~$1$ on
$\bigoplus_n C[-1]^{\ot n}$ such that $d^2=0$.
 Since a derivation of $\bigoplus_n C[-1]^{\ot n}$ is uniquely
determined by its restriction to $C[-1]\simeq C[-1]^{\ot 1}\subset
\bigoplus_n C[-1]^{\ot n}$, a noncounital curved $\Ainfty$\+coalgebra
structure on $C$ can be considered as a sequence of linear maps
$\mu_n\:C\rarrow C^{\ot n}$, \ $n=0$, $1$,~\dots\ of degree $2-n$.
 More precisely, define the maps~$\mu_n$ by the formula
$d(c) = \sum_{n=0}^\infty (-1)^{n+\sum_{j=1}^n (n-j)(|\mu_{n,j}(c)|+1)}
\mu_{n,1}(c)\ot\dsb\ot\mu_{n,n}(c)$, where $c\in C$ and
$\mu_n(c)=\mu_{n,1}(c)\ot\dsb\ot\mu_{n,n}(c)$
is a symbolic notation for the tensor $\mu_n(c)\in C^{\ot n}$.
 A \emph{convergence condition} must be satisfied: for any $c\in C$
one must have $\mu_n(c)=0$ for all but a finite number of
the degrees~$n$.
 Furthermore, the sequence of maps~$\mu_n$ must satisfy a sequence
of quadratic equations corresponding to the equation $d^2=0$ on
the derivation~$d$.

 A morphism of noncounital curved $\Ainfty$\+coalgebras $f\:C\rarrow D$
over a field~$k$ is, by the definition, a morphism of DG\+algebras
$\bigoplus_n C[-1]^{\ot n}\rarrow \bigoplus_n D[-1]^{\ot n}$ over~$k$.
 Since the graded algebra morphism from a graded tensor algebra
$\bigoplus_n C[-1]^{\ot n}$ is determined by its restriction to
$C[-1]\simeq C[-1]^{\ot 1}\subset\bigoplus_n C[-1]^{\ot n}$, a morphism
of noncounital curved $\Ainfty$\+coalgebras $f\:C\rarrow D$ can be
considered as a sequence of linear maps $f_n\:C\rarrow D^{\ot n}$, \
$n=0$, $1$,~\dots\ of degree $1-n$.
 More precisely, define the maps~$f_n$ by the formula $f(c) =
\sum_{n=0}^\infty (-1)^{n-1+\sum_{j=1}^n (n-j)(|f_{n,j}(c)|+1)}
f_{n,1}(c)\ot\dsb\ot f_{n,n}(c)$, where $c\in C$ and
$f_n(c)=f_{n,1}(c)\ot\dsb\ot f_{n,n}(c)\in D^{\ot n}$.
 A convergence condition must be satisfied: for any $c\in C$ one
must have $f_n(c)=0$ for all but a finite number of the degrees~$n$.
 Furthermore, the sequence of maps~$f_n$ must satisfy a sequence of
polynomial equations corresponding to the equation $d\circ f=f\circ d$
on the morphism~$f$.

 Let $C$ be a noncounital curved $\Ainfty$\+coalgebra over a field~$k$
and $M$ be a graded vector space over~$k$.
 A structure of \emph{noncounital left curved $\Ainfty$\+comodule}
over $C$ on $M$ is, by the definition, a structure of DG\+module over
the DG\+algebra $\bigoplus_n C[-1]^{\ot n}$ on the free graded left
module $\bigoplus_n C[-1]^{\ot n}\ot_kM$ over the graded algebra
$\bigoplus_n C[-1]^{\ot n}$.
 Analogously, a structure of \emph{noncounital right curved
$\Ainfty$\+comodule} over $C$ on a graded vector space $N$ is defined
as a structure of DG\+module over $\bigoplus_n C[-1]^{\ot n}$ on
the free graded right module $N\ot_k \bigoplus_n C[-1]^{\ot n}$.
 Since a derivation of a free graded module $\bigoplus_n C[-1]^{\ot n}
\ot_k M$ compatible with a given derivation of the graded algebra
$\bigoplus_n C[-1]^{\ot n}$ is determined by its restriction to
the subspace of generators $M\simeq C[-1]^{\ot 0}\ot_k M\subset
\bigoplus_n C[-1]^{\ot n}\ot_kM$, a noncounital left curved
$\Ainfty$\+comodule structure on $M$ can be considered as a sequence
of linear maps $\lambda_n\:M\rarrow C^{\ot n}\ot_kM$, \ $n=0$,
$1$,~\dots\ of degree $1-n$.
 More precisely, define the maps~$\lambda_n$ by the formula $d(x) =
\sum_{n=0}^\infty (-1)^{n+\sum_{j=-n}^{-1} (j+1)(|\lambda_{n,j}(x)|+1)}
\lambda_{n,-n}(x)\ot\dsb\ot\lambda_{n,-1}(x)\ot\lambda_{n,0}(x)$,
where $x\in M$ and $\lambda_n(x) = \lambda_{n,-n}(x)\ot\dsb\ot
\lambda_{n,-1}(x)\ot\lambda_{n,0}(x)\in C^{\ot n}\ot_kM$.
 A convergence condition must be satisfied: for any $x\in M$ one
must have $\lambda_n(x)=0$ for all but a finite number of
the degrees~$n$.
 Furthermore, the sequence of maps~$\lambda_n$ must satisfy a sequence
of nonhomogeneous quadratic equations corresponding to the equation
$d^2=0$ on the derivation~$d$ on $\bigoplus_n C[-1]^{\ot n}\ot_k M$.
 Analogously, a noncounital right curved $\Ainfty$\+comodule structure
on $N$ can be considered as a sequence of linear maps $\rho_n\:N\rarrow
N\ot_k C^{\ot n}$ defined by the formula $d(y) = \sum_{n=0}^\infty
(-1)^{n|\rho_{n,0}(y)|+\sum_{j=1}^n (n-j)(|\rho_{n,j}(y)|+1)}
\rho_{n,0}(y)\ot\rho_{n,1}(y)\ot\dsb\ot\rho_{n,n}(y)$, where $y\in N$
and $\rho_n(y) = \rho_{n,0}(y)\ot\rho_{n,1}(y)\ot\dsb\ot\rho_{n,n}(y)
\in N\ot_k C^{\ot n}$.

 The complex of morphisms between noncounital left curved
$\Ainfty$\+comodules $L$ and $M$ over a noncounital curved
$\Ainfty$\+coalgebra $C$ is, by the definition, the complex of
morphisms between left DG\+modules $\bigoplus_n C[-1]^{\ot n}\ot_kL$
and $\bigoplus_n C[-1]^{\ot n}\ot_kM$ over the DG\+algebra
$\bigoplus_n C[-1]^{\ot n}$.
 Analogously, the complex of morphisms between noncounital right
curved $\Ainfty$\+comodules $R$ and $N$ over $C$ is, by
the definition, the complex of morphisms between right DG\+modules
$R\ot_k\bigoplus_n C[-1]^{\ot n}$ and $N\ot_k\bigoplus_n C[-1]^{\ot n}$
over the DG\+algebra $\bigoplus_n C[-1]^{\ot n}$.
 A morphism of noncounital left curved $\Ainfty$\+comodules
$f\:L\rarrow M$ of degree~$i$ is the same that a sequence of linear
maps $f_n\:L\rarrow C^{\ot n}\ot_kM$, \ $n=0$, $1$,~\dots\ of
degree $i-n$ satisfying the convergence condition: for any $x\in L$
one must have $f_n(x)=0$ for all but a finite number of
the degrees~$n$.
 More precisely, define the maps~$f_n$ by the formula  $f(x) =
\sum_{n=0}^\infty (-1)^{n+\sum_{j=-n}^{-1} (j+1)(|f_{n,j}(x)|+1)}
f_{n,-n}(x)\ot\dsb\ot f_{n,-n}(x)\ot f_{n,0}(x)$, where $x\in L$ and
$f_n(x) = f_{n,-n}(x)\ot\dsb\ot f_{n,-1}(x)\ot f_{n,0}(x)\in
C^{\ot n}\ot_kM$.
 Any sequence of linear maps~$f_n$ satisfying the convergence
condition corresponds to a (not necessarily closed) morphism of
noncounital curved $\Ainfty$\+comodules~$f$.
 Analogously, a morphism of noncounital right curved
$\Ainfty$\+comodules $f\:R\rarrow N$ of degree~$i$ is the same that
a sequence of linear maps $f_n\:R\rarrow N\ot_k C^{\ot n}$, \ $n=0$,
$1$,~\dots\ of degree $i-n$ satisfying the convergence condition.
 More precisely, define the maps~$f_n$ by the formula $f(y) =
\sum_{n=0}^\infty (-1)^{n|f_{n,0}(y)|+\sum_{j=1}^n (n-j)(|f_{n,j}(y)|
+1)} f_{n,0}(y)\ot f_{n,1}(y)\ot\dsb\ot f_{n,n}(y)$, where $y\in R$
and $f_n(y) = f_{n,0}(y)\ot f_{n,1}(y)\ot\dsb\ot f_{n,n}(y)\in
N\ot_k C^{\ot n}$.

 Let $C$ be a noncounital curved $\Ainfty$\+coalgebra over a field~$k$
and $P$ be a graded vector space over~$k$.
 A structure of \emph{noncounital left curved $\Ainfty$\+contramodule}
over $C$ on $P$ is, by the definition, a structure of DG\+module
over the DG\+algebra $\bigoplus_n C[-1]^{\ot n}$ on the cofree graded
left module $\Hom_k(\bigoplus_n C[-1]^{\ot n},P)$ over the graded
algebra $\bigoplus_n C[-1]^{\ot n}$.
 The action of $\bigoplus_n C[-1]^{\ot n}$ in $\Hom_k(\bigoplus_n C[-1]
^{\ot n},P)$ is induced by the right action of
$\bigoplus_n C[-1]^{\ot n}$ in itself as explained in
\ref{injective-dg-modules} and~\ref{dg-mod-scalars}.
 Since a derivation of a cofree graded module $\Hom_k(\bigoplus_n
C[-1]^{\ot n},P)$ compatible with a given derivation of the graded
algebra $\bigoplus_n C[-1]^{\ot n}$ is determined by its composition
with the projection $\Hom_k(\bigoplus_n C[-1]^{\ot n},P)\rarrow P$
induced by the unit map $k\rarrow\bigoplus_n C[-1]^{\ot n}$,
a noncounital left curved $\Ainfty$\+contramodule structure on $P$
can be considered as a linear map $\pi\:\prod_{n=0}^\infty
\Hom_k(C^{\ot n},P)[n-1]\rarrow P$ of degree~$0$.
 More precisely, define the map $\pi$ by the rule that the image of
the element $d(g)$ under the projection to $P$ equals
$\pi((g_n)_{n=0}^\infty)$, where a map $g\:\bigoplus_n C[-1]^{\ot n}
\rarrow P$ and a sequence of maps $g_n\:C^{\ot n}\rarrow P$ are
related by the formula $g(c_1\ot\dsb\ot c_n) = (-1)^{n|g_n|+
\sum_{j=1}^n (n-j)(|c_j|+1)} g_n(c_1\ot\dsb\ot c_n)$ for $c_j\in C$.
 The map~$\pi$ must satisfy a system of nonhomogeneous quadratic
equations corresponding to the equation $d^2=0$ on
the derivation~$d$ on $\Hom_k(\bigoplus_n C[-1]^{\ot n},P)$.

{\hbadness=2300
 The complex of morphisms between noncounital left curved
$\Ainfty$\+contramodules $P$ and $Q$ over a noncounital curved
$\Ainfty$\+coalgebra $C$ is, by the definition, the complex of
morphisms between left DG\+modules $\Hom_k(\bigoplus_n C[-1]^{\ot n},
P)$ and $\Hom_k(\bigoplus_n C[-1]^{\ot n},Q)$ over the DG\+algebra
$\bigoplus_n C[-1]^{\ot n}$.
 A morphism of noncounital left curved $\Ainfty$\+contramodules
$f\:P\rarrow Q$ of degree~$i$ is the same that a linear map
$f_\sqcap\:\prod_{n=0}^\infty\Hom_k(C^{\ot n},P)[n] \rarrow Q$
of degree~$i$.
 More precisely, define the map $f_\sqcap$ by the rule that the image
of the element $f(g)$ under the projection to $Q$ equals
$f_\sqcap((g_n)_{n=0}^\infty)$, where a map $g\:\bigoplus_n C[-1]^
{\ot n}\rarrow P$ and a sequence of maps $g_n\:C^{\ot n}\rarrow P$ are
related by the above formula.
 Any linear map $f_\sqcap$ corresponds to a (not necessarily closed)
morphism of noncounital curved $\Ainfty$\+contramodules~$f$. \par}

\subsection{Strictly counital curved $\Ainfty$\+coalgebras}
\label{strictly-co1-curved-ainfty-coalgebras}
 Let $C$ be a noncounital curved $\Ainfty$\+coalgebra over a field~$k$.
 A homogeneous linear function $\eps\:C\rarrow k$ of degree~$0$ is
called a \emph{strict counit} if one has $\eps(\mu_{2,1}(c))\mu_{2,2}
(c) = c = \eps(\mu_{2,2}(c))\mu_{2,1}(c)$ and $\eps(\mu_{n,j}(c))
\mu_{n,1}(c)\ot\dsb\ot\mu_{n,j-1}(c)\ot\mu_{n,j+1}(c)\ot\dsb\ot
\mu_{n,n}(c) = 0$ for all $c\in C$, \ $1\le j\le n$, and $n\ne2$.
 A strict counit is unique if it exists.
 A \emph{strictly counital curved $\Ainfty$\+coalgebra} is
a noncounital curved $\Ainfty$\+coalgebra that admits a strict counit.
 A morphism of strictly counital curved $\Ainfty$\+coalgebras
$f\:C\rarrow D$ is a morphism of noncounital curved
$\Ainfty$\+coalgebras such that $\eps_D\circ f_1=\eps_C$ and
$\eps(f_{n,j}(c))f_{n,1}(c)\ot\dsb\ot f_{n,j-1}(c)\ot f_{n,j+1}(c)\ot
\dsb\ot f_{n,n}(c) = 0$ for all $c\in C$, \ $1\le j\le n$, and $n>1$.
 Notice that for a strictly counital curved $\Ainfty$\+coalgebra $C$
one has $\eps_C=0$ if and only if $C=0$.
 We will assume our strictly counital curved $\Ainfty$\+coalgebras
to have nonzero counits.

{\hfuzz=2pt\emergencystretch=0em
 A \emph{strictly counital left curved $\Ainfty$\+comodule} $M$ over
a strictly counital curved $\Ainfty$\+coalgebra $C$ is a noncounital
left curved $\Ainfty$\+comodule such that $\eps(\lambda_{1,-1}(x))
\lambda_{1,0}(x)\allowbreak = x$ and $\eps(\lambda_{n,j}(x))
\lambda_{n,-n}(x)\ot\dsb\ot\lambda_{n,j-1}(x)\ot\lambda_{n,j+1}(x)
\ot\dsb\ot\lambda_{n,-1}(x)\ot\lambda_{n,0}(x) = 0$ for all $x\in M$, \
$-n\le j\le -1$, and $n>1$.
 Analogously, a \emph{strictly counital right curved
$\Ainfty$\+comodule} $N$ over $C$ is a noncounital right curved
$\Ainfty$\+comodule such that $\eps(\rho_{1,1}(y))\rho_{1,0}(y) = y$
and $\eps(\rho_{n,j}(y))\rho_{n,0}(y)\ot\rho_{n,1}(y)\ot\dsb\ot
\rho_{n,j-1}(y)\ot\rho_{n,j+1}(y)\ot\dsb\ot\rho_{n,n}(y) = 0$ for all
$y\in N$, \ $1\le j\le n$, and $n>1$.
 Finally, a \emph{strictly counital left curved $\Ainfty$\+contramodule}
$P$ over $C$ is a noncounital left curved $\Ainfty$\+contramodule
satisfying the following condition.
 For any sequence of linear maps $g_n\:C^{\ot n}\rarrow P$, \ $n=0$,
$1$,~\dots\ of degree $s+n-1$ for which there exists a double sequence
of linear maps $g'_{n,j}\:C^{\ot n-1}\rarrow P$, \ $1\le j\le n$
such that $g_n(c_1\ot\dsb\ot c_n)=\sum_{j=1}^n\eps(c_j)
g'_{n,j}(c_1\ot\dsb\ot c_{j-1}\ot c_{j+1}\ot\dsb\ot c_n)$ for all~$n$
and $c_t\in C$ the equation $\pi((g_n)_{n=0}^\infty)=g'_{1,1}(1)$
must hold in $P^s$.\par}

 The complex of morphisms between strictly counital left curved
$\Ainfty$\+comodules $L$ and $M$ over a strictly counital left curved
$\Ainfty$\+coalgebra $C$ is the subcomplex of the complex of morphisms
between $L$ and $M$ as noncounital $\Ainfty$\+comodules consisting
of all morphisms $f\:L\rarrow M$ such that $\eps(f_{n,j}(x))
f_{n,-n}(x)\ot\dsb\ot f_{n,j-1}(x)\ot f_{n,j+1}(x)\ot\dsb\ot
f_{n,-1}(x) \ot f_{n,0}(x) = 0$ for all $x\in L$, \ $-n\le j\le -1$,
and $n>0$.
 Analogously, the complex of morphisms between strictly counital right
curved $\Ainfty$\+comodules $R$ and $N$ over $C$ is the subcomplex of
the complex of morphisms between $R$ and $N$ as noncounital
$\Ainfty$\+comodules consisting of all morphisms $f\:R\rarrow N$ 
such that $\eps(f_{n,j}(y))f_{n,0}(y)\ot f_{n,1}(y)\ot\dsb\ot
f_{n,j-1}(y)\ot f_{n,j+1}(y)\ot\dsb\ot f_{n,n}(y) = 0$ for all $y\in R$,
\ $1\le j\le n$, and $n>0$.
 Finally, the complex of morphisms between strictly counital left
curved $\Ainfty$\+contramodules $P$ and $Q$ over $C$ is the subcomplex
of the complex of morphisms between $P$ and $Q$ as noncounital
$\Ainfty$\+contramodules consisting of all morphisms $f\:P\rarrow Q$
satisfying the following condition.
 For any sequence of linear maps $g_n\:C^{\ot n}\rarrow P$, \ $n=0$,
$1$,~\dots\ of degree $s+n-i$ for which there exists a double sequence
of linear maps $g'_{n,j}\:C^{\ot j-1}\rarrow P$, \ $1\le j\le n$
such that $g_n(c_1\ot\dsb\ot c_n)=\sum_{j=1}^n\eps(c_j)
g'_{n,j}(c_1\ot\dsb\ot c_{j-1}\ot c_{j+1}\ot\dsb\ot c_n)$ for all~$n$
and $c_t\in C$ the equation $f_\sqcap((g_n)_{n=0}^\infty)=0$
must hold in $Q^s$, where $i$~is the degree of~$f$.

 Let $C$ be a noncounital curved $\Ainfty$\+coalgebra and $\eps_C\:
C\rarrow k$ be a homogeneous linear function of degree~$0$.
 Set $C_+=\ker\eps$.
 Then the graded tensor algebra $\bigoplus_n C_+[-1]^{\ot n}$ is
a subalgebra of the tensor algebra $\bigoplus_n C[-1]^{\ot n}$.
 Denote by $K_C$ the cokernel of the embedding $\bigoplus_n C_+[-1]
^{\ot n}\rarrow\bigoplus_n C[-1]^{\ot n}$ and by $\kappa_C\in K_C$
the element of $C/C_+[-1]\subset K_C$ for which $\eps(\kappa_C)=1$.
 Let $\theta_C\in\bigoplus_n C[-1]^{\ot n}$ be any element of
degree~$1$ whose image in $K_C$ is equal to~$\kappa_C$.
 Then the linear function $\eps_C\:C\rarrow k$ is a strict counit
if and only if the odd derivation $d'(a)=d(a)+[\theta_C,a]$
of degree~$1$ on the tensor algebra $\bigoplus_n C[-1]^{\ot n}$
preserves the subalgebra $\bigoplus_n C_+[-1]^{\ot n}$ and
the element $h=d(\theta_C)+\theta_C^2$ belongs to
$\bigoplus_n C_+[-1]^{\ot n}$.
 This condition does not depend on the choice of $\theta_C$.
 For strictly counital curved $\Ainfty$\+coalgebras $C$ and $D$,
a morphism of noncounital curved $\Ainfty$\+coalgebras $f\:C\rarrow D$
is a morphism of strictly counital curved $\Ainfty$\+coalgebras if
and only if $f(\bigoplus_n C_+[-1]^{\ot n})\subset\bigoplus_n
D_+[-1]^{\ot n}$ and $f(\kappa_C)=\kappa_D$. 

 Let $C$ be a strictly counital curved $\Ainfty$\+coalgebra and $M$
be a noncounital left curved $\Ainfty$\+comodule over~$C$.
 Then $M$ is a strictly counital curved $\Ainfty$\+comodule if and
only if the odd derivation $d'(z)=d(z)+\theta_C z$ of degree~$1$
on the free module $\bigoplus_n C[-1]^{\ot n}\ot_kM$ compatible with
the derivation $d'$ of the algebra $\bigoplus_n C[-1]^{\ot n}$
preserves the subspace $\bigoplus_n C_+[-1]^{\ot n}\ot_kM\subset
\bigoplus_n C[-1]^{\ot n}\ot_kM$.
 Analogously, a noncounital right curved $\Ainfty$\+comodule $N$
over $C$ is a strictly counital curved $\Ainfty$\+comodule if and
only if the odd derivation $d'(z) = d(z) - (-1)^{|z|}z\theta_C$
of degree~$1$ on the free module $N\ot_k\bigoplus_n C[-1]^{\ot n}$
compatible with the derivation $d'$ of the algebra $\bigoplus_n
C[-1]^{\ot n}$ preserves the subspace $N\ot_k\bigoplus_n C_+[-1]^{\ot n}
\subset N\ot_k\bigoplus_n C[-1]^{\ot n}$.
 Finally, a noncounital left curved $\Ainfty$\+contramodule $P$
over $C$ is a strictly counital curved $\Ainfty$\+contramodule if
and only if the odd derivation $d'(q)=d(q)+\theta_C q$ of
degree~$1$ on the cofree module $\Hom_k(\bigoplus_n C[-1]^{\ot n},P)$
compatible with the derivation $d'$ of the algebra
$\bigoplus_n C[-1]^{\ot n}$ preserves the subspace $\Hom_k(K_C,P)
\subset\Hom_k(\bigoplus_n C[-1]^{\ot n},P)$.

 For strictly counital left curved $\Ainfty$\+comodules $L$ and $M$
over a strictly counital curved $\Ainfty$\+coalgebra $C$, a (not
necessarily closed) morphism of noncounital $\Ainfty$\+comodules
$f\:L\rarrow M$ is a morphism of strictly counital $\Ainfty$\+comodules
if and only if one has $f(\bigoplus_n C_+[-1]^{\ot n}\ot_kL)\subset
\bigoplus_n C_+[-1]^{\ot n}\ot_kM$.
 Analogously, for strictly counital right curved $\Ainfty$\+comodules
$R$ and $N$ over $C$, a (not necessarily closed) morphism of noncounital
$\Ainfty$\+comodules $f\:R\rarrow N$ is a morphism of strictly counital
$\Ainfty$\+comodules if and only if one has $f(R\ot_k\bigoplus_n
C_+[-1]^{\ot n})\subset N\ot_k\bigoplus_n C_+[-1]^{\ot n}$.
 Finally, for strictly counital left curved $\Ainfty$\+contramodules
$P$ and $Q$ over $C$, a (not necessarily closed) morphism of noncounital
$\Ainfty$\+contramodules $f\:P\rarrow Q$ is a morphism of strictly
counital $\Ainfty$\+contramodules if and only if one has
$f(\Hom_k(K_C,P))\subset\Hom_k(K_C,Q)$.

 Let $C$ be a strictly counital curved $\Ainfty$\+coalgebra.
 Choose a homogeneous $k$\+linear section $w\:k\rarrow C$ of
the strict counit map $\eps\:C\rarrow k$.
 Define the element $\theta_C\in \bigoplus_n C[-1]^{\ot n}$ as
$\theta_C=w(1)\in C[-1]\subset \bigoplus_n C[-1]^{\ot n}$.
 Then $\theta_C$ is an element of degree~$1$ whose image in $K_C$
is equal to~$\kappa_C$.
 Let $d\:\bigoplus_n C_+[-1]^{\ot n}\rarrow\bigoplus_n C_+[-1]^{\ot n}$
be the restriction of the odd derivation $d'$ of $\bigoplus_n C[-1]^
{\ot n}$ defined by the above formula, and let $h\in\bigoplus C_+[-1]^
{\ot n}$ be the element defined above.
 Then $\Cb_w(C) = (\bigoplus_n C_+[-1]^{\ot n}\;d\;h)$ is
a CDG\+algebra.
 The CDG\+algebra $\Cb_w(C)$ is called the \emph{cobar-construction}
of a strictly counital curved $\Ainfty$\+coalgebra~$C$.

 Let $f\:C\rarrow D$ be a morphism of strictly counital curved
$\Ainfty$\+coalgebras.
 Let $w\:k\rarrow C$ and $w\:k\rarrow D$ be homogeneous $k$\+linear
sections, and let $\theta_C$ and $\theta_D$ be the corresponding
elements of degree~$1$ in the graded tensor algebras.
 The morphism of tensor algebras $f\:\bigoplus_n C[-1]^{\ot n}
\rarrow \bigoplus_n D[-1]^{\ot n}$ induces a morphism of graded tensor
algebras $\bigoplus_n C_+[-1]^{\ot n}\rarrow \bigoplus_n D_+[-1]^
{\ot n}$, which we will denote also by~$f$.
 The element $\eta_f = f(\theta_C) - \theta_D \in \bigoplus_n D[-1]
^{\ot n}$ has a zero image in $K_D$, so it belongs to $\bigoplus_n
D_+[-1]^{\ot n}$.
 Then the pair $(f,\eta_f)$ is a morphism of CDG\+algebras $\Cb_w(C)
\rarrow\Cb_w(D)$.
 Thus the cobar-construction $C\mpsto\Cb_w(C)$ is a functor from
the category of strictly counital curved $\Ainfty$\+coalgebras with
nonzero counits to the category of CDG\+algebras whose underlying
graded algebras are graded tensor algebras.
 One can easily see that this functor is an equivalence of categories.
 Alternatively, one can use any element~$\theta_C$ of degree~$1$
whose image in $K_C$ is equal to~$\kappa_C$ in the construction of
this equivalence of categories.

 Let $C$ be a strictly counital curved $\Ainfty$\+coalgebra,
$w\:k\rarrow C$ be a homogeneous $k$\+linear section, and $\theta_C
\in\bigoplus_n C[-1]^{\ot n}$ be the corresponding element of
degree~$1$.
 Let $M$ be a strictly counital curved $\Ainfty$\+comodule over~$C$.
 Set $d\:\bigoplus_n C_+[-1]^{\ot n}\ot_k M\rarrow\bigoplus_n
C_+[-1]^{\ot n}\ot_k M$ to be the restriction of the differential~$d'$
on $\bigoplus_n C[-1]^{\ot n}$ defined by the above formula.
 Then $\Cb_w(C,M) = (\bigoplus_n C_+[-1]^{\ot n}\ot_kM\;d)$ is
a left CDG\+module over the CDG\+algebra $\Cb_w(C)$.
 Analogously, for a strictly counital curved $\Ainfty$\+comodule $N$
over $C$ set $d\:N\ot_k\bigoplus_n C_+[-1]^{\ot n}\rarrow N\ot_k
\bigoplus_n C_+[-1]^{\ot n}$ to be the restriction of
the differential~$d'$ on $N\ot_k\bigoplus_n C[-1]^{\ot n}$
defined above.
 Then $\Cb_w(N,C) = (N\ot_k\bigoplus_n C_+[-1]^{\ot n}\;d)$ is
a right CDG\+module over the CDG\+algebra $\Cb_w(C)$.
 Finally, let $P$ be a strictly counital curved
$\Ainfty$\+contramodule over~$C$.
 Set $d\:\Hom_k(\bigoplus_n C_+[-1]^{\ot n},P)\rarrow
\Hom_k(\bigoplus_n C_+[-1]^{\ot n},P)$ to be the map induced by
the differential~$d'$ on $\Hom_k(\bigoplus_n C[-1]^{\ot n},P)$
defined above.
 Then $\Br^w(C,P) = (\Hom_k(\bigoplus_n C_+[-1]^{\ot n},P)\;d)$ is
a left CDG\+module over the CDG\+algebra $\Cb_w(C)$.

 To a (not necessarily closed) morphism of strictly counital left
curved $\Ainfty$\+comod\-ules $f\:L\rarrow M$ over $C$ one can assign
the induced map $\bigoplus_n C_+[-1]^{\ot n}\ot_k L\rarrow
\bigoplus_n C_+[-1]^{\ot n}\ot_k M$.
 So we obtain the DG\+functor $M\mpsto\Cb_w(C,M)$, which is
an equivalence between the DG\+category of strictly counital left
curved $\Ainfty$\+comodules over $C$ and the DG\+category of left
CDG\+modules over $\Cb_w(C)$ that are free as graded modules.
 Analogously, to a (not necessarily closed) morphism of strictly
counital right curved $\Ainfty$\+comodules $f\:R\rarrow N$ over $C$
one can assign the induced map $R\ot_k\bigoplus_n C_+[-1]^{\ot n}
\rarrow N\ot_k \bigoplus_n C_+[-1]^{\ot n}$.
 So we obtain the DG\+functor $N\mpsto\Cb_w(N,C)$, which is
an equivalence between the DG\+categories of strictly counital right
curved $\Ainfty$\+comodules over $C$ and right CDG\+modules over  
$\Cb_w(C)$ that are free as graded modules.
 Finally, to a (not necessarily closed) morphism of strictly counital
left curved $\Ainfty$\+contramodules $f\:P\rarrow Q$ over $C$ one
can assign the induced map $\Hom_k(\bigoplus_n C_+^{\ot n},P)
\rarrow\Hom_k(\bigoplus_n C_+^{\ot n},Q)$.
 So we obtain the DG\+functor $P\mpsto\Br^w(C,P)$, which is
an equivalence between the DG\+category of strictly counital left
curved $\Ainfty$\+contramodules over $C$ and the DG\+category of
left CDG\+modules over $\Cb_w(C)$ that are cofree as graded modules.

 Now let $C$ be a CDG\+coalgebra with a nonzero counit, $M$ be
a left CDG\+comodule over $C$, \ $N$ be a right CDG\+comodule over $C$,
and $P$ be a left CDG\+contramodule over~$C$.
 Let $w\:k\rarrow C$ be a homogeneous $k$\+linear section.
 Define a strictly counital curved $\Ainfty$\+coalgebra structure
on $C$ by the rules $\mu_0(c)=h(c)$, \ $\mu_1(c)=d(c)$, \ 
$\mu_2(c)=c_{(1)}\ot c_{(2)}$, and $\mu_n(c)=0$ for $n>2$.
 Define a structure of strictly counital left curved
$\Ainfty$\+comodule over $C$ on $M$ by the rules $\lambda_0(x)=d(x)$, \
$\lambda_1(x)= x_{(-1)}\ot x_{(0)}$, and $\lambda_n(x)=0$ for $n>1$,
where $x\in M$.
 Define a structure of strictly counital right curved
$\Ainfty$\+comodule over $C$ on $N$ by the rules $\rho_0(y)=d(y)$, \
$\rho_1(y)=y_{(0)}\ot y_{(1)}$, and $\rho_n(x)=0$ for $n>1$, where
$y\in N$.
 Finally, define a structure of strictly counital left curved
$\Ainfty$\+contramodule over $C$ on $P$ by the rule
$\pi((g_n)_{n=0}^\infty)=d(g_0)+\pi_P(g_1)$, where $d\:P\rarrow P$
is the differential on $P$ and $\pi_P\:\Hom_k(C,P)\rarrow P$ is
the contraaction map.
 Then the CDG\+algebra structure $\Cb_w(C)$ on the graded tensor
algebra $\bigoplus_n C[-1]^{\ot n}$ that was defined
in~\ref{bar-cobar-constr} coincides with the CDG\+algebra structure
$\Cb_w(C)$ constructed above, so our notation is consistent.
 The left CDG\+module structure $\Cb_w(C)\ot^{\tau_{C,w}}M$ on the free
graded module $\bigoplus_n C[-1]^{\ot n}\ot_kM$ that was defined
in~\ref{twisting-cochains-subsect} coincides with the left CDG\+module
structure $\Cb_w(C,M)$.
 The right CDG\+module structure $N\ot^{\tau_{C,w}}\Cb_w(C)$ on
the free graded module $N\ot_k\bigoplus_n C[-1]^{\ot n}$ coincides
with the right CDG\+module structure $\Cb_w(N,C)$.
 The left CDG\+module structure $\Hom^{\tau_{C,w}}(\Cb_w(C),P)$ on
the cofree graded module $\Hom_k(\bigoplus_n C[-1]^{\ot n},P)$
coincides with the left CDG\+module structure $\Br^w(C,P)$.

 A morphism of strictly counital curved $\Ainfty$\+coalgebras
$f\:C\rarrow D$ is called \emph{strict} if $f_n=0$ for all $n\ne1$.
 A \emph{coaugmented} strictly counital curved $\Ainfty$\+coalgebra
$C$ is a strictly counital curved $\Ainfty$\+coalgebra endowed with
a morphism of strictly counital curved $\Ainfty$\+coalgebras 
$k\rarrow C$, where the strictly counital curved $\Ainfty$\+coalgebra
structure on~$k$ comes from its structure of CDG\+coalgebra with
zero differential and curvature linear function.
 A coaugmented strictly counital curved $\Ainfty$\+coalgebra is
\emph{strictly coaugmented} if the coaugmentation morphism is strict.
 A morphism of coaugmented or strictly coaugmented strictly counital
curved $\Ainfty$\+coalgebras is a morphism of strictly counital curved
$\Ainfty$\+coalgebras forming a commutative diagram with
the coaugmentation morphisms.
 The categories of coaugmented strictly counital curved
$\Ainfty$\+coalgebras, strictly coaugmented strictly counital curved
$\Ainfty$\+coalgebras, and noncounital curved $\Ainfty$\+coalgebras
are equivalent.
 The DG\+category of strictly counital curved $\Ainfty$\+comodules
or $\Ainfty$\+contramodules over a coaugmented strictly counital
curved $\Ainfty$\+coalgebra $C$ is equivalent to the DG\+category of 
noncounital curved $\Ainfty$\+comodules or $\Ainfty$\+contramodules
over the corresponding noncounital curved $\Ainfty$\+coalgebra~$C$.
 If $C$ is a strictly coaugmented strictly counital curved
$\Ainfty$\+coalgebra and $w\:k\rarrow C$ is the coaugmentation map,
then the CDG\+algebra $\Cb_w(C)$ is in fact a DG\+algebra.

\subsection{Coderived category of curved $\Ainfty$\+comodules and
contraderived category of curved $\Ainfty$\+contramodules}
\label{curved-ainfty-co-contra-derived}
 Let $C$ be a strictly counital curved $\Ainfty$\+coal\-gebra over
a field~$k$.
 Let $w\:k\rarrow C$ be a homogeneous $k$\+linear section and
$B=\Cb_w(C)$ be the corresponding CDG\+algebra structure on
$\bigoplus_n C_+[-1]^{\ot n}$.

 The \emph{coderived category\/ $\sD^\co(C\comodl)$ of strictly
counital left curved $\Ainfty$\+comodules} over $C$ is defined
as the homotopy category of the DG\+category of strictly counital
left curved $\Ainfty$\+comodules over~$C$.
 The coderived category $\sD^\co(\comodr C)$ of strictly counital
right curved $\Ainfty$\+comodules over $C$ is defined in
the analogous way.
 The \emph{contraderived category\/ $\sD^\ctr(C\contra)$ of strictly
counital right curved $\Ainfty$\+contramodules} over $C$ is defined
as the homotopy category of the DG\+category of strictly counital
left curved $\Ainfty$\+contramodules over~$C$.

\begin{thm}
 The following five triangulated categories are naturally
equivalent: \par
\textup{(a)} the coderived category\/ $\sD^\co(C\comodl)$; \par
\textup{(b)} the contraderived category\/ $\sD^\ctr(C\contra)$; \par
\textup{(c)} the coderived category\/ $\sD^\co(B\modl)$; \par
\textup{(d)} the contraderived category\/ $\sD^\ctr(B\modl)$; \par
\textup{(e)} the absolute derived category\/ $\sD^\abs(B\modl)$. 
\end{thm}

\begin{proof}
 The isomorphism of triangulated categories~(c--e) is provided
by Theorem~\ref{finite-homol-dim-cdg-ring}(a), and the equivalence
of triangulated categories (a), (b), and~(e) is the assertion of
Theorem~\ref{finite-homol-dim-cdg-ring}(b) with projective and
injective graded modules replaced by free and cofree ones.
 It suffices to find for any left CDG\+module $M$ over $B$
a closed injection from $M$ to a CDG\+module $J$ such that
both $J$ and $J/M$ are cofree as graded $B$\+modules, and a closed
surjection onto $M$ from a CDG\+module $F$ such that both
CDG\+modules $M$ and $\ker(F\to M)$ are free as graded $B$\+modules.
 This can be easily accomplished with either of the constructions
of Theorem~\ref{finite-homol-dim-cdg-ring}
or Theorem~\ref{cdg-coalgebra-inj-proj-resolutions}.
 One only has to notice that for any graded module $M$ over a graded
tensor algebra $B$ the kernel of the map $B\ot_kM\rarrow M$ is a free
graded $B$\+module and the cokernel of the map $M\rarrow\Hom_k(B,M)$
is a cofree graded $B$\+module.
\end{proof}

 Let $C$ be a CDG\+coalgebra over~$k$; it can be considered as
a strictly counital curved $\Ainfty$\+coalgebra, and CDG\+comodules
and CDG\+contramodules over it can be considered as strictly counital
curved $\Ainfty$\+comodules and $\Ainfty$\+contramodules as explained
in~\ref{strictly-co1-curved-ainfty-coalgebras}.
 It follows from Theorem~\ref{nonconilpotent-duality} that
the coderived category of left CDG\+comodules over $C$ is equivalent
to the coderived category of strictly counital left curved
$\Ainfty$\+comodules and the contraderived category of left
CDG\+contramodules over $C$ is equivalent to the contraderived
category of strictly counital left curved $\Ainfty$\+contramodules,
so our notation is consistent.
 Thus the above Theorem provides the \emph{comodule-contramodule
correspondence for strictly counital curved $\Ainfty$\+coalgebras}.
 By Theorem~\ref{nonconilpotent-duality}(c), the comodule-contramodule
correspondence functors in the CDG\+coalgebra case agree with
the comodule-contramodule correspondence functors we have construced
in the strictly counital curved $\Ainfty$\+coalgebra case.

 For the definition of a conilpotent curved $\Ainfty$\+coalgebra
and the curved $\Ainfty$\+coal\-gebra analogues of some
assertions of Theorem~\ref{derived-category-ainfty-modules}, see
Remark~\ref{cofibrant-dg-alg}.

 The functor $\Cotor^C\:\sD^\co(\comodr C)\times\sD^\co(C\comodl)
\rarrow k\vect^\sgr$ is constructed by restricting the functor of
tensor product $\ot_B\:\Hot(\modr B)\times\Hot(B\modl)\rarrow
\Hot(k\vect)$ to the Cartesian product of the homotopy categories
of CDG\+mod\-ules that are free as graded modules.
 The functor $\Coext_C\:\sD^\co(C\comodl)^\op\times\sD^\ctr(C\contra)
\rarrow k\vect^\sgr$ is constructed by restricting the functor of
homomorphisms $\Hom_B\:\Hot(B\modl)^\op\times\Hot(B\modl)\rarrow
\Hot(k\vect)$ to the Cartesian product of the homotopy category
of CDG\+modules that are free as graded modules and the homotopy
category of CDG\+modules that are cofree as graded modules.
 The functor $\Ctrtor^C\:\sD^\co(\comodr C)\times\sD^\ctr(C\contra)
\rarrow k\vect^\sgr$ is constructed by restricting the functor of
tensor product $\ot_B\:\Hot(\modr B)\times\Hot(B\modl)\rarrow
\Hot(k\vect)$ to the Cartesian product of the homotopy category
of CDG\+modules that are free as graded modules and the homotopy
category of CDG\+modules that are cofree as graded modules.
 These definitions of $\Cotor^C$, \ $\Coext_C$, and $\Ctrtor^C$
agree with the definitions of the functors $\Cotor^C$, \ $\Coext_C$,
and $\Ctrtor^C$ for CDG\+coalgebras~$C$ by
Theorem~\ref{cotor-and-tor-coext-and-ext}.2.

\begin{rem}
 For any graded vector space $C$, consider the topological graded
tensor algebra $\prod_n C[-1]^{\ot n} = \varprojlim_n\bigoplus_{j=0}^n
C[-1]^{\ot j}$ (see Remark~\ref{cdg-coalgebra-inj-proj-resolutions}
for the relevant general definitions).
 One can define a noncounital (``uncurved'') $\Ainfty$\+coalgebra $C$
as a structure of augmented DG\+algebra with a continuous differential
on $\prod_n C[-1]^{\ot n}$.
 Such a structure is given by a sequence of linear maps $\mu_n\:
C\rarrow C^{\ot n}$, \ $n=1$, $2$,~\dots\ without any convergence
condition imposed on them (but satisfying a sequence of quadratic
equations corresponding to the equation $d^2=0$).
 A morphism of noncounital $\Ainfty$\+coalgebras is a continuous
morphism of the corresponding topological DG\+algebras (which
always preserves the augmentations).
 The notion of a noncounital $\Ainfty$\+coalgebra is neither more
nor less general than that of a noncounital curved $\Ainfty$\+coalgebra.
 Still, any noncounital curved $\Ainfty$\+coalgebra $C$ with $\mu_0=0$
can be considered as a noncounital $\Ainfty$\+coalgebra.
 A morphism $f$ of noncounital curved $\Ainfty$\+coalgebras with
$\mu_0=0$ can be considered as a morphism of noncounital
$\Ainfty$\+coalgebras provided that $f_0=0$.
 A noncounital left $\Ainfty$\+comodule $M$ over $C$ is a structure
of DG\+module with a continuous differential over the topological
DG\+algebra $\prod_n C[-1]^{\ot n}$ on the free topological graded
module $\prod_n C[-1]^{\ot n}\ot_k M$.
 Such a structure is given by a sequence of linear maps $\lambda_n\:
M\rarrow C^{\ot n}\ot_kM$, \ $n=0$, $1$,~\dots\ without any convergence
conditions imposed.
 A noncounital left $\Ainfty$\+contramodule $P$ over $C$ is
a structure of DG\+module over $\prod_n C[-1]^{\ot n}$ on the cofree
discrete graded module $\bigoplus_n\Hom_k(C[-1]^{\ot n},P)$ of
continuous homogeneous linear maps $\prod_n C[-1]^{\ot n}\rarrow P$,
where $P$ is discrete.
 Such a structure is given by a sequence of linear maps $\pi_n\:
\Hom_k(C^{\ot n},P)\rarrow P$, \ $n=0$, $1$,~\dots{}
 All the definitions of~\ref{non-co1-curved-ainfty-coalgebras}--%
\ref{strictly-co1-curved-ainfty-coalgebras} are applicable in this
situation, and all the results
of~\ref{strictly-co1-curved-ainfty-coalgebras} hold in this case.
 For a (strictly counital or noncounital) curved $\Ainfty$\+coalgebra
$C$ with $\mu_0=0$ there are forgetful functors from the 
DG\+categories of (strictly counital or noncounital) curved
$\Ainfty$\+comodules and $\Ainfty$\+contramodules to the corresponding
DG\+categories of uncurved $\Ainfty$\+comodules and
$\Ainfty$\+contramodules.
 Furthermore, let $C$ be an (uncurved) strictly counital
$\Ainfty$\+coalgebra.
 Define the derived categories of strictly counital
$\Ainfty$\+comodules and $\Ainfty$\+contramodules as the quotient
categories of the homotopy categories corresponding to the
DG\+categories of strictly counital $\Ainfty$\+comodules and
$\Ainfty$\+contramodules by the thick subcategories formed by all
the $\Ainfty$\+comodules and $\Ainfty$\+contramodules that are
acyclic with respect to $\lambda_0$ and~$\pi_0$.
 Then one can define the functors $\Cotor^{C,I}$, \ $\Coext_C^I$, \
$\Ctrtor^{C,I}$, \ $\Ext_C^I$, and $\Ext^{C,I}$ on the Cartesian
products of the derived categories of strictly counital
$\Ainfty$\+comodules and $\Ainfty$\+contramodules by applying
the functors of topological tensor product and continuous
homomorphisms over $\prod_n C_+[-1]^{\ot n}$ to the corresponding
(topological or discrete) CDG\+modules over $\prod_n C_+[-1]^{\ot n}$.
 These functors are even preserved by the restrictions of scalars
corresponding to quasi-isomorphisms of strictly counital
$\Ainfty$\+coalgebras (i.~e., morphisms $f$ such that $f_1$ is
a quasi-isomorphism of complexes with respect to~$\mu_1$).
 All these assertions follow from the result of~\cite{EM1}, just
as in the proof of Theorem~\ref{cotor-coext-first-kind}.
 In the case of a strictly counital $\Ainfty$\+coalgebra coming
from a DG\+coalgebra $C$, the above functors agree with the derived
functors $\Cotor^{C,I}$,\. $\Coext_C^I$, etc., defined
in~\ref{cotor-coext-first-kind}.
\end{rem}

\Section{Model Categories of DG-Modules, \\*
CDG-Comodules, and CDG-Contramodules}

By a \emph{model category} we mean a model category in the sense
of Hovey~\cite{Hov}.

\subsection{Two DG\+module model category structures}
 Let $A=(A,d)$ be a DG\+ring and $Z^0\sDG(A\modl)$ be the abelian
category of left DG\+modules over $A$ and closed morphisms between them.
 Denote by $Z^0\sDG(A\modl)_\proj$ and $Z^0\sDG(A\modl_\proj)$
the two full subcategories of $Z^0\sDG(A\modl)$ consisting of all
the projective DG\+modules in the sense of~\ref{projective-dg-modules}
and all the DG\+modules that are projective as graded $A$\+modules,
respectively. 
 Analogously, denote by $Z^0\sDG(A\modl)_\inj$ and
$Z^0\sDG(A\modl_\inj)$ the two full subcategories of $Z^0\sDG(A\modl)$
consisting of all the injective DG\+modules in the sense
of~\ref{injective-dg-modules} and all the DG\+modules that
are injective as graded $A$\+modules.

 It follows from the constructions in the proofs of
Theorems~\ref{projective-dg-modules}--\ref{injective-dg-modules}
that every object of $Z^0\sDG(A\modl)_\proj$ is homotopy equivalent
to an object of $Z^0\sDG(A\modl)_\proj\cap Z^0\sDG(A\modl_\proj)$
and every object of $Z^0\sDG(A\modl)_\inj$ is homotopy equivalent
to an object of $Z^0\sDG(A\modl)_\inj\cap Z^0\sDG(A\modl_\inj)$.

\begin{thm}
\textup{(a)} There exists a model category structure on the category
$Z^0\sDG(A\modl)$ with the following properties.
 A morphism is a weak equivalence if and only if it is
a quasi-isomorphism.
 A morphism is a cofibration if and only if it is injective and its
cokernel belongs to $Z^0\sDG(A\modl)_\proj\cap Z^0\sDG(A\modl_\proj)$.
 A morphism is a fibration if and only if it is surjective.
 An object is cofibrant if and only if it belongs to
$Z^0\sDG(A\modl)_\proj\cap Z^0\sDG(A\modl_\proj)$.
 All objects are fibrant. \par
\textup{(b)} There exists a model category structure on the category
$Z^0\sDG(A\modl)$ with the following properties.
 A morphism is a weak equivalence if and only if it is
a quasi-isomorphism.
 A morphism is a cofibration if and only if it is injective.
 A morphism is a fibration if and only if it is surjective and its
kernel belongs to $Z^0\sDG(A\modl)_\inj\cap Z^0\sDG(A\modl_\inj)$.
 All objects are cofibrant.
 An object is fibrant if and only if it belongs to
$Z^0\sDG(A\modl)_\inj\cap Z^0\sDG(A\modl_\inj)$.
\end{thm}

 We will call the model structure of part~(a) of Theorem
the \emph{projective model structure} and the model structure of
part~(b) the \emph{injective model structure} on the category
of DG\+modules $Z^0\sDG(A\modl)$.

\begin{proof}
 We will prove part~(a); the proof of~(b) is dual.
 It is clear that all limits and colimits exist in the abelian
category $Z^0\sDG(A\modl)$.
 The two-out-of-three axiom for weak equivalences is obvious, as
is the retraction axiom for weak equivalences and fibrations.
 The retraction axiom for cofibrations holds since the subcategory
$Z^0\sDG(A\modl)_\proj\cap Z^0\sDG(A\modl_\proj)\subset
Z^0\sDG(A\modl)$ is closed under direct summands.
 To prove the lifting properties, use~\cite[Lemma~9.1.1]{P}.
 One has to check that the group $\Ext^1(E,K)$ computed in the abelian
category $Z^0\sDG(A\modl)$ vanishes whenever $K\in Z^0\sDG(A\modl)
_\proj\cap Z^0\sDG(A\modl_\proj)$ and either of the DG\+modules $E$
and $K$ is acyclic.
 Indeed, let $E\rarrow M\rarrow K$ be such an extension.
 Since $K\in Z^0\sDG(A\modl_\proj)$, the extension of graded
$A$\+modules $E\rarrow M\rarrow K$ splits.
 So our extension of DG\+modules is given by a closed morphism
$K\rarrow E$.
 Since $K\in Z^0\sDG(A\modl)_\proj$ and one of the DG\+modules
$E$ and $K$ is acyclic, this closed morphism is homotopic to zero.
 It remains to construct the functorial factorizations.
 For any complex of abelian groups $N$, let $N'\rarrow N$ be
a functorial closed surjective map of complexes of abelian groups
such that the induced map of cohomology groups is surjective and $N'$
is a complex of free abelian groups with free cohomology groups.
 E.~g., one can take $N'$ to be the direct sum of $\boZ x$ over all
nonzero homogeneous elements $x$ of~$N$. 
 Let $f\:L\rarrow M$ be a closed morphism of DG\+modules over~$A$.
 Then $L\rarrow L\oplus A\ot_{\boZ}M'\rarrow M$ is a decomposition
of~$f$ into a cofibration followed by a fibration and $L\rarrow L\oplus
\cone(\id_A)[-1]\ot_{\boZ}M'\rarrow M$ is a decomposition of~$f$ into
a trivial cofibration followed by a fibration.
 To construct a decomposition of~$f$ into a cofibration followed by
a trivial fibration, start with its decomposition $L\rarrow E\rarrow M$
into a cofibration followed by a fibration.
 Let $M_1$ be the kernel of the closed morphism of DG\+modules
$E\rarrow M$; set $E_1=A\ot_{\boZ}M_1'$.
 Denote by $M_2$ the kernel of the closed morphism $E_1\rarrow E$, etc.
 Let $K$ be the total DG\+module of the complex of DG\+modules
$\dsb\rarrow E_2\rarrow E_1\rarrow E$ formed by taking infinite direct
sums.
 Then $L\rarrow K$ is a cofibration and $K\rarrow M$ is a trivial
fibration for the reasons explained in the proof of
Theorem~\ref{projective-dg-modules}.
\end{proof}

\subsection{CDG\+comodule and CDG\+contramodule model category
structures}  \label{cdg-co-contra-model-categories}
 Let $C$ be a CDG\+coalgebra over a field~$k$.
 Let $Z^0\sDG(C\comodl)$ and $Z^0\sDG(C\contra)$ be the abelian
categories of left CDG\+comodules and left CDG\+contramodules over~$C$
with closed morphisms between them.
 Denote by $Z^0\sDG(C\comodl_\inj)$ and $Z^0\sDG(C\contra_\proj)$
the full subcategories of $Z^0\sDG(C\comodl)$ and $Z^0\sDG(C\contra)$
formed by all the CDG\+comodules that are injective as graded
$C$\+comodules and all the CDG\+contramodules that are projective
as graded $C$\+contramodules.

\begin{thm} {\hbadness=6500
\textup{(a)} There exists a model category structure on the category
$Z^0\sDG(C\comodl)$ with the following properties.
 A morphism is a weak equivalence if and only if its cone is
a coacyclic CDG\+comodule over~$C$.
 A morphism is a cofibration if and only if it is injective.
 A morphism is a fibration if and only if it is surjective and its
kernel belongs to $Z^0\sDG(C\comodl_\inj)$.
 All objects are cofibrant.
 An object is fibrant if and only if it belongs to
$Z^0\sDG(C\comodl_\inj)$. \par}
\textup{(b)} There exists a model category structure on the category
$Z^0\sDG(C\contra)$ with the following properties.
 A morphism is a weak equivalence if and only if its cone is
a contraacyclic CDG\+contramodule over~$C$.
 A morphism is a cofibration if and only if it is injective and its
cokernel belongs to $Z^0\sDG(C\contra_\proj)$.
 A morphism is a fibration if and only if it is surjective.
 An object is cofibrant if and only if it belongs to
$Z^0\sDG(C\contra_\proj)$.
 All objects are fibrant.
\end{thm}

\begin{proof}
 We will prove part~(a).
 All limits and colimits exist in the category $Z^0\sDG(C\comodl)$,
since it is an abelian category with infinite direct sums and products.
 The two-out-of-three axiom for weak equivalences holds since
coacyclic CDG\+comodules form a triangulated subcategory of
$\Hot(C\comodl)$.
 The retraction axiom for weak equivalences and cofibrations is clear,
and for fibrations it follows from the fact that the subcategory
$Z^0\sDG(C\comodl_\inj)\subset Z^0\sDG(C\comodl)$ is closed under
direct summands.
 To check the lifting properties, we use the same Lemma that in
the previous proof together with
Theorem~\ref{cdg-coalgebra-inj-proj-resolutions}(a).
 Finally, let us construct the functorial factorizations.
 To decompose a closed morphism $f\:L\rarrow M$ into a cofibration
followed by a fibration, find an injective closed morphism from $L$
into an a CDG\+comodule $J$ that is injective as a graded $C$\+comodule.
 This can be done by using either the comodule version of the functor
$G^-$ from the proof of Theorem~\ref{finite-homol-dim-cdg-ring}, or
the construction from the proof of
Theorem~\ref{cdg-coalgebra-inj-proj-resolutions}.
 Then the morphism~$f$ can be presented as the composition
$L\rarrow J\oplus M\rarrow M$.
 To obtain a decomposition of $f$ into a cofibration followed by
a trivial fibration, one needs to make sure that the CDG\+comodule
$J$ is contractible.
 If one uses the functor $G^-$, it suffices to notice that its image
consists entirely of contractible CDG\+comodules, and if one uses
the construction from the proof of 
Theorem~\ref{cdg-coalgebra-inj-proj-resolutions}, one has to replace
$J$ with $\cone(\id_J)$.
 To construct a decomposition of~$f$ into a trivial cofibration
followed by a fibration, start with its decomposition $L\rarrow E
\rarrow M$ into a cofibration followed by a fibration.
 Let $L_{-1}$ be the cokernel of the morphism $L\rarrow E$ and
$L_{-1}\rarrow E_{-1}$ be a closed morphism with a coacyclic cone
from $L_{-1}$ to a CDG\+comodule $E_{-1}\in Z^0\sDG(C\comodl_\inj)$,
which can be constructed in the way of either
Theorem~\ref{noetherian-cdg-ring-case} or
Theorem~\ref{cdg-coalgebra-inj-proj-resolutions}.
 Set $K=\cone(E\rarrow E_{-1})[-1]$; then the morphism $L\rarrow K$
is a trivial cofibration and the morphism $K\rarrow M$ is a fibration.
\end{proof}

\begin{rem}
 Let $C$ be a DG\+coalgebra.
 The model category structure of the first kind on the abelian
category of DG\+comodules over $C$ is defined as follows.
 Let $Z^0\sDG(C\comodl)_\inj\subset Z^0\sDG(C\comodl)$ be the full
subcategory consisting of all injective DG\+comodules over $C$
in the sense of~\ref{dg-coalgebra-resolutions}.
 Since any coacyclic DG\+comodule is acyclic, every object of
$Z^0\sDG(C\comodl)_\inj$ is homotopy equivalent to an object of
$Z^0\sDG(C\comodl)_\inj\cap Z^0\sDG(C\comodl_\inj)$.
 Let weak equivalences in $Z^0\sDG(C\comodl)$ be quasi-isomorphisms,
cofibrations be injective morphisms, and fibrations be surjective
morphisms whose kernels belong to $Z^0\sDG(C\comodl)_\inj\cap
Z^0\sDG(C\comodl_\inj)$.
 One proves that this is a model category structure in the same way
as in the above Theorem, using the facts that $\Acycl(C\comodl)$
and $Z^0\sDG(C\comodl)_\inj\cap Z^0\sDG(C\comodl_\inj)$ form 
a semiorthogonal decomposition of $\Hot(C\comodl)$, any injective
closed morphism from an object of $Z^0\sDG(C\comodl_\inj)$ into
an object of $Z^0\sDG(C\comodl)$ splits as an injective morphism of
graded comodules over~$C$, and any object of $Z^0\sDG(C\comodl)$
admits a closed injection into an object of $Z^0\sDG(C\comodl_\inj)$.
 The model category structure of the first kind on the abelian
category of DG\+contramodules $Z^0\sDG(C\contra)$ over $C$ is
defined in the dual way.
\end{rem}

\subsection{Finite homological dimension CDG\+module case} 
\label{cdg-module-model-categories}
 Let $B$ be a CDG\+ring and $Z^0\sDG(B\modl)$ be the abelian category
of CDG\+modules over $B$ and closed morphisms between them.
 Denote by $Z^0\sDG(B\modl_\proj)$ and $Z^0\sDG(B\modl_\inj)$
the full subcategories of $Z^0\sDG(B\modl)$ formed by all
the CDG\+modules that are projective as graded $B$\+modules and
injective as graded $B$\+modules.

 Assume that the graded ring $B^\#$ has a finite left homological
dimension.

\begin{thm}
\textup{(a)} There exists a model category structure on the category
$Z^0\sDG(B\modl)$ with the following properties.
 A morphism is a weak equivalence if and only if its cone is
absolutely acyclic.
 A morphism is a cofibration if and only if it is injective and
its cokernel belongs to $Z^0\sDG(B\modl_\proj)$.
 A morphism is a fibration if and only if it is surjective.
 An object is cofibrant if and only if it belongs to
$Z^0\sDG(B\modl_\proj)$.
 All objects are fibrant. \par
\textup{(b)} There exists a model category structure on the category
$Z^0\sDG(B\modl)$ with the following properties.
 A morphism is a weak equivalence if and only if its cone is
absolutely acyclic.
 A morphism is a cofibration if and only if it is injective.
 A morphism is a fibration if and only if it is surjective and
its kernel belongs to $Z^0\sDG(B\modl_\inj)$.
 All objects are cofibrant.
 An object is fibrant if and only if it belongs to
$Z^0\sDG(B\modl_\inj)$.
\end{thm}

 We will call the model structure of part~(a) of Theorem
the \emph{projective model structure of the second kind} and the model
structure of part~(b) the \emph{injective model structure of
the second kind} on the category of CDG\+modules $Z^0\sDG(B\modl)$.

\begin{proof}
 Analogous to the proof of Theorem~\ref{cdg-co-contra-model-categories}.
 Let us spell out some details for part~(a).
 To check the lifting properties, use the same Lemma from~\cite{P}
together with Theorem~\ref{cdg-mod-orthogonality}.
 To construct the functorial factorizations, choose any functor
assigning to a graded left $B^\#$\+module $N$ a surjective morphism
onto it from a projective graded left $B^\#$\+module $F(N)$.
 E.~g., one can take $F(N)$ to be the direct sum on $B^\#x$ over all
homogeneous elements $x\in N$.
 To decompose a closed morphism $f\:L\rarrow M$ into a trivial
cofibration followed by a fibration, consider the surjective closed
morphism onto $M$ from the CDG\+module $P=G^+(F(M^\#))$.
 The CDG\+module $P$ is projective as a graded $B$\+module and
contractible, so $L\rarrow L\oplus P\rarrow M$ is the desired
factorization.
 To obtain a decomposition of~$f$ into a cofibration followed by
a trivial fibration, start with a decomposition $L\rarrow E\rarrow M$
of~$f$ into a cofibration followed by a fibration.
 Set $M_1=\ker(E\rarrow M)$ and choose a closed morphism
$E_1\rarrow M_1$ with an absolutely acyclic cone from a CDG\+module
$E_1\in Z^0\sDG(B\modl_\proj)$ into~$M_1$.
 Such a morphism exists by Theorem~\ref{finite-homol-dim-cdg-ring}(b).
 Set $K=\cone(E_1\rarrow E)$; then $L\rarrow K$ is a cofibration
and $K\rarrow M$ is a trivial fibration.
\end{proof}

\begin{rem}
 For a CDG\+ring $B$ such that the graded ring $B^\#$ satisfies
the condition~($*$) of~\ref{noetherian-cdg-ring-case}, one can define
the injective model structure of the second kind on the category
of CDG\+modules $Z^0\sDG(B\modl)$.
 In this model structure, weak equivalences are morphisms with
coacyclic cones.
 Cofibrations are injective morphisms, and fibrations are surjective
morphisms whose kernels belong to $Z^0\sDG(B\modl_\inj)$.
 Analogously, for a CDG\+ring $B$ such that the graded ring $B^\#$
satisfies the condition~($**$) of~\ref{coherent-cdg-ring-case}, one
can define the projective model structure of the second kind on
the category $Z^0\sDG(B\modl)$.
 In this model structure, weak equivalences are morphisms with
contraacyclic cones.
 Cofibrations are injective morphisms whose cokernels belong to
$Z^0\sDG(B\modl_\proj)$, and fibrations are surjective morphisms.
 The proofs are analogous to the above and based on
Theorems~\ref{noetherian-cdg-ring-case}--\ref{coherent-cdg-ring-case}.
\end{rem}

\subsection{Quillen equivalences}
 For a definition of Quillen adjunctions and equivalences,
see~\cite{Hov}.
 Below we list the important Quillen equivalences and adjunctions
arising from the constructions of Sections~1--6.
 In order to convey the information about the directions of 
the adjoint pairs, we will always mention the left adjoint functor
first and the right adjoint functor second.
 Also we will mention the category that is the source of the left
adjoint functor first, and the category that is the source of
the right adjoint functor second.
 No proofs are given, as they are all very straightforward.

 For any DG\+ring $A$, the projective and injective model category
structures on the category of DG\+modules $Z^0\sDG(A\modl)$ are
Quillen equivalent.
 The equivalence is provided by the adjoint pair of identity
functors between $Z^0\sDG(A\modl)$ and itself, where the identity
functor from $Z^0\sDG(A\modl)$ with its projective model structure
to $Z^0\sDG(A\modl)$ with its injective model structure is
considered as the left adjoint functor, and the identity functor
in the other direction is considered as the right adjoint functor.

 Analogously, for any CDG\+ring $B$ such that the graded ring $B^\#$
has finite left homological dimension, the projective and injective
model category structures on the category of CDG\+modules
$Z^0\sDG(B\modl)$ are Quillen equivalent.
 The equivalence is provided by the adjoint pair of identity
functors, where the identity functor from $Z^0\sDG(B\modl)$ with
its projective model structure of the second kind to $Z^0\sDG(B\modl)$
with its injective model structure of the second kind is considered
as the left adjoint functor, and the identity functor in the other
direction is considered as the right adjoint functor.

 The comodule-contramodule correspondence for a CDG\+coalgebra $C$
over a field~$k$ can also be understood as a Quillen equivalence.
 The pair of adjoint functors $\Phi_C\:Z^0\sDG(C\contra)\rarrow
Z^0\sDG(C\comodl)$ and $\Psi_C\:Z^0\sDG(C\contra)\rarrow
Z^0\sDG(C\comodl)$ is a Quillen equivalence between the model
category of left CDG\+contramodules over $C$ and the model category
of left CDG\+comodules over~$C$.

 Let $\tau\:C\rarrow A$ be an acyclic twisting cochain between
a conilpotent CDG\+coal\-gebra $C$ and a DG\+algebra~$A$.
 Then the adjoint pair of Koszul duality functors
$N\mpsto A\ot^\tau\!\.N$ and $M\mpsto C\ot^\tau\!\.M$ is a Quillen
equivalence between the model category $Z^0\sDG(C\comodl)$ of left
CDG\+comodules over $C$ and the model category $Z^0\sDG(A\modl)$
of left DG\+modules over $A$, with the projective model
structure on the latter.
 The adjoint pair of Koszul duality functors $P\mpsto\Hom^\tau(C,P)$
and $Q\mpsto\Hom^\tau(A,Q)$ is a Quillen equivalnce between
the model category $Z^0\sDG(A\modl)$ of left DG\+modules over $A$
and the model category $Z^0\sDG(C\contra)$ of left CDG\+contramodules
over $C$, with the injective model structure on the former.

 Analogously, let $C$ be a CDG\+coalgebra and $\tau=\tau_{C,w}\:
C\rarrow\Cb_w(C)$ be the natural twisting cochain.
  Then the adjoint pair of Koszul duality functors
$N\mpsto B\ot^\tau\!\.N$ and $M\mpsto C\ot^\tau\!\.M$ is a Quillen
equivalence between the model category $Z^0\sDG(C\comodl)$ of left
CDG\+comodules over $C$ and the model category $Z^0\sDG(A\modl)$
of left CDG\+modules over $B$, with the projective model
structure on the latter.
 The adjoint pair of Koszul duality functors $P\mpsto\Hom^\tau(C,P)$
and $Q\mpsto\Hom^\tau(B,Q)$ is a Quillen equivalnce between
the model category $Z^0\sDG(B\modl)$ of left DG\+modules over $B$
and the model category $Z^0\sDG(C\contra)$ of left CDG\+contramodules
over $C$, with the injective model structure on the former.

 Notice that in each of the two Koszul duality situations above,
the conilpotent and the nonconilpotent one, we have four model
categories, two of them on the DG\+algebra or CDG\+algebra side
and two on the CDG\+coalgebra side.
 The four Quillen equivalences between these four model categories
form a circular diagram: if one looks at the direction of, e.~g.,
the left adjoint functor in each of the adjoint pairs, one finds that
these functors map $Z^0\sDG(A\modl)$ or $Z^0\sDG(B\modl)$ with the
projective model structure to the same category with the injective
model structure to $Z^0\sDG(C\contra)$ to $Z^0\sDG(C\comodl)$
and back to $Z^0\sDG(A\modl)$ or $Z^0\sDG(B\modl)$ with
the projective model structure.

 Let $f\:A\rarrow B$ be a morphism of DG\+algebras
(see~\ref{dg-mod-scalars}).
 Then the pair of adjoint functors $E_f\:Z^0\sDG(A\modl)\rarrow
Z^0\sDG(B\modl)$ and $R_f\:Z^0\sDG(B\modl)\rarrow Z^0\sDG(A\modl)$
is a Quillen adjunction between the model categories $Z^0\sDG(A\modl)$
and $Z^0\sDG(B\modl)$ with the projective model structures.
 The pair of adjoint functors $R_f\:Z^0\sDG(B\modl)\rarrow
Z^0\sDG(A\modl)$ and $E^f\:Z^0\sDG(A\modl)\rarrow Z^0\sDG(B\modl)$
is a Quillen adjunction between the model categories
$Z^0\sDG(B\modl)$ and $Z^0\sDG(A\modl)$ with the injective
model structures.

 Let $f\:C\rarrow D$ be a morphism of CDG\+coalgebras
(see~\ref{cdg-coalgebra-scalars}).
 Then the pair of adjoint functors $R_f\:Z^0\sDG(C\comodl)\rarrow
Z^0\sDG(D\comodl)$ and $E_f\:Z^0\sDG\allowbreak(D\comodl)\rarrow
Z^0\sDG(C\comodl)$ is a Quillen adjunction between the model
categories of CDG\+comodules over $C$ and CDG\+comodules over~$D$.
 The pair of adjoint functors $E^f\:Z^0\sDG(D\contra)\rarrow
Z^0\sDG(C\contra)$ and $R^f\:Z^0\sDG\allowbreak(C\contra)\rarrow
Z^0\sDG(D\contra)$ is a Quillen adjunction between the model
categories of CDG\+contramodules over $D$ and CDG\+contramodules
over~$C$.

\Section{Model Categories of DG-Algebras and CDG-Coalgebras}

\subsection{Model category of DG\+algebras} 
\label{dg-algebras-model-category}
 Let $k$ be a commutative ring.
 Denote by $k\alg_\dg$ the category of DG\+algebras over~$k$
and by $k\alg_\dg^\aug$ the category of augmented DG\+algebras
over~$k$, i.~e., DG\+algebras endowed with a DG\+algebra
morphism onto~$k$.

 Given a DG\+algebra $A$, we will consider DG\+algebras of
the special form $A\langle x_{n,\alpha}\rangle$.
 These are DG\+algebras that, as graded algebras, are obtained
by adjoining to $A$ a double-indexes family of free homogeneous
generators $x_{n,\alpha}$, where $n=1$, $2$,~\dots\ and
$\alpha$~belongs to some index set.
 The differential on $A\langle x_{n,\alpha}\rangle$ must satisfy
the conditions that $A$ is a DG\+subalgebra of
$A\langle x_{n,\alpha}\rangle$ and for any adjoined free
generator $x_{n,\alpha}$ the element $d(x_{n,\alpha})$ belongs
to the subalgebra generated by $A$ and $x_{m,\beta}$ with $m<n$.

 We will also consider DG\+algebras of the even more special form
$A\langle x_\omega, dx_\omega\rangle$.
 These are DG\+algebras that, as graded algebras, are obtained
by adjoining to $A$ the free homogeneous generators $x_\omega$
and $dx_\omega$; the differential on $A\langle x_\omega,
dx_\omega\rangle$ is defined in the way that the notation suggests.
 So $A\langle x_\omega, dx_\omega\rangle$ can be also defined
as the DG\+algebra freely generated by a DG\+algebra $A$
and the elements~$x_\omega$.

 The same notation applies to augmented DG\+algebras, with
the following obvious changes.
 When working with augmented DG\+algebras, one adjoins free
homogeneous generators $x_{n,\alpha}$ or $x_\omega$ annihilated
by the augmentation morphism, and requires that their differentials
be also annihilated by the augmentation morphism.

 The following result is due to Hinich and Jardine~\cite{Hin,Jar}.

\begin{thm}
\textup{(a)}
 There exists a model category structure on $k\alg_\dg$ with
the following properties.
 A morphism is a weak equivalence if and only if it is
a quasi-isomorphism.
 A morphism is a fibration if and only if it is surjective.
 A morphism is a cofibration if and only if it is a retract
of a morphism of the form $A\rarrow A\langle x_{n,\alpha}\rangle$.
 A morphism is a trivial cofibration if and only if it is
a retract of a morphism of the form $A\rarrow A\langle x_\omega,
dx_\omega\rangle$.
 All DG\+algebras are fibrant.
 A DG\+algebra is cofibrant if and only if it is a retract
of a DG\+algebra of the form $k\langle x_{n,\alpha}\rangle$. \par
\textup{(b)}
 The same assertion applies to the category of augmented DG\+algebras
$k\alg_\dg^\aug$.
\end{thm}

\begin{proof}
 It is easy to see that all limits and colimits exist in $k\alg_\dg$.
 All of them are preserved by the forgetful functor to the category of
graded algebras; all limits and filtered colimits are even preserved
by the forgetful functor to the category of complexes of vector spaces.
 The two-out-of-three axiom for weak equivalences and the retraction
axiom for all three classes of morphisms are obvious.
 It is straightforward to check that any morphism of the form
$A\rarrow A\langle x_\omega, dx_\omega\rangle$ is a quasi-isomorphism.
 To decompose a morphism of DG\+algebras $A\rarrow B$ as
$A\rarrow A\langle x_\omega, dx_\omega\rangle \rarrow B$ with
a surjective morphism $A\langle x_\omega, dx_\omega\rangle \rarrow B$,
it suffices to have generators $x_\omega$ corresponding to all
the homogeneous elements of~$B$.
 To decompose a morphism $A\rarrow B$ as $A\rarrow A\langle 
x_{n,\alpha}\rangle \rarrow B$ with a surjective quasi-isomorphism
$A\langle x_{n,\alpha}\rangle\rarrow B$, the following inductive
process is used.
 Let $x_{1,\alpha}$ correspond to all the homogeneous elements of $B$
annihilated by the differential, and $x_{2,\beta}$ correspond to
all the homogeneous elements of~$B$.
 Then the morphism $A\langle x_{1,\alpha}, x_{2,\beta}\rangle\rarrow B$
is surjective and the induced morphism on the cohomology is also
surjective.
 Furthermore, let $x_{3,\gamma}$ correspond to all homogeneous cocycles
in  $A\langle x_{1,\alpha}, x_{2,\beta}\rangle$ whose images are
coboundaries in $B$, let $x_{4,\delta}$ correspond to all homogeneous
cocycles in $A\langle x_{1,\alpha}, x_{2,\beta}, x_{3,\gamma}\rangle$
whose images are coboundaries in $B$, etc.
 The lifting property of the morphisms $A\rarrow A\langle x_\omega,
d x_\omega\rangle$ with respect to surjective morphisms of
DG\+algebras is obvious, and the lifting property of the morphisms
$A\rarrow A\langle x_{n,\alpha}\rangle$ with respect to surjective
quasi-isomorphisms is verified by induction on~$n$ and
a straightforward diagram chase.
 It remains to check that all the quasi-isomorphisms of
DG\+algebras that are retracts of morphisms $E\rarrow E\langle
x_{n,\alpha}\rangle$ are actually retracts of morphisms
$A\rarrow A\langle x_\omega, d x_\omega\rangle$.
 Let $f\:A\rarrow B$ be such a quasi-isomorphism; decompose it as
$A\rarrow A\langle x_\omega, dx_\omega\rangle\rarrow B$, where
$A\langle x_\omega, dx_\omega\rangle\rarrow B$ is a surjective
morphism.
 Then the latter morphism is also a quasi-isomorphism by the
two-out-of-three property of quasi-isomorphisms.
 Therefore, the morphism~$f$, being a retract of a morphism
$E\rarrow E\langle x_{n,\alpha}\rangle$, has the lifting property
with respect to the morphism $A\langle x_\omega, dx_\omega\rangle
\rarrow B$.
 It follows that $f$~is a retract of the morphism $A\rarrow
A\langle x_\omega,dx_\omega\rangle$.
\end{proof}

 Let $k$ be a field.
 For the purposes of nonaugmented Koszul duality one needs to
consider the full subcategory $k\alg_\dg^+$ of $k\alg_\dg$ formed
by all the DG\+algebras with nonzero units.
 This full subcategory is obtained from the category $k\alg_\dg$
by excluding its final object $A=0$.
 There are no morphisms from this final object to any objects
in $k\alg_\dg$ that are not final.
 The category $k\alg_\dg^+$ is a model category without limits,
i.~e., the classes of morphisms in $k\alg_\dg^+$ that are weak
equivalences, cofibrations, and fibrations in $k\alg_\dg$ satisfy
all the axioms of model category except for the existence of
limits and colimits.
 In fact, limits of all nonempty diagrams, all coproducts, and
all filtered colimits exist in $k\alg_\dg^+$.
 At the same time, $k\alg_\dg^+$ has no final object, and
fibered coproducts sometimes do not exist in it.

\subsection{Limits and colimits of CDG\+coalgebras}  \label{cdg-limits}
 Let $k$ be a field.
 We will use the notation and terminology
of~\ref{co-algebra-bar-duality}.
 The forgetful functor from the category of conilpotent or
coaugmented coalgebras to the category of graded vector spaces
assigns to a coaugmented DG\+coalgebra $(C,w)$ the graded
vector space~$C/w(k)$.
 The forgetful functor from the category of conilpotent or
coaugmented DG\+coalgebras to the category of complexes of vector
spaces is defined in the analogous way.

\begin{lem}
\textup{(a)} Limits of all nonempty diagrams, all coproducts, and
all filtered colimits exist in the category of conilpotent
CDG\+coalgebras $k\coalg_\cdg^\conilp$.
 All limits of diagrams with a final vertex, and consequently all
filtered limits, are preserved by the forgetful functor to
the category of conilpotent graded coalgebras.
 All coproducts and filtered colimits are preserved by the forgetful
functors to the categories of conilpotent graded coalgebras,
coaugmented graded coalgebras, and graded vector spaces.  \par
\textup{(b)} All limits and colimits exist in the category of
conilpotent DG\+coalgebras $k\coalg_\dg^\conilp$.
 All of them are preserved by the forgetful functor to the category
of conilpotent graded coalgebras.
 All colimits are preserved by the forgetful functor to the categories
of coagmented graded coalgebras and complexes of vector spaces.
\end{lem}

\begin{proof}
 Let us start with considering limits and colimits in the category of
coaugmented graded coalgebras $k\coalg^\coaug$.
 It is easy to see that colimits exist in $k\coalg^\coaug$ and commute
with the forgetful functor to the category of graded vector spaces.
 To construct limits in $k\coalg^\coaug$, it suffices to obtain finite
products, filtered limits, and equalizers.
 Finite products in $k\coalg^\coaug$ are coaugmented coalgebras
cofreely cogenerated by the coaugmented coalgebras being
multiplied~\cite{Ago}.
 Filtered limits and equalizers in $k\coalg^\coaug$ are graded vector
subspaces of the filtered limits and equalizers in the category of
graded vector spaces.
 In particular, the filtered limit of a diagram of coalgebras
$C_\alpha$ is the subspace of its limit as a diagram of graded
vector spaces $\varprojlim C_\alpha$ equal to the full preimage of
the tensor product $\varprojlim C_\alpha \ot_k\varprojlim C_\alpha
\subset \varprojlim C_\alpha\ot_k C_\alpha$ under the limit of
the comultiplication maps $\varprojlim C_\alpha\rarrow
\varprojlim C_\alpha\ot_k C_\alpha$.
 The same applies to limits and colimits in the category of graded
coalgebras $k\coalg$ (with its different forgetful functor).
 Colimits in the category of conilpotent coalgebras $k\coalg^\conilp$
are compatible with those in $k\coalg^\coaug$, and limits in
$k\coalg^\conilp$ are the maximal conilpotent subcoalgebras of
the limits in $k\coalg^\coaug$.
 A simpler approach is to construct directly finite products in
$k\coalg^\conilp$ as the cofreely cogenerated conilpotent coalgebras.

 This suffices to clarify part~(b), so let us turn to~(a).
 Any diagram with a final vertex in $k\coalg_\cdg^\conilp$ is
isomorphic to a diagram of strict morphisms, i.~e., morphisms
of the form $(f,0)$, so its limit in $k\coalg_\cdg^\conilp$ 
agrees with that in $k\coalg^\conilp$.
 Products of pairs of objects and equalizers in $k\coalg_\cdg^\conilp$ 
have to be constructed explicitly.
 Here it is instructive to start with the dual case of coproducts of
pairs and coequalizers in the category of CDG\+algebras $k\alg_\cdg$.
 Given two morphisms $(f',a')$ and $(f'',a'')\:B\rarrow A$ in
$k\alg_\cdg$, one constructs their coequalizer as the quotient
algebra of $A$ by the two-sided ideal generated by all the elements
$f'(b)-f''(b)$ and $a'-a''\in A$, where $b\in B$.
 This is a certain graded quotient algebra of the coequalizer of
$f'$ and~$f''$ in the category of graded algebras $k\alg$.
 Given two CDG\+algebras $A$ and $B$, one constructs their coproduct
$A\sqcup B$ as the graded algebra freely generated by $A$, \ $B$, and
an element~$c$ of degree~$1$.
 The restriction of the differential in $A\sqcup B$ to $A$ coincides
with~$d_A$ and its restriction to $B$ is given by the formula
$b\mpsto d_B(b)-[c,b]$ for $b\in B$.
 The action of the differential in $A\sqcup B$ on the element $c$
is obtained from the equation $h_A+dc+c^2=h_B$.
 The constructions in $k\coalg_\cdg^\conilp$ are dual.

 Any diagram of a filtered inductive limit in $k\coalg^\conilp_\cdg$ is
isomorphic to a diagram of strict morphisms.
 This follows from the vanishing of the first derived functor of
filtered colimit in $k\vect$ computed in terms of the standard
bar-construction.
 Colimits of all diagrams of strict morphisms exist in
$k\coalg^\conilp_\cdg$ and are preserved by the three
forgetful functors.
\end{proof}

 The category $k\coalg^\conilp_\cdg$ has no final object, and fibered
coproducts sometimes do not exist in it.
 The closest candidate for a final object in $k\coalg^\conilp_\cdg$
is the following CDG\+coalgebra $(D,d,h)$.
 The coalgebra $D$ is cofreely cogenerated by the linear function~$h$,
i.~e., it is the graded tensor coalgebra of a one-dimensional vector
space concentrated in degree~$-2$.
 The differential~$d$ is zero.
 For any CDG\+coalgebra $D$ there is a morphism of CDG\+coalgebras
$C\rarrow D$; however, such a morphism is not unique.
 The CDG\+coalgebra $D$ is Koszul dual to a DG\+algebra with zero
cohomology.

 Notice also that $(D,d,h)$ is the final object in the category of
conilpotent CDG\+coalgebras and strict morphisms between them.
 It follows from the above proof that all limits and colimits exist
in the latter category.

\subsection{Model category of CDG\+coalgebras}
 In this subsection we presume all coalgebras to be conilpotent.
 In particular, by the graded coalgebra cofreely cogenerated by
a graded vector space, or the graded coalgebra cofreely cogenerated
by another graded coalgebra and a graded vector space, we mean
the corresponding universal objects in the category of conilpotent
graded coalgebras.
 These in general differ from the analogous universal objects in
the categories of arbitrary (or coaugmented) graded coalgebras.
 The same applies to the DG\+coalgebra cofreely cogenerated by
another DG\+coalgebra and a complex of vector spaces.
 Notice that the latter construction transforms quasi-isomorphisms
in either of its arguments to quasi-isomorphisms.

 All increasing filtrations below are presumed to be cocomplete.

 For any category $\sC$, denote by $\sC_\fin$ the category obtained
by formal adjoining of a final object~$*$ to~$\sC$.
 By the definition, there is one morphism into~$*$ from any object
of $\sC_\fin$, and there are no morphisms from~$*$ into any object
of $\sC$ but $*$~itself.
 We refer to the category $\sC_\fin$ as the \emph{finalized}
category~$\sC$.

\begin{thm}
\textup{(a)} There is a model category structure on the finalized
category of conilpotent CDG\+coalgebras $k\coalg^\conilp_{\cdg,\fin}$
with the following properties.
 A morphism of CDG\+coalgebras is a weak equivalence if and only if
it belongs to the minimal class of morphisms containing the filtered
quasi-isomorphisms and satisfying the two-out-of-three axiom.
 Also a morphism of CDG\+coalgebras is a weak equivalence if and only if
it is a composition of retracts of filtered quasi-isomorphisms.
 A morphism of CDG\+coalgebras is a cofibration if and only if
the underlying morphism of graded coalgebras is injective.
 A morphism of CDG\+coalgebras $C\rarrow D$ is a fibration if and
only if the graded coalgebra $C^\#$ is cofreely cogenerated by
the graded coalgebra $D^\#$ and a graded vector space.
 A morphism of CDG\+coalgebras is a trivial cofibration if and only if
it belongs to the minimal class of morphisms which contains
the injective filtered quasi-isomorphisms strictly compatible with
the filtrations, is closed under the composition, and contains
a morphism~$g$ whenever it contains morphisms $f$ and~$fg$.
 Also a morphism of CDG\+coalgebras is a trivial cofibration if and
only if it is a retract of an injective filtered quasi-isomorphism
strictly compatible with the filtrations.
 A morphism of CDG\+coalgebras is a trivial fibration if and only if
it is a retract of a morphism of CDG\+coalgebras $C\rarrow D$
such that $C$ and $D$ admit increasing filtrations $F$ compatible
with the comultiplications and differentials such that $F_0C=w(k)=F_0D$
and the DG\+coalgebra\/ $\gr_FC$ is cofreely cogenerated by
the DG\+coalgebra\/ $\gr_FD$ and an acyclic complex of vector spaces.
 The only morphism from any CDG\+coalgebra to the object~$*$ is
a cofibration, but not a weak equivalence.
 All objects of $k\coalg^\conilp_{\cdg,\fin}$ are cofibrant.
 A CDG\+coalgebra is fibrant if and only if its underlying graded
coalgebra is a cofree conilpotent coalgebra (graded tensor coalgebra). 
\par
\textup{(b)} There is a model category structure of the category of
conilpotent DG\+coalgebras $k\coalg^\conilp_\dg$ with the following
properties.
 A morphism of DG\+coalgebras is a weak equivalence if and only if
it belongs to the minimal class of morphisms containing the filtered
quasi-isomorphisms and satisfying the two-out-of-three axiom.
 Also a morphism of DG\+coalgebras is a weak equivalence if and only if
it is a composition of retracts of filtered quasi-isomorphisms.
 A morphism of DG\+coalgebras is a cofibration if and only if it is
injective.
 A morphism of DG\+coalgebras $C\rarrow D$ is a fibration if and
only if the graded coalgebra $C^\#$ is cofreely cogenerated by
the graded coalgebra $D^\#$ and a graded vector space.
 A morphism of DG\+coalgebras is a trivial cofibration if and only if
it belongs to the minimal class of morphisms which contains
the injective filtered quasi-isomorphisms strictly compatible with
the filtrations, is closed under the composition, and contains
a morphism~$g$ whenever it contains morphisms $f$ and~$fg$.
 Also a morphism of DG\+coalgebras is a trivial cofibration if and
only if it is a retract of an injective filtered quasi-isomorphism
strictly compatible with the filtrations.
 A morphism of DG\+coalgebras is a trivial fibration if and only if
it is a retract of a morphism of DG\+coalgebras $C\rarrow D$
such that $C$ and $D$ admit increasing filtrations $F$ compatible
with the comultiplications and differentials such that $F_0C=w(k)=F_0D$
and the DG\+coalgebra\/ $\gr_FC$ is cofreely cogenerated by
the DG\+coalgebra\/ $\gr_FD$ and an acyclic complex of vector spaces.
 All objects of $k\coalg^\conilp_\dg$ are cofibrant.
 A DG\+coalgebra is fibrant if and only if its underlying graded
coalgebra is a cofree conilpotent coalgebra (graded tensor coalgebra).
\end{thm}

\begin{proof}
 We will prove part~(a).
 Existence of limits and colimits in $k\coalg^\conilp_{\cdg,\fin}$
follows from Lemma~\ref{cdg-limits}.
 It suffices to check that coequalizers exist in
$k\coalg^\conilp_{\cdg,\fin}$.

 In the rest of the proof we use the Koszul duality functors
$\Br_v$ and $\Cb_w$ together with
Theorem~\ref{dg-algebras-model-category}.
 It will follow from our argument that the weak equivalences,
fibrations, and trivial fibrations of CDG\+coalgebras can be also
characterized in the following ways.
 A morphism of CDG\+coalgebras is a weak equivalence if and only if
the functor $\Cb_w$ transforms it into a quasi-isomorphism.
 A morphism of CDG\+coalgebras is a fibration if and only if it is 
a retract of a morphism obtained by a base change from a morphism
obtained by applying the functor $\Br_v$ to a fibration of DG\+algebras.
 A morphism of CDG\+coalgebras is a trivial fibration if and only if
it is a retract of a morphism obtained by a base change from a morphism
obtained by applying the functor $\Br_v$ to a trivial fibration of
DG\+algebras.

 The proof is based on several Lemmas.

\begin{lem1}
\textup{(i)} Any morphism of CDG\+coalgebras $C\rarrow D$ can be
decomposed into an injective morphism of CDG\+coalgebras $C\rarrow E$
followed by a morphism $E\rarrow D$ obtained by a base change from
a morphism obtained by applying the functor $\Br_v$ to a surjective
quasi-isomorphism of DG\+algebras. 
 Furthermore, the CDG\+coalgebras $E$ and $D$ admit increasing
filtrations $F$ compatible with the comultiplications and differentials
such that $F_0E=k=F_0D$ and the DG\+coalgebra\/ $\gr_FE$ is cofreely
cogenerated by the DG\+coalgebra\/ $\gr_FD$ and an acyclic complex of
vector spaces. \par
\textup{(ii)} Any morphism of CDG\+coalgebras $C\rarrow D$ can be
can be decomposed into an injective filtered quasi-isomorphism
$C\rarrow E$ strictly compatible with the filtrations and
a morphism $E\rarrow D$ obtained by a base change from a morphism
obtained by applying the functor $\Br_v$ to a surjective morphism
of DG\+algebras.
\end{lem1}

\begin{proof}
 This Lemma is based on a construction of Hinich~\cite{Hin2}.
 To prove~(i), decompose the morphism of DG\+algebras $\Cb_w(C)
\rarrow\Cb_w(D)$ into an injective morphism $\Cb_w(C)\rarrow A$
followed by a surjective quasi-isomorphism $A\rarrow\Cb_w(D)$.
 Consider the induced morphism $\Br_v(A)\rarrow\Br_v\Cb_w(D)$
and set $E$ to be the fibered product of the CDG\+coalgebras
$\Br_v(A)$ and $D$ over $\Br_v\Cb_w(D)$.
 Then $C\rarrow E\rarrow D$ is the desired decomposition.
 Now let $E\rarrow D$ be the morphism of CDG\+coalgebras obtained
from the morphism $\Br_v(A)\rarrow\Br_v(B)$ induced by a morphism
of DG\+algebras $A\rarrow B$ by the base change with respect to
a morphism of CDG\+coalgebras $D\rarrow\Br_v(B)$. 
 Define an increasing filtration $F$ on the DG\+algebra $A$
by the rules $F_0A=k$, \ $F_1A=k\oplus \ker(A\to B)$, and $F_2A=A$.
 Denote also by $F$ the induced filtrations on the quotient
DG\+algebra $B$ of $A$ and the bar-constructions $\Br_v(A)$ and
$\Br_v(B)$. 
 Define an increasing filtration $F$ on the CDG\+coalgebra $D$
by the rules $F_{2n}D=F_{2n+1}D=G_nD$, where $G$ is the natural
increasing filtration defined in~\ref{conilpotent-cobar-duality}
(where it is denoted by~$F$).
 Denote by $F$ the induced filtration on the fibered product~$E$. 
 Then the DG\+coalgebra $\gr_FE$ is the fibered product of
the DG\+coalgebras $\Br_v\gr_FA$ and $\gr_FD$ over $\Br_v\gr_FB$,
and the DG\+coalgebra $\Br_v\gr_FA$ is cofreely cogenerated
by $\Br_v\gr_FB$ and an acyclic complex of vector spaces, hence
the morphism of DG\+algebras $\gr_FE\rarrow\gr_FD$ has the same
property.

 Let us prove part~(ii).
 Let $A$ be the DG\+algebra freely generated by the DG\+algebra
$\Cb_w(C)$ and the acyclic complex of vector spaces
$\id_{\Cb_w(D)}[-1]$.
 Then the morphism of DG\+algebras $\Cb_w(C)\rarrow\Cb_w(D)$
factorizes into a trivial cofibration $\Cb_w(C)\rarrow A$ followed
by a fibration $A\rarrow\Cb_w(C)$.
 Let $F$ denote the natural increasing filtrations on the conilpotent
CDG\+coalgebras $C$ and $D$ and the induced filtrations on
$\Cb_w(C)$, \ $\Cb_w(D)$, and~$A$.
 Then the morphism $\Cb_w(\gr_FC)\rarrow\gr_FA$ is an injective
quasi-isomorphism and the morphism $\gr_FA\rarrow\Cb_w(\gr_FD)$
is surjective.
 Set $E$ to be the fibered product of the CDG\+coalgebras $\Br_v(A)$
and $D$ over $\Br_v\Cb_w(D)$.
 Let us show that the decomposition $C\rarrow E\rarrow D$ has
the desired properties.
 Denote by $F$ the induced filtration on $E$ as a CDG\+subcoalgebra
of $\Br_v(A)$.
 The morphism $\gr_FC\rarrow\Br_v(\gr_FA)$ is an injective
quasi-isomorphism, and it only remains to check that the morphism
$\gr_FE\rarrow\Br_v(\gr_FA)$ is a quasi-isomorphism.
 Notice that the DG\+coalgebra $\gr_FE$ is the fibered product of
$\Br_v(\gr_FA)$ and $\gr_FD$ over $\Br_v\Cb_w(\gr_FD)$.

 Let $n$ denote the nonnegative grading induced from the indexing of
the filtration~$F$.
 Introduce a decreasing filtration $G$ on $\gr_FA$ by the rules
that $G^0\gr_FA=\gr_FA$ and $G^j\gr_FA$ is the sum of the components
of the ideal $\ker(\gr_FA\to\Cb_w(\gr_FD))$ situated in the grading
$n\ge j$.
 This filtration is locally finite with respect to the grading~$n$.
 Let $G$ denote the induced decreasing filtration on $\gr_FE$ as
a DG\+subcoalgebra of $\Br_v(\gr_FA)$.
 Then the DG\+coalgebra $\gr_G\gr_FE$ is the fibered product of
$\Br_v(\gr_G\gr_FA)$ and $\gr_FD$ over $\Br_v\Cb_w(\gr_FD)$.
 The DG\+algebra $\gr_G\gr_FA$ and the DG\+coalgebra $\gr_G\gr_FD$
have two nonnegative gradings $n$ and~$j$.
 Introduce an increasing filtration $H$ on the DG\+algebra
$\gr_G\gr_FA$ for which the component $H_t$ is the sum of all
the components $(\gr_G\gr_FA)_{n,j}$ with $j=0$ and $t\ge n$
or $j>0$ and $t\ge 1$.
 Once again, denote by $H$ the induced increasing filtration on
$\gr_G\gr_FE$ as a DG\+subcoalgebra of $\Br_v(\gr_G\gr_FA)$.
 Then the DG\+coalgebra $\gr_H\gr_G\gr_FE$ is the fibered product
of $\Br_v(\gr_H\gr_G\gr_FA)$ and $\gr_FD$ over $\Br_v\Cb_w(\gr_FD)$.
 The DG\+coalgebra $\Br_v(\gr_H\gr_G\gr_FA)$ is cofreely cogenerated
by $\Br_v\Cb_w(\gr_FD)$ and a complex of vector spaces, and
the DG\+coalgebra $\gr_H\gr_G\gr_FE$ is cofreely cogenerated by
$\gr_FD$ and the same complex of vector spaces.
 Since the morphism $\gr_FD\rarrow\Br_v\Cb_w(\gr_FD)$ is
a quasi-isomorphism, so is the morphism $\gr_H\gr_G\gr_FE
\rarrow\Br_v(\gr_H\gr_G\gr_FA)$.
\end{proof}

\begin{lem2}
\textup{(i)} Injective filtered quasi-isomorphisms of CDG\+coalgebras
strictly compatible with the filtrations have the lifting property
with respect to the morphisms of CDG\+coalgebras $C\rarrow D$ such that
the graded coalgebra $C^\#$ is cofreely cogenerated by
the graded coalgebra $D^\#$ and a graded vector space. \par
\textup{(ii)} The class of morphisms of conilpotent graded coalgebras
$C\rarrow D$ such that $C$ is cofreely generated by $D$ and
a graded vector space is closed under retracts.
 If the associated quotient morphism of a morphism~$f$ with respect
to an increasing filtration belongs to this class, then
the morphism~$f$ itself does.
\end{lem2}

\begin{proof}
 Part~(i): let $g\:X\rarrow Y$ be an injective morphism of
CDG\+coalgebras endowed with increasing filtrations $F$ making
$f$ a filtered quasi-isomorphism strictly compatible with
the filtrations, i.~e., $F_nX=g^{-1}(F_nY)$.
 Clearly, we can restrict ourselves to diagrams of strict morphisms
of CDG\+coalgebras, i.~e., morphisms of the form $(f,0)$.
 Suppose that a morphism $Y\rarrow D$ has been lifted to morphisms
$X\rarrow C$ and $F_{n-1}Y\rarrow C$ in compatible ways; we need
to extend the second morphism to $F_nY$ in a compatible way.
 The quotient space $F_nY/(F_{n-1}Y+g(F_nX))$ is an acyclic complex of
vector spaces.
 Choose a graded subspace $V'$ in this complex such that
the restriction of the differential to $V'$ is injective and 
the complex coincides with $V'+d(V')$.
 Let $V$ be any graded subspace in $F_nY$ which projects
isomorphically onto~$V'$.
 By the assumption about the morphism $C\rarrow D$, one can extend
the morphism $F_{n-1}Y\rarrow C$ to a graded coalgebra morphism
$F_{n-1}Y+V\rarrow C$ in a way compatible with the morphisms into~$D$.
 Then the condition of compatibility with the differentials allows
to extend this morphism in the unique way to a strict CDG\+coalgebra
morphism $F_{n-1}Y+V+dV\rarrow C$.
 Finally, one combines this morphism with the morphism $F_nX\rarrow C$
to obtain the desired CDG\+coalgebra morphism $F_nY\rarrow C$.

 Part~(ii): This class of morphisms of conilpotent graded coalgebras
is characterized by the conditions that the induced morphism
$\Ext^1_C(k,k)\rarrow\Ext^1_D(k,k)$ is surjective and
the morphism $\Ext^2_C(k,k)\rarrow\Ext^2_D(k,k)$ is 
an isomorphism.
 Notice also the natural isomorphism $\Ext^i_E(k,k)\simeq
\Cotor_{-i}^E(k,k)$ for a coaugmented coalgebra~$E$.
\end{proof}

\begin{lem3}
\textup{(i)} The functor\/ $\Cb_w$ maps injective morphisms of
CDG\+coalgebras to cofibrations of DG\+algebras. \par
\textup{(ii)} Injective morphisms of CDG\+coalgebras have the lifting
property with respect to morphisms obtained by a base change from
morphisms obtained by applying the functor $\Br_v$ to surjective
quasi-isomorphisms of DG\+algebras. \par
\textup{(iii)} Injective morphisms of CDG\+coalgebras that are
transformed into quasi-isomorphisms by the functor\/ $\Cb_w$ have
the lifting property with respect to morphisms obtained by a base
change from morphisms obtained by applying the functor\/ $\Br_v$
to surjective morphisms of DG\+algebras.
\end{lem3}

\begin{proof}
 Parts (ii) and~(iii) follow immediately from part~(i) and
the adjunction of functors $\Br_v$ and $\Cb_w$.
 To check~(i), notice that for an injective morphism of conilpotent
CDG\+coalgebras $X\rarrow Y$, the DG\+algebra $\Br_v(Y)$ has
the form $\Br_v(X)\langle x_{n,\alpha}\rangle$.
 Indeed, let $F$ denote the natural increasing filtration on $Y$.
 Then the DG\+algebra $\Br_v(X+F_nY)$ has the form
$\Br_v(X+F_{n-1}Y)\langle x_{1,\alpha}, x_{2,\beta}\rangle$.
\end{proof}

 Now we can finish the proof of Theorem.

 It follows from (the proof of) Theorem~\ref{co-algebra-bar-duality}
that the minimal class of morphisms of CDG\+algebras containing
the filtered quasi-isomorphisms and satisfying the two-out-of-three
axiom coincides with the class of morphisms that are transformed into
quasi-isomorphisms by the functor $\Cb_w$.
 Analogously, the minimal class of morphisms that contains
the injective filtered quasi-isomorphisms strictly compatible with
the filtrations, is closed under the composition, and contains
a morphism~$g$ whenever it contains morphisms $f$ and~$fg$ coincides
with the class of injective morphisms that are transformed into
quasi-isomorphisms by the functor $\Cb_w$.

 Let $X\rarrow Y$ be a morphism of the latter class.
 Using Lemma~1(ii), decompose it into an injective filtered
quasi-isomorphisms $X\rarrow E$ strictly compatible with
the filtrations, followed by a morphism $E\rarrow Y$ obtained by
a base change from a morphism obtained by applying the functor
$\Br_v$ to a surjective morphism of DG\+algebras.
 By Lemma~3(iii), the morphism $X\rarrow Y$ has the lifting property
with respect to the morphism $E\rarrow Y$, hence the morphism 
$X\rarrow Y$ is a retract of the morphism $X\rarrow E$.

 Let $C\rarrow D$ be a morphism which is transformed into 
a quasi-isomorphism by the functor $\Cb_w$ and such that
$C^\#$ is cofreely cogenerated by $D^\#$ and a graded vector space.
 Using Lemma~1(i), decompose it into an injective morphism
$C\rarrow E$ followed by a morphism $E\rarrow D$ obtained by a base
change from a morphism obtained by applying the functor $\Br_v$ to
a surjective quasi-isomorphism of DG\+algebras.
 The latter morphism is also a filtered quasi-isomorphism, and
even admits increasing filtrations $F$ on $E$ and $D$ such that
the DG\+coalgebra $\gr_FE$ is cofreely cogenerated by $\gr_FD$ and
an acyclic complex of vector spaces.
 Hence the injective morphism $C\rarrow E$ is also transformed into
a quasi-isomorphism by the functor $\Cb_w$, and therefore is
a retract of an injective filtered quasi-isomorphism strictly
compatible with the filtrations.
 By Lemma~2(i), the morphism $C\rarrow E$ has the lifting property
with respect to the morphism $C\rarrow D$, so the morphism
$C\rarrow D$ is a retract of the morphism $E\rarrow D$.

 Let $C\rarrow D$ be a morphism that is transformed into
a quasi-isomorphism by the functor $\Cb_w$.
 Using Lemma~1(i) or~(ii), one can decompose it into an injective
morphism $C\rarrow E$ and a morphism $E\rarrow D$ such that
$E^\#$ is cofreely cogenerated by $D^\#$ and a graded vector space,
and both morphisms $C\rarrow E$ and $E\rarrow D$ are transformed
into quasi-isomorphisms by the functor $\Cb_w$.
 Thus the morphism $C\rarrow D$ is the composition of two retracts of
filtered quasi-isomorphisms.

 The rest of the proof is straightforward.
 In addition to the above Lemmas, one has to use the analogue
of Lemma~2 designed to imply the lifting and retraction properties
of the morphisms $C\rarrow *$, where $C^\#$ is a cofree conilpotent
graded coalgebra.
\end{proof}

 Consider the finalized category of DG\+algebras with nonzero units
$k\alg^+_{\dg,\fin}$.
 It has a model category structure in which a morphism of DG\+algebras
is a weak equivalence, cofibration, or fibration if and only if it
belongs to the corresponding class of morphisms in $k\alg_\dg$, and
for any DG\+algebra $A$ the morphism $A\rarrow *$ is a fibration and
a cofibration, but not a weak equivalence.
 It follows from the above that the functors $\Cb_w$ and $\Br_v$
define a Quillen equivalence between the model categories
$k\coalg^\conilp_{\cdg,\fin}$ and $k\alg^+_{\dg,\fin}$. 
 Let us point out that the natural Quillen adjunction between
the categories $k\alg_\dg$ and $k\alg^+_{\dg,\fin}$ is \emph{not}
a Quillen equivalence.
 Analogously, the functors $\Cb_w$ and $\Br_v$ define a Quillen
equivalence between the model categories $k\coalg^\conilp_\dg$
and $k\alg^\aug_\dg$.

\subsection{Cofibrant DG\+algebras}  \label{cofibrant-dg-alg}
 The following result demonstrates the importance of DG\+algebras
that are cofibrant in the model category structure of
Theorem~\ref{dg-algebras-model-category}(a).
 Let $k$ be a commutative ring of finite homological dimension.

\begin{thm}
 Let $A$ be a cofibrant DG\+algebra over~$k$.
 Then all the four triangulated subcategories\/ $\Acycl^\abs(A\modl)$, \
$\Acycl^\co(A\modl)$, \ $\Acycl^\ctr(A\modl)$, and\/
$\Acycl(A\modl)\subset\Hot(A\modl)$ coincide.
 Consequently, the derived categories of the first and the
second kind\/ $\sD(A\modl)$ and\/ $\sD^\abs(A\modl)$ are isomorphic.
\end{thm}

\begin{proof}
 The graded ring $A^\#$ has the homological dimension exceeding
the homological dimension of~$k$ by at most~$1$, being a retract of
a free graded algebra over~$k$.
 Hence $\Acycl^\co(A\modl)=\Acycl^\abs(A\modl)=\Acycl^\ctr(A\modl)$
by Theorem~\ref{finite-homol-dim-cdg-ring}(a).
 Let us prove that $\Acycl^\abs(A)=\Acycl(A)$ for a DG\+algebra $A$
of the form $A=k\langle x_{n,\alpha}\rangle$; this is clearly
sufficient.
 Let $C_+[-1]$ denote the free graded $k$\+module spanned by
$x_{n,\alpha}$; then $C=k\oplus C_+$ can be considered as a (strictly
coaugmented strictly counital) curved $\Ainfty$\+coalgebra over~$k$
(see~\ref{strictly-co1-curved-ainfty-coalgebras}).
 Set $F_nC$ to be the linear span of $k$ and $x_{m,\beta}$ for
$m\le n$.
 Let $M$ be a left DG\+module over~$A$.
 Consider the surjective closed morphism of DG\+modules
$A\ot_kM\rarrow M$.
 Its kernel can be naturally identified with $A\ot_k C_+[-1]\ot_kM$,
so there is the induced DG\+module structure on this triple tensor
product.
 The cone of the morphism of DG\+modules $A\ot_k C_+[-1]\ot_k M
\rarrow A\ot_kM$ can be naturally identified with $A\ot_kC\ot_kM$,
so on the latter triple tensor product there is also a natural
DG\+module structure.
 Consider the filtration $F$ on $A\ot_k C\ot_kM$ induced by
the filtration $F$ on~$C$.
 This filtration is compatible with the differential and the quotient
DG\+modules $(A\ot_k F_nC\ot_kM)/(A\ot_k F_{n-1}C\ot_kM)$ are
isomorphic to the DG\+modules $A\ot_k F_nC/F_{n-1}C\ot_kM$ with
the differentials induced by the differentials on $A$ and~$M$.
 It follows that the DG\+module $A\ot_kC\ot_kM$ is projective in
the sense of~\ref{projective-dg-modules} whenever $M$ is
a projective complex of $k$\+modules.
 Furthermore, this DG\+module is coacyclic whenever $M$ is coacyclic
as a complex of $k$\+modules.
 Since the cone of the morphism $A\ot_kC\ot_kM\rarrow M$ is always
absolutely acyclic and the ring~$k$ has a finite homological
dimension, it follows that $M$ is absolutely acyclic whenever
it is acyclic.
\end{proof}

 Let $A$ be a DG\+ring for which the underlying graded ring
$A^\#$ has a finite left homological dimension.
 When do the triangulated subcategories $\Acycl(A)$ and
$\Acycl^\abs(A)$ coincide?
 This cannot happen too often, as the absolute derived categories
$\sD^\abs(A\modl)$ and $\sD^\abs(B\modl)$ are isomorphic for any
two DG\+rings $A$ and $B$ that are isomorphic as CDG\+rings.
 At the same time, the derived category $\sD(A\modl)$ only depends
on the quasi-isomorphism class of $A$.
 These are two very different and hardly ever compatible kinds of
functoriality.
 The quasi-isomorphism classes of CDG\+isomorphic DG\+rings $A$ and
$B$ can be entirely unrelated to each other; see Examples below.

 Nevertheless, there are the following cases, in addition to
the above Theorem and Corollary~\ref{koszul-cogenerators}.2.
 When either $A^i=0$ for all $i>0$, or $A^i=0$ for all $i<0$,
the ring $A^0$ is semisimple, and $A^1=0$, one has
$\Acycl(A)=\Acycl^\abs(A)$.
 This follows from Theorems~\ref{dg-ring-bounded-cases}.1(d)--
\ref{dg-ring-bounded-cases}.2(d) and~\ref{finite-homol-dim-cdg-ring}(a).
 Besides, using Theorem~\ref{finite-homol-dim-cdg-ring}(b) one
can check~\cite{KLN} that $\Acycl(A)=\Acycl^\abs(A)$ whenever
$A$ has a zero differential.
 All of this assumes that $A^\#$ has a finite left homological
dimension; without this assumption, there are only some partial
results obtained in~\ref{dg-ring-bounded-cases}.

\begin{exs}
 Let $V$ be a (totally) finite-dimensional complex of $k$\+vector
spaces and $A$ be the graded algebra of endomorphisms of~$V$ endowed
with the induced differential~$d$.
 Then the DG\+algebra $(A,d)$ is isomorphic to the DG\+algebra
$(A,0)$ in the category of CDG\+algebras over~$k$.
 In particular, when $V$ is acyclic but nonzero, $(A,d)$ is also
acyclic, while $(A,0)$ is not.
 So the derived category $\sD((A,d)\modl)$ vanishes, while
the derived category $\sD((A,0)\modl)$ is equivalent to $\sD(k\vect)$.
 At the same time, the absolute derived categories
$\sD^\abs((A,d)\modl)$ and $\sD^\abs((A,0)\modl)$ are isomorphic.
 Furthermore, let $(D',d',h')$ and $(D'',d'',h'')$ be two coaugmented
CDG\+coalgebras with coaugmentations~$w$; then the CDG\+coalgebra
$(D'\oplus D''\;d'+d''\;h'+h'')$ has two induced coaugmentations
$w'$ and~$w''$.
 The DG\+algebra $\Cb_{w'}(D'\oplus D'')$ is quasi-isomorphic to
$\Cb_w(D')$ and the DG\+algebra $\Cb_{w''}(D'\oplus D'')$ is
quasi-isomorphic to $\Cb_w(D'')$.
 At the same time, the DG\+algebras $\Cb_{w'}(D'\oplus D'')$ and
$\Cb_{w''}(D'\oplus D'')$ are isomorphic as CDG\+algebras over~$k$.
\end{exs}

\begin{rem}
 One can call a strictly coaugmented strictly counital curved
$\Ainfty$\+coalgebra $C$ over a field~$k$ conilpotent if
the DG\+algebra $A=\bigoplus_n C_+[-1]^{\ot n}$ is cofibrant.
 Unlike comodules over a conilpotent coalgebra, which are all
conilpotent in the appropriate sense, not every strictly counital
curved $\Ainfty$\+comodule over a conilpotent curved
$\Ainfty$\+coalgebra is conilpotent.
 For any strictly counital curved $\Ainfty$\+coalgebra $C$,
consider the CDG\+coalgebra $U=\Br_v(A)$; it can be called
the conilpotent coenveloping CDG\+coalgebra of~$C$.
 The DG\+category of conilpotent curved $\Ainfty$\+comodules
over $C$ and strict morphisms between them can be defined as
the DG\+category of CDG\+comodules over~$U$.
 The faithful DG\+functor $N\mpsto A\ot^{\tau_{A,v}}\!\.N$ provides
an embedding of the DG\+category of conilpotent curved
$\Ainfty$\+comodules over $C$ and strict morphisms between them
into the DG\+category of strictly counital curved
$\Ainfty$\+comodules over~$C$.
 Analogously, the DG\+category of contranilpotent curved
$\Ainfty$\+contra\-modules over $C$ and strict morphisms between them
can be defined as the DG\+category of CDG\+contra\-modules over~$U$.
 The DG\+functor $Q\mpsto\Hom^{\tau_{A,v}}(A,Q)$ maps
the DG\+category of contranilpotent curved $\Ainfty$\+contra\-modules
over $C$ and strict morphisms between them to the DG\+category of
strictly counital curved $\Ainfty$\+contra\-modules.
 By Theorem~\ref{bar-construction-duality} and the above Theorem,
for any conilpotent curved $\Ainfty$\+coalgebra $C$
the coderived category of conilpotent curved $\Ainfty$\+comodules
and strict morphisms between them is equivalent to the coderived
category of strictly counital curved $\Ainfty$\+comodules,
and the contraderived category of contranilpotent curved
$\Ainfty$\+contra\-modules and strict morphisms between them is
equivalent to the contraderived category of strictly counital
curved $\Ainfty$\+contra\-modules.
\end{rem}

\appendix
\Section{Homogeneous Koszul Duality}
\label{homogeneous-koszul-appendix}

\subsection{Covariant homogeneous duality}  \label{covariant-duality}
 We will consider DG\+algebras, DG\+coal\-gebras, DG\+modules,
DG\+comodules, DG\+contramodules endowed with an additional
$\boZ$\+valued grading, which will be called the \emph{internal
grading} and denoted by lower indices.
 The internal grading is always assumed to be preserved by
the differentials, $d(X_n)\subset X_n$.
 All morphisms of internally graded objects are presumed to
preserve the internal gradings.
 The other grading, raised by~$1$ by the differentials, is also
always present; it is called the \emph{cohomological grading} and
denoted by upper indices, as above in this paper.
 We will always use the cohomological grading only in all the sign
rules, so the notation $|z|$ is understood to refer to
the cohomological grading.
 The same applies to the notation $V\mpsto V[1]$, which is
interpreted as shifting the cohomological grading and leaving
the internal grading unchanged.

 Let $k$ be a field, and let $A$ be a DG\+algebra and $C$ be
a DG\+coalgebra over~$k$ endowed with internal gradings.
 Assume that the bigraded vector spaces $A_+=A/k$ and $C_+=\ker
(C\to k)$ are concentrated in the positive internal grading,
i.~e., $A_+=\bigoplus_{n>0}A_n$ and $C_+=\bigoplus_{n>0}C_n$,
where $A_n$ and $C_n$ are complexes of vector spaces.
 Then $A$ admits a unique augmentation $v\:A\rarrow k$ preserving
the internal grading and $C$ admits a unique coaugmentation
$w\:k\rarrow C$ with the same property.
 Notice that any positively internally graded DG\+coalgebra $C$
is conilpotent.
 The bar and cobar-constructions $\Br_v(A)$ and $\Cb_w(C)$
from~\ref{bar-cobar-constr} produce a DG\+coalgebra and
a DG\+algebra endowed with internal and cohomological gradings. 

 We will presume twisting cochains $\tau\:C\rarrow A$
(see~\ref{twisting-cochains-subsect}) to have the internal degree~$0$
and the cohomological degree~$1$; besides, $\tau$ is assumed to satisfy
the equations $v\circ\tau=\tau\circ w=0$.
 These conditions hold for the natural twisting cochains
$\tau_{A,v}\:\Br_v(A)\rarrow A$ and $\tau_{C,w}\:C\rarrow\Cb_w(C)$.
 There are natural bijective correspondences between the three
sets: the set of all morphisms of internally graded DG\+algebras
$\Cb_w(C)\rarrow A$, the set of all morphisms of internally graded
DG\+coalgebras $C\rarrow\Br_v(A)$, and the set of all twisting
cochains $\tau\:C\rarrow A$ as above.
 So the functor $\Cb_w$ is left adjoint to the functor $\Br_v$. 

{\hbadness=5400
\begin{thm1}
 Let $\tau$ be a twisting cochain between a (positively
internally graded) DG\+coalgebra $C$ and DG\+algebra~$A$.
 Then the following conditions are equivalent: \par
\textup{(a)} the morphism of DG\+algebras\/ $\Cb_w(C)\rarrow A$
is a quasi-isomorphism; \par
\textup{(b)} the morphism of DG\+coalgebras $C\rarrow\Br_v(A)$
is a quasi-isomorphism; \par
\textup{(c)} the complex $A\ot^\tau\!\.C$ is quasi-isomorphic
to~$k$; \par
\textup{(d)} the complex $C\ot^\tau\!\.A$ is quasi-isomorphic
to~$k$. \par
 Besides, the functors\/ $\Br_v$ and\/ $\Cb_w$ transform
quasi-isomorphisms of DG\+algebras to quasi-isomorphisms of
DG\+coalgebras and vice versa.
 Consequently, these functors induce an equivalence between
the categories of DG\+algebras and DG\+coalgebras with
inverted quasi-isomorphisms.
\end{thm1}}

\begin{proof}
 The complex $A\ot^{\tau_{A,v}}\!\Br_v(A)$ is the reduced bar
resolution of the left DG\+module~$k$ over $A$, so the standard
contracting homotopy induced by the unit element of $A$ proves
its acyclicity.
 The same argument applies to the complexes $\Br_v(A)\ot^{\tau_{A,v}}
\!\.A$, \ $C\ot^{\tau_{C,w}}\!\Cb_w(C)$, and
$\Cb_w(C)\ot^{\tau_{C,w}}\!\.C$.
 If $(C',A',\tau')\rarrow(C'',A'',\tau'')$ is a morphism of twisting
cochains for which the morphism of DG\+coalgebras $C'\rarrow C''$
and the morphism of DG\+algebras $A'\rarrow A''$ are
quasi-isomorphisms, then the induced morphisms of complexes
$A'\ot^{\tau'}\!\.C'\rarrow A''\ot^{\tau''}\!\.C''$ 
and $C'\ot^{\tau'}\!\.A'\rarrow C''\ot^{\tau''}\!\.A''$ are
quasi-isomorphisms.
 Indeed, every internal degree component of $A\ot^\tau\!\.C$ and
$C\ot^\tau\!\.A$ is obtained from the tensor products of
the internal degree components of $A$ and $C$ by a finite number of
shifts and cones, and the operations of tensor product, shift, and
cone preserve quasi-isomorphisms of complexes of vector spaces.
 A similar argument proves that the functors $\Br_v$ and $\Cb_w$
map quasi-isomorphisms to quasi-isomorphisms, since the internal
degree components of $\Br_v(A)$ and $\Cb_w(C)$ are obtained from
the internal degree components of $A$ and $C$ by finite iterations
of tensor products, shifts, and cones.
 Finally, one shows by induction in the internal degree that if
a morphism of twisting cochains $(C,A',\tau')\rarrow (C,A'',\tau'')$
induces a quasi-isomorphism of complexes $A'\ot^{\tau'}\!\.C\rarrow
A''\ot^{\tau''}\!\.C$ or $C\ot^{\tau'}\!\.A'\rarrow
C\ot^{\tau''}\!\.A''$ then the morphism of $A'\rarrow A''$ is
a quasi-isomorphism, and analogously for morphisms of twisting
cochains $(C',A,\tau')\rarrow(C'',A,\tau'')$.
\end{proof}

 As in~\ref{acyclic-twisting-cochain}, a twisting cochain~$\tau$
satisfying the equivalent conditions of Theorem~1 is said to
be \emph{acyclic}.
 Any twisting cochain~$\tau$ induces two pairs of adjoint functors
between the homotopy categories of internally graded DG\+modules,
DG\+comodules, and DG\+contramodules $\Hot(A\modl)$, \
$\Hot(C\comodl)$, and $\Hot(C\contra)$.
 The functor $M\mpsto C\ot^\tau\!\.M\:\Hot(A\modl)\rarrow
\Hot(C\comodl)$ is right adjoint to the functor $N\mpsto A\ot^\tau\!\.
N\:\Hot(C\comodl)\rarrow\Hot(A\modl)$ and the functor $P\mpsto
\Hom^\tau(C,P)\:\Hot(A\modl)\rarrow\Hot(C\contra)$ is left adjoint
to the functor $Q\mpsto\Hom^\tau(A,Q)\:\Hot(C\contra)\rarrow
\Hot(A\modl)$.

 Let $\Hot(A\modl^\up)$ and $\Hot(A\modl^\down)$ denote the homotopy
categories of (internally graded) left DG\+modules over $A$
concentrated in nonnegative and nonpositive internal degrees,
respectively.
 Let $\Hot(C\comodl^\up)$ denote the homotopy category of left
DG\+comodules over $C$ concentrated in nonnegative internal degrees
and $\Hot(C\contra^\down)$ denote the homotopy category of left
DG\+contramodules over $C$ concentrated in nonpositive internal
degrees.
 Let $\sD(A\modl^\up)$, \ $\sD(A\modl^\down)$, \ $\sD(C\comodl^\up)$,
and $\sD(C\contra^\down)$ denote the corresponding derived categories,
i.~e., the triangulated categories obtained from the respective
homotopy categories by inverting the quasi-isomorphisms.

\begin{thm2}
 Let $\tau\:C\rarrow A$ be an acyclic twisting cochain. Then \par
\textup{(a)} The functors $C\ot^\tau\!\.{-}:\Hot(A\modl^\up)
\rarrow\Hot(C\comodl^\up)$ and $A\ot^\tau\!\.{-}\:\allowbreak
\Hot(C\comodl^\up)\rarrow\Hot(A\modl^\up)$ induce mutually inverse
equivalences between the derived categories\/ $\sD(A\modl^\up)$
and\/ $\sD(C\comodl^\up)$. \par
\textup{(b)} The functors $\Hom^\tau(C,{-})\:\Hot(A\modl^\down)
\rarrow\Hot(C\contra^\down)$ and $\Hom^\tau\allowbreak(A,{-})\:
\Hot(C\contra^\down)\rarrow\Hot(A\modl^\down)$ induce mutually inverse
equivalences between the derived categories\/ $\sD(A\modl^\down)$
and\/ $\sD(C\contra^\down)$.
\end{thm2}

\begin{proof}
 One can deduce this Theorem from
Theorem~\ref{acyclic-twisting-cochain}, using the facts that
any acyclic DG\+comodule from $\Hot(C\comodl^\up)$ is coacyclic
and any acyclic DG\+contra\-module from $\Hot(C\contra^\down)$ is
contraacyclic.
 Indeed, any nonnegatively internally graded DG\+comodule $N$
over~$C$ is the inductive limit of DG\+comodules obtained by
shifts and cones from the internal degree components $N_n$
considered as DG\+comodules with the trivial comodule structure.
 Analogously, any nonpositively internally graded DG\+contramodule
$Q$ over $C$ is the projective limit of DG\+contramodules obtained
by shifts and cones from the internal degree components $Q_n$
considered as DG\+contramodules with the trivial contramodule
structure. 
 (Similarly one can show that any nonnegatively internally graded
DG\+module over $A$ that is projective as a bigraded $A^\#$\+module
is projective in the sense of~\ref{projective-dg-modules} and any
nonpositively internally graded DG\+module over $A$ that is
injective as a bigraded $A^\#$\+module is injective in the sense
of~\ref{injective-dg-modules}.)
 Alternatively, here is a direct proof.

 The functors $M\mpsto C\ot^\tau\!\.M$ and $N\mpsto A\ot^\tau\!\.N$
map acyclic DG\+modules to acyclic DG\+comodules and vice versa,
since the internal degree components of $C\ot^\tau\!\.M$ and
$A\ot^\tau\!\.N$ are obtained from the internal degree components
of $M$ and $N$ by tensoring them with the internal degree
components of $C$ and $A$ and taking shifts and cones a finite
number of times for each of the components of $C\ot^\tau\!\.M$
and $A\ot^\tau\!\.N$.
 To check that the adjunction morphism $A\ot^\tau C\ot^\tau M
\rarrow M$ is a quasi-isomorphism, notice that the integral degree
components of its kernel can be obtained from the positive internal
degree components of $A\ot^\tau C$ by tensoring them with
the internal degree components of $M$ and taking shifts and cones
a finite number of times.
 Analogously one shows that the adjunction morphism $N\rarrow
C\ot^\tau A\ot^\tau N$ is a quasi-isomorphism.
 The proof of part~(b) is similar, with the only difference that
one has to use exactness of the functor $\Hom_k$ instead of~$\ot_k$.
\end{proof}

 The mutually inverse functors $C\ot^\tau{-}$ and $A\ot^\tau{-}$
between the derived categories $\sD(A\modl^\up)$ and
$\sD(C\comodl^\up)$ transform DG\+modules with the trivial action
of $A$ into cofree DG\+comodules over $C$, free DG\+modules over~$A$
into DG\+comodules with the trivial action of~$C$.
 The mutually inverse functors $\Hom^\tau(C,{-})$ and $\Hom^\tau
(A,{-})$ between the derived categories $\sD(A\modl^\down)$ and
$\sD(C\contra^\down)$ transform DG\+modules with the trivial action
of $A$ into free DG\+contramodules over $C$, cofree DG\+modules
over~$A$ into DG\+contramodules with the trivial action of~$C$.

 As explained in~\ref{acyclic-twisting-cochain}, for any DG\+module
$M\in\sD(A\modl^\up)$ the complex $C\ot^\tau\!\.M$ computes
$\Tor^A(k,M)$ and for any DG\+comodule $N\in\sD(C\comodl^\up)$
the complex $A\ot^\tau\!\.N$ computes $\Cotor^C(k,N)\simeq\Ext_C(k,N)$.
 For any DG\+module $P\in\sD(A\modl^\down)$ the complex
$\Hom^\tau(C,P)$ computes $\Ext_A(k,P)$ and for any DG\+contramodule
$Q\in\sD(C\contra^\down)$ the complex $\Hom^\tau(A,Q)$ computes
$\Coext_C(k,Q)\simeq\Ctrtor^C(k,Q)$.
 The DG\+coalgebra $C$ itself, considered as a complex, computes
$\Tor^A(k,k)$, and the DG\+algebra $A$, considered as a complex,
computes $\Cotor^C(k,k)\simeq\Ext_C(k,k)$.

\begin{cor}
 Let $f\:A'\rarrow A''$ be a morphism of DG\+algebras and
$g\:C'\rarrow C''$ be a morphism of DG\+coalgebras.
 Then the functor $\boI R_f\:\sD(A''\modl^\up)\rarrow\sD(A'\modl^\up)$
is an equivalence of triangulated categories if and only if
the functor $\boI R_f\:\sD(A''\modl^\down)\rarrow\sD(A'\modl^\down)$
is an equivalence of triangulated categories and if and only if
$f$~is a quasi-isomorphism.
 The functor $\boI R_g\:\sD(C'\comodl^\up)\rarrow\sD(C''\comodl^\up)$ is
an equivalence of triangulated categories if and only if
the functor $\boI R^g\:\sD(C'\contra^\down)\rarrow\sD(C''\contra^\down)$
is an equivalence of triangulated categories and if and only if
$g$~is a quasi-isomorphism.
\end{cor}

\begin{proof}
 The facts that whenever $f$ or $g$ is a quasi-isomorphism the related
restriction-of-scalars functors are equivalences of triangulated
categories can be deduced from Theorems~\ref{dg-mod-scalars}
and~\ref{cdg-coalgebra-scalars}, since any quasi-isomorphism of
positively internally graded DG-coalgebras is a filtered
quasi-isomorphism.
 To prove the inverse implications, it suffices to apply the adjoint
functors $\boL E_f$, \ $\boR E^f$, \ $\boR E_g$, and $\boL E^g$ to
the DG\+modules $A$ and $\Hom_k(A,k)$, the DG\+comodule $C$, and
the DG\+contramodule $\Hom_k(C,k)$.
 All assertions of Corollary also follow from Theorems~1 and~2.
\end{proof}

 Let $\oA=\bigoplus_{n=0}^\infty \oA_n$ and $\oC=\bigoplus_{n=0}^\infty
\oC_n$ be a Koszul graded algebra and a Koszul graded coalgebra
quadratic dual to each other~\cite{P0}.
 Let $C$ be the internally graded DG\+coalgebra with the components
$C_n^j=0$ for all $j\ne 0$ and $C_n^0=\oC_n$, and let $A$ be
the internally graded DG\+algebra with the components $A_n^j=0$ for
$j\ne n$ and $A_n^n=\oA_n$, the multiplication in $A$ and
the comultiplication in $C$ being given by the multiplication in
$\oA$ and the comultiplication in $\oC$, and the differentials in
$A$ and $C$ being zero.
 Let $\tau\:C\rarrow A$ be the twisting cochain that is given by
the isomorphism $\oC_1\simeq\oA_1$ in the internal degree~$1$ and
vanishes in all the other internal degrees.
 By the definition of quadratic duality and Koszulity,
$\tau$~is acyclic.

 The homotopy category $\Hot(C\comodl)$ and the derived category
$\sD(C\comodl)$ can be identified in an obvious way with the homotopy
category $\Hot(\oC\comodl)$ and the derived category $\sD(\oC\comodl)$
of the abelian category of graded $\oC$\+comodules.
 Analogously, the homotopy category $\Hot(C\contra)$ and the derived
category $\sD(C\contra)$ can be identified with the homotopy category
$\Hot(\oC\contra)$ and the derived category $\sD(C\contra)$ of
the abelian category of graded $\oC$\+contramodules.
 Finally, an obvious transformation of the bigrading allows to identify
the homotopy category $\Hot(A\modl)$ and the derived category
$\sD(A\modl)$ with the homotopy category $\Hot(\oA\modl)$ and
the derived category $\sD(\oA\modl)$ of the abelian category of
graded $\oA$\+modules.

 Hence we obtain the equivalences of derived categories
$\sD(\oA\modl^\up)\simeq\sD(\oC\comodl^\up)$ and $\sD(\oA\modl^\down)
\simeq\sD(\oC\contra^\down)$, where $\oA\modl^\up$, \ $\oA\modl^\down$,
\ $\oC\comodl^\up$, and $\oC\contra^\down$ denote the abelian
categories of nonnegatively graded $\oA$\+modules, nonpositively
graded $\oA$\+modules, nonnegatively graded $\oC$\+comodules, and
nonpositively graded $\oC$\+contramodules, respectively.
 Notice that these are unbounded derived categories of the abelian
categories of the categories of graded modules, comodules, and
contramodules satisfying no finiteness conditions.
 The only restriction imposed on the complexes of modules, comodules,
and contramodules is that of the nonpositivity or nonnegativity of
the internal grading.

\begin{rem}
 The results of Theorems~1--2 can be generalized to the following
situation.
 Let $\sE$ be an (associative, noncommutative) tensor DG\+category
(with shifts and cones) and $\sA\subset H^0(\sE)$ be a full
triangulated subcategory that is a two-sided tensor ideal.
 Consider graded algebra objects $A$ and graded coalgebra objects 
$C$ over $\sE$ such that $A_n=0=C_n$ for $n<0$ and both $A_0$ and
$C_0$ are the unit object of~$\sE$.
 Let us call a morphism of algebras $A'\rarrow A''$ or a morphism
of coalgebras $C'\rarrow C''$ a quasi-isomorphism if its cone
belongs to $\sA$ in every degree.
 Then the analogue of Theorem~1 holds for such graded algebras and
coalgebras.
 Furthermore, let $\sM$ be a left module DG\+category (with shifts
and cones) over the tensor DG\+category $\sE$ and let
$\sB\subset H^0(\sM)$ be a full triangulated subcategory such that
$H^0(\sE)\ot \sB\subset\sB$ and $\sA\ot H^0(\sM)\subset\sB$.
 Consider nonnegatively graded module objects $M$ and comodule
objects $N$ in $\sM$ over the graded algebras $A$ and graded
coalgebras~$C$.
 Let us call a morphism of modules $M'\rarrow M''$ or a morphism
of comodules $N'\rarrow N''$ a quasi-isomorphism if its cone
belongs to $\sB$ in every degree.
 Then the analogue of Theorem~2 holds for such graded modules and
comodules.
 In particular, one can take $\sE$ to be the tensor DG\+category of
DG\+bimodules over a DG\+ring $R$ that are flat as right DG\+modules,
$\sM$ to be the left module DG\+category of left DG\+modules over $R$,
and $\sA$ and $\sB$ to be the triangulated subcategories of acyclic
DG\+(bi)modules.
 One can also take $\sE$ to be the tensor DG\+category of
CDG\+bicomodules over a CDG\+coalgebra $C$ that are injective as
right graded comodules, $\sM$ to be the left module DG\+category of
left CDG\+comodules over $C$, and $\sA$ and $\sB$ to be
the triangulated subcategories of coacyclic CDG\+(bi)comodules.
 One can also use right module DG\+categories opposite to
the DG\+categories of left DG\+modules or left CDG\+contramodules
with respect to the action functors $\Hom_R$ or $\Cohom_C$,
assuming DG\+bimodules to be projective as left DG\+modules or
CDG\+bicomodules to be injective as left graded comodules.
\end{rem}

\subsection{Contravariant duality}
 Let $R$ be a (noncommutative) ring and let $A$ and $B$ be two
internally graded DG\+rings such that $A_n=0=B_{-n}$ for $n<0$, \ 
$A_0^j=0=B_0^j$ for $j\ne 0$, and $A_0^0=R=B_0^0$.
 Assume further that all the internal degree components $A_n$ are
bounded complexes of finitely generated projective left $R$\+modules
and all the components $B_n$ are bounded complexes of finitely
generated projective right $R$\+modules.
 We will consider nonnegatively internally graded left DG\+modules
$M$ over $A$ and nonpositively internally graded right DG\+modules
$N$ over $B$, and we will assume all their internal degree components
$M_n$ and $N_n$ to be bounded complexes of finitely generated
projective $R$\+modules.
 Let $\Hot(A\modl^\up)$ and $\Hot(\modrdown B)$ denote the homotopy
categories of such DG\+modules, and $\sD(A\modl^\up)$ and
$\sD(\modrdown B)$ denote the corresponding derived categories.

 Given two internally graded DG\+rings $A$ and $B$ as above, construct
two opposite DG\+rings $T^l$ and $T^r$ in the following way.
 Let $T^l=T^r$ be the subcomplex in $\prod_n B_{-n}\ot_R A_n$
consisting of all the elements $t$ satisfying the equation
$rt=tr$ for all $r\in R$.
 Define the multiplication in $T^l$ by the formula $(b'\ot a')
(b''\ot a'') = (-1)^{|a'|(|b''|+|a''|)} b'b''\ot a''a'$ and
the multiplication in $T^r$ by the formula $(b'\ot a')(b''\ot a'')
= (-1)^{|b''|(|b'|+|a'|)}b''b'\ot a'a''$.
 It is easy to see that these multiplications are well-defined,
associative and compatible with the differential, and that $T^l$ and
$T^r$ are two opposite DG\+rings in the sense of~\ref{second-tor-ext}.

 By a \emph{twisting cochain} for the DG\+rings $A$ and $B$ we
mean an element $\tau\in\prod_{n=1}^\infty B_{-n}\ot_R A_n$ of
(cohomological) degree~$1$ satisfying the equation $\tau^2+d\tau=0$
in $T^l$, or equivalently, the equation $\tau^2-d\tau=0$ in $T^r$.

 For a DG\+ring $A$, consider the internally and cohomologically
graded $R$\+$R$\+bimodule $\Hom_R(A_+,R)[-1]$, where $A_+=A/R$.
 There in a natural ``trace'' element $\tau_A\in\prod_{n=1}^\infty
\Hom_R(A_n,R)\ot_R A_n$.
 Consider the ring $B=\Cb^l(A)$ freely generated by the bimodule
$\Hom_R(A_+,R)[-1]$ over~$R$.
 Then $\tau_A$ is an element of degree~$1$ in the cohomologically
graded ring~$T^l$.
 Endow $B=\Cb^l(A)$ with the unique differential for which
$\tau_A^2+d\tau_A=0$ in~$T^l$.
 Then $B$ is an internally graded DG\+ring and $\tau_A$ is a twisting
cochain for $A$ and~$B$.
 Analogously, for a DG\+ring $B$ consider the internally and
cohomologically graded $R$\+$R$\+bimodule $\Hom_{R^\rop}(B_+,R)[-1]$,
where $B_+=B/R$.
 There is a natural ``trace'' element $\tau_B\in\prod_{n=1}^\infty
B_{-n}\ot_R\Hom_{R^\rop}(B_{-n},R)$.
 Consider the ring $A=\Cb^r(B)$ freely generated by the bimodule
$\Hom_{R^\rop}(B_+,R)[-1]$ over~$R$.
 Then $\tau_B$ is an element of degree~$1$ in the cohomologically
graded ring~$T^l$.
 Endow $A=\Cb^r(B)$ with the unique differential for which
$\tau_B^2+d\tau_B=0$ in~$T^l$.
 Then $A$ is an internally graded DG\+ring and $\tau_A$ is a twisting
cochain for $A$ and~$B$.

 There are natural bijective correspodences between three sets:
the set of morphisms of internally graded DG\+rings $\Cb^l(A)\rarrow B$,
the set of morphisms of internally graded DG\+rings $\Cb^r(B)\rarrow A$,
and the set of all twisting cochains $\tau$ for $A$ and~$B$.
 So the contravariant functors $\Cb^l$ and $\Cb^r$ are left adjoint
to each other.

 Let $\tau$ be a twisting cochain for $A$ and $B$.
 Given a DG\+module $M$ over $A$, construct a DG\+module
$\Hom_\tau(M,B)$ over $B$ in the following way.
 As an internally and cohomologically graded right $B$\+module,
$\Hom_\tau(M,B)$ coincides with $\Hom_R(M,B)$, where the action
of $B$ is given by the formula $(fb)(x)=(-1)^{|b||x|}f(x)b$.
 On the complex $\Hom_R(M,B)$ there is a structure of left
DG\+module over the DG\+ring $T^l$ given by the formula
$(b\ot a)(f)(x)=(-1)^{|a||f|}bf(ax)$ for $x\in M$ and
$f\in\Hom_R(M,B)$.
 The left action of $T^l$ on $\Hom_R(M,B)$ commutes with the right
action of~$B$.
 Thus one can twist the differential~$d$ on $\Hom_R(M,B)$ by replacing
it with $d+\tau$, where $\tau$ acts on $\Hom_R(M,B)$ as an element
of~$T^l$.
 So we obtain the DG\+module $\Hom_\tau(M,B)$ over~$B$.

 Given a DG\+module $N$ over $B$, construct a DG\+module
$\Hom(N,A)_\tau$ over $A$ as follows.
 As an internally and cohomologically graded left $A$\+module,
$\Hom(N,A)_\tau$ coincides with $\Hom_{R^\rop}(N,A)$.
 On the complex $\Hom_{R^\rop}(N,A)$ there is a structure of left
DG\+module over the DG\+ring $T^l$ given by the formula
$(b\ot a)(g)(y)=(-1)^{(|b|+|a|)(|g|+|y|)}g(yb)a$ for $y\in N$
and $g\in \Hom_{R^\rop}(N,A)$.
 The left action of $T^l$ on $\Hom_{R^\rop}(N,A)$ supercommutes with
the left action of~$A$.
 Thus one can twist the differential~$d$ on $\Hom_{R^\rop}(N,A)$ by
replacing it with $d+\tau$, where $\tau$ acts as an element of~$T^l$.
 So we obtain the DG\+module $\Hom(N,A)_\tau$ over~$A$.

\begin{thm1}
 Let $\tau$ be a twisting cochain between internally graded DG\+rings
$A$ and~$B$.
 Then the following conditions are equivalent: \par
\textup{(a)} the morphism of DG\+rings $\Cb^l(A)\rarrow B$
is a quasi-isomorphism; \par
\textup{(b)} the morphism of DG\+rings $\Cb^r(B)\rarrow A$
is a quasi-isomorphism; \par
\textup{(c)} the complex $\Hom_\tau(A,B)$ is quasi-isomorphic
to~$R$; \par
\textup{(d)} the complex $\Hom(B,A)_\tau$ is quasi-isomorphic
to~$R$. \par
 Besides, the functors\/ $\Cb^l$ and\/ $\Cb^r$ transform
quasi-isomorphisms of internally graded DG\+rings to quasi-isomorphisms.
 Consequently, they induce an anti-equivalence between the categories
of DG\+rings of the $A$ and $B$ kind with inverted quasi-isomorphisms.
\end{thm1}

\begin{proof}
 The complex $\Hom_{\tau_A}(A,\Cb^l(A))$ is the reduced
cobar-resolution of the left DG\+module $R$ over $A$ relative to
$R$, and the complex $\Hom(\Cb^l(A),A)_{\tau_A}$ is the reduced
bar-resolution of the left DG\+module $R$ over $A$ relative
to~$R$.
 So these complexes are quasi-isomorphic to~$R$.
 The rest of the proof is analogous to that of
Theorem~\ref{covariant-duality}.1.
 It uses the exactness and reflexivity properties of the tensor
product, $\Hom$, and duality for bounded complexes of finitely
generated projective modules.
\end{proof}

 A twisting cochain $\tau$ satisfying the equivalent conditions of
Theorem~1 is said to be \emph{acyclic}.
 Any twisting cochain $\tau$ induces a pair of left adjoint
contravariant functors $M\mpsto\Hom_\tau(M,B)\:\Hot(A\modl^\up)^\op
\rarrow\Hot(\modrdown B)$ and $N\mpsto\Hom(N,A)_\tau\:\Hot(\modrdown B)
^\op\rarrow\Hot(A\modl^\up)$.

\begin{thm2}
 Let $\tau$ be an acyclic twisting cochain for internally graded
DG\+rings $A$ and~$B$.
 Then the functors\/ $\Hom_\tau({-},B)$ and\/ $\Hom({-},A)_\tau$ induce
mutually inverse anti-equivalences between the derived categories\/
$\sD(A\modl^\up)$ and\/ $\sD(\modrdown B)$.
\end{thm2}

\begin{proof}
 The functor $\Hom_\tau({-},B)$ preserves acyclicity of DG\+modules,
since the internal degree components of $\Hom_\tau(M,B)$ are obtained
from the internal degree components of $M$ by applying the functor
$\Hom_R$ into the internal degree components of $B$ and taking shifts
and cones.
 To check that the adjunction morphism $\Hom(\Hom_\tau(M,B),A)_\tau
\rarrow M$ is a quasi-isomorphism, notice that the internal degree
components of its kernel can be obtained from the positive internal
degree components of $\Hom(B,A)_\tau$ by applying the functor~$\ot_R$
with the internal degree components of $M$ and taking shifts and cones.

 In fact, stronger assertions hold: every DG\+module in
$\Hot(A\modl^\up)$ that is projective as a graded $A$\+module is
projective in the sense of~\ref{projective-dg-modules}, and
consequently the functor $\Hom_\tau({-},B)$ maps acyclic DG\+modules
to contractible ones. 
\end{proof}

 The mutually inverse functors $\Hom_\tau({-},B)$ and $\Hom({-},A)_\tau$
transform the trivial DG\+module $R$ over $A$ into the free DG\+module
$B$ over $B$ and the free DG\+module $A$ over $A$ into
the trivial DG\+module $R$ over~$B$.

 For any DG\+module $M$ over $A$, the complex $\Hom_\tau(M,B)$
computes $\Ext_A(M,R)$.
 Indeed, $\Hom(\Hom_\tau(M,B),A)_\tau$ is a projective DG\+module
over~$A$ quasi-isomorphic to~$M$, and $\Hom_\tau(M,B)\simeq
\Hom_A(\Hom(\Hom_\tau(M,B),A)_\tau,R)$.
 Analogously, for any DG\+module $N$ over $B$, the complex
$\Hom(N,A)_\tau$ computes $\Ext_{B^\rop}(N,R)$.
 The DG\+ring $B$ itself computes the Yoneda Ext-ring
$\Ext_A(R,R)^\rop$, and the DG\+ring $A$ computes $\Ext_B(R,R)$.

\begin{cor}
 Let $f\:A'\rarrow A''$ be a morphism of nonnegatively internally
graded DG\+rings as above.
 Then the functor $\boI R_f\:\sD(A''\modl^\up)\rarrow\sD(A'\modl^\up)$
is an equivalence of triangulated categories if and only if
$f$~is a quasi-isomorphism.
\end{cor}

\begin{proof}
 This can be deduced either from Theorem~\ref{dg-mod-scalars}, or
from Theorems~1 and~2.
\end{proof}

\begin{rem}
 In the case when $R$ has a finite (left and right) homological
dimension there are several other versions of the above duality.
 E.~g., one can consider DG\+rings $A$ and $B$ whose internal
degree components $A_n$ and $B_n$ are bounded complexes of finitely
generated projective $R$\+modules, and DG\+modules $M$ and $N$ whose
components $M_n$ and $N_n$ are unbounded complexes of finitely
generated projective $R$\+modules.
 One can also consider DG\+rings $A$ and DG\+modules $M$ whose
components are complexes of finitely generated projective left
$R$\+modules bounded from above, and DG\+rings $B$ and DG\+modules
$N$ whose components are complexes of finitely generated projective
right $R$\+modules bounded from below.
 The above Theorems 1 and~2 hold in these situations.
 Assuming additionally that $R$ is (left and right) Noetherian one can
even drop the projectivity assumption on the bigrading components of
the modules in all of these results.
 This is so because the derived categories of modules whose
components $M_n$ and $N_n$ are (appropriately bounded) complexes
of finitely generated $R$\+modules are equivalent to the derived
categories of modules whose components are (accordingly bounded)
complexes of finitely generated projective $R$\+modules.
 (Cf.~\cite{MR}.)
\end{rem}

\Section{$\D$--$\Omega$ Duality}

\subsection{Duality functors}  \label{d-omega-duality-functors}
 Let $X$ be a scheme; denote by $\O_X$ its structure sheaf.
 A \emph{quasi-coherent CDG\+algebra} $\B$ over $X$ is a graded
quasi-coherent $\O_X$\+algebra $\B=\bigoplus_i \B^i$ endowed with
the following structure.
 For each affine open subscheme $U\subset X$ there is a structure
of CDG\+ring on the graded ring $\B(U)$, i.~e., a (not necessarily
$\O_X$\+linear) odd derivation~$d\:\B(U)\rarrow\B(U)$ of degree~$1$
and an element $h\in\B^2(U)$ satisfying the usual equations are given.
 For each pair of embedded affine subschemes $U\subset V\subset X$
an element $a_{UV}\in \B^1(U)$ is fixed such that the pair
$(\rho_{UV},a_{UV})$ is a morphism of CDG\+rings $\B(V)\rarrow\B(U)$,
where $\rho_{UV}$ is the restriction map in the sheaf~$\B$.
 For any three embedded affine subschemes
$U\subset V\subset W\subset X$, the morphism $(\rho_{UW},a_{UW})$
is equal to the composition of the morphisms $(\rho_{VW},a_{VW})$
and $(\rho_{UV},a_{UV})$, i.~e., the equation $a_{UW}=a_{UV}+
\rho_{UV}(a_{VW})$ holds.

 A morphism of quasi-coherent CDG\+algebras $\B\rarrow\A$ over~$X$
is a morphism of graded quasi-coherent $\O_X$\+algebras together with
a family of morphisms of CDG\+rings $\B(U)\rarrow\A(U)$ defined
for all affine open subschemes $U\subset X$ and satisfying the obvious
compatibility conditions.
 A \emph{quasi-coherent} (left or right) \emph{CDG\+module} $\N$ over
a quasi-coherent CDG\+algebra $\B$ is an $\O_X$\+quasi-coherent
sheaf of graded modules over $\B$ together with a family of
differentials $d\:\N(U)\rarrow\N(U)$ of degree~$1$ on the graded
$\B(U)$\+modules $\N(U)$ defined for all affine open subschemes
$U\subset X$ which make $\N(U)$ into CDG\+modules over
the CDG\+rings $\B(U)$ and satisfy the obvious compatibilities
with respect to the restriction morphisms $\N(V)\rarrow\N(U)$ and
$\B(V)\rarrow\B(U)$.

 Quasi-coherent left CDG\+modules over a fixed quasi-coherent
CDG\+algebra $\B$ form a DG\+category denoted by $\sDG(\B\modl)$
and quasi-coherent right CDG\+modules form the DG\+category
$\sDG(\modr\B)$; the corresponding homotopy categories are
denoted by $\Hot(\B\modl)$ and $\Hot(\modr\B)$.
 When $X$ is affine, the category of quasi-coherent CDG\+algebras
$\B$ over $X$ is equivalent to the category of CDG\+rings $(B,d,h)$
for which the graded ring $B$ is an $\O(X)$\+algebra;
the categories $\sDG(\B\modl)$ and $\sDG(\modr\B)$ are equivalent
to the DG\+categories $\sDG(B\modl)$ and $\sDG(\modr B)$.

{\hbadness=1200
 The coderived categories $\sD^\co(\B\modl)$
and $\sD^\co(\modr\B)$ are defined in the usual way as the quotient
categories of the homotopy categories $\Hot(\B\modl)$ and
$\Hot(\modr\B)$ by their minimal triangulated subcategories
containing the total quasi-coherent CDG\+modules of the exact triples
in the abelian category of quasi-coherent CDG\+modules (with closed
morphisms between them) and closed under infinite direct sums.
 Notice that we do \emph{not} define contraderived categories of
quasi-coherent CDG\+modules, because infinite products of
quasi-coherent sheaves are not well-behaved.
 A reader who finds the above definitions too sketchy may wish to
consult Section~\ref{cdg-module-section} for a detailed discussion
of the affine case. \par}

\medskip
 From now on, let $X$ be a smooth algebraic variety (a smooth
separated scheme of finite type) over a field~$k$.
 Let $\E$ be an algebraic vector bundle (locally free sheaf)
over~$X$.
 We denote by $\D_{X,\E}$ the sheaf of differential operators
on $X$ acting in the sections of~$\E$.
 These words have a unique meaning in the case of a field~$k$
of characteristic~$0$, and in the finite characteristic case we
are interested in the ``crystalline'' differential operators.
 More precisely, let $\END(\E)$ be the sheaf of endomorphisms
of~$\E$ and $F_1\D_{X,\E}$ be the sheaf of differential operators
of order at most~$1$ acting in the sections of~$\E$.
 Then we define $\D_{X,\E}$ as the sheaf of rings on $X$ generated
by the bimodule $F_1\D_{X,\E}\supset\END(\E)$ over the quasi-coherent
sheaf of $\O_X$\+algebras $\END(\E)$ with the relations $uv-vu=[u,v]$
for any two local sections $u$ and $v$ of $F_1\D_{X,\E}$ such that
the image of at least one of them in $F_1\D_{X,E}/\END(\E)\simeq
\END(\E)\ot_{\O_X}\!\.\T_X$ belongs to $\T_X\subset\END(\E)\ot_{\O_X}
\!\.\T_X$, where $\T_X$ denotes the sheaf of vector fields on~$X$.

 Assume that $X$ is affine.
 Then there exists a global connection $\nabla_\E$ in the vector
bundle~$\E$.
 Denote by $\End(\E)$ and $D_{X,\E}$ the global section rings
of the sheaves of rings $\END(\E)$ and $\D_{X,\E}$ on~$X$.
 Let $\Omega(X,\End(\E))=\Omega(X)\ot_{\O(X)}\End(\E)$
denote the graded ring of global differential forms on $X$ with
coefficients in $\End(\E)$.
 The connection $\nabla_\E$ induces a connection $\nabla_{\END(\E)}$
in the vector bundle $\END(\E)$, so there is the de Rham
differential~$d$ on $\Omega(X,\End(\E))$ related to
the connection~$\nabla_{\END(\E)}$.
 Let $h\in\Omega^2(X,\End(\E))$ be the curvature of
the connection~$\nabla_\E$.
 Then the differential~$d$ and the curvature element~$h$ define
a structure of CDG\+ring on the graded ring $\Omega(X,\End(\E))$.
 The de Rham differential on the graded
$\Omega(X,\End(\E))$\+module $K(X,\E)=\Omega(X,\End(\E))
\ot_{\End(\E)}D_{X,\E}\simeq\Omega(X)\ot_{\O(X)}D_{X,\E}$ makes it
a left CDG\+module over $\Omega(X,\End(\E))$ with a commuting
structure of a right $D_{X,\E}$\+module.

 Now let $X$ be an arbitrary smooth algebraic variety.
 Denote by $\Omega_X(\END(\E))=\Omega_X\ot_{\O_X}\END(\E)$
the sheaf of differential forms on $X$ with coefficients in
the vector bundle $\END(\E)$; it is a graded quasi-coherent
$\O_X$\+algebra.
 Choosing connections $\nabla_{\E|_U}$ on the restrictions of $\E$
on all the affine open subschemes $U\subset X$, one can define
a structure of quasi-coherent CDG\+algebra on $\Omega_X(\END(\E))$.
 The graded $\Omega_X(\END(\E))$\+module $\K(X,\E) =
\Omega_X(\END(\E))\ot_{\END(\E)}\D_{X,\E}\simeq
\Omega_X\ot_{\O_X}\D_{X,\E}$ becomes a quasi-coherent left
CDG\+module over $\Omega_X(\END(\E))$ with a commuting structure
of a (quasi-coherent) right $\D_{X,\E}$\+module.

 To any complex of right $\D_{X,\E}$\+modules $\M$ one can assign
the quasi-coherent right CDG\+module $\HOM_{\D_{X,\E}^\rop}
(\K(X,\E),\M)\simeq\HOM_{\END(\E)^\rop}(\Omega_X(\END(\E)),\M)$
over $\Omega_X(\END(\E))$.
 Conversely, to any quasi-coherent right CDG\+module $\N$ over
$\Omega_X(\END\allowbreak(\E))$ one can assign the complex of right
$\D_{X,\E}$\+modules $\N\ot_{\Omega_X(\END(\E))}\K(X,\E)\simeq
\N\ot_{\END(\E)}\D_{X,\E}$.
 This defines a pair of adjoint functors between the homotopy
categories $\Hot(\modr\Omega_X(\END(\E)))\rarrow\Hot(\modr\D_{X,\E})$
and $\Hot(\modr\D_{X,\E})\rarrow\Hot(\modr\Omega_X(\END(\E)))$.

 Whenever $X$ is affine, there is also a pair of adjoint functors
between the left module categories.
 To any complex of left $D_{X,\E}$\+modules $P$ one can assign
the left CDG\+module $K(X,\E)\ot_{D_{X,\E}}P\simeq\Omega(X,\End(\E))
\ot_{\End(\E)}P$ over $\Omega(X,\End(\E))$.
 Conversely, to any left CDG\+module $Q$ over $\Omega(X,\End(\E))$
one can assign the complex of left $D_{X,\E}$\+modules $\Hom_
{\Omega(X,\End(\E))}(K(X,\E),Q)\simeq\Hom_{\End(\E)}(D_{X,\E}\;Q)$.
 This defines a pair of adjoint functors between the homotopy
categories $\Hot(D_{X,\E}\modl)\rarrow\Hot(\Omega(X,\End(\E))\modl)$
and $\Hot(\Omega(X,\End(\E))\modl)\rarrow\Hot(D_{X,\E}\modl)$.

\subsection{Duality theorem}  \label{d-omega-duality-theorem}
 Let $X$ be a smooth algebraic variety over a field~$k$ and $\E$
be an algebraic vector bundle over~$X$.
 Denote by $\Hot(\modr\D_{X,\E})$ and $\sD(\modr\D_{X,\E})$
the (unbounded) homotopy and derived categories of the abelian
category of (quasi-coherent) right $\D_{X,\E}$\+modules.

{\hbadness=3500
\begin{thm}
 \textup{(a)} The adjoint functors\/
$\HOM_{\END(\E)^\rop}(\Omega_X(\END(\E)),{-})\:\Hot(\modr\D_{X,\E})
\allowbreak\rarrow \Hot(\modr\Omega_X(\END(\E)))$ and\/
${-}\ot_{\END(\E)}\D_{X,\E}\:\Hot(\modr\Omega_X(\END(\E)))\rarrow 
\Hot(\modr\D_{X,\E})$ induce functors\/ $\sD(\modr\D_{X,\E})\rarrow
\sD^\co(\modr\Omega_X(\END(\E)))$ and\/ $\sD^\co(\modr\Omega_X
(\END(\E)))\rarrow\sD(\modr\D_{X,\E})$, which are mutually inverse
equivalences of triangulated categories. \par
 \textup{(b)} Assume that $X$ is affine. Then the adjoint functors\/
$\Omega(X,\End(\E))\ot_{\End(\E)}{-}\:
\Hot(D_{X,\E}\modl)\rarrow\Hot(\Omega(X,\End(\E))\modl)$ and\/
$\Hom_{\End(\E)}(D_{X,\E}\;{-})\:\allowbreak\Hot(\Omega(X,\End(\E))
\modl)\rarrow\Hot(D_{X,\E}\modl)$ induce functors\/ $\sD(D_{X,\E}\modl)
\rarrow\sD^\ctr(\Omega(X,\End(\E))\modl)$ and\/ $\sD^\ctr(\Omega(X,
\End(\E))\modl)\rarrow\sD(D_{X,\E}\modl)$, which are mutually inverse
equivalences of triangulated categories.
\end{thm}}

\begin{proof}
 We will prove part~(a); the proof of part~(b) is analogous up to
duality.
 For simplicity of notation, denote by $\P_{X,\E}$
the graded $\END(\E)$\+$\END(\E)$\+bimodule $\bigoplus_i
\END(\E)\ot_{\O_X}\Lambda^i\T_X$ of polyvector fields on $X$
with coefficients in $\END(\E)$; then $\HOM_{\END(\E)^\rop}
(\Omega_X(\END(\E)),\M)\simeq\M\ot_{\END(\E)}\P_{X,\E}$ for
any right $\END(\E)$\+module $\M$.
 We will make use of the decreasing filtration $F$ of
$\Omega_X(\END(\E))$ defined by the rule $F^i\Omega_X(\END(\E))=
\bigoplus_{j\ge i}\Omega_X^j(\END(\E))$, the dual increasing
filtration $F$ of $\P_{X,\E}$ given by the rule $F_i\P_{X,\E}=
\bigoplus_{j\le i}\P^j_{X,\E}$, where $\P^j_{X,\E}=\END(\E)\ot_{\O_X}
\Lambda^j\T_X$, and the increasing filtration $F$ of the sheaf
$\D_{X,\E}$ by the order of differential operators.
 To prove that the induced functors exist, one can notice
that the category of right $\D_{X,\E}$\+modules has a finite
homological dimension, so $\sD(\modr\D_{X,\E})\simeq\sD^\co
(\modr\D_{X,\E})$, and the functors ${-}\ot_{\END(\E)}\P_{X,\E}$
and ${-}\ot_{\END(\E)}\D_{X,\E}$ transform exact triples of
complexes of right $\D_{X,\E}$\+modules into exact triples of
quasi-coherent right CDG\+modules over $\Omega_X(\END(\E))$ and
vice versa.
 Alternatively, for any complex of right $\D_{X,\E}$\+modules
$\M$ consider the increasing filtration of the quasi-coherent
CDG\+module $\M\ot_{\END(\E)}\P_{X,\E}$ induced by the increasing
filtration $F$ of $\P_{X,\E}$.
 It is a filtration by quasi-coherent CDG\+submodules and
the associated quotient quasi-coherent CDG\+modules
$\M\ot_{\END(\E)}\P^i_{X,\E}$ only depend on the $\END(\E)$\+module
structures on the components of $\M$, so it suffices to know that
the category of right $\END(\E)$\+modules has a finite homological
dimension and consequently any acyclic complex of right
$\END(\E)$\+modules is coacyclic.

 Furthermore, for any complex of right $\D_{X,\E}$\+modules $\M$
we have to show that the morphism of complexes of right
$\D_{X,\E}$\+modules $\M\ot_{\END(\E)}\P_{X,\E}\ot_{\END(\E)}\D_{X,\E}
\rarrow\M$ is a quasi-isomorphism.
 Consider the increasing filtration $F$ of $\M\ot_{\END(\E)}\P_{X,\E}
\ot_{\END(\E)}\D_{X,\E}$ by complexes of $\END(\E)$\+submodules
induced by the increasing filtrations $F$ of $\P_{X,\E}$ and
$\D_{X,\E}$; then the associated quotient complex of this filtration
is the tensor product of the complex of right $\END(\E)$\+modules $\M$
with the Koszul complex of $\END(\E)$\+$\END(\E)$\+bimodules
$\P_{X,\E}\ot_{\END(\E)}\gr_F\D_{X,\E}$.
 Since the positive grading componenents of the Koszul complex are
finite acyclic complexes of flat left $\END(\E)$\+modules, we are done.
 Finally, let $\N$ be a quasi-coherent right CDG\+module over
$\Omega_X(\END(\E))$; let us prove that the cone of the morphism
of quasi-coherent right CDG\+modules $\N\rarrow\N\ot_{\END(\E)}
\D_{X,\E}\ot_{\END(\E)}\P_{X,\E}$ is coacyclic.
 Let us reduce the question to the case when $\N$ is a complex of
right $\END(\E)$\+modules with a trivial CDG\+module structure, i.~e.,
$\N\cdot\Omega^i_X(\END(\E))=0$ for $i>0$.
 Here it suffices to consider the finite decreasing filtration $G$
on $\N$ defined by the rules $G^i\N=\N\cdot F^i\Omega_X(\END(\E))$
and the induced filtration on $\N\ot_{\END(\E)}\D_{X,\E}
\ot_{\END(\E)}\P_{X,\E}$.
 Now when $\N$ is just a complex of right $\END(\E)$\+modules,
the increasing filtration $F$ of $\N\ot_{\END(\E)}\D_{X,\E}
\ot_{\END(\E)}\P_{X,\E}$ induced by the increasing filtrations $F$
of $\D_{X,\E}$ and $\P_{X,\E}$ is a filtration by quasi-coherent
CDG\+submodules, the associated quotient quasi-coherent CDG\+module
is also just a complex of right $\END(\E)$\+modules, and it remains
to use the acyclicity and flatness of the Koszul complex again.
\end{proof}

\subsection{Coderived and contraderived categories}
 Let $X$ be a smooth affine algebraic variety of dimension~$d$
over a field~$k$ and $\E$ be an algebraic vector bundle over~$X$.
 Consider also the algebraic vector bundle $\Omega^d\E^*=
\Omega_X^d\ot_{\O_X}\E^*$ over $X$, where $\E^*$ is the dual
vector bundle to~$\E$ and $\Omega_X^d$ is the line bundle of
differential forms of top degree.
 It is well-known that the ring $D_{X,\E}^\rop$ opposite to $D_{X,\E}$
is naturally isomorphic to $D_{X,\Omega^d\E^*}$, so right
$D_{X,\Omega^d\E^*}$\+modules can be considered as left
$D_{X,\E}$\+modules.

 On the other hand, there is a natural equivalence between
the DG\+category of right CDG\+modules over the CDG\+ring
$\Omega(X,\End(\Omega^d\E^*))$ corresponding to the vector
bundle $\Omega^d\E^*$ and the DG\+category left CDG\+modules
over the CDG\+ring $\Omega(X,\End(\E))$ corresponding to
the vector bundle~$\E$.
 This equivalence assigns to a right CDG\+module $N$ over
$\Omega(X,\End(\Omega^d\E^*))$ the left CDG\+module
$\Omega^d(X)\ot_{\O(X)}N$ over $\Omega(X,\End(\E))$.
 The corresponding two equivalences of homotopy categories transform
the functor $\Hom_{\End(\Omega^d\E^*)^\rop}(\Omega(X,\End(\Omega^d\E^*))
,{-})\:\Hot(\modr D_{X,\Omega^d\E^*})\rarrow\Hot(\modr\Omega(X,
\End(\Omega^d\E^*)))$ into the functor $\Omega(X,\End(\E))\ot_{\End(\E)}
{-}\:\Hot(D_{X,\E}\modl)\allowbreak\rarrow\Hot(\Omega(X,\End(\E))\modl)$.

 Thus by Theorem~\ref{d-omega-duality-theorem} the same functor
$\Omega(X,\End(\E))\ot_{\End(\E)}{-}$ induces equivalences of
the derived category $\sD(D_{X,\E}\modl)$ with both the coderived
category $\sD^\co(\Omega(X,\End(\E))\modl)$ and the contraderived
category $\sD^\ctr(\Omega(X,\End(\E))\modl)$ of CDG\+modules.
 The next Theorem provides a more explicit equivalence between
these coderived and contraderived categories making a commutative
diagram with the above two equivalences.

{\hbadness=10000
\begin{thm}
 There is a natural equivalence of triangulated categories\/
$\sD^\co(\Omega(X,\End(\E))\modl)\simeq
\sD^\ctr(\Omega(X,\End(\E))\modl)$.
\end{thm}}

\begin{proof}
 The proof is similar to that of
Theorem~\ref{gorenstein-cdg-ring-case}, with the following changes
(cf.~\cite[Theorem~5.5]{P}).
 In place of graded modules of finite injective (projective)
dimension, one considers the graded modules over the graded
ring $\Omega(X,\End(\E))$ that are induced from graded
modules over $\End(\E)$, or equivalently, coinduced from
graded modules over $\End(\E)$.
 These are the graded modules of the form $\Omega(X,\End(\E))
\ot_{\End(\E)}P$, or equivalently, of the form
$\Hom_{\End(\E)}(\Omega(X,\End(\E)),M)$, where $P$ and $M$
are graded left $\End(\E)$\+modules.
 The quotient category of the homotopy category of left CDG\+modules
over $\Omega(X,\End(\E))$ whose underlying graded modules are
induced (coinduced) from graded $\End(\E)$\+modules by the minimal
triangulated subcategory containing the exact triples of left
CDG\+modules that as exact triples of graded modules are induced
(coinduced) from exact triples of graded $\End(\E)$\+modules is
equivalent to both $\sD^\co(\Omega(X,\End(\E))\modl)$ and
$\sD^\ctr(\Omega(X,\End(\E))\modl)$.
\end{proof}

\subsection{Filtered $\D$\+modules}  \label{filtered-d-modules}
 Let $X$ be a smooth algebraic variety over a field~$k$ and
$\E$ be an algebraic vector bundle over~$X$.
 Let $F$ denote the increasing filtration of the sheaf of
differential operators $\D_{X,\E}$ by the order of
differential operators.
 Consider the category $\modrfu\D_{X,\E}$ of right
$\D_{X,\E}$\+modules $\M$ endowed with a cocomplete increasing
filtration $F$ (by quasi-coherent $\O_X$\+submodules) compatible
with the filtration $F$ on $\D_{X,\E}$ and such that $F_{-1}\M=0$.
 The category $\modrfu\D_{X,\E}$ has a natural exact category
structure; a triple of $\D_{X,\E}$\+modules is exact if and only
if all its filtration components are exact triples of $\O_X$\+modules.
 Denote by $\Hot(\modrfu\D_{X,\E})$ and $\sD(\modrfu\D_{X,\E})$
the (unbounded) homotopy and derived categories of this exact category.
 When $X$ is affine, one can also consider the exact category
$D_{X,\E}\modl_\fil^\down$ of left $D_{X,\E}$\+modules endowed with
a complete filtration $F$ indexed by nonpositive integers,
$\dsb\subset F_{-2}M\subset F_{-1}M\subset F_0M=M$, and compatible
with the filtration $F$ on $D_{X,\E}$.
 The homotopy and derived categories of this exact category will be
denoted by $\Hot(D_{X,\E}\modl_\fil^\down)$ and
$\sD(D_{X,\E}\modl_\fil^\down)$.

 Assuming that $X$ is affine, consider the graded ring
$\Omega(X,\End(\E))\sptilde$ constructed in~\ref{dg-categories-constr}
(now we are not interested in its DG\+ring structure).
 For any $X$, the graded rings $\Omega(U,\End(\E|_U))$ for affine
subschemes $U\subset X$ glue together naturally forming a sheaf
of graded rings $\Omega(X,\END(\E))\sptilde$ over~$X$.
 The sheaf $\Omega(X,\END(\E))\sptilde$ is nonpositively graded and
its zero-degree component is the sheaf of rings $\END(\E)$.
 Notice that $\Omega(X,\END(\E))\sptilde$ is not a quasi-coherent
$\O_X$\+algebra, as the image of $\O_X$ is not contained in its center.
 The abelian category of quasi-coherent sheaves of graded (left or
right) modules over the sheaf of graded rings $\Omega(X,\END(\E))
\sptilde$ is isomorphic to the abelian category of quasi-coherent
(left or right) CDG\+modules over the quasi-coherent CDG\+algebra
$\Omega(X,\END(\E))$ and closed morphisms between them.

 We will find it convenient to consider $\Omega(X,\END(\E))\sptilde$ as
a sheaf of bigraded rings with an internal and a cohomological grading,
or a sheaf of internally graded DG\+rings with a zero differential,
as explained in Appendix~\ref{homogeneous-koszul-appendix}.
 So we place $\Omega(X,\END(\E))\sptilde$ in the nonpositive internal
grading~$n$ and the nonnegative cohomological grading $i=-n$ running
from $i=0$ to $i=\dim X+1$.
 The obvious bigrading transformation identifies sheaves of
internally graded DG\+modules over $\Omega(X,\END(X))\sptilde$
considered as a sheaf of internally graded DG\+rings with complexes
of sheaves of graded modules over the sheaf of graded rings
$\Omega(X,\END(X))\sptilde$.
 We will prefer the DG\+module language.

 Denote by $\Hot(\modrup\Omega(X,\END(\E))\sptilde)$ the homotopy
category of sheaves of nonnegatively internally graded 
$\O_X$\+quasi-coherent right DG\+modules over
$\Omega(X,\END(\E))\sptilde$ and by
$\sD(\modrup\Omega(X,\END(\E))\sptilde)$ the corresponding derived
category, i.~e., the quotient category of the homotopy category by
the thick subcategory of sheaves of DG\+modules that are acyclic
as complexes of sheaves in every internal grading.
 When $X$ is affine, consider also the homotopy category
$\Hot(\Omega(X,\End(\E))\sptilde\modl^\down)$ of nonpositively
internally graded left DG\+modules over the internally graded
DG\+ring $\Omega(X,\End(\E))\sptilde$ and the corresponding
derived category $\sD(\Omega(X,\End(\E))\sptilde\modl^\down)$.
 These homotopy and derived categories are isomorphic to
the (unbounded) homotopy and derived categories of the corresponding
abelian categories of quasi-coherent (nonnegatively) graded sheaves
or (nonpositively) graded modules.

 The tensor product $\L(X,\E)=\Omega(X,\END(\E))\sptilde\ot_{\END(\E)}
\D_{X,\E}$ has a natural structure of a sheaf of (internally ungraded)
left DG\+modules over $\Omega(X,\END(\E))\sptilde$ with the differential
defined by the formula $d_\L(\phi\ot p) = (-1)^{|\phi|}\phi d_\K(1\ot p)
- (-1)^{|\phi|}\phi\delta\ot p$, where $d_\K$ denotes the differential
of the CDG\+module $\K(X,\E)=\Omega(X,\END(\E))\ot_{\END(\E)}\D_{X,\E}$,
while $\phi\in\Omega(X,\END(\E))\sptilde$ and $p\in\D_{X,\E}$.
 One defines the filtration $F$ on $\L(X,\E)$ as the tensor product of
the increasing filtration $F$ on $\D_{X,\E}$ and the filtration $F$
on $\Omega(X,\END(\E))\sptilde$ indexed by nonpositive integers and
associated with the internal grading of $\Omega(X,\END(\E))\sptilde$.
 The latter filtration is given by the rule $F_j\Omega(X,\END(\E))
\sptilde = \bigoplus_{n\le j}\Omega(X,\END(\E))_n\sptilde$, where
$n=-\dim X-1$,~\dots,~$0$.
 This endowes $\L(X,\E)$ with a structure of a complex of filtered
right $\D_{X,\E}$\+modules with the right action of $\D_{X,\E}$
commuting with the left action of $\Omega(X,\END(\E))\sptilde$.
 For $X$ affine, we denote by $L(X,\E)$ the DG\+module of
global sections of the sheaf $\L(X,\E)$.

{\hbadness=2000
 Using the sheaf of DG\+modules $\L(X,\E)$, one can construct a pair
of adjoint functors between the homotopy categories
$\Hot(\modrfu\D_{X,\E})$ and $\Hot(\modrup\Omega(X,\END(\E))\sptilde$.
 To a complex of filtered right $\D_{X,\E}$\+modules $\M$ one assigns
the sheaf of internally graded right DG\+modules
$\HOM_{\.\bigoplus_j F_j\D_{X,\E}^\rop}(\bigoplus_j F_j\L(X,\E)\;
\bigoplus_j F_j\M)\simeq\HOM_{\END(\E)^\rop}(\Omega(X,\END(\E))
\sptilde\;\bigoplus_j F_j\M)$, where the direct sums of the filtration
components $\bigoplus_j F_j\M$ etc.\ are considered as (sheaves of)
internally graded modules in the internal grading~$j$.
 To a sheaf of internally graded right DG\+modules $\N$ over
$\Omega(X,\END(\E))\sptilde$ one assigns the complex of filtered
right $\D_{X,\E}$\+modules $(\bigoplus_n\N_n)
\ot_{\Omega(X,\END(\E))}\L(X,\E)\simeq(\bigoplus_n\N_n)\ot_{\END(\E)}
\D_{X,\E}$, where $\bigoplus_n\N$ is the notation for the sheaf of
internally ungraded DG\+modules corresponding to~$\N$, and the tensor
products are endowed with the filtration~$F$ induced by
the filtration $F$ on $\L(X,\E)$ or $\D_{X,\E}$ and the increasing
filtration $F$ on $\bigoplus_n\N$ corresponding to the grading of~$\N$.
\par}

 When $X$ is affine, one can also use the DG\+module $L(X,\E)$ to
construct a pair of adjoint functors between the homotopy categories
$\Hot(D_{X,\E}\modl_\fil^\down)$ and $\Hot(\Omega(X,\End(\E))\sptilde
\modl^\down)$.
 To a complex of filtered left $D_{X,\E}$\+modules $M$ one assigns
the internally graded left DG\+module
$\bigoplus_j F_jL(X,\E)\ot_{\,\bigoplus_jF_jD_{X,\E}}
\bigoplus_jM/F_{j-1}M\simeq\Omega(X,\End(\E))\sptilde\ot_{\End(\E)}
\bigoplus_jM/F_{j-1}M$, where the direct sum of the filtration
quotients $\bigoplus_j M/F_{j-1}M$ is considered as an internally
graded $\bigoplus_j F_jD_{X,\E}$\+module with the internal grading~$j$.
 To an internally graded left DG\+module $N$ over
$\Omega(X,\End(\E))\sptilde$ one assigns the complex of filtered
left $D_{X,\E}$\+modules $\Hom_{\Omega(X,\End(\E))\sptilde}
(L(X,\E)\;\prod_n N_n)\simeq\Hom_{\End(\E)}(D_{X,\E}\;\prod_n N_n)$,
where $\prod_n N_n$ is the internally ungraded DG\+module constructed
by taking infinite products of components with a fixed cohomological
grading, and the modules $\Hom$ are endowed with the filtration $F$
induced by the filtration $F$ on $L(X,\E)$ or $D_{X,\E}$ and
the complete filtration $F$ on $\prod_n N_n$, indexed by nonpositive
integers and coming from the grading of~$N$.

{\hbadness=1600
\begin{thm}
 \textup{(a)} The adjoint functors $\M\mpsto\HOM_{\END(\E)^\rop}
(\Omega(X,\END(\E))\sptilde\;\bigoplus_j F_j\M)$ and $\N\mpsto
(\bigoplus_n\N_n)\ot_{\END(\E)}\D_{X,\E}$ between the homotopy
categories $\Hot(\modrfu\D_{X,\E})$ and\/ $\Hot(\modrup
\Omega(X,\END(\E))\sptilde)$ induce functors\/ $\sD(\modrfu\D_{X,\E})
\rarrow\sD(\modrup\allowbreak\Omega(X,\END(\E))\sptilde)$ and\/
$\sD(\modrup\Omega(X,\END(\E))\sptilde)\rarrow\sD(\modrfu\D_{X,\E})$,
which are mutually inverse equivalences of triangulated categories. \par
 \textup{(b)} Assume that $X$ is affine.
 Then the adjoint functors $M\mpsto\Omega(X,\End(\E))\sptilde
\ot\allowbreak_{\End(\E)}\bigoplus_jM/F_{j-1}M$ and $N\mpsto
\Hom_{\End(\E)}(D_{X,\E}\;\prod_n N_n)$ between the homotopy
categories\/ $\Hot(D_{X,\E}\modl_\fil^\down)$ and\/
$\Hot(\Omega(X,\End(\E))\sptilde\modl^\down)$ induce functors\/
$\sD(D_{X,\E}\modl_\fil^\down)\rarrow\sD(\Omega(X,\End(\E))\sptilde
\modl^\down)$ and\/ $\sD(\Omega(X,\End(\E))\sptilde\modl^\down)
\rarrow\sD(D_{X,\E}\modl_\fil^\down)$, which are mutually
inverse equivalences of triangulated categories.
\end{thm}}

\begin{proof}
 This is essentially a particular case of the generalization of
Theorem~\ref{covariant-duality}.2 described in
Remark~\ref{covariant-duality}.
 First notice that the derived category of the exact category of
filtered $\D_{X,\E}$\+modules with the filtrations indexed by
nonnegative integers is equivalent to the derived category of
the abelian category of (quasi-coherent) sheaves of nonnegatively graded
modules over the sheaf of graded rings $\bigoplus_j F_j\D_{X,\E}$,
the equivalence being given by the functor $\M\mpsto\bigoplus_jF_j\M$.
 Analogously, the derived category of the exact category of
filtered $D_{X,\E}$\+modules with the filtrations indexed by
nonpositive integers is equivalent to the derived category of
the abelian category of nonpositively graded modules over
the graded ring $\bigoplus_j F_jD_{X,\E}$, the equivalence
being given by the functor $M\mpsto\bigoplus_j M/F_{j-1}M$.
 Furthermore, as in the proof of Theorem~\ref{d-omega-duality-theorem},
denote by $\P_{X,\E}\sptilde$ the graded
$\END(\E)$\+$\END(\E)$\+bimodule $\HOM_{\END(\E)^\rop}
(\Omega(X,\END(\E))\sptilde,\END(\E))$.
 To check that the components of positive internal grading of
the complexes $\P_{X,\E}\sptilde\ot_{\END(\E)}\bigoplus_jF_j\D_{X,\E}$
and $(\bigoplus_jF_j\D_{X,\E})\ot_{\END(\E)}\P_{X,\E}\sptilde$ are
acyclic, consider the decreasing filtrations $G$ on these complexes
induced by the filtration $G$ on $\bigoplus_jF_j\D_{X,\E}$ given by
the rule $G^tF_j\D_{X,\E}=F_{j-t}\D_{X,\E}$ for $t\ge0$ and
the filtration $G$ on $\P_{X,\E}\sptilde$ given by the rules
$G^0\P_{X,\E}\sptilde=\P_{X,\E}\sptilde$, \ $G^1\P_{X,\E}\sptilde=
\ker(\P_{X,\E}\sptilde\to\P_{X,\E})$.
 The associated graded complexes to these filtrations are the tensor
products over~$k$ of the Koszul complexes $\P_{X,\E}\ot_{\END(\E)}
\gr_F\D_{X,\E}$ and $\gr_F\D_{X,\E}\ot_{\END(\E)}\P_{X,\E}$ with
the Koszul complex for the symmetric algebra and the exterior
coalgebra in one variable.
 (Notice that, in particular, the first-degree component of
$\P_{X,\E}\sptilde$ is isomorphic to $F_1\D_{X,\E}\sptilde$.)
 It remains to use the exactness of the tensor products of complexes
of $\END(\E)$\+modules with finite complexes of flat
$\END(\E)$\+modules and homomorphisms from finite complexes of 
projective $\End(\E)$\+modules into complexes of $\End(\E)$\+modules.
\end{proof}

 Notice that one has $\sD(\modrup\Omega(X,\END(\E))\sptilde)=
\sD^\co(\modrup\Omega(X,\END(\E))\sptilde)$, i.~e., an
$\O_X$\+quasi-coherent sheaf of nonnegatively internally graded
right DG\+modules over $\Omega(X,\END(\E))\sptilde$ is coacyclic
whenever it is acyclic.
 Analogously, one has $\sD(\Omega(X,\End(\E))\sptilde\modl^\down)
=\sD^\ctr(\Omega(X,\End(\E))\sptilde\modl^\down)$ when $X$ is affine.
 One proves this using the filtrations of (sheaves of) internally
graded DG\+modules associated with the internal gradings, as
explained in the proof of Theorem~\ref{covariant-duality}.2,
together with finiteness of the homological dimension of
the category of quasi-coherent sheaves over~$X$.

 The functor of forgetting the filtration, acting from the derived
category of filtered right $\D_{X,\E}$\+modules with filtrations
indexed by nonnegative integers to the derived category of right
$\D_{X,\E}$\+modules, corresponds to the following functor
$\sD(\modrup\Omega(X,\END(\E))\sptilde)\rarrow\sD^\co(\modr
\Omega(X,\END(\E)))$.
 Given an $\O_X$\+quasi-coherent sheaf $\N$ of nonnegatively graded
right DG\+modules over $\Omega(X,\END(\E))\sptilde$, one assigns
the quasi-coherent right CDG\+module $\bigoplus_n\N_n$ over
$\Omega(X,\END(\E))$ with the action of $\Omega(X,\END(\E))$ induced
by the action of $\Omega(X,\END(\E))\sptilde$ in $\N$ and
the differential $d=\delta+d'$, where $d'$ denotes the differential
in the DG\+module $\N$.
 When $X$ is affine, the functor of forgetting the filtration, acting
from the derived category of filtered left $D_{X,\E}$\+modules with
filtrations indexed by nonpositive integers to the derived category
of left $D_{X,\E}$\+modules, corresponds to the following functor
$\sD(\Omega(X,\End(\E))\sptilde\modl^\down)\rarrow\sD^\ctr
(\Omega(X,\End(\E))\modl)$.
 Given a nonpositively graded left DG\+module $N$ over
$\Omega(X,\End(\E))\sptilde$, one constructs the left CDG\+module
$\prod_n N_n$ over $\Omega(X,\End(\E))$ with the action of
$\Omega(X,\End(\E))$ induced by the action of $\Omega(X,\End(\E))
\sptilde$ and the differential $d=\delta+d'$.

\begin{rem}
 For an affine variety~$X$, the derived categories of filtered
$D_{X,\E}$\+modules with filtrations indexed by the integers
(see Example~\ref{functors-phi-and-psi}) can be also described as
certain exotic derived categories of (the abelian categories of)
graded modules over $\Omega(X,\End(\E))\sptilde$.
 Namely, define the semiderived category (cf.~\cite{P}) of graded
right modules over $\Omega(X,\End(\E))\sptilde$ (with the grading
indexed by the integers) as the quotient category of the homotopy
category of graded right modules over $\Omega(X,\End(\E))\sptilde$
by the thick subcategory formed by all complexes of graded modules
that are coacyclic as complexes of graded modules over
$\Omega(X,\End(\E))\subset\Omega(X,\End(\E))\sptilde$.
 Then the derived category of the exact category of filtered right
$D_{X,\E}$\+modules with complete and cocomplete filtrations indexed
by the integers (and compatible with the filtration $F$ on $D_{X,\E}$)
is equivalent to the semiderived category of graded right modules
over $\Omega(X,\End(\E))\sptilde$.
 Analogously, define the semiderived category of graded left modules
over $\Omega(X,\End(\E))\sptilde$ as the quotient
category of the homotopy category of graded left modules over
$\Omega(X,\End(\E))\sptilde$ by the thick subcategory formed by
all complexes that are contraacyclic as complexes of graded modules
over $\Omega(X,\End(\E))$.
 Then the derived category of the exact category of filtered left
$D_{X,\E}$\+modules with complete and cocomplete filtrations indexed
by the integers is equivalent to the semiderived category of graded
left modules over $\Omega(X,\End(\E))\sptilde$.
 The functors providing these equivalences of categories are defined
as follows.
 The functors assigning complexes of graded modules over
$\Omega(X,\End(\E))\sptilde$ to complexes of filtered
$D_{X,\E}$\+modules are constructed exactly as in
the nonnegative/nonpositive grading/filtration case above.
 When one constructs the functors assigning complexes of filtered
modules over $D_{X,\E}$ to complexes of graded modules over
$\Omega(X,\End(\E))\sptilde$, one has to pass to the completion
(for right modules) or the cocompletion (for left modules) with
respect to the filtration of $D_{X,\E}$\+modules after performing
the procedures described above for the nonnegative/nonpositive
grading/filtration case.
 The functors of passing to the associated graded
$\gr_FD_{X,\E}$\+modules on the derived categories of filtered
$D_{X,\E}$\+modules correspond to the forgetful functors assigning
graded modules over $\Omega(X,\End(\E))$ to graded modules over
$\Omega(X,\End(\E))\sptilde$ (in both the left and right module
situations).
 Here the derived category of graded $\gr_FD_{X,\E}$\+modules
is identified with the coderived and contraderived categories of
graded modules over $\Omega(X,\End(\E))$.
 The analogous results hold for arbitrary nonhomogeneous Koszul
rings over base rings of finite homological dimension (see~\cite{P})
in place of the filtered ring of differential operators.
\end{rem}

\subsection{Coherent $\D$\+modules}  \label{coherent-d-modules}
 Let $X$ be a smooth algebraic variety over a field~$k$ and
$\E$ be an algebraic vector bundle over~$X$.
 Let $\sD^\b(\modrcoh\D_{X,\E})$ denote the bounded derived category
of coherent (locally finitely generated) right $\D_{X,\E}$\+modules.
 Furthermore, let $\sD^\b(\modrfc\D_{X,\E})$ denote the bounded
derived category of the exact category of filtered right
$\D_{X,\E}$\+modules with locally finitely generated filtrations,
i.~e., coherent right $\D_{X,\E}$\+modules $\M$ with filtrations $F$
by coherent $\O_X$\+submodules compatible with the filtration $F$
on $\D_{X,\E}$ and such that $F_n\M=F_{n-1}\M\cdot F_1\D_{X,\E}$
for large~$n$.
 By~\cite[Exercise~III.6.8]{Har} there are enough vector bundles
on $X$, i.~e., any coherent sheaf on $X$ is the quotient sheaf of
a locally free sheaf (of finite rank).

 Let $\sD^\abs(\modrcoh\Omega_X(\END(\E)))$ denote the absolute
derived category of $\O_X$\+coherent right CDG\+modules
over the quasi-coherent CDG\+algebra $\Omega_X(\END(\E))$.
 Here a quasi-coherent right CDG\+module $\N$ over $\Omega_X(\END(\E))$
is said to be $\O_X$\+coherent if it has a finite number of nonzero
graded components only and all these graded components are coherent
$\O_X$\+modules.
 The absolute derived category $\sD^\abs(\modrcoh\Omega_X(\END(\E)))$
is defined as the quotient category of the homotopy category of
$\O_X$\+coherent right CDG\+modules over $\Omega_X(\END(\E))$ by
its minimal thick subcategory containing the total CDG\+modules of
exact triples of $\O_X$\+coherent right CDG\+modules over
$\Omega_X(\END(\E))$.
 Finally, let $\sD(\modrgc\Omega_X(\END(\E))\sptilde)$ denote
the derived category of $\O_X$\+coherent sheaves of internally
graded right DG\+modules over $\Omega_X(\END(\E))\sptilde$,
the $\O_X$\+coherence being defined as above.

 Although this is clearly not necessary for the next Theorem to be
valid, it will be convenient for us to presume our filtrations and
internal gradings to be indexed by nonnegative integers, and
incorporate this assumption into the definitions of
$\sD^\b(\modrfc\D_{X,\E})$ and $\sD(\modrgc\Omega_X(\END(\E))\sptilde)$.

{\hbadness=1100
\begin{thm}
\textup{(a)} The triangulated categories\/ $\sD^\b(\modrcoh\D_{X,\E})$
and\/ $\sD^\abs(\modrcoh\allowbreak\Omega_X(\END(\E)))$ are naturally
equivalent. \par
\textup{(b)} The triangulated categories\/ $\sD^\b(\modrfc\D_{X,\E})$
and\/ $\sD(\modrgc\Omega_X(\END(\E))\sptilde)$ are naturally
equivalent.
\end{thm}}

\begin{proof}
 There is a natural fully faithful functor $\sD^\b(\modrcoh\D_{X,\E})
\rarrow\sD(\modr\D_{X,\E})$.
 The natural functor $\sD^\abs(\modrcoh\Omega_X(\END(\E)))\rarrow
\sD^\co(\modr\Omega_X(\END(\E)))$ is also fully faithful, as one can
show in the way of the proof of Theorem~\ref{fin-gen-cdg-mod}.1.
 The natural functors $\sD^\b(\modrfc\D_{X,\E})\rarrow\
\sD(\modrfu\D_{X,\E})$ and $\sD(\modrgc\Omega(X,\END(\E))\sptilde)
\rarrow\sD(\modrup\Omega(X,\END(\E))\sptilde)$ are also clearly
fully faithful.
 We only have to show that the subcategories we are interested in
correspond to each other under the equivalences of categories
$\sD(\modr\D_{X,\E})\simeq\sD^\co(\modr\Omega_X(\END(\E)))$ and
$\sD(\modrfu\D_{X,\E})\simeq\sD(\modrup\Omega_X(\END(\E))\sptilde)$.
 It follows immediately from the constructions of the duality
functors that they map $\sD^\abs(\modrcoh\Omega_X(\END(\E)))$ into
$\sD^\b(\modrcoh\D_{X,\E})$ and
$\sD(\modrgc\Omega_X(\END(\E))\sptilde)$ into 
$\sD^\b(\modrfc\D_{X,\E})$.

 Furthermore, for any vector bundle $\F$ over $X$ the right
$\D_{X,\E}$\+module $\F\ot_{\O_X}\D_{X,\E}$ corresponds to
the trivial $\O_X$\+coherent CDG\+module $\F\ot_{\O_X}\END(\E)$ over
$\Omega(X,\END(\E))$.
 Analogously, for any $n\ge 0$ the filtered right $\D_{X,\E}$\+module
$\F(n)\ot_{\O_X}\D_{X,\E}$ with the filtration $F$ induced by
the filtration $F$ on $\D_{X,\E}$ and the filtration $F$ on $\F(n)=\F$
given by the rules $F_{n-1}\F(n)=0$, \ $F_n\F(n)=\F(n)$ corresponds to
the $\O_X$\+coherent sheaf of internally graded DG\+modules
$\F(n)\ot_{\O_X}\END(\E)$ over $\Omega(X,\END(\E))\sptilde$ living
in the internal degree~$n$ and the cohomological degree~$0$.
 So it remains to check that the filtered right $\D_{X,\E}$\+modules
$\F(n)\ot_{\O_X}\D_{X,\E}$ generate $\sD^\b(\modrfc\D_{X,\E})$,
as it will then follow that the right $\D_{X,\E}$\+modules 
$\F\ot_{\O_X}\D_{X,\E}$ generate $\sD^\b(\modrcoh\D_{X,\E})$.
 
 It is important that the Rees algebra $\bigoplus_j F_j\D_{X,\E}$
is Noetherian and has a finite homological dimension.
 Given a coherent filtered right $\D_{X,\E}$\+module $\M$, construct
its resolution $\dsb\rarrow\Q_1\rarrow\Q_0\rarrow\M$ consisting of
filtered $\D_{X,\E}$\+modules isomorphic to finite directs sums of
filtered $\D_{X,\E}$\+modules of the form $\F(n)\ot_{\O_X}\D_{X,\E}$.
 For $d$~large enough, the embedding $\K\rarrow\Q_{d-1}$ of
the image of the morphism $\Q_d\rarrow\Q_{d-1}$ will split locally
in~$X$ as an embedding of filtered right $\D_{X,\E}$\+modules.
 It follows that $F_0\K\cdot\D_{X,\E}\subset\K$ will be locally
a direct summand of $F_0\Q_{d-1}\cdot\D_{X,\E}$ and $\K/(F_0\K\cdot
\D_{X,\E})$ will be a direct summand of $\Q_{d-1}/(F_0\Q_{d-1}\cdot
\D_{X,\E})$ as filtered right $\D_{X,\E}$\+modules, so one can prove
that $\K$ is a direct sum of filtered $\D_{X,\E}$\+modules of
the form $\F(n)\ot_{\O_X}\D_{X,\E}$ arguing by induction.

 Another approach is to notice that all four our ``bounded-coherent''
triangulated categories consist of compact objects in
the corresponding larger triangulated categories and generate them
as triangulated categories with infinite direct sums.
 It follows immediately by a result of~\cite{Neem1} that
the subcategories of compact objects in the larger triangulated
categories are obtained by adjoining to the ``bounded-coherent''
categories the images of all idempotent endomorphisms.
 So the duality functors identify the desired ``bounded-coherent''
categories up to adjoining the images of idempotents.
 Three of the four ``bounded-coherent'' categories are bounded
derived categories of Noetherian abelian categories, so they are
closed under the images of idempotents.
 To check that the absolute derived category of $\O_X$\+coherent
CDG\+modules $\sD^\abs(\modrcoh\Omega(X,\END(\E)))$ also contains
images of its idempotents, one has to use the forgetful functor
$\sD^\b(\modrfc\D_{X,\E})\rarrow\sD^\b(\modrcoh\D_{X,\E})$ and
the corresponding ``totalization'' functor
$\sD(\modrgc\Omega(X,\END(\E))\sptilde)\rarrow
\sD^\abs(\modrcoh\Omega(X,\END(\E)))$.
\end{proof}

\begin{rem}
 All the results of this appendix are applicable to any filtered
sheaf of rings $(\D,F)$ over a separated Noetherian scheme $X$
such that the associated graded sheaf of rings $\gr_F\D$ is
a quasi-coherent graded $\O_X$\+algebra isomorphic to the tensor
product of the algebra of endomorphisms of a vector bundle $\E$
and the symmetric algebra of a vector bundle~$\L$.
 To such a filtered algebra one can assign a quasi-coherent
CDG\+algebra structure on the tensor product of the algebra of
endomorphisms of $\E$ and the exterior algebra of $\L^*$,
construct the duality functors, and then repeat verbatim
the arguments above.
 In particular, this applies to any sheaf of twisted differential
operators~\cite[section~2]{BB}, the enveloping algebra
corresponding to a central extension of a Lie algebroid
(cf.\ Example~\ref{koszul-generators}), etc.
 In fact, the only conditions one has to impose on the associated
graded algebra to $(\D,F)$ are the coherence and local freeness of
the components, finite homological dimension, Noetherianness,
and Koszulity.
 Even more general results can be found in~\cite[sections~0.4
and~11.8]{P}.
\end{rem}

\end{document}